\documentclass[12pt,a4paper]{amsart}

 \usepackage{framed}
\usepackage[T1]{fontenc}
\input xypic
\usepackage{amsmath}
\usepackage{amssymb}
\usepackage[usenames]{color}
\usepackage{amssymb}
\usepackage{latexsym}
 \usepackage{graphicx}
\usepackage{psfrag}

\newcommand{\noextension}[1]{} 

\newcommand{\finalpaper}[1]{} 
\newcommand{\arxivpreprint}{} 

\newcommand{\motivation}{} 

\newcommand{\modified}{}

\newcommand{\extension}{}

\newcommand{\MTEXT}{}

\newcommand{\MATTITEXT} {} 

\newcommand{\observation}[1]{}
\newcommand{\generalizations}[1]{}

\newcommand{\ryksi}{r_2}
\newcommand{\rkaksi}{r_1}

\def\cirt{\overline }
\def\Hsymbol{H}
\def\tsparameter{s}
\newcommand{\symbolP}{v}
\newcommand{\symbolQ}{w}

\newcommand{\cell}{\ell}
\newcommand{\hattuM}{M}
\newcommand{\Sclo}{\mathcal S^{cl}}
\newcommand{\Scle}{\mathcal S_e}

\newcommand{\xxi}{\xi}
\newcommand{\zzeta}{\zeta}

\renewcommand{\div}{\hbox{div}}

\newcommand{\Ric}{\hbox{Ric}}
\newcommand{\Ein}{\hbox{Ein}}
\def \Box {$\square$}

\newcommand{\mltext}{} 

\newcommand{\newtext}{}
\newcommand{\tobecheckedtext}{}

\newcommand{\HOX}[1]{}

\newcommand{\hiddenfootnote}[1]{}

\newcommand{\rhoepsilon}{{\rho}}
\newcommand{\bsequence}{{\bf b}}
\def\X{{\mathcal X}}
\def\Y{{\mathcal Y}}
\def\U{{\mathcal U}}
\def\G{{\mathcal G}}

\def\S{{\mathcal S}}
\def\T{{\mathcal T}}
\def\P{{\mathcal P}}
\def\qP {{\mathcal R}}

\def\pointear{{\bf e}}

\renewcommand{\H}{{\mathbb H}} 
\newcommand{\R}{{\mathbb R}} 
\newcommand{\D}{{\cal D}} 
\newcommand{\I}{{\cal I}} 
\newcommand{\be}{{\mathcal E}} 
 
\newcommand{\A}{{\cal A}}
\newcommand{\B}{{\cal B}}  
  
\newcommand{\K}{{\cal K}}  

\newcommand{\C}{{\mathbb C}} 
\newcommand{\N}{{\mathbb N}} 
\newcommand{\cal}{\mathcal }

\renewcommand{\L}{{\mathcal L}} 
\newcommand{\V}{{\cal V}} 
\newcommand{\W}{{\cal W}}

\def\hat{\widehat}
\def\tilde{\widetilde}
\def \bfo {\begin {eqnarray*} }
\def \efo {\end {eqnarray*} }
\def \ba {\begin {eqnarray*} }
\def \ea {\end {eqnarray*} }
\def \beq {\begin {eqnarray}}
\def \eeq {\end {eqnarray}}
\def \supp {\hbox{supp}\,}
\def \dim{\hbox{dim}\,}

\def \re {\hbox{Re}\,}
\def \im {\hbox{Im}\,}

\def \WF {\hbox{WF}\,}

\def \dist {\hbox{dist}}
\def\diag{\hbox{diag }}
\def \det {\hbox{det}}

\def\bra{\langle}
\def\cet{\rangle}

\def \e {\varepsilon}
\def \p {\partial}
\def \a {\alpha}
\renewcommand{\b}{\beta}

\def\M{{\mathcal M}}
\def\F{{\mathcal F}}

\def\Z{{\mathbb Z}}

\newtheorem{definition}{Definition}[section] 
\newtheorem{theorem}[definition]{Theorem} 
\newtheorem{lemma}[definition]{Lemma} 
 
\newtheorem{proposition}[definition]{Proposition} 
\newtheorem{corollary}[definition]{Corollary}

\begin{document}
\title[Inverse problems in spacetime I]
{Inverse problems in spacetime I:
Inverse problems for Einstein equations
\extension{\\ -- Extended preprint version}
}
\date{}
\author{Yaroslav Kurylev, Matti Lassas,  Gunther Uhlmann}
\address{Yaroslav Kurylev, UCL; Matti Lassas,  University of Helsinki;
Gunther Uhlmann, University of Washington,  and 
 University of Helsinki.  {\rm y.kurylev@ucl.ac.uk, Matti.Lassas@helsinki.fi,  
  gunther@math.washington.edu}}

\email{}

\maketitle


{\bf Abstract:}
{\it  
We consider inverse problems for the coupled Einstein equations  and
the matter field equations 
 on a 4-dimensional  globally hyperbolic
Lorentzian manifold $(M,g)$. 
We give a positive answer to the question: Do the active measurements, done in a neighborhood 
$U\subset M$
of a freely falling observed $\mu=\mu([s_-,s_+])$, determine the conformal structure of the spacetime in the minimal causal diamond-type set  $V_g=J_g^+(\mu(s_-))\cap J_g^-(\mu(s_+))\subset M$ containing $\mu$?

More precisely, we consider the Einstein equations  coupled with the scalar field equations
and study the system
$\Ein(g)=T$, $T=T(g,\phi)+\F_1$, and $\square_g\phi-\mathcal V^\prime(\phi)=\F_2$, where 
the sources $\F=(\F_1,\F_2)$ correspond to perturbations of 
the physical fields which we control.
The sources $\F$ need to be such that the fields $(g,\phi,\F)$  are solutions of this system and satisfy the
conservation law $\nabla_jT^{jk}=0$.
Let $(\hat g,\hat \phi)$ be the background fields corresponding to the vanishing source $\F$.
We prove that the observation of the solutions $(g,\phi)$
in the set $U$ 
corresponding to sufficiently small sources  $\F$ supported 
in  $U$ determine $V_{\hat g}$  as a differentiable manifold and the conformal structure of the metric $\hat g$ in the domain 
$V_{\hat g}$. 
The methods developed here have
potential to be applied to a large class of inverse problems for non-linear hyperbolic
equations  encountered e.g.\ in various practical imaging problems.
}



\tableofcontents




\section{Introduction and main results}
We consider inverse problems for the non-linear Einstein equations   coupled with matter field
equations. In this paper, we consider for the matter fields 
 the simplest possible model,
the  scalar field equations and  study the perturbations of a globally hyperbolic Lorent\-zian manifold $(M,\hat g)$ of dimension
$(1+3)$, {where 
the metric signature of
$\hat g$ is $(-,+,+,+)$.}

Roughly speaking,  we study the following problem: Can an observer 
in a space-time 
determine
the structure of the surrounding space-time by doing measurements near its
world line. More precisely, when 
$U_{\hat g}$ is a neighborhood of a time-like geodesic $\hat \mu$, we assume
that 
we can control sources supported in an open neighborhood
$W_{\hat g}\subset U_{\hat g}$ of $\hat \mu$ and measure the physical fields in the set $U_{\hat g}$. We ask,
can the properties of the metric (the metric itself or its conformal class)
be determined in a suitable larger 
set $J(p^-,p^+)$, $p^\pm=\hat \mu (s_\pm)$ that is not contained in the set $U_{\hat g}$, see Fig.\ 1(Left). 
This paper considers  inverse problems for active measurements and the corresponding
problem for passive measurements is studied in the second part of this paper,
\cite{Paper-part-2}. 

\begin{center}

\psfrag{1}{}
\psfrag{2}{$U_{\hat g}$}
\psfrag{3}{\hspace{-.2cm}$J({p^-},{p^+})$}
\psfrag{4} {}
\psfrag{5}{\hspace{-.2cm}$\mu_{z,\eta}$}
\includegraphics[height=5.5cm]{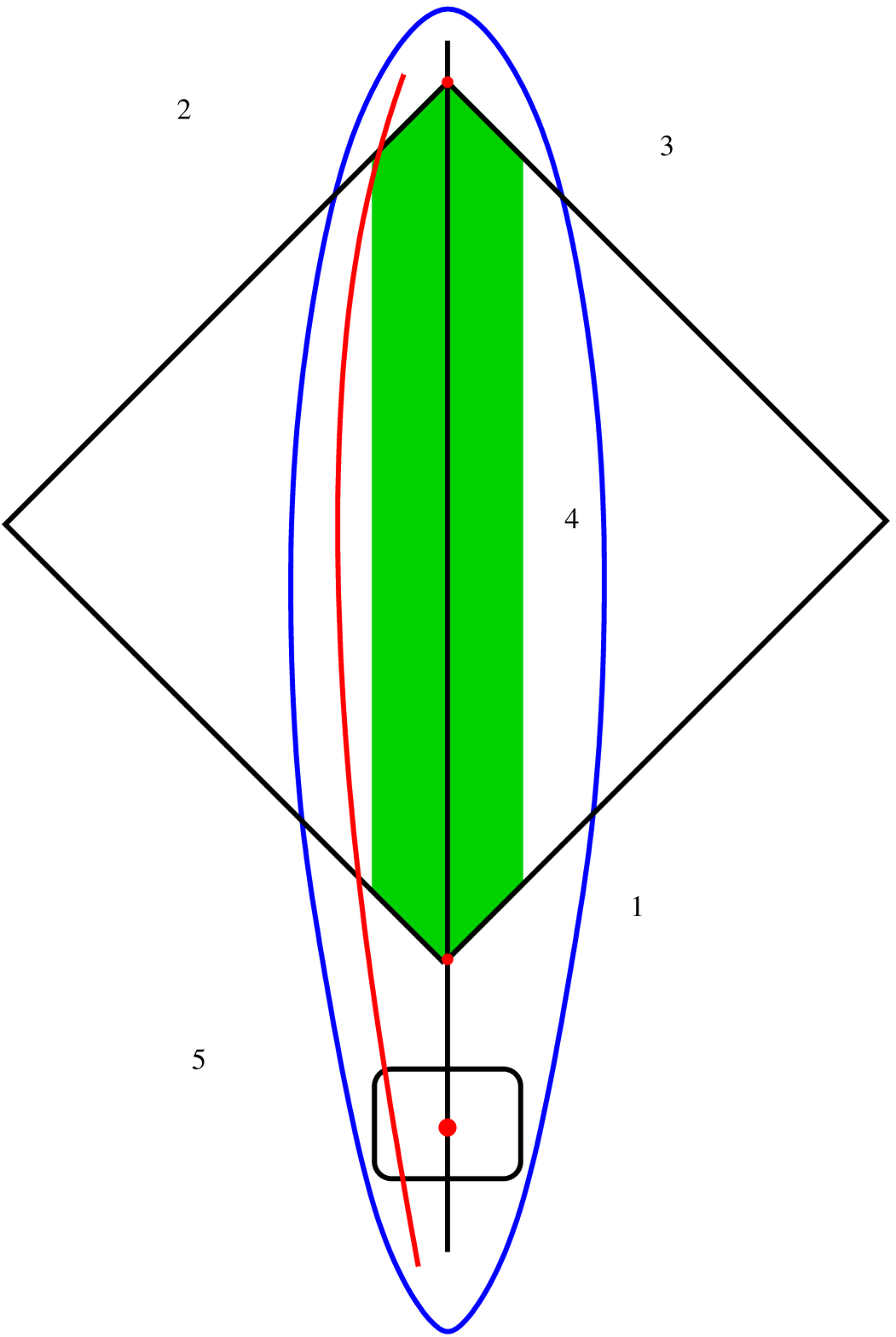}
\psfrag{1}{$y_0$}
\psfrag{2}{$\Sigma$}
\psfrag{3}{$\Sigma_1$}
\psfrag{4}{$y^\prime$}
\includegraphics[height=3.5cm]{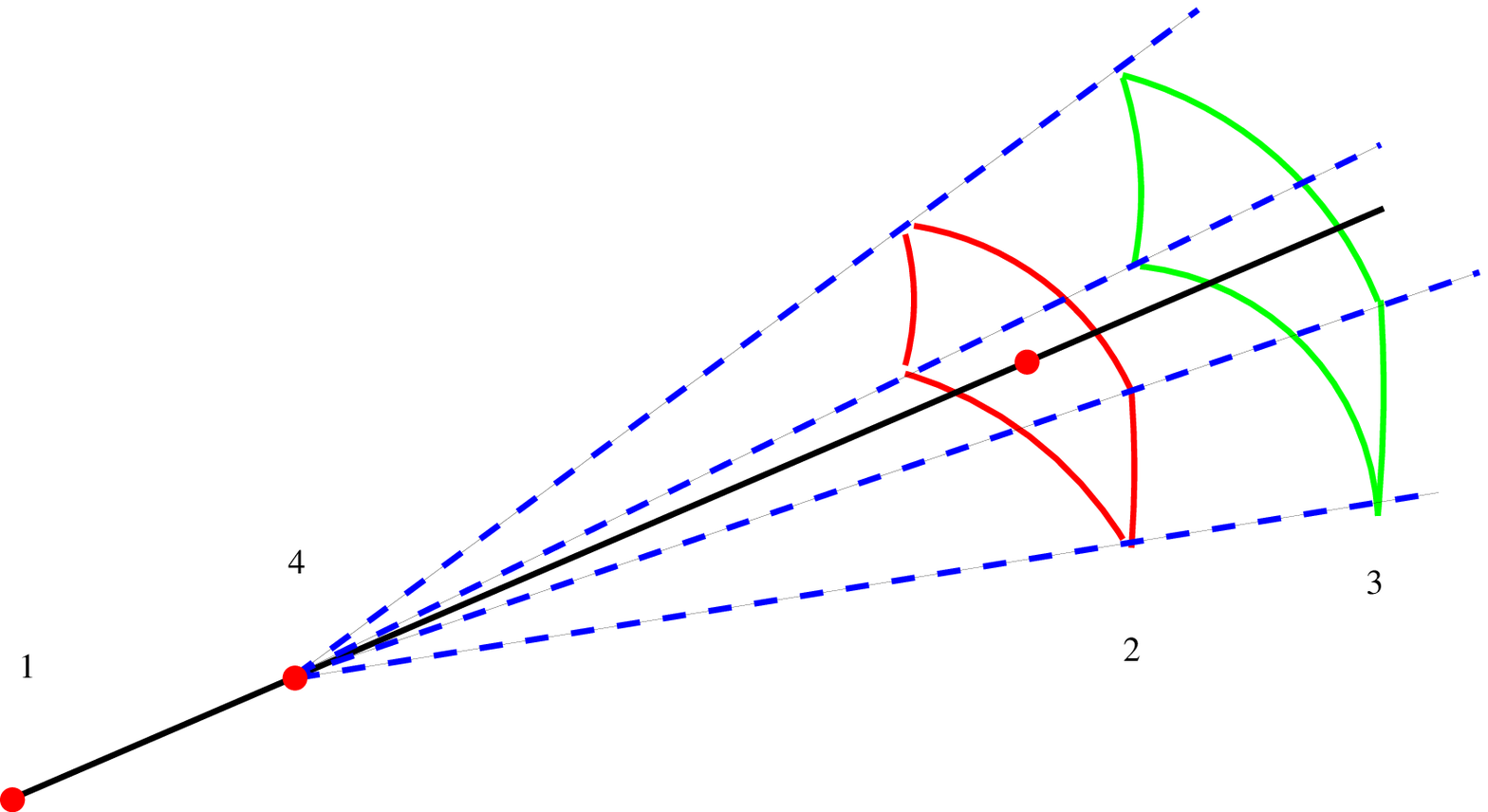}
\end{center}
{\it FIGURE 1. 
{\bf Left:} This is a schematic figure in  $\R^{1+1}$.
  The black vertical line is the freely falling
observer $\hat \mu([-1,1])$.
The rounded black square
is $\pi(\mathcal U_{z_0,\eta_0})$ that is is a neighborhood of $z_0$,
and the red curve  passing through $z\in\pi(\mathcal U_{z_0,\eta_0})$ is the
time-like geodesic $\mu_{z,\eta}([-1,1])$. 
The boundary of the domain $U_{\hat g}$ where we observe waves is shown on blue.
The green area is the set $W_{\hat g}\subset U_{\hat g}$
where sources are supported, and the black ``diamond'' is the set
$J(p^-,p^+)=J_{\hat g}^+({p^-})\cap J^-_{\hat g}({p^+})$.

{\bf Right:} 
This is a schematic figure in the space $\R^{3}$.
It describes the location of a {distorted plane wave (or a piece of a spherical wave)} $\dot u$ at different 
time moments. This wave
propagates near the geodesic $\gamma_{x_0,\zeta_0}((0,\infty))\subset \R^{1+3}$,
$x_0=(y_0,t_0)$
 and is singular on
 a subset of a light cone emanated from  $x^\prime=(y^\prime,t^\prime)$. 
  The piece of the distorted plane wave is sent from the surface $\Sigma\subset \R^3$,  it starts to propagate, and at a later time its singular support is the surface $\Sigma_1$.%
}

\subsubsection{Notations}\label{sec:notations 1}



Let $(M,g)$ be a $C^\infty$-smooth $(1+3)$-dim\-ens\-ion\-al 
time-orientable Lor\-entz\-ian manifold.
For $x,y\in M$ we say that $x$ is in the chronological past of $y$ and denote $x\ll y$ if $x\not =y$ and there is a time-like path from $x$ to $y$.
If $x\not=y$ and there is a causal path from $x$ to $y$,
we say that $x$ is  in the causal past of $y$ and  denote $x<y$.
If 
$x< y$ or $x=y$ we denote $x\leq y$.
The chronological future $I^+(p)$ of $p\in M$
consist of all points $x\in M$ such that $p\ll x$,
and the causal future $J^+(p)$ of $p$
consist of all points $x\in M$ such that $p\leq x$.
One defines similarly  the chronological past $I^-(p)$ of $p$
 and
the causal past $J^-(p)$ of $p$.
For a set $A$ we denote $J^\pm(A)=\cup_{p\in A}J^\pm(p)$.
We also denote $J(p,q):=J^+(p)\cap J^-(q)$ and $I(p,q):=I^+(p)\cap I^-(q)$.
If we need to emphasize the metric $g$ which is used to define the causality, 
we denote  $J^\pm(p)$ by $J^\pm_g(p)$ etc.

Let $\gamma_{x,\xi}(t)=\gamma^g_{x,\xi}(t)=\exp_x(t\xi)$ denote a geodesics in $(M,g)$. 
The projection from the tangent bundle $TM$ to the base point of a vector is denoted by
$\pi:TM\to M$. 
Let $L_xM$ denote the light-like directions of $T_xM$, 
and $L_x^+M$ and $L_x^-M$  denote the future and past pointing light-like vectors,
respectively. 
We also denote $\L^+_{g}(x)=\exp_x(L^+_xM)\cup\{x\}$ the union of the image of the future light-cone
in the exponential map of $(M,g)$ and the point $x$.

By  \cite{Bernal}, an open time-orientable Lorenzian manifold $(M,g)$ is globally hyperbolic
if and only if there are no closed causal paths in $M$ and
for all $q^-,q^+\in M$ such that $q^-<q^+$ the set 
$J(q^-,q^+)\subset M$ is compact.
We assume throughout the paper that $(M,g)$ is globally hyperbolic. 


When $g$ is a Lorentzian metric, having eigenvalues $\lambda_j(x)$ and eigenvectors 
$v_j(x)$ in some local coordinates, we will use also the corresponding
Riemannian metric, denoted by $g^+$ which has the  
eigenvalues $|\lambda_j(x)|$ and the eigenvectors 
$v_j(x)$ in the same local coordinates.
Let $B_{g^+}(x,r)=\{y\in M;\ d_{g^+}(x,y)<r\}$.

\subsubsection{Perturbations of a global hyperbolic metric} 
\label{subsubsec: Gloabal hyperbolicity}
Let  $(\hattuM ,\hat g)$ be a $C^\infty$-smooth globally hyperbolic Lorentzian
manifold. We will call $\hat g$ the  background metric 
 on $M$ and consider
its small perturbations. 
A Lorentzian 
metric $g_1$ dominates the metric $g_2$, if all vectors $\xi$
that are light-like or time-like with respect to the metric $g_2$ are
time-like with respect to the metric $g_1$, and in this case we denote $g_2<g_1$.
As  $(\hattuM ,\hat g)$ is globally hyperbolic, it follows from \cite{Geroch} that there is a Lorentzian metric $\tilde g$ such that  $(\hattuM ,\tilde g)$
is   globally hyperbolic and $\hat g<\tilde g$. One can
assume 
that the metric $\tilde g$ is smooth. We use the positive definite Riemannian metric $\hat g^+$ to define
norms in  the spaces $C^k_b(\hattuM )$ of functions with bounded $k$ derivatives 
and the  Sobolev spaces $H^s(\hattuM )$.

By \cite{Bernal2}, the globally hyperbolic manifold $(\hattuM ,\tilde g)$ has
an isometry $\Phi$  
 to the smooth product manifold $(\R\times N,\tilde  h)$,
where $N$
is a 3-dimensional manifold and the metric $\tilde  h$ can be written as
$\tilde  h=-\beta(t,y) dt ^2+\kappa(t,y)$ where $\beta:\R\times  N\to (0,\infty)$ is a smooth function and 
$\kappa(t,\cdotp)$ is a Riemannian metric on $ N$ depending smoothly
on $t\in \R$, and the submanifolds $\{t^{\prime}\}\times N$ are $C^\infty$-smooth
Cauchy surfaces for all $t^{\prime}\in \R$. We define the smooth time function
${\bf t}:\hattuM \to \R$ by setting
 ${\bf t}(x)=t$ if $\Phi(x)\in \{t\}\times N$.
Let us next identify these isometric manifolds,
that is, we denote $\hattuM =\R\times N$.
%

For $t\in \R$, let $\hattuM (t)=(-\infty,t)\times N$ and, for a fixed $t_0>0$ and $t_1>t_0$,
let $\hattuM _j=\hattuM (t_j)$, $j=0,1$.
Let $r_0>0$ be sufficiently small and $\V(r_0)$ be the set of metrics $g$ on $\hattuM _1= (-\infty,t_1)\times N$,  which  $C^8_b(\hattuM _1)$-distance
to 
 $\hat g$ is less that $r_0$ and 
coincide with $\hat g$ in $M(0)= (-\infty,0)\times N$. 


\subsubsection{Observation domain $U$}
For $g\in \V(r_0)$, let $\mu_g:[-1,1]\to M_1$ be a freely falling observer, that is, a time-like geodesic on $(M,g)$. Let $-1<s_{-3}<s_{-2}<s_{-1}<s_{+1}<s_{+2}<s_{+3}<1$ be such that $p^-=\mu_g(s_{-1})
\in \{0\}\times N$ and let $p^+=\mu_g(s_{+1})$. Below, we  denote $s_\pm=s_{\pm 1}$ and $\hat \mu=\mu_{\hat g}$.

When  $z_0=\hat \mu(s_{-2})\in M(0)$ and $\eta_0=\p_s \hat \mu(s_{-2})$,  we denote by  $\U_{z_0,\eta_0}(h)$  
the open $h$-neighborhood
of $(z_0,\eta_0)$ in the Sasaki metric of $(TM,\hat g^+)$.  
We use below a small parameter $\hat h>0$.
For $(z,\eta)\in \U_{z_0,\eta_0}(2\hat h)$
we define on $(M,g)$  a freely falling observer 
$\mu_{g,z,\eta}:[-1,1]\to M$,
such that   $\mu_{g,z,\eta}(s_{-2})=z$, and $\p_s\mu_{g,z,\eta}(s_{-2})=\eta$.
We assume that $\hat h$ is so small that $\pi(\U_{z_0,\eta_0}(2\hat h))\subset 
M(0)$
and for all $g\in \V(r_0)$ and $(z,\eta)\in \U_{z_0,\eta_0}(2\hat h)$ the geodesic $\mu_{g,z,\eta}([-1,1])\subset M$
is well defined and time-like and satisfies 
\beq\label{kaava D}
\mu_{g,z,\eta}(s_{-j-1})\in I_{g}^-(\mu_{g,z_0,\eta_0}(s_{-j})),\ \
\mu_{g,z,\eta}(s_{j+1})\in I_{g}^+(\mu_{g,z_0,\eta_0}(s_{j})),\hspace{-1cm}
\eeq
for $j=1,2$.
We denote, see Fig.\ 1(Left), 
$ \U_{z_0,\eta_0}= \U_{z_0,\eta_0}(\hat h)$ and 
\beq\label{eq: Def Wg with hat}
U_g\hspace{-1mm}=\hspace{-4mm}\bigcup_{(z,\eta)\in \U_{z_0,\eta_0}} \hspace{-4mm}\mu_{g,z,\eta}([-1,1]), \quad \hat U=U_{\hat g}.
\eeq 
\subsection{Formulation of the inverse problem}

\subsubsection{Inverse problems for non-linear wave equations}
The solution of inverse problems is often done by constructing the coefficients
of the equations using invariant methods, e.g.\ using travel time coordinates. Thus in the topical studies of the subject, inverse problems are  
 formulated invariantly, that is, on manifolds, see e.g.\ \cite {AKKLT,BelKur,dsFKS,Eskin, GsB, GST,LeU}.
Many physical models lead to non-linear differential equations.  In small perturbations, 
these equations can be approximated by linear equations,
and most of the previous results on hyperbolic inverse problems in the multi-dimensional case  concern 
linear models. Moreover, the existing uniqueness results are limited to
the time-independent or real-analytic coefficients \cite{AKKLT,Bel1,BelKur,Eskin,KKL} as 
these results are based on Tataru's unique 
continuation principle \cite{Tataru1,Tataru2}. Such unique continuation results 
have   been shown to fail for general metric
tensors which are not analytic in the time variable \cite{Alinhac}.
Even some linear inverse problem are not uniquely
solvable.  In fact, the counterexamples for these problems have
been used in the so-called transformation optics. This has led to  
models for  fixed frequency invisibility cloaks, see e.g.\ 
\cite{GKLU1} and references therein.
These applications give one more motivation 
to study inverse problems.

Earlier studies on inverse problems  for  non-linear  equations  have  concerned  
  parabolic  equations
 \cite {Isakov1}, elliptic equations \cite {IsakovNachman,Kang,Sun,SunU},
and 1-dimensional hyperbolic equations \cite {Nakamura1}. 
The present paper differs from the earlier studies in   
 that in our approach we do not consider the non-linearity as a perturbation, which effect
is small with special solutions, 
 but as a tool that helps us to solve the inverse problem. Indeed, the non-linearity
 makes it possible to solve a non-linear inverse problem which linearized version is not
 yet solved. This is the key novel feature of the paper.

\subsubsection{Einstein equations}
\MTEXT{Below, we use the Einstein summation convention. The roman indexes
$i,j,k$ etc.\ run usually over indexes of spacetime variables as the greek letters are 
reserved to other indexes in sums.}
The Einstein tensor of a Lorentzian metric $g=g_{jk}(x)$
 is
\ba
\Ein_{jk}(g)=\Ric_{jk}(g)-\frac 12 (g^{pq}\,\Ric_{pq} (g))
g_{jk}.
\ea
Here, $\Ric_{pq} (g)$ is the Ricci curvature of the metric $g$.
We define the divergence of a 2-covariant tensor $T_{jk}$ to be 
$(\div_g T)_k=\nabla_n(g^{nj}T_{jk})$.

Let us consider the
Einstein equations  in the presence of matter,
\beq\label{Einmat1} 
& &\Ein_{jk}(g)=T_{jk},\\
& &\div_{g}T=0,\label{Einmat2}
\eeq
for a Lorentzian metric $g$
and a stress-energy tensor $T$
related to the distribution of mass and energy.
We recall that by Bianchi's identity $\div_{g}(\Ein(g))=0$ and thus
the equation (\ref{Einmat2}),  called the conservation law for the stress-energy tensor, follows automatically from (\ref{Einmat1}). 
 \subsubsection{Reduced Einstein tensor}


%
%
%

Let $m\geq 5$, $t_1>t_0>0$ and  $g^\prime\in \V(r_0)$ be a $C^m$-smooth metric
that satisfy the Einstein equations  $\Ein(g^\prime)=T^\prime$
on $\hattuM (t_1)$.  When $r_0$ above is small enough, there is
a diffeomorphism $f:\hattuM (t_1)\to f(\hattuM (t_1))\subset \hattuM $ that is  a 
$(g^{\prime},\hat g)$-wave map $f:(\hattuM (t_1),g^\prime)\to (\hattuM,\hat g)$
and satisfies $\hattuM (t_0)\subset f(\hattuM (t_1))$.
Here, $f:(\hattuM (t_1),g^\prime)\to ( \hattuM,\hat g)$ is a wave map,
 see \cite[Sec.\ VI.7.2 and App.\ III, Thm. 4.2] {ChBook}, if
 \beq
\label{C-problem 1 pre}& &\square_{g^{\prime},\hat g} f=0\quad\hbox{in } \hattuM (t_1),\\
\label{C-problem 2 pre}& &f=Id,\quad \hbox{in  } (-\infty,0)\times N,
\eeq
where
$$\square_{g^{\prime},\hat g} f
= 
(g^{\prime})^{  jk}(\frac \p{\p x^j}\frac \p{\p x^k} f^A -\Gamma^{{\prime} n}_{jk}(x)\frac \p{\p x^n}f^A+
\hat \Gamma^A_{BC}(f(x))
\frac \p{\p x^j}f^B\,\frac \p{\p x^k}f^C),$$
and $\hat \Gamma^A_{BC}$ denotes the Christoffel symbols of metric $\hat g$
and  $\Gamma^{ {\prime} j}_{kl}$ are the Christoffel symbols of metric $g^{\prime}$,
%
%
see \cite[formula (VI.7.32)]{ChBook}.
The wave map has the property that $g= f_*g^{\prime}$  satisfies
  $\Ein(g)= \Ein_{\hat g}(g)$,
  where 
$\Ein_{\hat g}(g)$ is the $\hat g$-reduced Einstein tensor, 
see  also formula (\ref{Reduced Einstein tensor}) below,
\ba
(\Ein_{\hat g} g)_{pq}=-\frac 12 g^{jk}\hat \nabla_j\hat \nabla_k  g_{pq}
+\frac 14( g^{nm}g^{jk}\hat \nabla_j\hat \nabla_k  g_{nm})g_{pq}
+P_{pq}(g,\hat \nabla g),
\ea
where
$\hat \nabla_j$ is the covariant differentation with respect to the metric $\hat g$
and $P_{pq}$ is a polynomial function of 
$ g_{nm}$, $ g^{nm}$, and $ \hat \nabla_jg_{nm}$ with coefficients
depending on the metric $\hat g_{nm}$ and its derivatives. Considering
the  wave map $f$ as a transformation of coordinates, we see that $g=f_*g^{\prime}$ and $T=f_*T^{\prime}$  satisfy the $\hat g$-reduced Einstein equations
\beq\label{eq: EE in correct coordinates}
\Ein_{\hat g}(g)= T\quad\hbox{on }\hattuM (t_0).
\eeq
In the literature, the above is often
stated by saying that 
the reduced
Einstein equations  (\ref{eq: EE in correct coordinates})  is
the Einstein equations   written
with the wave-gauge corresponding to the metric $\hat g$. The equation (\ref{eq: EE in correct coordinates}) is a quasi-linear hyperbolic system of equations for $g_{jk}$.
We emphasize that
a solution of the reduced Einstein equations  can be a solution
of the original Einstein equations  only if the stress energy
tensor satisfies the conservation law $\nabla^{g}_jT^{jk}=0$.
It is usual also to assume that  the energy density is non-negative. For instance,
the weak energy condition requires that $T_{jk}X^jX^k\geq 0$ for all time-like
vectors $X$.  Next, we couple the Einstein equations  with matter fields  and
 formulate the direct problem for  the $\hat  g$-reduced Einstein equations.

\subsubsection{The initial value problem with sources} 

We consider metric and physical fields on a Lorentzian manifold
$(M,g)$. This is an informal discussion.
We aim to study an inverse problem with active measurements.
As measurements cannot be implemented in Vacuum (as
the Einstein equations  is uniquely solvable with fixed initial data),
we have to add matter fields in the model. We 
consider the coupled system of the Einstein equations
and the equations for $L$ scalar fields  $ \phi
=(\phi_\ell)_{\ell=1}^L$ with some sources $\F^1$ and $\F^2$.

Let $\hat g$ and $\hat \phi
=(\hat \phi_\ell)_{\ell=1}^L$  be $C^\infty$-background fields  on $M$.
Consider  
\beq\label{eq: adaptive model with no source}
& &\Ein_{\hat g}(g) =T,\quad T_{jk}={\bf T}_{jk}(g,\phi)+\F^1_{jk},\quad
\hbox{in }M_0=(-\infty,t_0)\times N,
\\ \nonumber
& &{\bf T}_{jk}(g,\phi)=\bigg(\sum_{\ell=1}^L(\p_j\phi_\ell \,\p_k\phi_\ell 
-\frac 12 g_{jk}g^{pq}\p_p\phi_\ell \,\p_q\phi_\ell\bigg)-
V(\phi)g_{jk},
\\ \nonumber
& &\square_g\phi_\ell - V^\prime_\ell(\phi)=\F^2_\ell,\quad \ell=1,2,3,\dots,L,\\
& &\nonumber g=\hat g\hbox{ and $\phi_\ell=\hat \phi_\ell$ in $M_0\setminus
J^+_g(p^-)$,}
\eeq
where $\F^1$ and $\F^2$ are supported in $U_g^+\cap J^+_g(p^-)$,
\MTEXT{$V\in C^\infty(\R^L)$, and $V^\prime_\ell(s)=\frac {\p }{\p s_\ell}V(s)$,
$s=(s_1,s_2,\dots,s_L)$. A typical model is 
$V(s)= \sum_{\ell=1}^L\frac 12 m^2s_\ell^2$. Above,
$\square_g\phi = |\det(g)|^{-\frac 12}\p_p 
(|\det(g)|^{\frac12}
g^{pq}\p_q\phi)$.
} 
We assume that  the background fields $\hat g$ and  $\hat \phi$
satisfy the equations (\ref{eq: adaptive model with no source}) with $\F^1=0$ and 
$\F^2=0$. Note that above $J^+_{g}(p^-)\cap M_0\subset 
J^+_{\tilde g}(p^-)$ when $g\in \V(r_0)$ and $r_0$  is small enough.
To obtain a physically meaningful model,
we need to assume that the physical conservation law in relativity,
\beq\label{conservation law0}
\nabla_p(g^{pk} T_{kj})=0,\quad\hbox
{for $j=1,2,3,4$, where }  T_{kj}={\bf T}_{kj}(g,\phi)+\F^1_{kj}\hspace{-1cm}
\eeq
 is satisfied. Here $\nabla=\nabla^g$ is the connection corresponding to $g$. 
As will be noted in Subsection
 \ref{sssec: relation},  the reduced Einstein tensor $\Ein_{\hat g}(g)$
is equal to the Einstein tensor $\Ein(g)$ when
  $(g,\phi)$ satisfies the system (\ref{eq: adaptive model with no source})
 and the conservation law (\ref{conservation law0}).
We mainly need local existence results\footnote{In this paper we do not use optimal smoothness for the solutions
in classical  $C^k$ spaces or Sobolev space $W^{k,p}$ but just suitable smoothness for which
the non-linear wave equations can be easily analyzed using $L^2$-based Sobolev spaces.} for the system (\ref{eq: adaptive model with no source}). The global existence problem for the related systems 
 has recently attracted much interest in the mathematical
community and many important results been obtained, see e.g.\
\cite{Ch-K,Dafermos2,K1,K3,Li1,Li2}.

We  encounter above the difficulty that
the source $\F=(\F^1,\F^2)$ in (\ref{eq: adaptive model with no source}) has to satisfy the condition (\ref{conservation law0})
that depends on the solution $g$ of  (\ref{eq: adaptive model with no source}).
This makes the formulation of active measurements in relativity difficult.
In Appendix \noextension{A}\extension{C}
 we consider a model where
the source term $\F^1$ corresponds to e.g.\ fluid fields consisting of particles
whose 4-velocity vectors are controlled and  $\F^2$ contains a term corresponding 
to a secondary source function that adapts
the changes of $g,\phi,$ and $\F^1$ so that the physical conservation law (\ref{conservation law0})
is satisfied. This model is considered in detail in 
\cite{Paper-inpreparation}. However, in this
 paper we replace the adaptive source functions by a general assumption of microlocal linearization stability that does not fix the physical model for the source fields $\F$.


%

{

\subsubsection{Definition of measurements}

For $r>0$ let, see Fig.\ 1(Left),
\beq\label{observer neighborhood with hat}
W_g(r)\hspace{-1mm}=\hspace{-4mm}\bigcup_{s_-<s<s_+-r} \hspace{-4mm}I_{g}( \mu_g(s),\mu_g(s+r)),
\eeq
and let $r_1>0$ be so small that  $W_{g}(2r_1)\subset U_{g}$
for all $g\in \V(r_0)$. 
%
We denote  $W_{g}=W_{g}(r_1)$. 
We use  Fermi-type coordinates:
Let $Z_j(s)$, $j=1,2,3,4$
be a parallel frame of linearly independent time-like vectors at $\mu_g(s)$ such
that $Z_1(s)=\dot\mu_g(s)$. Let $\Phi_g:(t_j)_{j=1}^4\mapsto
\exp_{\mu_g(t_1)}(\sum_{j=2}^4 t_jZ_j(t_1))$.  
We assume that $r_1>0$ is so small that 
$\Psi_{g}=\Phi_{g}^{-1}$ defines coordinates in $W_{g}(2r_1)$.
We define 
the norm-like functions  
\ba
\mathcal N_{\hat g}^{(k)}(g)\hspace{-1mm}=\hspace{-1mm}
\|(\Psi_{g})_*g\hspace{-1mm}-\hspace{-1mm}(\Psi_{\hat g})_*\hat g\|_{C^{k}_b(\overline {\Psi_{\hat g}(W_{\hat g})})},
\ \ \mathcal N^{(k)}(\F)\hspace{-1mm}=\hspace{-1mm}
\|(\Psi_{g})_*\F\|_{C^{k}_b(\overline {\Psi_{\hat g}(W_{\hat g})})},
\ea
where $k\in \N$,
that measures the $C^{k}$ distance of $g$ from $\hat g$ and $\F$ from zero in
the Fermi-type coordinates. As  we have assumed 
that the background metric $\hat g$ and the field $\hat \phi$
are $C^\infty$-smooth, we can  consider as smooth sources
as we wish. Thus we use below  smoothness assumptions on the sources 
that are far from the optimal ones.
\medskip

Let us  define 
the \MTEXT{source-observation 4-tuples} corresponding to me\-as\-ur\-em\-ents in  $U_{g}$ with sources $\F$ 
 supported in $W_{g}$. 
Let  $\e>0$ and 
$k_0\geq 8$ and define 
\beq\nonumber
{\cal D}(\hat g,\hat \phi,\e)=\{[(U_g,
g|_{U_g},\phi|_{U_g},\F|_{U_g})]&;&(g,\phi,\F)\hbox{ are $C^{k_0+3}$-smooth 
}\\ 
& &  \label{eq: main data}
\hspace{-4cm}\hbox{solutions of  
(\ref{eq: adaptive model with no source}) and (\ref{conservation law0})  with }\F\in C^{k_0+3}_0(W_{g};\B^L),
\\  \nonumber
& &\hspace{-6.5cm}\hbox{$J_g^+(\supp(\F))\cap J_g^-(\supp(\F))\subset 
W_g, $ $\mathcal N^{({ {k_0}})}(\F)<\e$, $\mathcal N^{({{k_0}})}_{\hat g}(g)<\e$}\}.
\eeq
Above, the sources $\F$ above are considered as sections of the bundle $\B^L$,
where $\B^L$ is a vector bundle on $\hattuM $ that is the product bundle of
the bundle of symmetric $(0,2)$-tensors and the trivial vector bundle with the fiber $\R^L$.
Also, $[(U_g,
g,\phi,\F)]$ denotes the equivalence class of
all
Lorenztian manifolds $(U^{\prime},g^{\prime})$ 
and  functions $\phi^{\prime}=(\phi^{\prime}_\ell)_{\ell=1}^L$ and the  tensors $\F^{\prime}$ defined on 
a $C^\infty$-smooth manifold $U^{\prime}$,
such that there is $C^\infty$-smooth 
diffeomorphism $\Psi:U^{\prime}\to U_g$ satisfying 
$\Psi_*g^{\prime}=g$, $\Psi_*\phi_\ell^{\prime}=\phi_\ell$,  and $\Psi_*\F^\prime=\F$.

\MTEXT{In many inverse problems one considers a Dirichlet-to-Neumann map or, equivalently to that,
the Cauchy data set that is the graph of the Dirichlet-to-Neumann map.
Similarly, the source-observation 4-tuples ${\cal D}(\hat g,\hat \phi,\e)$ could be considered
as  graph of a ``source-to-field'' map but due to the conservation law the 
source-to-field map could be defined only on a subset of sources supported in $U_g$.
To avoid the difficulties related to the fact that we do not have a good characterization
for this subset, nor do we know the wave map coordinates in $U_{g}$, we do not define
 a ``source-to-field'' map but use the data set
  ${\cal D}(\hat g,\hat \phi,\e)$.}

\extension{We will analyze the smoothness objects in 
${\cal D}(\hat g,\hat \phi,\e)$ in different coordinates.
Observe that when $g^{\prime}\in C^{k_0+3}$
and $\F^{\prime}\in C^{k_0+3}$, the  Fermi coordinates $\Psi_{g^{\prime}}$ are
$ C^{k_0+1}$-smooth and thus $(\Psi_g^{\prime})_*g^{\prime}$ is 
$ C^{k_0+1}$-smooth. Moreover, the $(g^{\prime},\hat g)$-wave map $f$
is $C^{k_0}$-smooth and thus the metric $g^{\prime}$ 
and the source $\F^{\prime}$ in the wave map coordinates,
that is, $f_*g^{\prime}$ and $f_*\F^{\prime}$, are 
$C^{k_0-1}$-smooth.  However, to consider local existence
results for the Einstein-scalar field equations, we need to consider 
the case when the norm of
the source is small. We do this next.

Observe that when $\mathcal N^{({ {k_0}})}(\F)<\e$, $\mathcal N^{({{k_0}})}_{\hat g}(g)<\e$
we can locally solve 
the $(g,\hat g)$-wave map $\tilde f$
is $C^{k_0-3}$-smooth and thus the metric $g$ 
and the source $\F$ in the wave map coordinates,
that is, $\tilde f_*g$ and $\tilde f_*\F$, are 
$C^{k_0-4}$-smooth in $W_{\tilde f_*g}$. Since $\F$ vanishes outside the domain
$W_g$, we see that $\tilde f_*\F$ vanishes outside $W_{\tilde f_*g}$
and thus  $\tilde f_*\F$ is $C^{k_0-4}$-smooth  in the set $M_0$.
 As  $k_0\geq 8$, this implies that we can later obtain local existence
of the Einstein-scalar field equations (\ref{eq: adaptive model with no source})
when $\mathcal N^{({ {k_0}})}(\F)<\e$, $\mathcal N^{({{k_0}})}_{\hat g}(g)<\e$
and $\e$ is small enough.
}

\MTEXT{Note that $[(U_{\hat g},\hat g,\hat \phi,0)]$ is the only element in  ${\cal D}(\hat g,\hat \phi,\e)$
for which the $\F$-component is zero. Thus  
the collection $ {\cal D}(\hat g,\hat \phi,\e)$ determines the isometry type of $(U_{\hat g},\hat g)$.}

\subsubsection{Linearized equations}
{We need also to consider the linearized version of the equations 
(\ref{eq: adaptive model with no source}) that  have the form (in  local coordinates)
\beq\label{linearized eq: adaptive model with no source}
& &\square_{\hat g} \dot g_{jk}+A_{jk}(\dot g,\dot \phi,\p \dot g,\p \dot \phi)=f^1,\quad
\hbox{in }M_0,\\
\nonumber
& &\square_{\hat g}\dot \phi_\ell +
B_\ell(\dot g,\dot \phi,\p \dot g,\p \dot \phi)
=f^2,\quad \ell=1,2,3,\dots,L,
\eeq
where $A_{jk}$ and $B_\ell$ are first order linear differential
operators which coefficients depend on $\hat g$ and $\hat \phi$.
When $g_\e$ and $\phi_\e$ are solutions of (\ref{eq: adaptive model with no source})
with source $\F_\e$ depending smoothly on $\e\in \R$ such that 
$(g_\e,\phi_\e,\F_\e)|_{\e=0}=(\hat g,\hat \phi,0)$,
 then $(\dot g,\dot \phi,f)=(\p_\e g_\e,\p_\e \phi_\e,\p_\e \F_\e)|_{\e=0}$ solve (\ref{linearized eq: adaptive model with no source}).
}

{Let us consider the concept of the {\it linearization stability (LS)} for the source problems, cf.\ 
\cite{Brill1,Ch-Deser,L1,Girbau}, and references therein:
Let
$s_0>4$ and  consider 
a $C^{s_0+4}$-smooth  source $f=(f^1,f^2)$ that is supported in $W_{\hat g}$ 
and satisfies
the linearized conservation law 
\beq\label{eq: lineariz. conse. law PRE}
& &\hspace{-1.5cm}
{
\hspace{-1.5cm}\frac 12 \hat g^{pk}\hat \nabla_p f^1_{kj}+
\sum_{\ell=1}^L f^2_\ell \, \p_j\hat \phi_\ell
=0,\quad j=1,2,3,4.\hspace{-1cm}
}
\eeq
Let  $(\dot g,\dot \phi)$  be the solution of the linearized
Einstein equations  (\ref{linearized eq: adaptive model with no source}) with source $f$.
We say that  {\it $f$ has the LS-property} in $C^{s_0}(M_0)$
if  there are $\e_0>0$ and a family $\mathcal F_\e=(\mathcal F^1_\e,\mathcal F^2_\e)$
of sources, supported in $W_{g_\e}$ for all $\e\in [0,\e_0)$,
and functions $(g_\e,\phi_\e)$ 
that all depend smoothly on $\e\in [0,\e_0)$ in $C^{s_0}(M_0)$
such that 
\beq\label{eq: LSp}
& &\hbox{$(g_\e,\phi_\e)$  satisfies the 
equations 
(\ref{eq: adaptive model with no source}) and the conservation law (\ref{conservation law0})},\\
& &\hbox{
$(g_\e,\phi_\e,\F_\e)|_{\e=0}=(\hat g,\hat \phi,0)$,
and
$(\dot g,\dot \phi,f)=(\p_\e g_\e,\p_\e \phi_\e,\p_\e \F_\e)|_{\e=0}$}.
\nonumber
\eeq
In this case, we say that 
$f=(f^1,f^2)$  has the LS-property with the family  $\F_\e$, $\e\in [0,\e_0)$.


\medskip

Note that above (\ref{eq: lineariz. conse. law PRE}) is obtained
by linearization of the  conservation law (\ref{conservation law0}).

The above linearization stability, concerning the local existence of the solutions in $M_0=(-\infty,t_0)\times N$,
is valid under quite general conditions, see \cite{Brill2}.
Below we will require that for some sources $f$, supported in a neighborhood $V$ of a point $y$,
we can find functions $\F_\e$ that are also supported in the set $V$.
The conditions when this happen are considered below\noextension{ and in Appendix A and \cite{Paper-inpreparation}}.

\motivation{
Next, we consider sources that are conormal distributions.
When $Y\subset U_{\hat g}$
is a 2-dimensional space-like submanifold,
consider local coordinates defined 
in  $V\subset \hattuM _0$ such that  
$Y\cap V\subset \{x\in \R^4;\ x^jb^1_j=0,\  x^jb^2_j=0\}$,
where $b^1_j,b^2_j\in \R$. Next we slightly abuse the notation by identifying $x\in V$ with its
coordinates $X(x)\in \R^4$. We denote  $f\in \I^{n}(Y)$, $n\in \R$, if in the
above local coordinates, $f$ can be written as
\beq\label{def eq: b b-prime}
f(x^1,x^2,x^3,x^4)=
\re \int_{\R^2}e^{i(\theta_1b^1_m+\theta_2b^2_m)x^m}\sigma_f(x,\theta_1,\theta_2)\,
d\theta_1d\theta_2,\hspace{-2cm}
\eeq
where 
 $\sigma_f(x,\theta)\in S_{0,1}^{n}(V;\R^2)$, $\theta=(\theta_1,\theta_2)$ 
is a classical symbol.  A function  $c(x,\theta)$ that is $n$-positive homogeneous in $\theta$,
i.e.,
$c(x,s\theta)=s^nc(x,\theta)$ for $s>0$, is the principal symbol of $f$ 
if there is $\phi\in C^\infty_0(\R^2)$ being 1
near zero such that
$\sigma_f(x,\theta)-(1-\phi(\theta))c(x,\theta)\in S_{0,1}^{n-1}(V;\R^2)$.
When $\eta=\theta_1b_1+\theta_2b_2\in N^*_xY$, we say that
$\tilde c(x,\eta)=c(x,\theta)$  is the value of the principal symbol of $f$ at $(x,\eta)\in N^*Y$.


 We need a condition that we call {\it microlocal linearization stability}: 
%
\medskip  

\noindent
\MTEXT{
    {\bf Assumption  $\mu$-LS (Microlocal linearization stability)}: 
\modified{\it Let $n_0\in \Z_-$, $Y\subset U_{\hat g}$
be a 2-dimensional space-like submanifold, $V\subset U_{\hat g}$ 
  an open local coordinate neighborhood of $y\in Y$ with 
the coordinates $X:V\to \R^4$, $X^j(x)=x^j$ such that 
 $X(Y\cap V)\subset \{x\in \R^4;\ x^jb^1_j=0,\  x^jb^2_j=0\}$.
 Let, in addition, $(y,\eta)\in N^*Y$ be a light-like covector,
 $\mathcal W\subset N^*Y$ be a conic neighborhood of $(y,\eta)$,
$(c_{jk})_{j,k=1}^4$ be a symmetric $(0,2)$ tensor at $y$
that satisfies
\beq\label{mu linearized conservation law for symbols}
& &\hat g^{lk}(y)\eta_l
c_{kj}
=0,\quad \hbox{for all } j=1,2,3,4,
\eeq
and $ (d_{\ell})_{\ell=1}^L\in \R^L$.
Then, for any $n\in \Z_-,$ $n\leq n_0$ there are $f^1_{jk}\in \I^{n}(Y)$, $(j,k)\in \{1,2,3,4\}^2$,
and $f^2 _\ell \in \I^{n}(Y)$, $\ell=1,2,\dots,L$,
 supported in $V$ with symbols that are in $S^{-\infty}$ outside the neighborhood $\mathcal W$ of
 $(y,\eta)$. The  principal symbols of $f^1_{jk}(x)$ and $f^2 _\ell(x)$ in the $X$-coordinates, denoted by $\tilde f^1_{jk}$ and
  $\tilde f^2_\ell$, respectively, are at $(y,\eta)$ equal to 
$\tilde f^1_{jk}(y,\eta)=c_{jk}$ and  $\tilde f^2_\ell(y,\eta)=d_{\ell}$. 
Moreover,  the source
 $f=(f^1,f^2)$  satisfies the linearized conservation law (\ref{eq: lineariz. conse. law PRE})
and  
 $f$  has \modified{the LS property (\ref{eq: LSp})
 in $C^{s_1}(M_0)$, $s_1\geq 13$, 
  with a family  $\F_\e$, $\e\in [0,\e_0)$
  such that  $\F_\e$ are supported in $V$.}}}
\subsection{Main results}

Our  main result is the following uniqueness theorem for the inverse problem for
the Einstein-scalar field equations:    
}

    \medskip

 \begin{theorem}\label{alternative main thm Einstein} Let  $(\hattuM^{(1)} ,\hat g^{(1)})$ and $(\hattuM^{(2)} ,\hat g^{(2)})$ be globally hyperbolic
manifolds and $(\hat g^{(j)},\hat \phi^{(j)})$, $j=1,2$ satisfy 
Einstein-scalar field equations (\ref{eq: adaptive model with no source}) with vanishing sources $\F^1=0$
and $\F^2=0$.
Also,
 assume that there are neighborhoods $U_{\hat g^{(j)}}$, $j=1,2$ of the time-like geodesics $\hat \mu_j\subset \hattuM^{(j)}$ where the  
 Assumption $\mu$-LS is valid.
Moreover,
 assume that,  for some $\e>0$, the sets ${\cal D}(\hat g^{(1)},\hat \phi^{(1)},\e)$ and
${\cal D}(\hat g^{(2)},\hat \phi^{(2)},\e)$ are the same, that is, the measurements done in
$U_{\hat g^{(1)}}$ and $U_{\hat g^{(2)}}$ coincide. Let  $p_j^-=\hat \mu_j(s_-)$ and $p_j^+=\hat \mu_j(s_+)$.
\motivation{Then there is a  diffeomorphism $\Psi:I_{\hat g^{(1)}}(p^-_1,p^+_1)\to I_{\hat g^{(2)}}(p^-_2,p^+_2)$ such that the metric $\Psi^*\hat g^{(2)}$ is conformal to
$\hat g^{(1)}$.}
 \end{theorem}

\medskip

 The  theorem above says  that the data ${\cal D}(\hat g,\hat \phi,\e)$ 
determine uniquely the manifold $I_{\hat g}({p^-},{p^+})\subset M$ and 
the conformal type of $\hat g$ in $I_{\hat g}({p^-},{p^+})$. 
Reconstruction of the  conformal structure of the manifold provides naturally less 
information than  
finding the whole metric structure, but the conformal structure is crucial for many questions of analysis and physics, see e.g. \cite {CB-I-P1,GrZ}.
Roughly speaking,
the above result means that if the manifold $(\hattuM _0,\hat g)$ is unknown, then
the source-to-observation pairs corresponding to
freely falling sources which are near a
freely falling observer $\mu_{\hat g}$ and the measurements of 
the metric tensor and the scalar fields in a neighborhood $U_{\hat g}$ of  $\mu_{\hat g}$,
 determine the  metric tensor up to conformal transformation
in the set $I_{\hat g}({p^-},{p^+})$.
\noextension{ 
 In this paper we present the proof of this result,
but mention for the convenience of the reader that extended versions of 
some technical computations discussed briefly  can be
found in the preprint \cite{preprint}.}



\begin{center}
\arxivpreprint{$  $\hspace{-10cm}
\includegraphics[height=4.0cm]{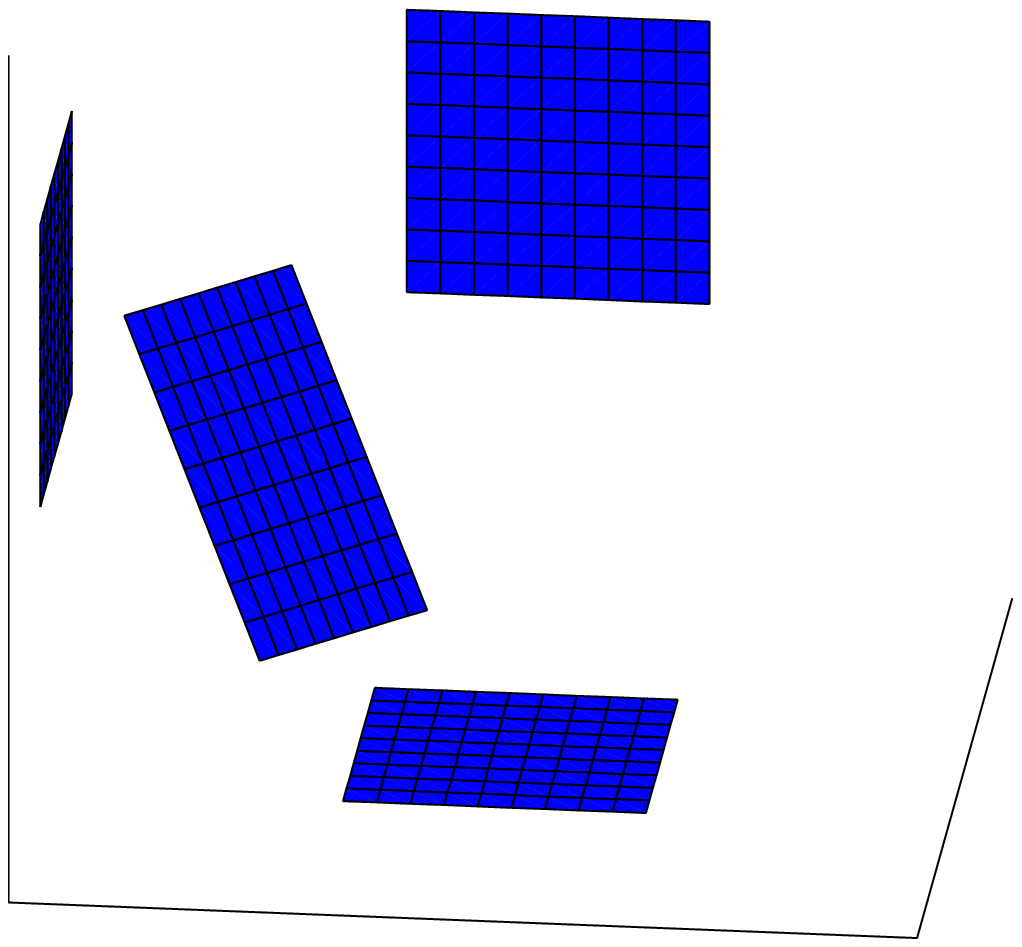}\quad 
\includegraphics[height=4.0cm]{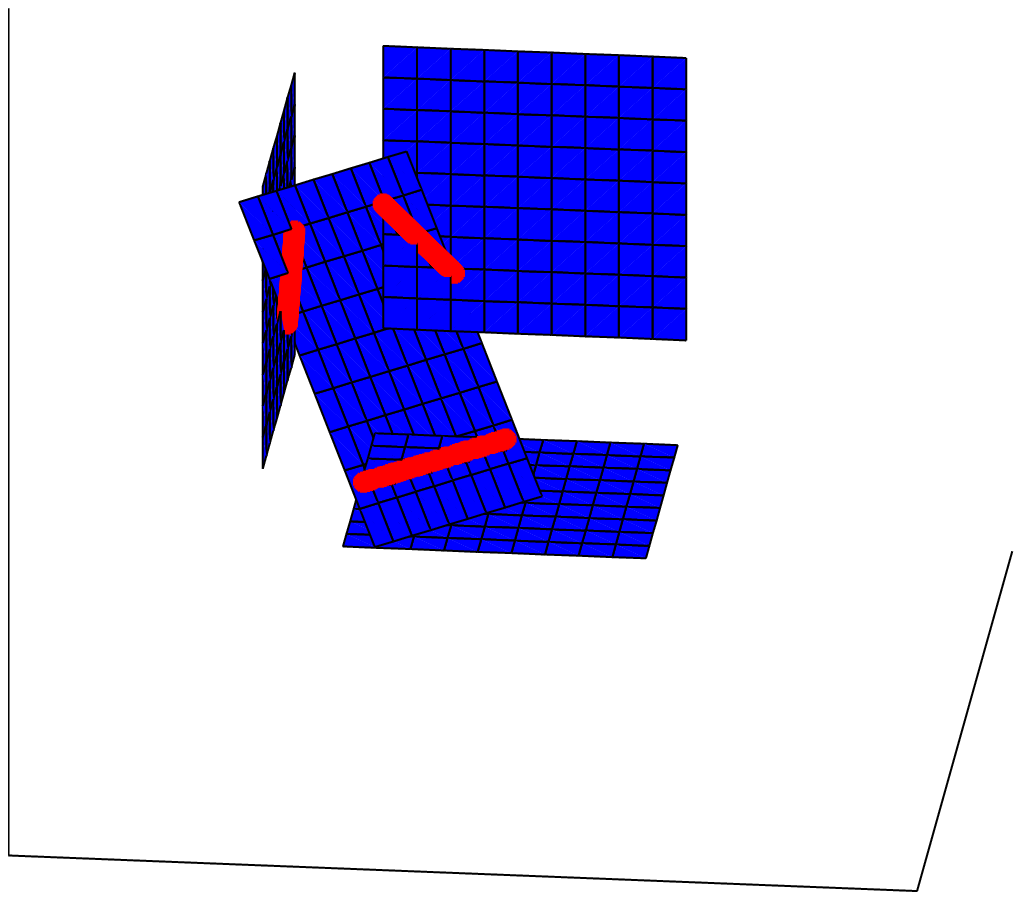}\quad
\includegraphics[height=4.0cm]{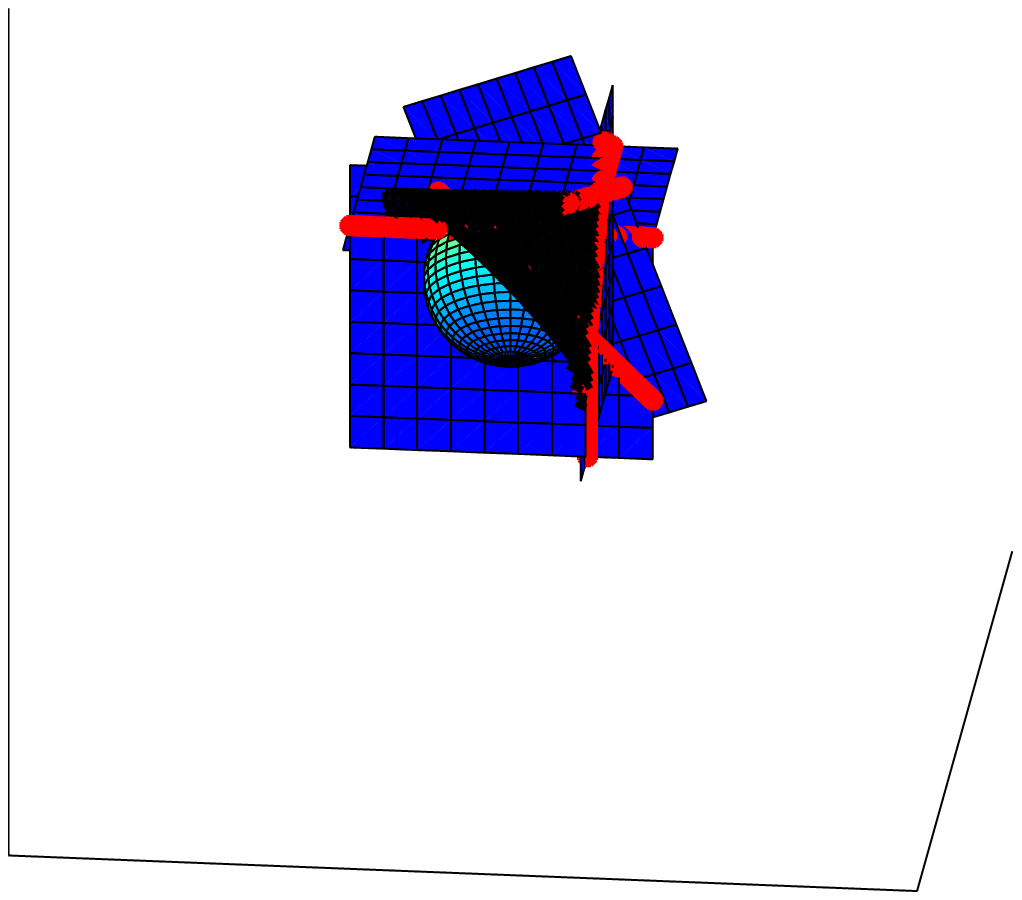}\quad
\includegraphics[height=4.0cm]{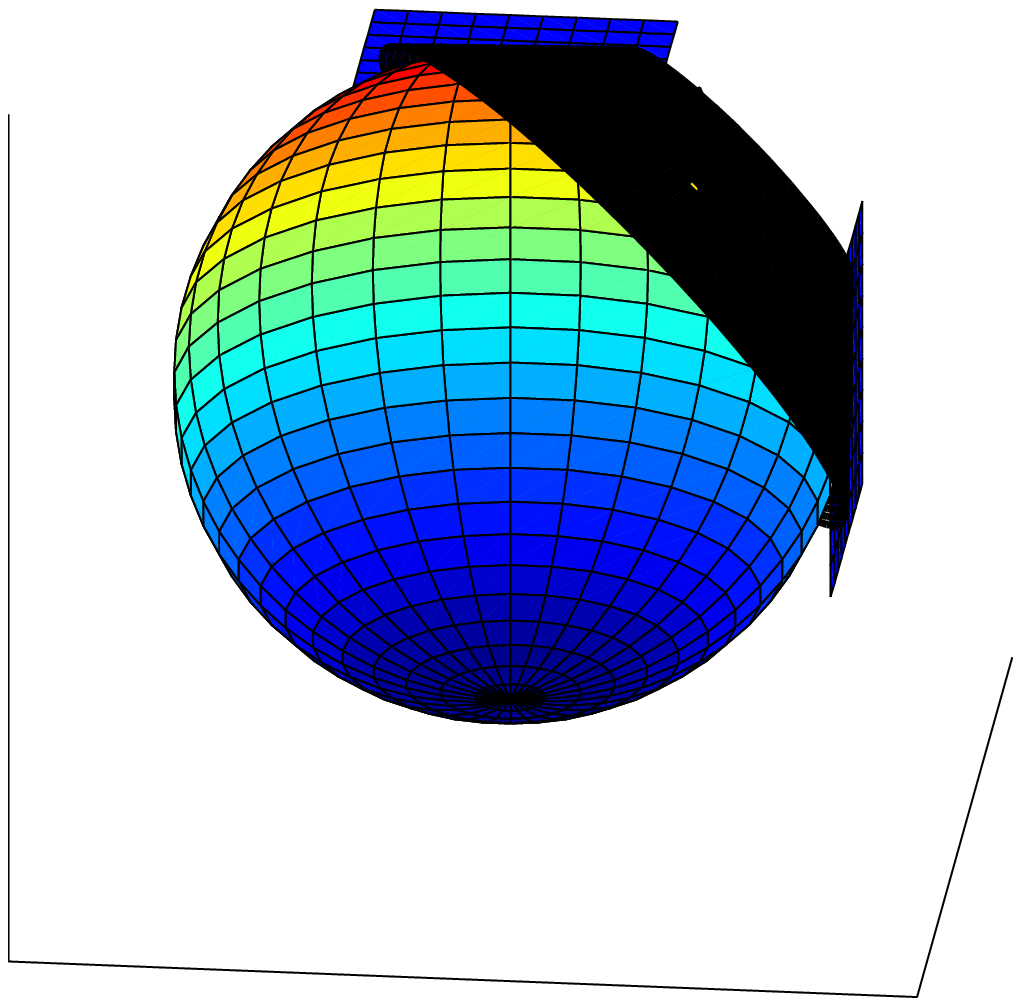}\hspace{-10cm}$  $}
\end{center}
{\it FIGURE 2. Four plane waves propagate in space. 
When the planes intersect, the non-linearity of the hyperbolic system  produces
new waves.
The four figures show
the waves before the interaction of the waves start, when 2-wave interactions have  started, when all  4 waves have just interacted, and later
after the interaction. {\bf Left:} Plane waves before  interacting.
 {\bf Middle left:}
 The 2-wave interactions (red line segments) appear  but do not cause
new propagating singularities. {\bf Middle right and Right:}  All plane waves have intersected and new 
 waves have appeared.The 3-wave interactions cause 
new conic waves (black surface). Only one such wave is shown in the figure. The 4-wave interaction causes
a point source in spacetime that sends a spherical wave to all future light-like 
directions. This spherical wave is essential in our considerations. For an animation
on these interactions, see the supplementary video.} 
\medskip

\extension{
\begin{center}
$  $\hspace{-10cm}
\includegraphics[height=4.0cm]{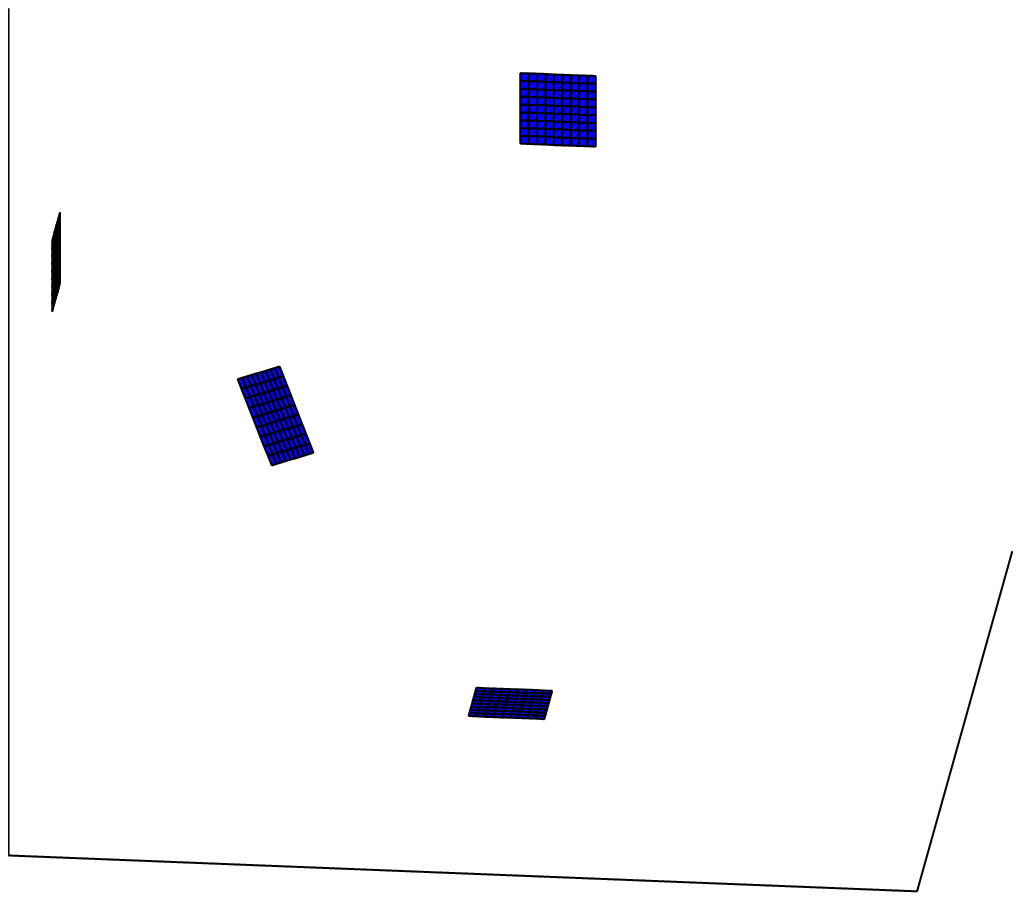}\quad 
\includegraphics[height=4.0cm]{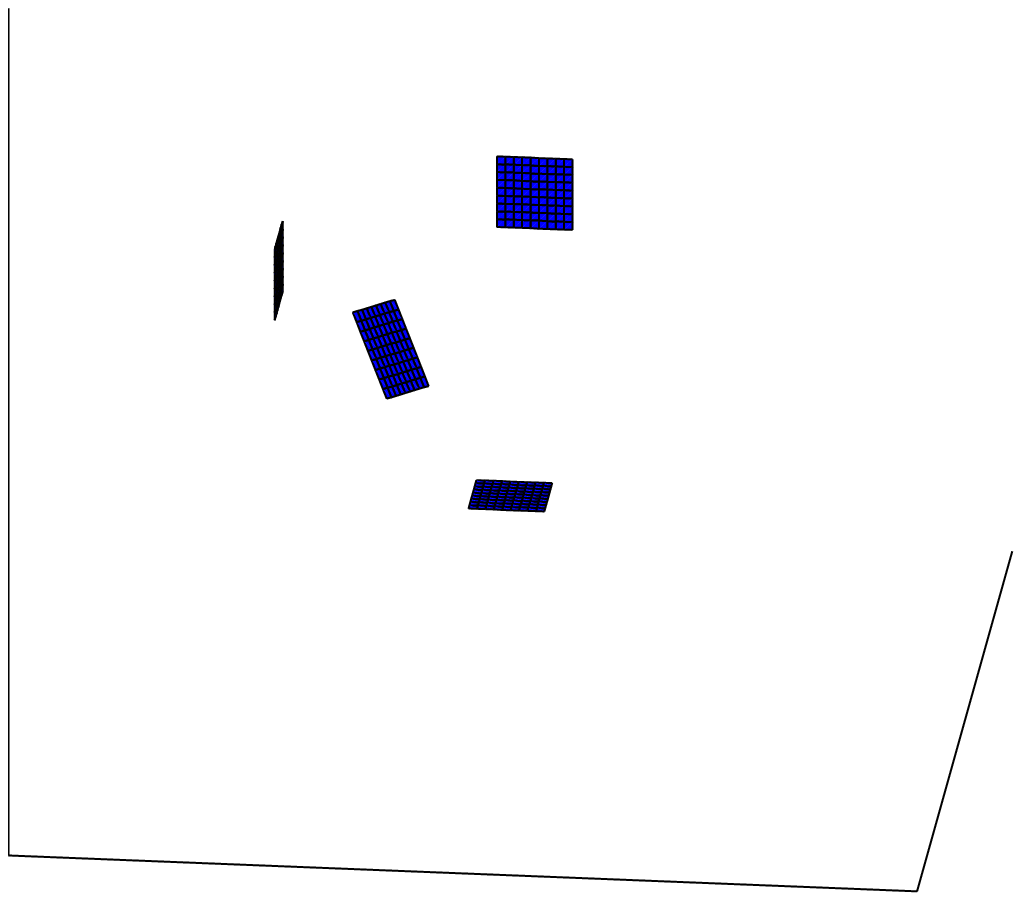}\quad
\includegraphics[height=4.0cm]{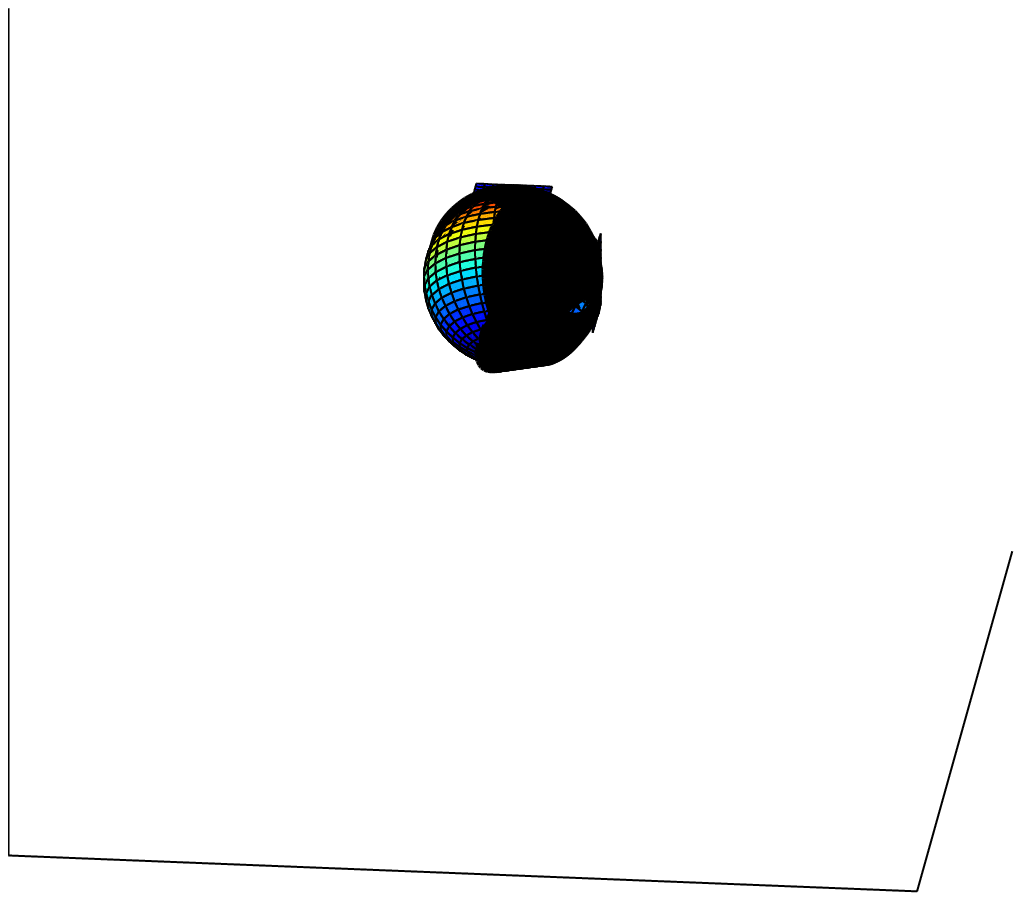}\quad
\includegraphics[height=4.0cm]{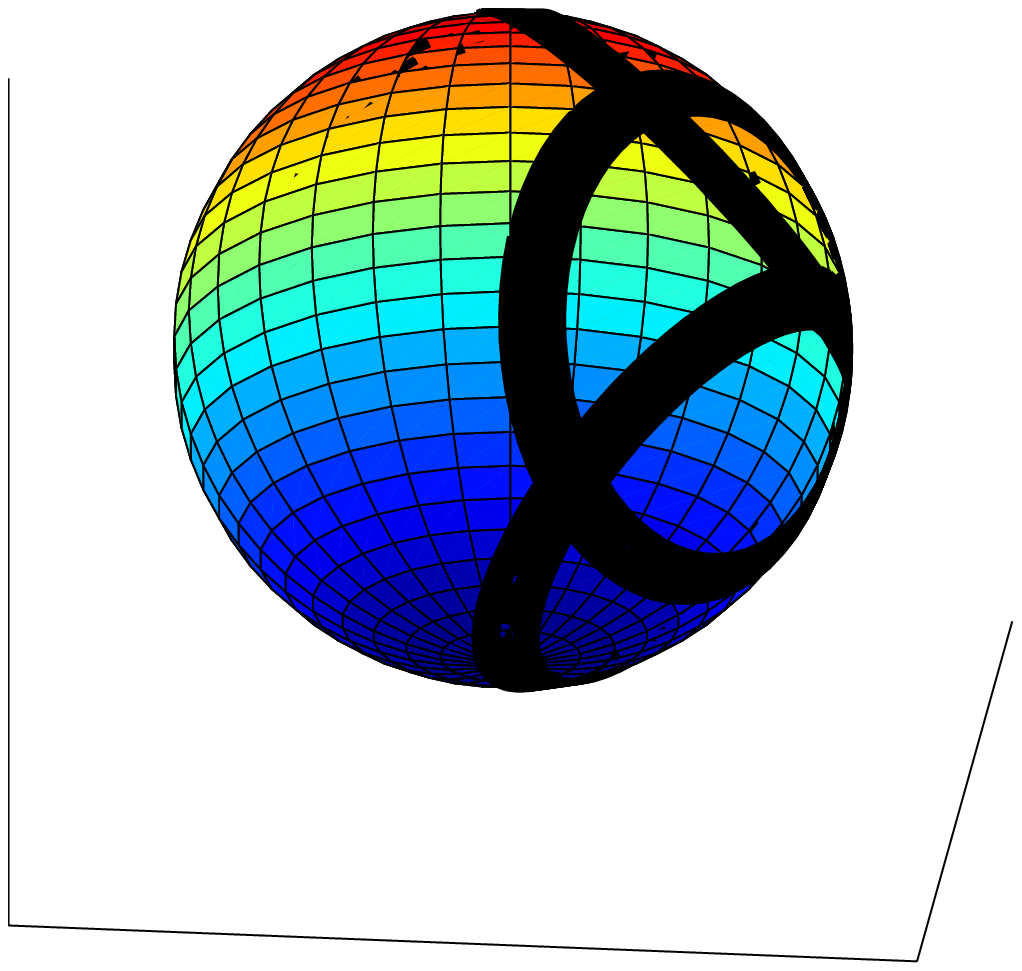}
\hspace{-10cm}$  $
\end{center}
{\it FIGURE 2B. The Figure 2 with pieces of plane waves that have
smaller width. The time moments are the same as in Figure 2, but now
all conic waves produced by 3-interactions (black cones) are visible.
Note that there are four conic waves but one of those is not visible as it is
behind the sphere.}
\medskip
}

{\bf The outline of the proof of Theorem \ref{alternative main thm Einstein}:} We consider 4 sources that send 
distorted plane waves from $W_{\hat g}$
that propagate near geodesics $\gamma_{x_j,\xi_j}$, see Fig.\ 1(Right).
Due to the non-linearity of the Einstein equations, these waves interact
and may produce a point source in the space-time, see Fig.\ 2 on the interaction
of waves. All four waves interact if the geodesics  $\gamma_{x_j,\xi_j}$ 
intersect at a single point $q$ of the space-time, see Fig.\ 3(Left).
We show in Sec.\ 3, that if the intersection of these geodesics happens before the
conjugate points of the geodesics  (i.e.\ the caustics of the waves), then there 
are some sources  satisfying (\ref{eq: lineariz. conse. law PRE}) such that the
produced spherical plane  wave has a non-vanishing singularity on the future light cone
emanating from $q$,  see Fig.\ 3(Right). Thus 
we can observe in $U_{\hat g}$
the set of the earliest light observations of
 the point $q$,  $\be_{U_{\hat g}}(q)$. By varying
the starting directions $({x_j,\xi_j})$ of the geodesics we can observe
the  collection  of sets of the earliest light observations 
for all points $q\in I(p^-,p^+)$. This data
determine uniquely the topological, differentiable and the conformal structures
on  $I(p^-,p^+)$, as is shown in the second part of this paper,  \cite{Paper-part-2}.
\MTEXT{In the proof we have to deal with several technical difficulties: First, 
the wave produced by the 4th order interaction consists of many terms which
could cancel each other. In Sec.\ 3, we  show that the principal symbol of this wave,
considered in the wave map coordinates,
does not vanish in a generic situation.
Second, we do not know the wave map coordinates in $U_{\hat g}$
and thus we have to consider observations in  normal coordinates. This change of coordinates
is a gauge transform where some of the produced
singularities may vanish. However, in Sec.\ 4 we show  that some singularities persists and can be observed. Third,
caustics may produce singularities whose interactions are difficult to analyze. We avoid
this  by using global Lorentzian geometry (in Sec.\ 2) and show that no caustics affect in
the earliest observations. We use this in Sec.\ 5 to give a step-by-step construction of the diamond set
$J^+(p^-)\cap J^-(p^+)$.}  

\begin{center}
\psfrag{1}{\hspace{-2mm}$x_1$}
\psfrag{2}{\hspace{-2mm}$x_2$}
\psfrag{3}{$q$}
\includegraphics[height=4.5cm]{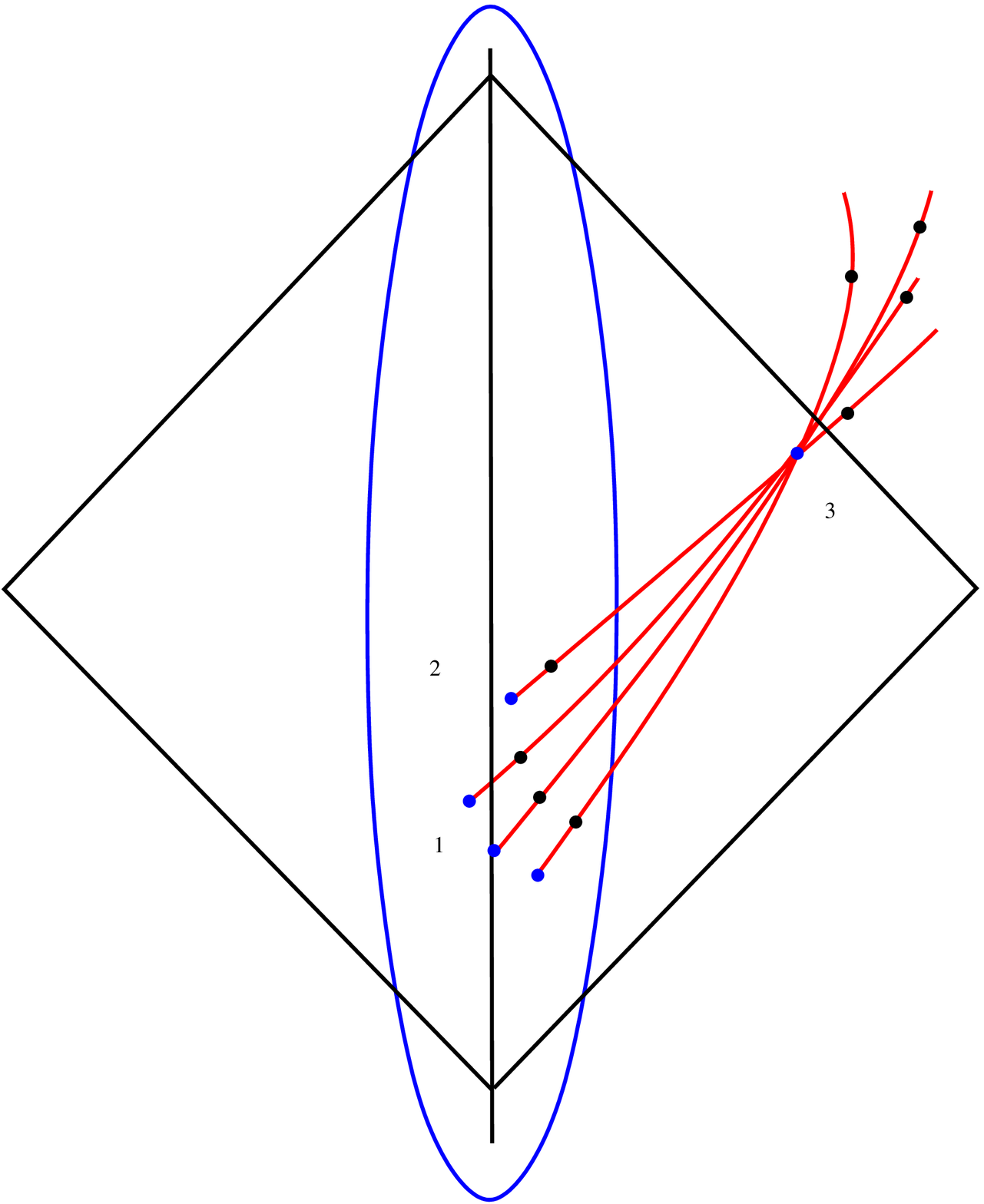}\quad\quad
\psfrag{1}{$q$}
\psfrag{2}{\hspace{-2mm}$x$}
\includegraphics[height=4.5cm]{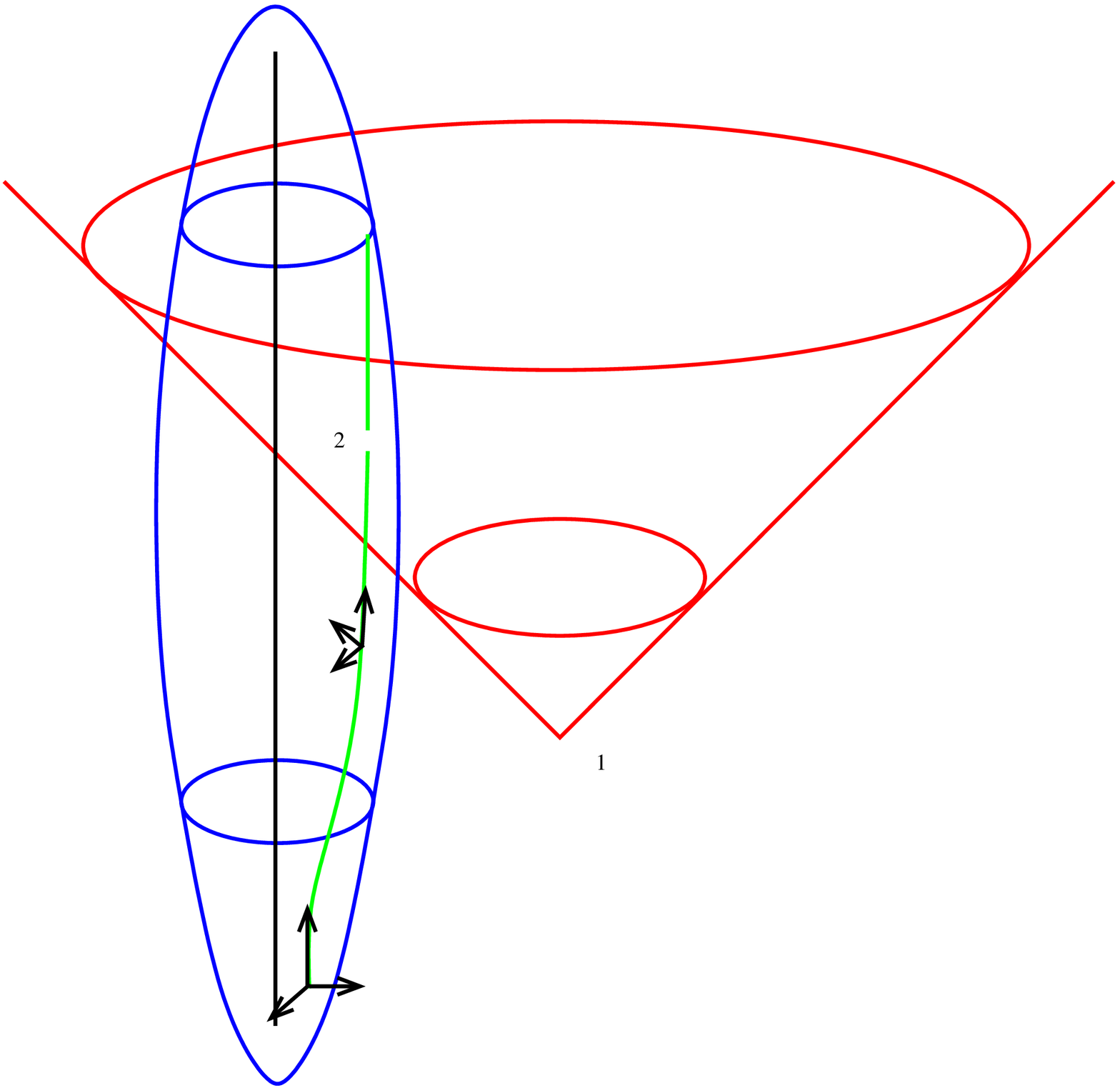}
\end{center}
{\it FIGURE 3. {\bf Left:} 
The four light-like geodesics $\gamma_{x_j,\xi_j}([0,\infty))$, $j=1,2,3,4$ starting at
the blue points $x_j$
intersect at $q$ before the first cut points of $\gamma_{x_j,\xi_j}([t_0,\infty))$,
denoted by black points. The points $\gamma_{x_j,\xi_j}(t_0)$ are also shown as black points.
{\bf Right:}
 The future light cone $\mathcal L^+_{\hat g}(q)$ emanating from the point $q$ is  shown
as a red cone. The set of the light observation points $\mathcal P_{U_{\hat g}}(q)$, see Definition \ref{def. O_U},  is the intersection
of the set $\mathcal L^+_{\hat g}(q)$ and the set $U_{\hat g}$.
The green curve is the geodesic
 $\mu=\mu_{\hat g,z,\eta}$.
This geodesic intersects the future light cone $\mathcal L^+_{\hat g}(q)$
at the point $x$.
The black vectors are the frame $(Z_j)$ that is obtained using parallel translation along
the geodesic $\mu$. Near the intersection point $x$ we detect singularities in normal coordinates
centered at $x$ and associated to the frame $(Z_j)$. 
}
\medskip
\vfill

\MTEXT{We want to point out that by the main theorem, if we have
two non-conformal  spacetimes, a generic  
measurement gives different results
on these manifolds. In particular, this implies that perfect 
space-time cloaking, in sense of light rays, see \cite{Fridman,McCall}, with a smooth metric  in a globally hyperbolic universe  is not possible.}

The assumptions of Theorem \ref{alternative main thm Einstein} are valid in many cases.  For instance, consider the a case when the background fields vary sufficiently:
 \medskip

\noindent
 {\bf Condition A}:
Assume that at any $x\in \overline U_{\hat g}$ there is
a permutation $\sigma:\{1,2,\dots,L\}\to \{1,2,\dots,L\}$, denoted $\sigma_x$, such that the 
$5\times 5$ matrix $[ B_{jk}^\sigma(\hat \phi(x),\nabla \hat \phi(x))]_{j,k\leq 5}$
is invertible, where
\ba
& &B_{j\ell}^\sigma(\hat \phi(x),\nabla \hat \phi(x))=\frac \p {\p x^j}  \hat \phi_{\sigma(\ell)}(x),\quad 
\hbox{for $j\leq 4$, $\ell=1,2,3,4,5$,} \\
& &B_{jk}^\sigma(\hat \phi(x),\nabla \hat \phi(x))
=\hat \phi_{\sigma(\ell)}(x),\quad \hbox{for $j=5$, $\ell=1,2,3,4,5$.}
\ea
  \medskip
 
 When the  Condition A is valid, also the condition $\mu$-LS is valid,
 see
 Appendix \noextension{A and \cite{Paper-inpreparation}}\extension{C}. 
  Very roughly speaking, Condition A means that the background fields vary so
 much that one could implement a measurement with some suitable sources $\F$ by taking
the needed ``energy'' from the varying $\phi_\ell$ fields.

 Theorem \ref{alternative main thm Einstein} can in some cases be improved
 so that also  the conformal factor of the metric tensor can be reconstructed.
Indeed, Theorem \ref{alternative main thm Einstein} and 
Corollary 1.3 of \cite{Paper-part-2} imply 
that if $V\subset I_{\hat g}(p^-,p^+)$ is Vacuum, i.e., Ricci-flat,
and all points $x\in V$ 
can be connected by a curve $\alpha\subset V^{int}$ to points of $U_{\hat g}$. Then
under the assumptions of Theorem \ref{alternative main thm Einstein},  the
whole metric tensor $g$ in $V$ can be reconstructed.

Theorem \ref{alternative main thm Einstein} deals with an inverse problem for ``near field''  measurements. 
We remark that for inverse problems for  linear equations the measurements of the Dirichlet-to-Neumann map or the ``near field''
measurements  are
 equivalent to scattering or ``far field'' information \cite{Ber}. Analogous considerations
 for non-linear equations have not yet been done but are plausible. On
 related inverse scattering problems, see \cite{GrZ,MsBV}.

Also, one can ask if 
one can make an approximate image of the space-time doing
 only one measurement.  In general, in many
inverse problems several measurements can be packed together to 
one measurement.
For instance, for the wave equation with a  time-independent simple metric
this is done in \cite{HLO}. Similarly, Theorem \ref{alternative main thm Einstein}  and its proof
make it possible to do approximate reconstructions in a suitable class of
manifolds with only one measurement, see Remark 5.1. 

\motivation{
The techniques considered in this paper can be used also to study
inverse problems for non-linear hyperbolic systems encountered in applications.
For instance, in medical imaging, in the
  the recently developed Ultrasound Elastography imaging technique
 the elastic material parameters are reconstructed by sending 
 (s-polarized) elastic waves that are imaged using (p-polarized) elastic waves,
 see  
e.g.\ \cite{Hoskins,McLaughlin1}.
This imaging method uses  interaction of waves and is based on the non-linearity of the system.
}

\section{Geometry of the observation times and the cut points}
%
%
%
%
%

\extension{
By \cite{Bernal}, a globally hyperbolic manifold, as defined in the
introduction, satisfies the  strong causality condition:
\beq\label{strong claim}
& &\hbox{For every  $z\in M$ and every neighborhood $V\subset M$ of $z$ there is 
}\\
\nonumber & &\hbox{a neighborhood $V^{\prime}\subset M$ of $z$ that if $x,y\in V^{\prime}$ and  $\alpha\subset M$
is }
\\
\nonumber & &\hbox{a causal path connecting $x$ to $y$ then $\alpha\subset V$.}
\eeq
}

\extension{

We use the following simple result:
\medskip

\noindent
{\bf Lemma 2.A.1} {\it Let $z\in M$.
Then
there  is a neighborhood $V$ of $z$ so that
\smallskip 

\noindent
(i) If the geodesics $\gamma_{y,\eta}([0,s])\subset V$ and $\gamma_{y,\eta^{\prime}}([0,s^{\prime}])\subset V$,
$s,s^{\prime}>0$
satisfy $\gamma_{y,\eta}(s)=\gamma_{y,\eta^{\prime}}(s^{\prime})$,
then $\eta= c\eta^{\prime}$  and $s^\prime=cs$ with some  $c>0$. 
\smallskip 

\noindent
(ii) For any $y\in V$, $\eta\in T_yM\setminus 0$ there is $s>0$ such that $\gamma_{x,\eta}(s)\not \in V$.
}
\medskip

\noindent
{\bf Proof.} The property (i) follows from 
\cite[Prop.\ 5.7]{ONeill}.  Making  $s>0$ so small that $\overline V\subset B_{g^+}(z,\rho)$ with 
a sufficiently small $\rho$, the claim (ii) follows from \cite[Lem.\ 14.13]{ONeill}. 
\hfill \Box \medskip
 }
 


Let us consider points $x, y\in M$. If  $x< y$,
we define  the time separation function $\tau(x,y)\in [0,\infty)$
to be the supremum of the lengths 
$
L(\alpha)=\int_0^1 \sqrt{-g(\dot\alpha(s),\dot\alpha(s))}\,ds
$
of the piecewise smooth
causal paths $\alpha:[0,1]\to M$ from $x$ to $y$. If  the condition $x< y$ does not
hold, we define $\tau(x,y)=0$. 

Since $M$ is globally hyperbolic,
the time separation function $(x,y)\mapsto \tau(x,y)$  is continuous in  $M\times M$
by \cite[Lem.\ 14.21]{ONeill} and  the sets $J^\pm(q)$ are closed
 by \cite[Lem.\ 14.22]{ONeill}.
Also, any points $x, y\in M$, $x< y$
 can be connected by a causal geodesic whose
 length is $\tau(x,y)$ by \cite[Prop.\ 14.19]{ONeill}.   


When $(x,\xi)$ is a  light-like vector, we define  $\T(x,\xi)$ to
be the length of
the maximal interval on which $\gamma_{x,\xi}:[0,\T(x,\xi))\to M$ is defined.
Below,  to simplify notations, we sometimes use the notation $\gamma_{x,\xi}([0,\infty))$ for the geodesic
 $\gamma_{x,\xi}([0,\T(x,\xi))$.

When $(x,\xi_+)$ is a  future pointing light-like vector,
and $(x,\xi_-)$ is a  past pointing light-like vector,
we define the modified cut locus functions, c.f.\ \cite[Def.\ 9.32]{Beem},
\beq\label{eq: max time}
& &\rho_g(x,\xi_+)=\sup\{s\in [0,\T(x,\xi_+));\ \tau(x,\gamma_{x,\xi_+}(s))=0\},\\
\nonumber
& &\rho_g(x,\xi_-)=\sup\{s\in [0,\T(x,\xi_-));\ \tau(\gamma_{x,\xi_-}(s),x)=0\}.
\eeq
The point $\gamma_{x,\xi}(s)|_{s=\rho(x,\xi)}$ is called the cut point
on the geodesic  $\gamma_{x,\xi}$.

Using \cite[Thm. 9.33]{Beem}, we see that the function $\rho_g(x,\xi)$ is lower semi-continuous on a globally hyperbolic
Lorentzian manifold $(M,g)$.

%

Below, in this section, we consider the manifold $(M,\hat g)$ and 
 denote by $\gamma_{x,\xi}$
the  geodesics of $(M,\hat g)$ and $\rho(x,\xi)=\rho_{\hat g}(x,\xi)$. Also, we
denote $\mu_{{\hat g},z_0,\eta_0}=\hat \mu$, $p^\pm=\mu_{\hat g}(s_\pm)$, 
 $p_{+2}={\hat \mu}(s_{+2})$, and
 $p_{-2}={\hat \mu}(s_{-2})$.
Recall that by  (\ref{kaava D}),
$\mu_{{\hat g},z_0,\eta_0}(s_{\pm j})\in I^\mp(\mu_{{\hat g},z,\eta}(s_{\pm (j+1)})),$ $j=1,2$,
for all $(z,\eta)\in \U_{z_0,\eta_0}$. 

\begin{definition}\label{def: f functions}
Let $\mu=\mu_{\hat g,z,\eta}$, $(z,\eta)\in \U_{z_0,\eta_0}$. 
For $x\in J^+(\mu(-1))\cap J^-(\mu(+1))$ we define $f_\mu^\pm(x)\in [-1,1]$
by setting 
\ba
& &
 f_\mu^+(x) =\inf (\{s\in (-1,1);\tau(x,\mu(s))>0\}\cup\{1\}),
\\
& &f_\mu^-(x) =\sup (\{s\in (-1,1);\ \tau(\mu(s),x)>0\}\cup\{-1\}).
\ea 
\end{definition}

We need the following simple properties of these functions.

\begin{lemma} \label{B: lemma} 
Let $\mu=\mu_{{\hat g},z,\eta}$, $(z,\eta)\in \U_{z_0,\eta_0}$,
and   $x\in  J^-(p_{+2})\cap J^+(p_{-2})$.

\smallskip 

\noindent
(i) The function $s\mapsto \tau(x,\mu(s))$ is non-decreasing on the interval $s\in [-1,1]$
and strictly increasing on $s\in [f_\mu^+(x),1]$.
\smallskip 

\noindent
(ii) We have that 
 $s_{-3}< f_\mu^+(x)<s_{+3}$.
\smallskip 

\noindent (iii)  Let $y=\mu(f_\mu^+(x))$. Then $\tau(x,y)=0$.
Also, if $x\not\in \mu$, there is a light-like geodesic $\gamma([0,s])$ in $M$ from $x$ to  $y$ with no conjugate points
on $\gamma([0,s))$.

\smallskip 

\noindent (iv)
The maps $f^+_\mu:J^-(p_{+2})\cap J^+(p_{-2})\to [-1,1]$  is  continuous. 
%
%
\extension{
The analogous results hold for $f^-_\mu:J^-(p_{+2})\cap J^+(p_{-2})\to [-1,1]$. 

\noindent (v)  For $q\in J^-(p^+)\setminus I^-(p^-)$ the map $F:\U_{z_0,\eta_0}\to \R;$
$F(z,\eta)=f^+_{\mu(z,\eta)}(q)$ is continuous.}
\smallskip

\noextension{The analogous results hold for $f^-_\mu:J^-(p_{+2})\cap J^+(p_{-2})\to [-1,1]$.}
\end{lemma}

\noindent
{\bf Proof.} (i)
and (ii) follows from the definition of $f^+_\mu$
and the fact that $p_{\pm 2}\in I^\mp(\mu(s_{\pm 3}))$ by (\ref{kaava D}). Claim
(iii) follows from   \cite[Lem.\ 10.51]{ONeill}.

(iv)  Assume that $x_j\to x$ in $J^-(p_{+2})\cap  J^+(p_{-2})$ as $j\to\infty$.
Let
$s_j=f_\mu^+(x_j)$ and $s=f_\mu^+(x)$.  As  $\tau$ is continuous,
 for any $\e>0$ we have 
$\lim_{j\to\infty} \tau(x_j,\mu(s+\e))=\tau(x,\mu(s+\e))>0$ and thus
for $j$ large enough 
$s_j\leq s+\e.$ Thus $\limsup_{j\to \infty}s_j\leq s$.
Assume next that $\liminf_{j\to \infty}s_j=\tilde s< s$ and denote $\e=\tau(\mu(\tilde s),\mu(s))>0$. Then
by the reverse triangle inequality, \cite[Lem.\ 14.16]{ONeill},
$\liminf_{j\to\infty}\tau(x_j,\mu(s))\geq
\liminf_{j\to\infty}\tau(\mu(s_j),\mu(s))
\geq \e$, and as $\tau$ is continuous in $M\times M$,
we obtain $\tau(x,\mu(s))\geq \e$, which is not possible since $s=f_\mu^+(x)$.
Hence $s_j\to s$ as $j\to \infty$, proving (iv). The analogous results for $f_\mu^-$
follow similarly.
\extension{

(v) Observe that as $J^+(q)$ is a closed set,  $F(z,\eta)$ is equal to the smallest value $s\in [-1,1]$
such that $\mu_{z,\eta}(s)\in J^+(q)$.
Let $(z_j,\eta_j)\to (z,\eta)$ in $(TM,g^+)$ 
as $j\to \infty$ and $s_j=F(z_j,\eta_j)$
and $\underline s=\liminf_{j\to \infty}s_j$.
As  the map $(z,\eta,s)\mapsto \mu_{z,\eta}(s)$
is continuous, we see that  for a suitable subsequence
 $\mu_{z,\eta}(\underline s)=\lim_{k\to \infty} \mu_{z_{j_k},\eta_{j_k}}(s_{j_k})\in J^+(q)$
 and hence $F(z,\eta)\leq \underline s=\liminf_{j\to \infty}F(z_j,\eta_j)$.
 This shows that $F$ is lower-semicontinuous.
 
 On the other hand, let $\overline s=F(z,\eta)$. As 
 $\mu_{z,\eta}$ is a time-like geodesic,  for any $\e\in (0,1-\overline s)$ we have
 $\tau(q,\mu_{z,\eta}(\overline s+\e))>0$. Since $\tau$
 and the map $(z,\eta,s)\mapsto \mu_{z,\eta}(s)$
are continuous, we see that there is $j_0$ such that 
 if $j>j_0$ then $\tau(q,\mu_{z_j,\eta_j}(\overline s+\e))>0$.
 Hence, $F(z_j,\eta_j)\leq \overline s+\e$. Thus
 $\limsup_{j\to \infty}F(z_j,\eta_j)\leq \overline s+\e$,
 and as $\e>0$ can be chosen to be arbitrarily small, we have  $\limsup_{j\to \infty}F(z_j,\eta_j)\leq \overline s=
 F(z,\eta)$. Thus $F$ is also upper-semicontinuous. This proves (v).}
\hfill \Box \medskip

Let $W\subset M$. We define the earliest points of set $W$ on 
the curve $\mu_{z,\eta}=\mu_{z,\eta}([-1,1])$, and in the set $U$, respectively, to be 
\beq\label{Earliest element sets}
& &\pointear_{z,\eta}(W)\hspace{-1mm}=\hspace{-1mm}\{\mu_{z,\eta}(\inf\{s\in [-1,1];\ \mu_{z,\eta}(s)\in W\})\},
\hbox{ if }\mu_{z,\eta}\cap W\not =\emptyset, \hspace{-10mm}\\
\nonumber 
& &\pointear_{z,\eta}(W)\hspace{-1mm}=\hspace{-1mm}\emptyset,
\quad\hbox{if }\mu_{z,\eta}\cap W =\emptyset.
\eeq

\extension{
\medskip

\noindent
{\bf Lemma  2.A.2} 
{\it Let $K\subset M$ be a compact set. Then there is $R_1>0$ such that
if $\gamma_{y,\theta}([0,l])\subset K$ is a light-like geodesic
with $\|\theta\|_{g^+}=1$, then $l\leq R_1$. In the case when $K=J(p^-,p^+)$,
with $q^-,q^+\in M$ we have $\gamma_{y,\theta}(t)\not \in J(q^-,q^+)$ for $t>R_1$.
}
\medskip

The proof of this lemma is standard, but we include it for the convenience of the reader.
\medskip

\noindent
{\bf Proof.}
Assume that there are no such $R_1$.
Then
there are geodesics $\gamma_{y_j,\theta_j}([0,l_j])\subset K$, $j\in \Z_+$
such that $\|\theta_j\|_{g^+}=1$ and $l_j\to \infty$ as $j\to \infty$. 
Let us choose a subsequence $(y_j,\theta_j)$ which
converges to some point $(y,\theta)$ in $(TM,g^+)$. 
As  $\theta_j$ are light-like, also $\theta$ is light-like.

Then, we observe that for all $R_0>0$ the functions $t\mapsto \gamma_{y_j,\theta_j}(s),$
converge in $C^1([0,R_0];M)$ to $s\mapsto \gamma_{y,\theta}(t),$ as $j\to \infty$. 
As  $\gamma_{y_j,\theta_j}([0,l_j])\subset K$ for all $j$,
we see that $\gamma_{y,\theta}([0,R_0])\subset K$ for all $R_0>0$.
Let $z_n=\gamma_{y,\theta}(n)$, $n\in \Z_+$. As  $K$ is compact,
we see that there is a subsequence $z_{n_k}$ which converges
to a point $z$ as $n_k\to \infty$. Let now $V\subset M$ be a small convex neighborhood
of $z$ such that each geodesic starting from $V$ exits the set $V$ (cf.\ Lemma 
2.A.1). 
Let $V^{\prime}\subset M$ be a neighborhood
of $z$  so that the strong causality condition (\ref{strong claim})
is satisfied for $V$ and $V^{\prime}$. Then we see that there is $k_0$ such that
if $k\geq k_0$ then $z_{n_k}\in V^{\prime}$, implying that  
$\gamma_{y,\theta}([n_{k_0},\infty))\subset V$. 
This is a contradiction and thus the claimed $R_1>0$ exists.

Finally, in the case when $K=J(q^-,q^+)$, $q^-,q^+\in M$ we see that
if $q(s)=\gamma_{y,\theta}(s) \in K$ for some $s>R_1$, then for all $\tilde s\in [0,s]$
we have 
$q(\tilde s)\leq q(s) \leq q^+$ and $q_-\leq q(0)\leq q(\tilde s)$. Thus 
$q(\tilde s)\in K$ for all $\tilde s\in [0,s]$. As  $t>R_1$, this is not possible by
the above reasoning, and thus the last assertion follows.
\hfill \Box \medskip
}

\hiddenfootnote{ 
SLAVA ASKED THAT WE INCLUDE THE FOLLOWING LEMMA TO SIMPLIFY THE PROOF. NOTE THAT THIS LEMMA HAS BEEN INCLUDED IN A PROOF IN THE "COMBINING RESULTS" SECTION.\begin{lemma}\label{lemma:Slava1}
Let  $x_1\in I^-(p^+)\setminus J^-(p^-)$.  Let $\gamma_1([0,t_1])$ and $\gamma_2([0,t_2])$  be two future going light-like geodesic with 
$\gamma_1(0)=\gamma_2(0)=x_1$ and that there are no $c>0$ such that
 $\dot \gamma_1(0)=c\dot \gamma_2(0)$. Assume that 
$x_2=\gamma_1(s_1)=\gamma_2(s_2)\in I^-(p^+)$ with $s_1>0$. 
Let $z_j=\pointear_\mu(x_j)\in \mu$, $j=1,2$. Then either $z_1\ll z_2$ or $x_2=z_1$.
\end {lemma}

\noindent
{\bf Proof.}  As  $s_1>0$, the fact that there are no closed causal curves implies that $s_2>0$ and
that $z_1\leq z_2$. As $z_1,z_2\in \mu$ and $\mu$  is time-like curve, we then have either $z_1\ll z_2$
or $z_1=z_2$. If the former one is true, the claim holds. In the case when $z_1=z_2$,
and let us assume, opposite to our claim, that $z_1\not =x_2$. 

Then $\tau(x_1,z_1)=0$ and $\tau(x_2,z_2)=0$. As $z_1=z_2\not =x_2$, this 
implies that there is 
a light-like geodesics, denoted $\gamma_{x_2,\xi}([0,a])$, from $x_2$ to $z_1$.
Moreover, let
$\eta_j$, $j=1,2$ be a curve from $x_1$ to $z_1$ which 
is the union of  $\gamma_j([0,s_j])$
from $x_1$ to $x_2$ and 
the light-like geodesics $\gamma_{x_2,\xi}([0,a])$, from $x_2$ to $z_1$. 
Then either $\eta_1$ or $\eta_2$ is a causal curve  
which is not a light-light geodesic from $x_1$ to $z_1$. Hence by 
\cite[Prop.\ 10.46, 10.51]{ONeill} there is a time-like curve connecting 
 $x_1$ to $z_1$ and thus $\tau(x_1,z_1)>0$. However, this contradicts  the equation $\tau(x_1,z_1)=0$,
 and hence $z_1=z_2\not =x_2$ is not possible. This proves the claim.
\hfill \Box}

Below, we use  for a pair $(x,\xi)\in L^+\hattuM $  the notation 
\beq\label{eq: x(h) notation}
& &(x(h),\xi(h))=(\gamma_{x,\xi}(h),\dot \gamma_{x,\xi}(h)).
\eeq


 
 \MATTITEXT{
 
\begin{center}
\psfrag{1}{${\hat x}$}
\psfrag{2}{$x$}
\psfrag{3}{$z$}
\psfrag{6}{\hspace{-1mm}$q_1$}
\psfrag{5}{$q_0$}
\psfrag{4}{\hspace{-3mm}$p_2$}
\psfrag{7}{\hspace{-1mm}$x_1$}
\psfrag{8}{$\gamma_{{x},\xi}$}
\psfrag{9}{$\gamma_{{\hat x},\hat \zeta}$}
\includegraphics[height=5.5cm]{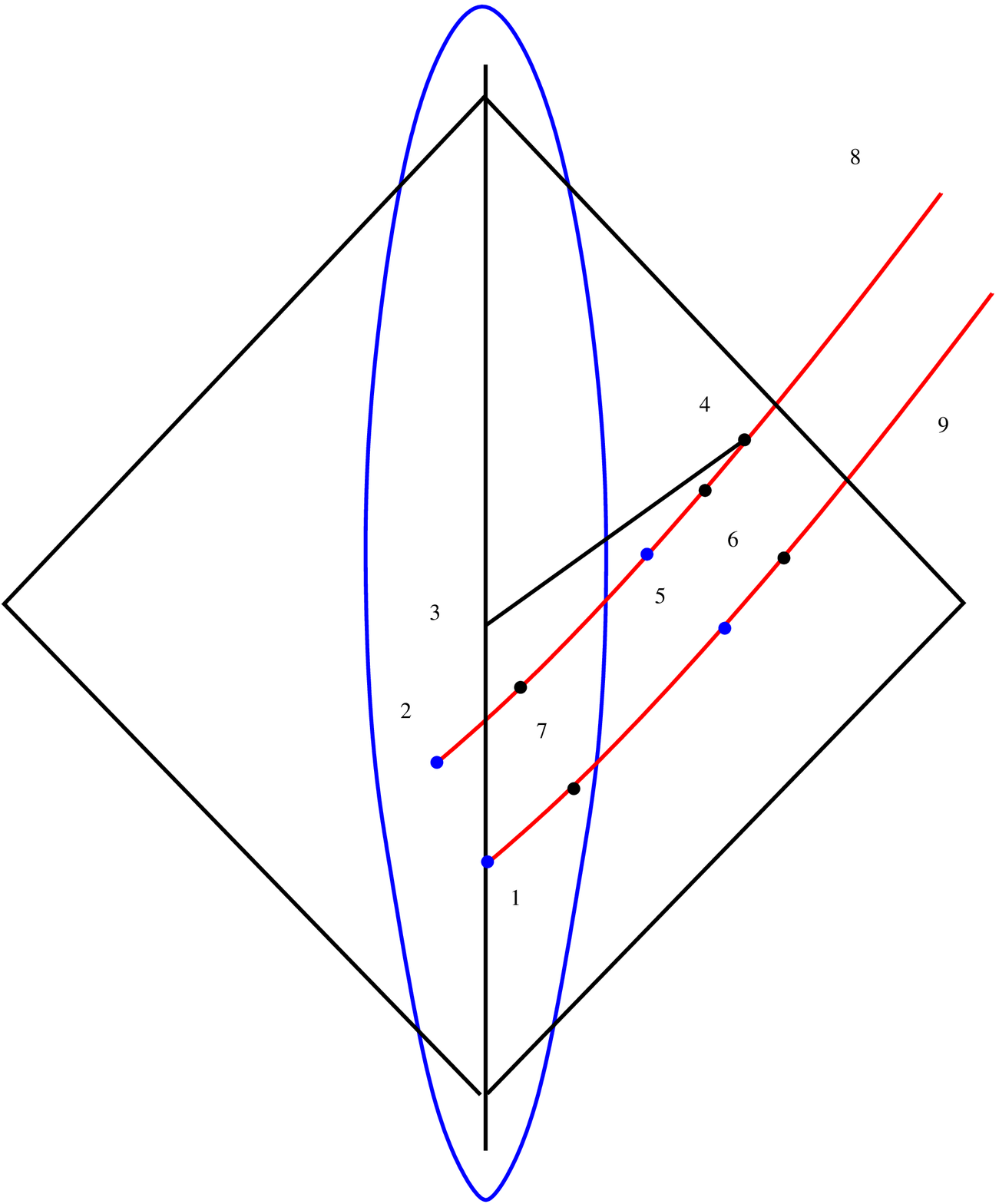}
\psfrag{1}{$x_1$}
\psfrag{2}{\hspace{-1mm}$x_2$}
\psfrag{3}{$q$}
\psfrag{4}{\hspace{-0mm}$p_1$}
\psfrag{5}{\hspace{-2mm}$p_2$}
\psfrag{6}{}
\psfrag{7}{$x_5$}
\includegraphics[height=5.5cm]{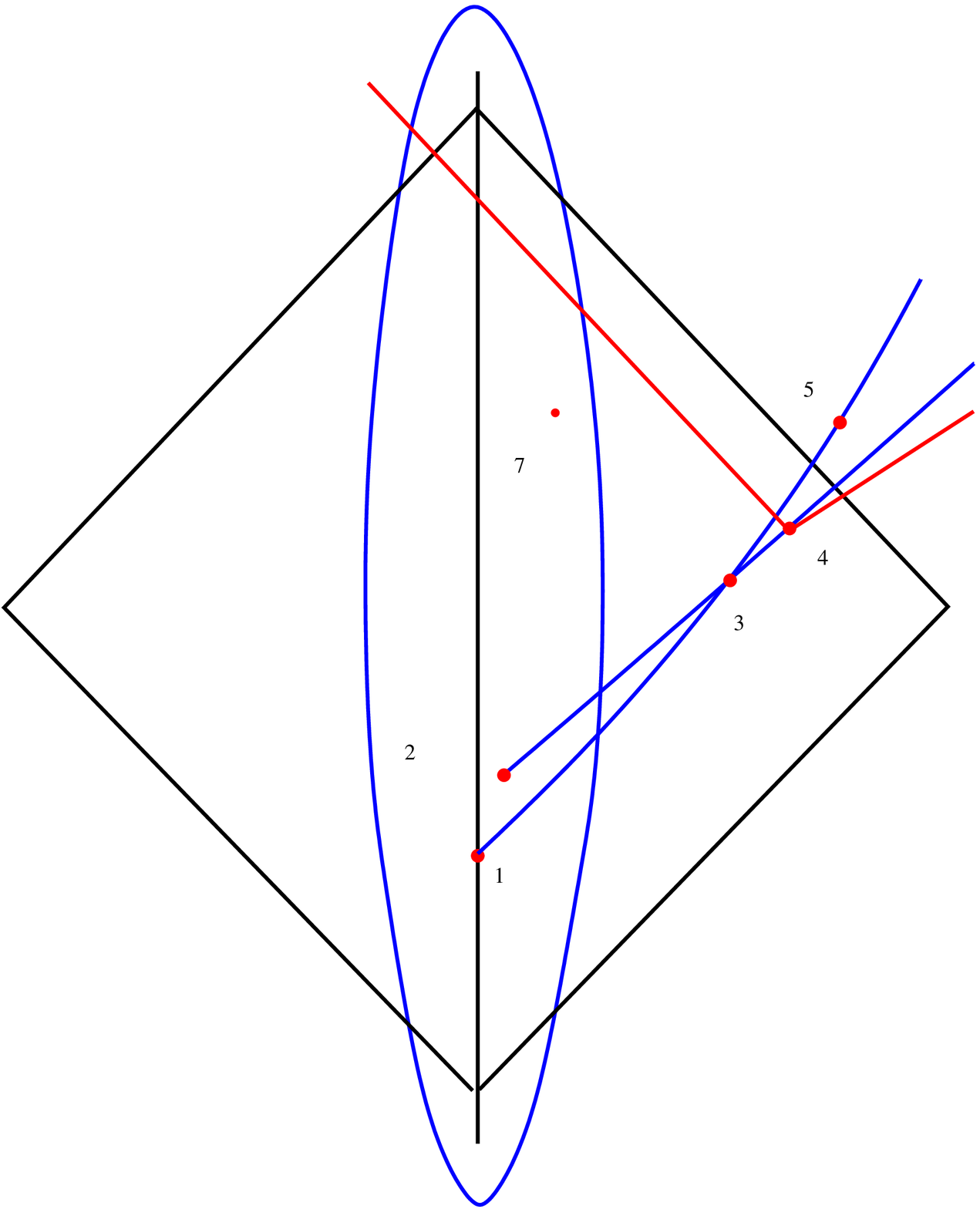}
\end{center}


\noindent
{\it FIGURE 4. {\bf Left.} Figure shows the situation in Lemma \ref{lem: detect conjugate 0}. The point ${\hat x}=\hat \mu(r_1)$ is on the time-like
geodesic $ \hat \mu$ shown as a black line. 
The black diamond is the set $J_{{\hat g}}( p^-,p^+)$,
 $(x,\xi)$ is a light-like direction close to $({\hat x},\hat\zeta)$,
 and $x_1=\gamma_{{x},\xi}(t_0)=x(t_0)$.
The points $q_0=\gamma_{x,\xi}(\rho(x,\xi))$ and $q_1=
\gamma_{x(t_0),\xi(t_0)}(\rho(x(t_0),\xi(t_0)))$ are the  first cut point on  
$\gamma_{{x},\xi}$ corresponding to the points $x$ and $x_1$,
respectively. The blue and black points on $\gamma_{{\hat x},\hat \zeta}$ are the corresponding
cut points on 
 $\gamma_{{\hat x},\hat \zeta}$.
Also, $z=\hat \mu(r_2)$, 
  $r_2=f^-_{\hat \mu}(p_2)$.
{\bf Right:}
The figure shows the configuration in formulas (\ref{eq: summary of assumptions 1}) and  (\ref{eq: summary of assumptions 2}).
We send light-like geodesics $\gamma_{x_j,\xi_j}([t_0,\infty))$ from $x_j$,
$j=1,2,3,4.$ The boundary $\p \mathcal V((\vec x,\vec \xi), t_0)$ is denoted by red line segments. We assume the these  geodesics intersect
at the point $q$ before their 
 first cut points $p_j$.
}

\medskip
\motivation{ Later, we will consider  wave packets sent from 
a point $x$ that propagate
near a geodesic $\gamma_{x,\xi}([0,\infty))$. These waves may have singularities near the conjugate
points of the geodesic and due to this we analyze next how the conjugate points
move along a geodesic when the initial point of the geodesic is moved
from $x$ to $\gamma_{x,\xi}(t_0)$. }

 \begin{lemma} \label{lem: detect conjugate 0} 
   There are 
$\vartheta_1,\kappa_1,\kappa_2>0$ such that 
for  all  $\hat x=\hat \mu(r_1)$ with $r_1\in [s_-,s_+]$,
${\hat \zeta}\in L^+_{\hat x}M$, $\|{\hat \zeta}\|_{{\hat g}^+}=1$,
$t_0\in [\kappa_1,4\kappa_1]$,  
and 
$(x,\xi)\in L^+M$ satisfying $d_{{\hat g}^+}(({\hat x},{\hat \zeta}),(x,\xi))\leq \vartheta_1$
the following holds:

(i)  $0<t\leq 5\kappa_1$, then $\gamma_{x,\xi}(t)\in U_{\hat g}$ and  $f^-_{\hat \mu}(\gamma_{\hat x,{\hat \zeta}}(t))=r_1$.



(ii) 
Assume that $t_2
\in [t_0+\rho(\gamma_{x,\xi}(t_0),\dot \gamma_{x,\xi}(t_0)),\T(x,\xi))$  and  $p_2=\gamma_{x,\xi}(t_2)\in J^-({\hat \mu}(s_{+2}))$. 
Then    
 $r_2=f^-_{\hat \mu}(p_2)$ satisfies $r_2-r_1>2\kappa_2$.
\end{lemma}

Note  that above in (ii) we can choose $t_2= t_0+\rho(\gamma_{x,\xi}(t_0),\dot \gamma_{x,\xi}(t_0))$ in which case
$p_2$ is the first cut point $q_1$ of $\gamma_{x,\xi}([t_0,\infty))$, see Fig.\ 4(Left).

{\tobecheckedtext 

\noindent
{{\bf Proof.} 
  Let $B=\{({\hat x},{\hat \zeta})\in L^+M;\ {\hat x}\in {\hat \mu}([s_-,s_+]),\ \|{\hat \zeta}\|_{  {\hat g}^+}=1\}$.
Since $B$ is compact, the positive and lower semi-continuous function  $\rho(x,\xi)$ obtains
its minimum on $B$. Hence 
we see that (i)  holds when $\kappa_1\in (0, \frac 15 \inf\{\rho({\hat x},{\hat \zeta});\ ({\hat x},{\hat \zeta})\in B\})$ is small enough.

(ii) Let $K$ be the compact set
 $K=\{(x,\xi)\in L^+M;\ d_{{\hat g}^+}(({x},{\xi}),B)\leq \vartheta_1\}$. Also, let
$T_{+}(x,\xi)=\sup \{t\geq 0\ ;\ \gamma_{x,\xi}(t)\in J^-(p_{+2})\}$
and \ba
K_0=\{(x,\xi)\in K; \ {\rho(x(\kappa_1),\xi(\kappa_1))+\kappa_1
} 
\leq T_{+}(x,\xi)\},
\quad K_1=K\setminus K_0.
\ea 
Using  \cite[Lem.\ 14.13]{ONeill}, we see that $T_{+}(x,\xi)$ is bounded in $K$.
 {Note that for $t_0\geq \kappa_1$ and $a>t_0$ the geodesic
$\gamma_{x,\xi}([t_0,a])$ can have a cut point only if $\gamma_{x,\xi}([\kappa_1,a])$
has a cut point and thus
$t_0+\rho(x(t_0),\xi(t_0))\geq \kappa_1+\rho(x(\kappa_1),\xi(\kappa_1))$.
}
If $K_0=\emptyset$,
 the claim is valid as the condition $p_2\in J^-(p_{+2})$  does not hold for any $(x,\xi)\in K_1$.
 Thus it is enough consider the case when  $K_0\not=\emptyset$.

We can also assume that
 $\vartheta_1>0$ is so small that for all  $(x,\xi)\in K$ we have $f^-_{\hat \mu}(x)>s_{-2}$.
Then, by Lemma \ref{B: lemma},  the map $L: G_0=\{(x,\xi,t)\in K\times \R_+;\ 
{\rho(x(\kappa_1),\xi(\kappa_1))+\kappa_1\leq } t\leq T_{+}(x,\xi)\}\to \R,$ defined by
$
L(x,\xi,t)= f^-_{\hat \mu}(\gamma_{x,\xi}(t))- f^-_{\hat \mu}(x),
$ 
is continuous.
Since $\rho(x,\xi)$ is lower semi-continuous
and $T_{+}(x,\xi)$ is upper semi-continuous  and bounded 
we have  that the sets $K_0$  and $G_0$ are compact. 

{For $(x,\xi,t)\in G_0$,
the geodesic $\gamma_{x,\xi}([\kappa_1,t])$ has a cut point in which case
 we see, for $y=\gamma_{x,\xi}(t)$, 
that 
$\tau(x,y)>0$. 
Thus, for $z_1=\hat \mu(f^-_{\hat \mu}(x))$, we have $\tau(z_1,y)\geq
\tau(z_1,x)+\tau(x,y)\geq \tau(x,y)>0$. This shows that
 $L(x,\xi,t)>0$.
Since $G_0$ is compact and $L$ is continuous and strictly positive, $\e_1:=\inf\{L(x,\xi,t);\
(x,\xi,t)\in G_0\}>0$.

}

As  $f^-_{\hat \mu}$ is continuous and
$\hat \mu([-1,1])$ is compact, we have that, by making $\vartheta_1$ smaller if necessary,
we can assume that if ${\hat x}\in \hat \mu$ and $d_{{\hat g}^+}(x,{\hat x})\leq \vartheta_1$ then
 $|f^-_{\hat \mu}(x)-f^-_{\hat \mu}({\hat x})|<\e_1/2$. 
%
Let $\kappa_2=\e_1/4$.
Then, $\rho(x(\kappa_1),\xi(\kappa_1))+\kappa_1<t_2<\mathcal T_+(x,\xi)$ so that  $r_2=f^-_{\hat \mu}(p_2)$ and  $r_3=f^-_{\hat \mu}(x)$ satisfies 
$r_2-r_3\geq \e_1$ and 
$r_2-r_1>\e_1/2$.
This proves the claim.}
\hfill \Box \medskip
}

Note that for proving the unique solvability of the inverse problem we 
need to consider  two manifolds, $(M^{(1)},\hat g^{(1)})$ and $(M^{(2)},\hat g^{(2)})$
 with the same  data. For these manifolds, we can choose $\vartheta_1,\kappa_j$ so that they are
the same for both manifolds.
\subsubsection{Geometric results on  the light observation sets}

Let  us first consider $M$ with a fixed metric $\hat g$. Denote below in this subsection
 $U=U_{\hat g}$.
See (\ref{Earliest element sets}) for the notations we use.
%

\begin{definition}\label{def. O_U}
We define   the light-observation  set of the point $q\in M$ to be $\P_U(q)=(\L_{\hat g}^+(q)\cup \{q\})\cap U=
 \{\gamma_{q,\eta}(r)\in M;\ r\geq 0,\ \eta\in L_q^{+}M, 
\ \gamma_{q,\eta}(r)\in U\}$, see Fig.\ 3(Right). 
\MTEXT{The set of the earliest light observations of $q$ is 
 $\be_U(q)=\bigcup_{(z,\eta)\in \U_{z_0,\eta_0}}\pointear_{z,\eta}(\P_U(q))$, that is,}
\ba
\be_U(q)=\{\gamma_{q,\eta}(r)\in M;\ r\in [0,\rho(q,\eta)],\ \eta\in L_q^{+}M, 
\ \gamma_{q,\eta}(r)\in U\}\subset \P_U(q).
\ea 
\end{definition}

 Below, 
 when $X$ is a set, let  $P(X)=2^X=\{Z;\ Z\subset X\}$ denote the power set of $X$.
When $\Phi:U_1\to U_2$ is a map, we 
say that the power set extension of $\Phi$ is the map
$\tilde \Phi: 2^{U_1}\to 2^{U_2}$ given by
$\tilde \Phi(U') =\{\Phi(z); \, z \in U'  \}$ for $U' \subset U$.
  We need the following theorem proven in \cite{Paper-part-2} with  $V=I_{\hat g}^+(p^-)\cap I_{\hat g}^-(p^+)\subset 
  I_{\hat g}^-(\mu_{\hat g}(s_{+2}))\setminus I_{\hat g}^-(\mu_{\hat g}(s_{-2})).$
%
%
%

\begin{theorem}\label{main thm paper part 2}
Let $(M_j,{\hat g}_j)$, $j=1,2$ be two open, $C^\infty$-smooth, globally hyperbolic 
Lorentzian manifolds of dimension $(1+3)$ and let $p^+_j,p^-_j\in M_j$ be the points of a 
  time-like geodesic $\mu_{{\hat g}_j}([-1,1])\subset M_j$,
$p_j^ \pm=\mu_{{\hat g}_j}(s_\pm)$. Let  $U_j\subset M_j$ be a neighborhood
of $\mu_{{\hat g}_j}([s_-,s_+])$ and
 $V_j = I_{{\hat g}_j}^-(p^+_j)\cap  I_{{\hat g}_j}^+(p^-_j) \subset M_j$. Then
\medskip

\noindent
(i)  The map $\be_{U_1}:q\mapsto  \be_{U_1}(q)$ in injective in $V_1$.

\medskip 
\noindent
(ii) Let us denote by
 $
\be_{U_j}(V_j)=\{\be_{U_j}(q);\ q\in V_j\}\subset 2^{U_j}
$
the collections of  the
sets of the earliest light observations 
 on the manifold $(M_j,{\hat g}_j)$ of the points in the set $V_j$.  
Assume that there is a conformal diffeomorphism
$\Phi:U_1\to U_2$ such that $\Phi(\mu_1(s))=\mu_2(s)$, $s\in [s_-,s_+]$ 
and the power set extension $\tilde \Phi$ of $\Phi$ satisfies
\ba
\tilde \Phi(\be_{U_1}(V_1))= \be_{U_2}(V_2).
\ea
Then there is a diffeomorphism $\Psi:V_1\to V_2$
such the metric $\Psi^*{\hat g}_2$ is conformal to ${\hat g}_1$ and
$\Psi|_{V_1\cap U_1}=\Phi$.
\end{theorem}

\section{Analysis of the Einstein equations in wave coordinates}

\subsection{Asymptotic analysis of the reduced Einstein equations}
\subsubsection{The reduced Einstein tensor}
\label{sssec: relation}

For the  Einstein equations, we will consider a smooth background metric $\hat g$ on $\hattuM $
and the smooth metric $\tilde g$ for which $\hat g<\tilde g$ and $(\hattuM ,\tilde g)$
is globally hyperbolic. We also use the notations defined in Section \ref{sec:notations 1}. 
In particular, we identify $M=\R\times N$ and consider the metric tensor $g$ on $\hattuM _0=
(-\infty,t_0)\times N$, $t_0>0$ that
coincide with $\hat g$ in $(-\infty,0)\times N$. Recall also that we consider
a freely falling observer $\hat \mu=\mu_{\hat g}:[-1,1]\to \hattuM _0$
for which $\hat \mu(s_-)={p^-}\in  [0,t_0]\times N$. 
We denote 
 $L_x^+M_0=L^+_x(M_0,\hat g)$ and $L^+M_0=L^+(M_0,\hat g)$, and
the cut locus function on $(M_0,\hat g)$ by $\rho(x,\xi)=\rho_{\hat g}(x,\xi)$.
Also, recall that  $\hat U=U_{\hat g}$ the neighborhood of the geodesic $\hat \mu=\mu_{\hat g}$. We denote by $\gamma_{x,\xi}(t)$ the geodesics of $(M_0,\hat g)$.
%

Following \cite{FM} we recall that 
\beq\label{q-formula2copy}
\Ric_{jk}(g)&=& \Ric_{jk}^{(h)}(g)
+\frac 12 (g_{j q}\frac{\p \Gamma^q}{\p x^{k}}+g_{k q}\frac{\p \Gamma^q}{\p x^{j}})
\eeq
where  $\Gamma^q=g^{mn}\Gamma^q_{mn}$,
\beq\label{q-formula2copyB}
& &\hspace{-1cm}\Ric_{jk}^{(h)}(g)=
-\frac 12 g^{pq}\frac{\p^2 g_{jk}}{\p x^p\p x^q}+ P_{jk},
\\ \nonumber
& &\hspace{-2cm}P_{jk}=  
g^{ab}g_{pq}\Gamma^p_{j b} \Gamma^q_{k a}+ 
\frac 12(\frac{\p g_{jk}}{\p x^a}\Gamma^a  
+ \nonumber
g_{k l}  \Gamma^l _{ab}g^{a q}g^{bd}  \frac{\p g_{qd}}{\p x^j}+
g_{j l} \Gamma^l _{ab}g^{a q}g^{bd}  \frac{\p g_{qd}}{\p x^k}).\hspace{-2cm}
\eeq
Note that $P_{jk}$ is a polynomial of $g_{pq}$ and $g^{pq}$ and first derivatives of $g_{pq}$.

The $\hat g$-reduced Einstein tensor 
 $\Ein_{\hat g} (g)$ and   Ricci tensor
  $\Ric_{\hat g} (g)$  
are 
\beq& &
\label{Reduced Ric tensor}
(\Ric_{\hat g} (g))_{jk}=\Ric_{jk} g-\frac 12 (g_{jn} \hat \nabla _k F^n+ g_{kn} \hat \nabla _j F^n)
\\ \nonumber
& &\quad\quad\quad\quad\quad=\Ric_{jk}^{(h)}(g)
+\frac 12 (g_{j q}\frac{\p }{\p x^{k}}(g^{ab}\hat \Gamma^q_{ab})+g_{k q}\frac{\p }{\p x^{j}}(g^{ab}\hat \Gamma^q_{ab})),
\\
\label{Reduced Einstein tensor}
& &(\Ein_{\hat g} (g))_{jk}=(\Ric_{\hat g} (g))_{jk}-\frac 12 (g^{ab}(\Ric_{\hat g} g)_{ab})g_{jk},
\eeq
where $ F^n$ are the harmonicity functions given by
\beq\label{Harmonicity condition AAA}
F^n=\Gamma^n-\hat \Gamma^n,\quad\hbox{where }
\Gamma^n=g^{jk}\Gamma^n_{jk},\quad
\hat \Gamma^n=g^{jk}\hat \Gamma^n_{jk},
\eeq
where $\Gamma^n_{jk}$ and $\hat \Gamma^n_{jk}$ are the Christoffel symbols
for $g$ and $\hat g$, respectively.
The harmonicity functions $F^n$ of the solution $(g,\phi)$ of 
 the equations (\ref{eq: adaptive model with no source}) vanish when
 the conservation law (\ref{conservation law0}) is valid, see  \cite[eq.\ (14.8)]{Ringstrom}. Thus by (\ref{Reduced Ric tensor}), the conservation law
 implies that the solutions of the reduced Einstein equations (\ref{eq: EE in correct coordinates}) satisfy
 of the Einstein equations  (\ref{Einmat1}).


%

%

\subsubsection{Local existence of solutions}\label{subsec: Direct problem}

{\MATTITEXT
{Let  we consider the solutions $(g,\phi)$ of the equations  (\ref{eq: adaptive model with no source})
with  source $\F$.
To consider their local existence, let us  denote $u:=(g,\phi)-(\hat g,\hat \phi)$.

It follows from  by \cite[Cor.\ A.5.4]{BGP}
that $\K_j=J^+_{\tilde g}({p^-})\cap \overline M_j$ is compact.
Since $\hat  g<\tilde g$,
we see that if $r_0$ above is small enough,
for all $g\in \V(r_0)$, see subsection \ref{subsubsec: Gloabal hyperbolicity},
we have $g|_{\K_1}<\tilde g|_{\K_1}$. In particular, we have  $J^+_{g}(p^-)\cap \hattuM _1\subset 
J^+_{\tilde g}(p^-)$.  
%


%
%
Let us assume that 
 $\F$ is small enough in the norm  $C^4_b(M_0)$ and that it is  supported in a compact set  $\K=J_{\tilde g}({p^-})\cap  [0,t_0]\times N\subset \overline M_0$.
Then
 we can write the equations (\ref{eq: adaptive model with no source}) for $u$  in the form
 \beq\label{eq: notation for hyperbolic system 1}
& &P_{g(u)}(u)=\F,\quad x\in \hattuM _0,\\
& & \nonumber u=0\hbox{ in $(-\infty,0)\times N$, where}\\
\nonumber
& &P_{g(u)}(u):=g^{jk}(x;u)\p_j\p_k  u(x)+H(x,u(x),\p u(x)). 
\eeq
Here, the notation $g^{jk}(x;u)$ is used to indicate that the metric depends on the solution $u$.
More precisely, as the metric and the scalar field are $(g,\phi)=u+(\hat g,\hat \phi)$, we have
$(g^{jk}(x;u))_{j,k=1}^4=(g_{jk}(x))^{-1}$.
Moreover, above   $(x,v,w)\mapsto H(x,v,w)$ is a smooth function which is a second order polynomial in $w$ 
with coefficients being smooth functions of $v$, $\hat g$, and the derivatives of $\hat g$, \cite{Taylor3}.
{Note that when the norm of $\F$ in  $C^4_b(M_0)$  is small enough,  we have $\supp(u)\cap M_0\subset \K$.
We note that one could also consider non-compactly supported sources or initial data,
see \cite {CIP}. Also, the scalar field-Einstein system can be considered
with much less regularity that is done below, see \cite {CB-I-P2, CB-I-P3}.} 
\observation{
\medskip

{\bf Remark 3.1.}
We note that  the Laplace-Beltrami
operator can be written  for 
as $\square_g \psi_\ell =g^{jk}\p_j\p_k \phi_\ell-g^{pq} \Gamma^n \p_n \phi_\ell=
g^{jk}\p_j\p_k \phi_\ell- \Gamma^n \p_n \phi_\ell$
and thus  in the $(g,\hat g)$-wave map coordinates we have
\beq\label{eq: wave operator in wave gauge}
\square_g \phi_\ell=g^{jk}\p_j\p_k \phi_\ell-g^{pq}\hat \Gamma^n_{pq} \p_n \phi_\ell.
\eeq
Thus, the scalar field equation $\square_g \phi_\ell+m\phi_\ell=0$
does not involve derivatives of $g$.  As this can be the case in
 (\ref{eq: notation for hyperbolic system 1}),
the system (\ref{eq: notation for hyperbolic system 1}) is
a slight generalization of (\ref{eq: adaptive model with no source}).
\medskip
}



Let $s_0\geq 4$ be an even integer. 
Below we will consider the solutions $u=(g-\hat g,\phi-\hat \phi)$ and the sources $\F$ as sections
of the bundle $\B^L$ on $M_0$.
We will consider these functions as elements
of the section-valued Sobolev spaces $H^s(M_0;\B^L)$    etc.
Below, we
omit the bundle $\B^L$ in these notations and denote
 $H^s(M_0;\B^L)=H^s(M_0)$. We use the same convention
for the spaces
\ba
E^s=\bigcap_{j=0}^s C^j([0,t_0];H^{s-j}(N)),\quad s\in \N.
\ea
Note that $E^s\subset C^p([0,t_0]\times N)$ when $0\leq p<s-2$.
\MTEXT{Local existence results for (\ref {eq: notation for hyperbolic system 1})
follow from the standard techniques
for quasi-linear equations developed e.g.\ in 
 \cite{HKM} or  \cite{Kato1975},
or \cite[Section 9
]{Ringstrom}. These yield that 
when $\F$ is supported in the compact set $\K$ and $\|\F\|_{E^{s_0}}<c_0$, where
$c_0>0$  is  small  enough,  
there  exists  a unique function $u$ satisfying equation (\ref{eq: notation for hyperbolic system 1}) 
on $M_0$  with the source $\F$. Moreover, \noextension{$\|u\|_{{E^{s_0}}}\leq
C_1\|\F\|_{E^{s_0}}.$}
\extension{\beq\label{eq: Lip estim}
\|u\|_{{E^{s_0}}}\leq
C_1\|\F\|_{E^{s_0}}.
\eeq}
\noextension{For a detailed analysis, see Appendix B in \cite{preprint}.}\extension{(For details, see  Appendix B).} 
}

\subsubsection{Asymptotic expansion for  the non-linear wave equation}\label{subset: AENWE}


Let us  consider a small parameter $\e>0$ and the sources
 $\F=\F_\e$, depending smoothly on $\e\in [0,\e_0)$, with $\F_\e|_{\e=0}= 0$,  $\p_\e\F_\e|_{\e=0}= {\bf f}$ and $ {\bf f}=( {\bf f}^1, {\bf f}^2)$, for
 which the equations (\ref{eq: notation for hyperbolic system 1}) have
 a solution $u_\e$. Denote $\vec h=(h_{(j)})_{j=1}^4$, where $h_{(j)}=\p_\e^j\F_\e|_{\e=0}$
 so that ${\bf f}=h_{(1)}$.
Below, we always assume that $\F_\e$ is supported in $\K$ and
$\F_\e\in E^s$, where $s\geq s_0+10$ is an odd integer.
We consider the solution $u=u_\e$ of (\ref{eq: notation for hyperbolic system 1})
with $\F=\F_\e$
 and write it in the form
\beq\label{eq: epsilon expansion}
 u_\e(x)=\sum_{j=1}^4\e^j w^j(x)+w^{res}(x,\e).
 \eeq
%
To obtain the equations for $ w^j$, we use the
 representation    (\ref{Reduced Einstein tensor}) for the $\hat g$-reduced  Einstein tensor.
Below, we use the notation $ w^j=(( w^j)_{pq})_{p,q=1}^4,(( w^j)_\ell)_{\ell=1}^L)$
where $(( w^j)_{pq})_{p,q=1}^4$ is the $g$-component of $ w^j$ and 
$(( w^j)_\ell)_{\ell=1}^L$
is the $\phi$-component of $ w^j$.
\MTEXT{Below, we use also the notation 
where the components of $w=((g_{pq})_{p,q=1}^4,(\phi_\ell)_{\ell=1}^L)$ are re-enumerated so that $w$ is represented 
as a $(10+L)$-dimensional vector, i.e., we write
 $w=(w_m)_{m=1}^{10+L}$ (cf.\ Voigt notation).}
\extension{
Using formulas (\ref{q-formula2copyB}), (\ref{Reduced Ric tensor}), and (\ref{Reduced Einstein tensor}),
the solution  $u_\e(x)$ of the equation (\ref{eq: notation for hyperbolic system 1})
can be written in the form
\ba
u_\e&=& \e  w^1+\e^2  w^2+\e^3  w^3+\e^4  w^4+O(\e^5),\\
 w^1&=&{\bf Q}h_{(1)},\\
 w^2&=&{\bf Q}(A[ w^1, w^1])+{\bf Q}h_{(2)},\\
 w^3&=&2{\bf Q}(A[ w^1,{\bf Q}(A[ w^1, w^1])])+{\bf Q}(B[  w^1, w^1, w^1])+{\bf Q}h_{(3)},\\
 w^4
&=&{\bf Q}(A[ {\bf Q}(A[ w^1, w^1]),{\bf Q}(A[ w^1, w^1])\\
& &+4{\bf Q}(A[ w^1,{\bf Q}(A[ w^1,{\bf Q}(A[ w^1, w^1])])])+2{\bf Q}(A[ w^1,{\bf Q}(B[ w^1, w^1, w^1])])\\
& &+3{\bf Q}(B[ w^1, w^1,{\bf Q}(A[ w^1, w^1])])+{\bf Q}(C[ w^1, w^1, w^1, w^1])+{\bf Q}h_{(4)},
\ea
where ${\bf Q}={\bf Q}_{\hat g}=(\square_{\hat g} +V)^{-1}$ is the causal inverse of
the linearized Einstein equation,
$A$ is the notation for a generic 2nd  order multilinear operator in $u$
and  derivatives of order 2 obtained via pointwise multiplication of derivatives of $u$,  $B$ is a 3rd order multilinear operator in $u$
and  derivatives of order 2,  and
$C$ is a 4th order multilinear operator in $u$
and  derivatives of order 2, so that the sums of the orders of the
derivatives appearing in the terms $A,B$, and $C$ is at most two.
In a more explicit  way, $A$, $B$, and $C$ are of the form
\ba
& &A[v_1,v_2]=\sum_{|\a|+|\beta|\leq 2}a_{\a\b pq}(x)(\p_x^\a v^p_1(x))\,\cdotp(\p_x^\beta v^q_2(x)),\\
& &B[v_1,v_2,v_3]=\sum_{|\a|+|\beta|+|\a^{\prime}|\leq 2}b_{\a\beta\a^{\prime}pqr}(x)(\p_x^\a v_1^p(x))\,
\cdotp(\p_x^\beta v_2^q(x))
\,\cdotp(\p_x^{\a^{\prime}} v_3^r(x)),\\
& &C[v_1,v_2,v_3,v_4]=\\
& &\sum_{|\a|+|\beta|+|\a^{\prime}|+|\beta^{\prime}|\leq 2}c_{\a\beta\a^{\prime}\beta^{\prime}pqrs}(x)(\p_x^\a v_1^p(x))\,\cdotp(\D^\beta v_2^q(x))
\,\cdotp(\p_x^{\a^{\prime}} v^r_3(x))\,\cdotp(\p_x^{\beta^{\prime}} v^s_4(x)),
\ea
where the components of $v_k=((g_{ab}^{(k)}),\phi^\ell_{(k)})$ are represented in form $v_k=(v_k^p)_{p=1}^{10+L}$ (cf.\ Voigt notation).
 We emphasize that above the total order of derivatives
is always less or equal to two.
We can write the above formula without using multilinear forms:}
We have
 that $ w^j$, $j=1,2,3,4$ are
 given 
by
   \beq\nonumber 
 w^j&=&(g^j,\phi^j)={\bf Q}_{\hat g}\mathcal H^j,\quad j=1,2,3,4,\hbox{ where }\\
\nonumber
 \mathcal H^1&=&h_{(1)},\\ \nonumber
  \mathcal H^2&=&(2\hat g^{jp}w^1_{pq}\hat g^{qk}\p_j\p_k   w^1,0)
+\A^{(2)}( w^1,\p  w^1)+h_{(2)},\\
  \mathcal H^3&=& (\mathcal G_3,0)
  +\A^{(3)}( w^1,\p  w^1, w^2,\p  w^2)+h_{(3)},
 \label{w1-w3}
\\
\nonumber
\mathcal G_3&=&
-6\hat g^{jl}w^1_{li}\hat g^{ip}w^1_{pq}\hat g^{qk}\p_j\p_k   w^1+\\
& &\ \ +
3\hat g^{jp}w^2_{pq}\hat g^{qk}\p_j\p_k   w^1+3\hat g^{jp}w^1_{pq}\hat g^{qk}\p_k\p_j   w^2,
\nonumber
\eeq
and
\beq  \nonumber
\hspace{-10mm}  \mathcal H^4&=& (\mathcal G_4,0)
+\A^{(4)}( w^1,\p  w^1, w^2,\p  w^2, w^3,\p  w^3)+h_{(4)},\\
\nonumber
\hspace{-10mm} \mathcal G_4 &=&24\hat g^{js}w^1_{sr}\hat g^{rl}w^1_{li}\hat g^{ip} w^1_{pq}\hat g^{qk}\p_k\p_j   w^1
+6\hat g^{jp}w^2_{pq}\hat g^{qk}\p_k\p_j   w^2+\\
& &\ \  \label{w4}
-18\hat g^{jl}w^1_{li}\hat g^{ip} w^2_{pq}\hat g^{qk}\p_k\p_j   w^1
-12\hat g^{jl}w^1_{li}\hat g^{ip} w^1_{pq}\hat g^{qk}\p_k\p_j   w^2+
\hspace{-15mm}\\
& &\ \ +
3\hat g^{jp}w^3_{pq}\hat g^{qk}\p_k\p_j   w^1+
3\hat g^{jp}w^1_{pq}\hat g^{qk}\p_k\p_j   w^3.
  \nonumber
\eeq
 Moreover,
${\bf Q}_{\hat g}=(\square_{\hat g}+V(x,D))^{-1}$
is the causal inverse of the operator $\square_{\hat g}+V(x,D)$
where $V(x,D)$ is a first order differential operator with coefficients
depending on $\hat g$ and its derivatives
 and
 $\A^{(\a)}$, $\a=2,3,4$ denotes a sum of a multilinear operators  of orders $m$,
 $2\leq m\leq \a$ 
 having at a point $x$ the 
representation
\beq\label{eq: mixed source terms}
& &\hspace{1cm}(\A^{(\a)}(v^1,\p v^1,v^2,\p v^2,v^3,\p v^3))(x)
\\ &=& 
\nonumber
\sum_{} \bigg(a^{(\a)}_{abcijkP_1P_2P_3pq}(x)\,
 (v^1_a(x))^i(v^2_b(x))^j (v^3_c(x))^k \cdotp \\
 \nonumber
 & &\quad \quad\quad \cdotp P_1 ( \p v^1(x))P_2 ( \p  v^2(x)) P_3 ( \p v^3(x))  
 \bigg)\hspace{-1cm} 
\eeq
where  
$(v^1_a(x))^i$ denotes the $i$-th power of $a$-th component of $v^1(x)$
and  the sum is taken over the indexes $a,b,c,p,q,n,$ and integers
$ i,j,k$.
The homogeneous monomials $P_d(y)=y^{\beta_d}$,
$\beta_d=(b_1,b_2,\dots,b_{4(10+L)})\in \N^{4(10+L)}$, $d=1,2,3$
 having orders $|\beta_d|$, respectively,  
 where $P_d(y)=1$ for $d>\alpha$,
 and \beq
\label{Eq: F condition1a}
& &i+2j+3k+|\beta_1|+2|\beta_2|+3|\beta_3|
= \alpha,\\
& &\hbox{$|\beta_1|+|\beta_2|+|\beta_3|\leq 2$.}
\label{Eq: F condition3}\eeq
Here, by (\ref{Eq: F condition1a}), the term $\A^{(\a)}$ produces
a term of order $O(\e^\alpha)$ when $v^j=w^j$
and
condition (\ref{Eq: F condition3}) means
that  $\A^{(\a)}(v^1,\p v^1,v^2,\p v^2,v^3,\p v^3)$
contain only terms where the sum of the powers of derivatives
of $v^1,v^2$, and $v^3$ is at most two. 

 

 By \cite[App.\ III, Thm.\ 3.7]{ChBook}, or alternatively, the proof of \cite [Lemma 2.6]{HKM} adapted for manifolds, we see that
the estimate  $\| {\bf Q}_{\hat g}H\|_{E^{s_1+1}}\leq C_{s_1}\|H\|_{E^{s_1}}$
 holds  for all   
  $H\in E^{s_1}$, $s_1\in \Z_+$  
 that are supported in $\K_0=J^+_{\tilde g}(p^-)\cap \overline \hattuM _0$.
Note that we are interested only on the local solvability of the Einstein
equations. 
  
 Consider next an even integer $s\geq s_0+4$ and $\F_{ \e}\in E^{s}$
that depends smoothly on $ \e$.
Then, by defining $ w^j$ via the equations (\ref{w1-w3})-(\ref{w4}) with $h_{(j)}\in E^{s}$, $j\leq 4$ 
and using results of  \cite{HKM} we obtain that in (\ref{eq: epsilon expansion}) we have
%
$ w^j=\p_\e^j u_\e|_{\e=0}\in E^{s+2-j}$, $j=1,2,3,4$
and  
$  \|w^{res}(\,\cdotp,\e)\|_{E^{s-4}}\leq C\e^5$. 
%
%
%
\extension{
 Let us explain the details of this. Recall that $s\geq s_0+4$ is an even integer and $h_{(j)}\in E^{s}$.  First, we have $w^1={\bf Q}h_{(1)} \in  E^{s+1}$.  
 By defining $ w^j$ via the above equations with $h_{(j)}\in E^{s}$,
 we see that in the equation for $w^j$ we have non-linear terms
 where the first derivative of $w^j$ is multiplied by the first and zeroth 
 derivatives of $w^k$, $k<j$ and  
 non-linear terms
 where  $w^j$ is multiplied by the second, first, and zeroth 
 derivatives of $w^k$, $k<j$. Moreover, we obtain a source term
 where the second, first, and zeroth 
 derivatives of $w^k$, $k<j$ are multiplied by each other.
 In other words, denoting $W^j=(w^k)_{k=1}^j$,
 we have
 \ba
& & (\square_{\hat g}+V)w^j+A^n(W^j,\p W^j)\p_n w^j+\\
& & \quad  +
 B(W^j,\p W^j,\p^2W^j)w^j=H(W^j,\p W^j,\p^2W^j),\quad\hbox{in }M_0,\\
 & &\supp(w_j)\subset \K,
 \ea
 Using \cite [Thm.\ I]{HKM}, we see that the solution  exists in $E^{s_0}$.
 Let us next
 consider the regularity of the solution.
  Assuming that $W^j\in  E^{s+2-j+1}$, we have that 
 $H(W^j,\p W^j,\p^2W^j)\in  E^{s+2-j-1}$ and that the coefficient 
 $ B(W^j,\p W^j,\p^2W^j)$  of $w^j$ in the equation satisfy
 $ B(W^j,\p W^j,\p^2W^j)\in E^{s+2-j-1}$. Using the 
 fact that $\| {\bf Q}_{\hat g}H\|_{E^{s_1+1}}\leq C_{s_1}\|H\|_{E^{s_1}}$ with
 $s_1\leq s+1-j$, we see that $w^j\in E^{s+2-j}$. Starting from the case $j=1$ 
 we can repeat this argument for $j=2,3,4$.
 
 We have from above that 
  $ w^j\in E^{s+2-j}$,
for  $j=1,2,3,4$. Thus, by using Taylor expansion of the coefficients
  in the equation (\ref{eq: notation for hyperbolic system 1})
 we see that the approximate 4th order expansion
 $u^{app}_\e=\e  w^1+\e^2  w^2+\e^3  w^3+\e^4  w^4$
 satisfies an equation of the form
 \beq\label{eq: notation for hyperbolic system 1b}
& &P_{g(u^{app}_\e)}(u^{app}_\e)
\hspace{-1mm}=
\F(\e) \hspace{-1mm}+\hspace{-1mm}H^{res}(\,\cdotp,\e),\ \ x\in \hattuM _0,
\hspace{-5mm}\\
& & \nonumber \supp(u^{app}_\e)\subset \K,
\eeq
%
%
%
such that 
 $$\|H^{res}(\,\cdotp,\e)\|_{E^{s-4}}
 \leq c_1\e^5.  $$ 
 Using the Lipschitz continuity of the solution
 of the equation (\ref{eq: notation for hyperbolic system 1b}) with respect to
 the source term, see Appendix B,
 we have that
 there are $c_1,c_2>0$ such that for all $0<\e<c_1$ the function
 $w^{res}(x,\e)=u_\e(x)-u^{app}_\e(x)$ satisfies 
  \ba
 \|u_\e-u^{app}_\e\|_{E^{s-4}}
 \leq c_2\e^5.
 \ea
Thus in (\ref{eq: epsilon expansion}) we have
$  \|w^{res}(\,\cdotp,\e)\|_{E^{s-4}}\leq C\e^5$ and
$ w^j=\p_\e^j u_\e|_{\e=0}\in E^{s-4}$, $j=1,2,3,4.$
%
%
%
}

 \subsection{Distorted plane wave solutions for the linearized equations}
 \subsubsection{Lagrangian distributions}\label{Lagrangian distributions}
Let us recall the definition of the conormal and Lagrangian distributions that we will use below.
Let $X$ be a manifold of dimension $n$
and $\Lambda\subset T^*X\setminus \{0\}$ be
a Lagrangian submanifold.
Let  
$\phi(x,\theta)$, $\theta\in \R^N$ be a non-degenerate
phase function that locally parametrizes $\Lambda$.
We say that a distribution $u\in {\cal D}^{\prime}(X)$  
is a Lagrangian distribution associated with $\Lambda$ and denote
$u\in \I^m(X;\Lambda)$, if in local coordinates $u$ can be represented 
 as an oscillatory integral,
\beq\label{def: conormal}
u(x)=\int_{\R^N}e^{i\phi(x,\theta)}a(x,\theta)\,d\theta,
\eeq
where $a(x,\theta)\in S^{m+n/4-N/2}(X; \R^N)$,
see  \cite{GU1,H4,MU1}.

In particular, when $S\subset X$ is a submanifold,  its conormal bundle
$N^*S=\{(x,\xi)\in T^*X\setminus \{0\};\ x\in S,\ \xi\perp T_xS\}$ is a Lagrangian submanifold.
If $u$ is a  Lagrangian distribution  associated with $\Lambda_1$
where $\Lambda_1=N^*S$, we say that $u$ is a conormal distribution.

Let us next consider the case when $X=\R^n$ and
let
$(x^1,x^2,\dots,x^n)=(x^{\prime},x^{\prime\prime},x^{\prime\prime\prime})$ 
be the Euclidean coordinates with $x^{\prime}=(x_1,\dots,x_{d_1})$,
  $x^{\prime\prime}=(x_{d_1+1},\dots,x_{d_1+d_2})$,
 $x^{\prime\prime\prime}=(x_{d_1+d_2+1},\dots,x_{n})$. 
 If $S_1=\{x^{\prime}=0\}\subset \R^n$, $\Lambda_1=N^*S_1$
then $u\in \I^m(X;\Lambda_1)$ can be represented by (\ref{def: conormal})
with $N=d_1$ and  $\phi(x,\theta)=x^{\prime}\cdotp \theta$.
 
Next we recall the definition of $\I^{p,l}(X;\Lambda_1,\Lambda_2)$, the space of
the distributions $u$ in ${\cal D}^{\prime}(X)$ associated to two
cleanly intersecting Lagrangian manifolds $\Lambda_1,\Lambda_2\subset T^*X
\setminus \{0\}$, see \cite {CGP,GU1,MU1}. These classes have been widely
used in the study of inverse problems, see \cite {CS2,FG}.  
Let us start with the case when $X=\R^n$.

Let $S_1,S_2\subset \R^n$ be the linear subspaces of codimensions
$d_1$ and $d_1+d_2$, respectively,
$S_2\subset S_1$, given by
$S_1=\{x^{\prime}=0\}$, $S_2=\{x^{\prime}=x^{\prime\prime}=0\}$.
Let us denote $\Lambda_1=N^*S_1,$  $\Lambda_2=N^*S_2$. Then
$u\in \I^{p,l}(\R^n;N^*S_1,N^*S_2)$ if and only if
\ba
u(x)=\int_{\R^{d_1+d_2}}e^{i (x^{\prime}\cdotp \theta^{\prime}+
x^{\prime\prime}\cdotp \theta^{\prime\prime})}a(x,\theta^{\prime},\theta^{\prime\prime})\,d\theta^{\prime}d\theta^{\prime\prime},
\ea
where the symbol $
a(x,\theta^{\prime},\theta^{\prime\prime})$ belongs in the product type symbol class
$S^{\mu_1,\mu_2}(\R^n; (\R^{d_1}\setminus 0)\times \R^{d_2})$
that is the space of function  $a\in C^{\infty}(\R^n\times \R^{d_1}\times \R^{d_2})$ that satisfy
\beq\label{product symbols}
|\p_x^\gamma\p_{\theta^{\prime}}^\alpha \p_{\theta^{\prime\prime}}^\beta a(x,\theta^{\prime},\theta^{\prime\prime})|
\leq C_{\alpha\beta\gamma K}(1+|\theta^{\prime}|+|\theta^{\prime\prime}|)^{\mu_1-|\alpha|}
(1+|\theta^{\prime\prime}|)^{\mu_2-|\beta|}\hspace{-1.7cm}
\eeq
for all $x\in K$, multi-indexes  $\alpha,\beta,\gamma$, and
compact sets $K \subset \R^n$. Above,
$\mu_1=p+l-d_1/2+n/4$ and
$\mu_2=-l-d_2/2$.

When $X$ is a manifold of dimension $n$ and 
$\Lambda_1,\Lambda_2\subset T^*X
\setminus \{0\}$ are two
cleanly intersecting Lagrangian manifolds, 
we define the class $\I^{p,l}(X;\Lambda_1,\Lambda_2)\subset \mathcal D^\prime(X)$ to consist
of locally finite sums of distributions of the form $u=Au_0$, where 
$u_0\in \I^{p,l}(\R^n;N^*S_1,N^*S_2)$ and $S_1,S_2\subset \R^n$ are the linear subspace of codimensions
$d_1$ and $d_1+d_2$, respectively, such that
$S_2\subset S_1$, 
and $A$ is a Fourier integral
operator of order zero with a canonical relation $\Sigma$ for which
$ \Sigma\circ (N^*S_1)^\prime\subset \Lambda_1^\prime$ and 
$\Sigma\circ (N^*S_2)^\prime\subset \Lambda_2^\prime$. 
Here, for $\Lambda\subset T^*X$ we denote $ \Lambda^\prime=\{(x,-\xi)\in T^*X;\ (x,\xi)\in \Lambda\}$,
and for $\Sigma\subset T^*X\times  T^*X$ we denote 
$\Sigma^\prime=\{(x,\xi,y,-\eta);\ (x,\xi,y,\eta)\in \Sigma\}.$

In most cases, below $X=M$. We denote then
  $\I^{p}(M;\Lambda_1)=\I^{p}(\Lambda_1)$
and 
 $\I^{p,l}(M;\Lambda_1,\Lambda_2)=\I^{p,l}(\Lambda_1,\Lambda_2)$.
Also, $\I(\Lambda_1)=\cup_{p\in \R}\I^{p}(\Lambda_1)$.

By \cite{GU1,MU1}, microlocally away from $\Lambda_1$ and $\Lambda_0$,
\beq\label{microlocally away}
& &\I^{p,l}(\Lambda_0,\Lambda_1)\subset \I^{p+l}(\Lambda_0\setminus \Lambda_1)\quad
\hbox{and}\quad \I^{p,l}(\Lambda_0,\Lambda_1)\subset \I^{p}(\Lambda_1\setminus \Lambda_0),
\eeq
respectively.
Thus the principal symbol of $u\in \I^{p,l}(\Lambda_0,\Lambda_1)$
is well defined on $\Lambda_0\setminus \Lambda_1$ and  $\Lambda_1\setminus \Lambda_0$. We denote  $\I(\Lambda_0,\Lambda_1)=\cup_{p,q\in \R}\I^{p,q}(\Lambda_0,\Lambda_1)$.

Below, when $\Lambda_j=N^*S_j,$ $j=1,2$ are conormal bundles of smooth cleanly
intersecting submanifolds
$S_j\subset M$ of codimension $m_j$, where $\dim(M)=n$,
 we
use the traditional notations,
\beq\label{eq: traditional}\quad\quad
\I^\mu(S_1)=\I^{\mu+m_1/2-n/4}(N^*S_1),\quad \I^{\mu_1,\mu_2}(S_1,S_2)
=\I^{p,l}(N^*S_1,N^*S_2),\hspace{-1cm}
\eeq
where $p=\mu_1+\mu_2+m_1/2-n/4$ 
 and $l=-\mu_2-m_2/2$,
and call such distributions the conormal distributions  associated to $S_1$ or 
product type  conormal distributions  associated to
$S_1$ and $S_2$, respectively.  {By \cite{GU1},
$\I^\mu(X;S_1)\subset L^p_{loc}(X)$ for $\mu< -m_1(p-1)/p$, $1\leq p<\infty$.}

For the
wave operator $\square_{g}$ on the globally hyperbolic manifold $(M,{g})$,
$\hbox{Char}\,(\square_{g})$ is the set of light-like co-vectors with respect to $g$.
For $(x,\xi)\in \hbox{Char}\,(\square_{g})$, $\Theta_{x,\xi}$
denotes the bicharacteristic of $\square_{g}$. Then,
$(y,\eta)\in \Theta_{x,\xi}$ if and only if there is $t\in\R$ such
that
for $a=\eta^\sharp$ and $b=\xi^\sharp$ we have
 $(y,a)=(\gamma^{g}_{x,b}(t),\dot \gamma^{g}_{x,b}(t))$
where $\gamma^{g}_{x,b}$ is a light-like geodesic with respect to the metric $g$
with the initial data $(x,b)\in LM$. Here, we use notations $(\xi^\sharp)^j=g^{jk} \xi_k$ and $(b^\sharp)_j=g_{jk} b^k$.

Let $P=\square_{g}+B^0+B^j\p_j$, where $B^0$  is a scalar function and $B_j$ is a vector
field.
Then $P$ is a classical pseudodifferential operator of real
principal type and order $m=2$  on $M$,
and  \cite{MU1}, see also \cite{K2},
$P$ has a parametrix $Q\in \I^{p,l}(\Delta^\prime_{T^*M},\Lambda_P)$, 
$p=\frac 12-m$, $l=-\frac 12$ ,
where $\Delta_{T^*M}=N^*(\{(x,x);\ x\in M\})$ and $\Lambda_g\subset T^*M\times T^*M$
is the Lagrangian manifold associated to the canonical relation of  the operator $P$,
that is, 
\beq\label{eq: Lambda g}
\Lambda_g=\{(x,\xi,y,-\eta);\ (x,\xi)\in\hbox{Char}\,(P),\ (y,\eta)\in \Theta_{x,\xi}\},
\eeq
where $\Theta_{x,\xi}\subset T^*M$ is the bicharacteristic of $P$ containing $(x,\xi)$.  
 When
 $(M,g)$ is a globally hyperbolic manifold, the operator  
$P$ has a causal inverse operator, see e.g.\  \cite[Thm.\ 3.2.11]{BGP}.
We denote it by  $P^{-1}$ and by \cite{MU1}, we have
$P^{-1} \in \I^{-3/2,-1/2}(\Delta^\prime_{T^*M},\Lambda_g)$.
We will repeatedly use the fact (see \cite[Prop. 2.1]{GU1}) that 
if $F\in \I^{p}(\Lambda_0)$ and $\Lambda_0$ intersects Char$(P)$ transversally
so that all bicharacterestics of $P$ intersect $\Lambda_0$ only finitely many
times, then $(\square_{g}+B^0+B^j\p_j )^{-1}F\in \I^{p-3/2,-1/2}(\Lambda_0,\Lambda_1)$
where $\Lambda_1^{\prime}=\Lambda_g\circ \Lambda_0^{\prime}$ is called the flowout from $\Lambda_0$ on Char$(P)$,
that is,
\ba
\Lambda_1=\{(x,-\xi);\ 
(x,\xi,y,-\eta)\in \Lambda_g,\ (y,\eta)\in \Lambda_0 \}.
\ea


 \subsubsection{The linearized Einstein equations  and 
 the linearized conservation law
 }

We will below consider  sources $\F=\e {\bf f}(x)$ and  solution $u_\e$ satisfying
(\ref{eq: notation for hyperbolic system 1}),
where ${\bf f}=({\bf f}^{(1)},{\bf f} ^{(2)})$.

We consider the linearized Einstein equations  and the 
linearized wave $w^1=\p_\e u_\e|_{\e=0}$ in (\ref{eq: epsilon expansion}) that we
denote by $u^{(1)}=w^1.$ It satisfies  
 the linearized Einstein equations  (\ref{linearized eq: adaptive model with no source}) that we write as
\beq\label{lin wave eq}
& &\square_{\hat g}  u^{(1)}+V(x,\p_x) u^{(1)}={\bf f},
\eeq
where $v\mapsto V(x,\p_x)v$ is a linear first order partial differential operator with
coefficients depending on  $\hat g$ and its derivatives.

Assume that $Y\subset \hattuM _0$
is a 2-dimensional space-like submanifold and consider local coordinates defined 
in  $V\subset \hattuM _0$. Moreover, assume that 
in these local coordinates $Y\cap V\subset \{x\in \R^4;\ x^jb_j=0,\  x^jb^\prime_j=0\}$,
where $b^\prime_j\in \R$  and let ${\bf f}=({\bf f}^{(1)},{\bf f}^{(2)})\in \I^{n+1}(Y)$, $n\leq n_0=-17$,  be defined by 
\beq\label{eq: b b-prime AAA}
{\bf f}(x^1,x^2,x^3,x^4)=
\re \int_{\R^2}e^{i(\theta_1b_m+\theta_2b^\prime_m)x^m}\sigma_{\bf f}(x,\theta_1,\theta_2)\,
d\theta_1d\theta_2.\hspace{-2cm}
\eeq
\MTEXT{Here, we assume that 
 $\sigma_{\bf f}(x,\theta)$, $\theta=(\theta_1,\theta_2)$ 
is a $\B^L$-valued classical symbol and we denote  the principal symbol of ${\bf f}$ by
 $c(x,\theta)$, or component-wise,
$((c^{(1)}_{jk}(x,\theta))_{j,k=1}^4,(c^{(2)}_{\ell}(x,\theta))_{\ell=1}^L)$. When 
$x\in Y$ and $\xi=(\theta_1b_m+\theta_2b^\prime_m)dx^m$ so that $(x,\xi)\in N^*Y$,
we denote the value of the principal symbol of ${\bf f}$ at $(x,\xi)$ by
$\tilde c(x,\xi)=c(x,\theta)$, that is component-wise,
$\tilde c^{(1)}_{jk}(x,\xi)=c^{(1)}_{jk}(x,\theta)$
 and $\tilde c^{(2)}_{\ell}(x,\xi)=c^{(2)}_{\ell}(x,\theta)$.
We say that this is the principal symbol} of ${\bf f}$ at $(x,\xi)$, associated to the phase function $\phi(x,\theta_1,\theta_2)=(\theta_1b_m+\theta_2b^\prime_m)x^m$.
The above defined principal symbols can be defined invariantly, see \cite{GuU1}. 

%

%
 
We will below consider what happens when  ${\bf f}=({\bf f}^{(1)},{\bf f}^{(2)})\in \I^{n+1}(Y)$ satisfies 
the {\it linearized conservation
law} (\ref{eq: lineariz. conse. law PRE}). {\motivation 
Roughly speaking, these four linear conditions
imply that  the principal symbol of the source ${\bf f}$ satisfies four linear conditions.
Furthermore, the  linearized conservation law implies that also
the linearized wave $u^{(1)}$ produced by  ${\bf f}$
satisfies four linear conditions that we call the linearized harmonicity conditions,
and finally, the principal symbol of the wave $u^{(1)}$ has to satisfy four linear conditions.
Next we explain these conditions in detail.}

When  (\ref{eq: lineariz. conse. law PRE})  is valid, we have
\beq\label{eq: lineariz. conse. law symbols}
& &\hat g^{lk}\xi_l\tilde c^{(1)}_{kj} (x,\xi)
=0,\quad\hbox{for 
$j\leq 4$ and $\xi\in N^*_x Y$}.
\eeq
We say that this is the {\it linearized conservation law for the principal symbols}.
\extension{Note that $\I^{\mu}(Y)\subset C^s(M_0)$ when 
$s\leq -\mu-3$. We will later use such indexes $\mu$ so that
we can use $s=13$, cf.\ Assumption $\mu$-LS.}

\subsubsection{The harmonicity condition for the linearized solutions}
%

Assume that $(g,\phi)$ satisfy equations (\ref{eq: adaptive model with no source})
and  the conservation law (\ref{conservation law0})
is valid. 
The conservation law (\ref{conservation law0}) and 
the $\hat g$-reduced Einstein equations (\ref{eq: adaptive model with no source}) imply,
see e.g.\ \cite{ChBook,Ringstrom}, 
 that the harmonicity functions $\Gamma^j=g^{nm} \Gamma^j_{nm}$
 satisfy 
 \beq\label{harmonicity condition}
g^{nm} \Gamma^j_{nm}= g^{nm}\hat \Gamma^j_{nm}.
\eeq
Next we 
denote $u^{(1)}=(g^1,\phi^1)=(\dot g,\dot\phi)$, see (\ref{w1-w3}), 
and discuss the implications  of (\ref{harmonicity condition}) for the metric
component $\dot g$ of the solution of
the linearized Einstein equations. 

}

}

\extension{
We  do next calculations  in local coordinates of $M_0$ and 
denote $\p_k=\frac\p{\p x^k}$.
Direct calculations show that 
$
h^{jk}=g^{jk}\sqrt{-\det(g)}
$
satisfies $\p_k h^{kq}= - \Gamma^q_{kn}h^{nk}$.
Then (\ref{harmonicity condition})
implies that 
\beq\label{harmonic condition - alternative2}
\p_k h^{kq}=
-\hat \Gamma^q_{kn}h^{nk}.
\eeq
We call (\ref{harmonic condition - alternative2}) the {\it  harmonicity condition} for  the metric $g$.

Assume now that  $g_\e$ and $\phi_\e$ satisfy (\ref{eq: adaptive model with no source})
with source $\F=\e f$
 where $\e>0$
 is a small parameter.
 We define $h_\e^{jk}=g_\e^{jk}\sqrt{-\det(g_\e)}$
 and denote $\dot g_{jk}=\p_\e (g_\e)_{jk}|_{\e=0}$,
$\dot g^{jk}=\p_\e (g_\e)^{jk}|_{\e=0}$,
and $\dot h^{jk}=\p_\e h_\e^{jk}|_{\e=0}$.

%
%
%

The equation (\ref{harmonic condition - alternative2})
yields then\footnote{The treatment on this de Donder-type gauge condition is known in
the folklore of the field. For a similar gauge condition to (\ref{harmonic condition for dots1}) 
in harmonic coordinates, see \cite[pages 6 and 250]{Maggiore},  or  
\cite[formulas 107.5, 108.7, 108.8]{Landau}, or \cite[p.\ 229-230]{HE}.
}
\beq\label{harmonic condition for dots1}
\p_k\dot h^{kq}=-\hat \Gamma^q_{kn}\dot h^{nk}.
\eeq
A direct computation  shows that 
\ba
\dot h^{ab}
&=&(-\det(\hat g))^{1/2}\kappa^{ab},
\ea
where $\kappa^{ab}=\dot g^{ab}-\frac 12 \hat g^{ab}\hat g_{qp}\dot g^{pq}$.
Thus (\ref{harmonic condition for dots1}) gives
\beq\label{extra 2a}
\p_a((-\det(\hat g))^{1/2}\kappa^{ab})
=
-\hat \Gamma^b_{ac}(-\det(\hat g))^{1/2}\kappa^{ac}
\eeq
that implies 
$
\p_a\kappa^{ab}
+\kappa^{nb}\hat\Gamma^a_{an}+\kappa^{an}\hat \Gamma^b_{an}=0,
$
or equivalently, 
\beq\label{extra 2b}
\hat \nabla_a\kappa^{ab}=0.
\eeq
We call (\ref{extra 2b}) the {\it  linearized harmonicity condition} for $g$. 
Writing this for $\dot g$, we obtain
\beq\label{harmonicity condition AAA}
& &\hspace{-1cm}-\hat g^{an}\p_a \dot g_{nj}+\frac 12 \hat g^{pq}  \p_j\dot g_{pq}=m^{pq}_j\dot g_{pq}
\eeq
where $m_j$  depend on $\hat g_{pq}$ and its derivatives.
On similar conditions for the polarization tensor, see \cite[form.\ (9.58) and example  9.5.a, p. 416]{Padmanabhan}.
}
\subsubsection{Properties of the principal symbols of the waves}

Let $K\subset \hattuM _0$ be a light-like submanifold of dimension 3 that in
 local coordinates  $X:V\to \R^4$, $x^k=X^k(y)$ 
is given by $K\cap V\subset \{x\in \R^4;\ b_kx^k=0\}$, where $b_k\in \R$ are constants.
Assume that the solution $u^{(1)}=(\dot g,\dot \phi)$ of the 
linear wave equation (\ref{lin wave eq}) with the right hand
side vanishing  in $V$ is such that  $u^{(1)}\in \I^\mu (K)$ with $\mu\in \R$.
Below we use $\mu=n-\frac 12$  where $n\in \Z_-$, $n\leq n_0=-17.$
Let us write    $\dot g_{jk}$  as an oscillatory integral using 
 a phase function $\varphi(x,\theta)=b_kx^k\theta$, and
a classical symbol $\sigma_{\dot g_{jk}}(x,\theta)\in S^n_{cl} (\R^4,\R)$, 
\beq\label{eq: oscil}
\dot g_{jk}(x^1,x^2,x^3,x^4)=
\re \int_{\R}e^{i(\theta b_mx^m)}\sigma_{\dot g_{jk}}(x,\theta)\,
d\theta,
\eeq
where $n=\mu+\frac 12$. 
We denote
the  (positively homogeneous) principal symbol of $\dot g_{jk}$ by 
$a_{jk}(x,\theta)$. 
When $x\in K$ and $\xi=\theta b_kdx^k$ so that $(x,\xi)\in N^*K$,
we denote the value of $a_{jk}$ at $(x,\theta)$ by 
 $\tilde a_{jk}(x,\xi)$, that is,  $\tilde a_{jk}(x,\xi)=a_{jk}(x,\theta)$.

%
%
Then, 
if $\dot  g_{jk}$ satisfies the linearized harmonicity condition (\ref{harmonicity condition}),
its principal symbol $\tilde a_{jk}(x,\xi)$ satisfies 
\beq\label{harmonicity condition for symbol}
& &\hspace{-1cm}-\hat g^{mn}(x)\xi_m v_{nj}+\frac 12  \xi_j(\hat g^{pq}(x) v_{pq})=0,
\quad v_{pq}=\tilde a_{pq}(x,\xi),
\eeq
where $j=1,2,3,4$ and  $\xi=\theta b_kdx^k\in N_x^*K$.
 If (\ref{harmonicity condition for symbol}) holds, we say that the {\it  harmonicity condition for the symbol} is satisfied for $ \tilde a(x,\xi)$
at $(x,\xi)\in N^*K$.

\subsubsection{Distorted plane waves satisfying a linear wave equation}


 Next we consider a  
 distorted plane  wave whose singular support is
concentrated near a  geodesic. \motivation{These waves, sketched in Fig.\ 1(Right),
propagate near the geodesic $\gamma_{x_0,\zeta_0}([t_0,\infty))$
and are singular on a surface $K(x_0,\zeta_0;t_0,s_0)$, defined below in (\ref{associated submanifold}), that is a subset
of the light cone $\L^+_{\hat g}(x^\prime)$, $x^{\prime}=\gamma_{x_0,\xi_0}(t_0)$. The parameter $s_0$ gives a ``width'' of
the wave packet and when $s_0\to 0$, its singular support
tends to the set $\gamma_{x_0,\zeta_0}([2t_0,\infty))$.
Next we will define these wave packets.}

%
%
%

We define the 3-submanifold
$K(x_0,\zeta_0;t_0,s_0)\subset M_0$  associated to  $(x_0,\zeta_0)\in L^{+}(M_0,\hat g)$, $x_0\in U_{\hat g}$ and  parameters $t_0,s_0\in \R_+ $
as 
\beq\label{associated submanifold}
& &K(x_0,\zeta_0;t_0,s_0)=\{\gamma_{x^{\prime},\eta}(t)\in M_0;\ \eta\in \W,\ t\in (0,\infty)\},\hspace{-1cm}
\eeq
where  $(x^{\prime},\zeta^{\prime})=(\gamma_{x_0,\zeta_0}(t_0),\dot \gamma_{x_0,\zeta_0}(t_0))$
and $\W\subset L^{+}_{x^{\prime}}(M_0,\hat g)$ is a  neighborhood of $\zeta^{\prime}$
consisting of vectors $\eta\in L^+_{x^{\prime}}(M_0)$ satisfying
$
\left \|\eta -\zeta^{\prime}\right \|_{\hat g^+}<s_0.
$
Note that $K(x_0,\zeta_0;t_0,s_0)\subset \L^+_{\hat g}(x^{\prime})$ is a subset
of the light cone starting at  $x^{\prime}=\gamma_{x_0,\xi_0}(t_0)$ and that it 
is singular
at  the point $x^{\prime}$. 
Let 
$S=\{x\in \hattuM _0;{\bf t}(x)={\bf t}(  \gamma_{x_0,\zeta_0}(2t_0) )\}$  
be a Cauchy surface which intersects 
$\gamma_{x_0,\zeta_0}(\R)$ transversally at the point $ \gamma_{x_0,\zeta_0}(2t_0)$.
When $t_0>0$  is  small enough,
$Y(x_0,\zeta_0;t_0,s_0)=S\cap K(x_0,\zeta_0;t_0,s_0)$ is a smooth 2-dimensional space-like surface that is
a subset of $U_{\hat g}$.

Let $\Lambda(x_0,\zeta_0;t_0,s_0)$ be the Lagrangian manifold that is the  flowout from  
$N^*Y(x_0,\zeta_0;t_0,s_0)\cap  N^*K(x_0,\zeta_0;t_0,s_0)$ on Char$(\square_{\hat g})$ in the future direction. \modified{When $K^{reg}\subset K= K(x_0,\zeta_0;t_0,s_0)$
is the set of points $x$ that have a neighborhood $W$  such that $K\cap W$ is   a smooth 3-dimensional submanifold,
we have $N^*K^{reg}\subset \Lambda(x_0,\zeta_0;t_0,s_0).$}
\MTEXT{Below,  we represent locally the elements $w\in \B_x$ in the fiber of the bundle
$\B$ 
as a $(10+L)$-dimensional vector,
 $w=(w_m)_{m=1}^{10+L}$.}

\begin{lemma}\label{lem: Lagrangian 1} 
Let $n\leq n_0=-17$ be an integer,  $t_0,s_0>0$, $Y=Y(x_0,\zeta_0;t_0,s_0)$,
$K=K(x_0,\zeta_0;t_0,s_0)$, $\Lambda_1=\Lambda(x_0,\zeta_0;t_0,s_0)$, and $(y,\xi)\in N^*Y\cap \Lambda_1$.
Assume that ${\bf f}=({\bf f}_1,{\bf f}_2)\in \I^{n+1}(Y)$,
  is a $\B^L$-valued conormal
distribution  that is supported in  a neighborhood 
$V\subset \hattuM _0$ of  $\gamma_{x_0,\zeta_0}\cap Y=\{\gamma_{x_0,\zeta_0}(2t_0)\}$ and has  a $\R^{10+L}$-valued 
classical symbol.
Denote 
the principal symbol of ${\bf f}$ by 
${{\tilde f}}(y,\xi)=({{\tilde f}}_k(y,\xi))_{k=1}^{10+L}
$,
and assume that the symbol of ${\bf f}$ vanishes near
the  light-like directions
in $N^*Y\setminus N^*K$.

Let $u^{(1)}=(\dot g,\dot\phi)$ be a solution of the 
linear wave equation (\ref{lin wave eq})
with the source ${\bf f}$. Then $u^{(1)}$,
considered as a vector valued Lagrangian distribution on the set $M_0\setminus Y$,  satisfies
$
u^{(1)}\in \I^{n-1/2} ( M_0\setminus Y;\Lambda_1), 
$
and its principal symbol  
$\tilde a(x,\eta)=(\tilde a_{j}(x,\eta))_{j=1}^{10+L}$ at $(x,\eta)\in  \Lambda_1$,  is given by 
\beq\label{eq: R propagation}\tilde a_j(x,\eta)=\sum_{k=1}^{10+L} R_j^k(x,\eta,y,\xi){{\tilde f}}_{k}(y,\xi)
,\hspace{-1cm}
\eeq
where the pairs $(y,\xi)$ and $(x,\eta)$ are on the same bicharacteristics of $\square_{\hat g}$, 
and $y\ll x$. Observe that  $((y,\xi),(x,\eta))\in \Lambda_{\hat g}^\prime$, and in addition,
$(y,\xi)\in N^*Y\cap N^*K$.
%
 Moreover, the matrix $(R_j^k(x,\eta,y,\xi))_{j,k=1}^{10+L}$
is invertible. 
\end{lemma}
We call the solution $u^{(1)}$   a distorted plane  wave that is associated
to the submanifold $K(x_0,\zeta_0;t_0,s_0)$. 
%

{\bf Proof.} 
It follows from  \cite{MU1}
that the causal inverse  of the scalar wave operator  $\square_{\hat g}+V(x,D)$,
where $V(x,D)$ is a 1st order differential operator, 
satisfies $(\square_{\hat g}+V(x,D))^{-1}\in
 \I^{-3/2,-1/2}(\Delta^\prime_{T^*M_0},\Lambda_{\hat g})$. Here,
$\Delta_{T^*M_0}$ is the conormal bundle of the diagonal of $M_0\times M_0$
and $\Lambda_{\hat g}$ is the flow-out of the canonical relation of $\square_{\hat g}$.
A geometric representation for its kernel is given in \cite{K2}.
An analogous result holds for
%
%
%
%
%
%
the matrix valued wave operator, $\square_{\hat g}I+V(x,D)$,
when  $V(x,D)$ is a 1st order differential operator, that is,
 $(\square_{\hat g}I+V(x,D))^{-1}\in \I^{-3/2,-1/2}(\Delta^\prime_{T^*M_0},\Lambda_{\hat g})$,
 see  \cite{MU1} and \cite{Dencker}.
By \cite[Prop.\ 2.1]{GU1}, this yields  $u^{(1)}\in \I^{n-1/2}(\Lambda_1)$
and the formula (\ref{eq: R propagation}) where
$R=(R_j^k(x,\eta,y,\xi))_{j,k=1}^{10+L}$ is obtained by solving a system
of ordinary differential equation along a bicharacteristic curve. Making  similar considerations
for the adjoint of the $(\square_{\hat g}I+V(x,D))^{-1}$, i.e., considering
the propagation of singularities using reversed causality, we see that
the matrix $R$ is invertible.
 \hfill \Box \medskip

\observation{{\bf Remark 3.2.}
Let us write the above principal symbol
${{\tilde f}}(y,\eta)\in \R^{10+L}$ as
${{\tilde f}}(y,\eta)=({{\tilde f}}^{\prime}(y,\eta),{{\tilde f}}^{\prime\prime}(y,\eta))\in \R^{10}\times
\R^L$, where  ${{\tilde f}}^{\prime}(y,\eta)$ corresponds
to the principal symbol of the  $g$-component of the wave
and ${{\tilde f}}^{\prime\prime}(y,\eta)$ the  $\phi$-component of the wave.
Recall that due to (\ref{eq: wave operator in wave gauge}),
the equations $\square_g \phi_\ell+m\phi_\ell=0$
does not involve derivatives of $g$. 
This implies that the  
matrix $R(y,\eta,x,\xi)=[R_j^k(y,\eta,x,\xi)]_{j,k=1}^{L+10}$ in (\ref{eq: R propagation}) is such that
if ${{\tilde f}}(y,\eta)=({{\tilde f}}^{\prime}(y,\eta),{{\tilde f}}^{\prime\prime}(y,\eta))$
is such that ${{\tilde f}}^{\prime\prime}(y,\eta)=0$, then
\beq
\label{new eq: R propagation}
\tilde a(y,\eta)=\left(\begin{array}{c}
\tilde a^{\prime}(y,\eta)\\
0
\end{array}\right)
=
R(y,\eta,x,\xi)\left(\begin{array}{c}
{\tilde f}^{\prime}(y,\eta) \\
{{\tilde f}}^{\prime\prime}(y,\eta)
\end{array}\right)
,\hspace{-1cm}
\eeq
that is, $\tilde a^{\prime\prime}(y,\eta)=0$. Roughly speaking,
this means that during the propagation along a bicharacteristic, 
the leading order singularities can
not move from the $g$-component of the wave to $\phi$-component.
\medskip
}

Below, let $(y,\xi)\in  N^*Y\cap \Lambda_1$ and $(x,\eta)\in T^*M_0$
be a light-like co-vector such that $(x,\eta)\in\Theta_{y,\xi}$, $x\not \in Y$ and, $y\ll x$.

Let $\B^L_x$  be the fiber of the bundle $\B^L$ at $x$ and 
 $\mathfrak S_{x,\eta}$ be the space of 
the elements in $\B^L_x$  satisfying
the harmonicity condition for the symbols (\ref{harmonicity condition for symbol})
at $({x,\eta})$.
{Let $(y,\xi)\in N^*Y$ and
$\mathfrak C_{y,\xi}$ be the set of elements $b$ in $\B^L_y$
that satisfy the  linearized conservation law for symbols, i.e., 
(\ref{eq: lineariz. conse. law symbols}).  

Let $n\leq n_0$ and $t_0,s_0>0$, $Y=Y(x_0,\zeta_0;t_0,s_0)$,
$K=K(x_0,\zeta_0;t_0,s_0)$, $\Lambda_1=\Lambda(x_0,\zeta_0;t_0,s_0)$,
and $b_0\in \mathfrak C_{x,\xi}$. By Condition $\mu$-SL,
there is  a conormal distribution ${\bf f}\in \I^{n+1}(Y)=\I^{n+1}(N^*Y)$
such that ${\bf f}$ satisfies 
the linearized conservation law
(\ref{eq: lineariz. conse. law PRE})
and the principal symbol
$\tilde f$ 
of ${\bf f}$, defined on $N^*Y$, \MTEXT{satisfies
  $\tilde f(y,\xi)=b_0$.} 
Moreover, by Condition $\mu$-SL there is a family of sources
$\F_\e$, $\e\in [0,\e_0)$ such that $\p_\e\F_\e|_{\e=0} ={\bf f}$
and a solution $u_\e+(\hat g,\hat \phi)$ of the Einstein equations
with the source $\F_\e$ that depend smoothly on $\e$
and $u_\e|_{\e=0}=0$. 
Then $\dot u=\p_\e u_\e|_{\e=0} \in \I^{n-1/2} ( M_0\setminus Y;\Lambda_1)$.

As  $\dot u=(\dot g,\dot\phi)$ satisfies the linearized harmonicity condition (\ref{harmonicity condition}),
the principal symbol $\tilde a(x,\eta)=(\tilde a_1(x,\eta),\tilde a_2(x,\eta))$
of $\dot u$ satisfies $\tilde a(x,\eta)\in \mathfrak S_{x,\eta}$.
This shows that the map $R=R(x,\eta,y,\xi)$, given by $R:{{\tilde f}}(y,\xi)\mapsto \tilde a(x,\eta)$
that is defined in Lemma \ref{lem: Lagrangian 1},
satisfies $R:\mathfrak C_{y,\xi}\to \mathfrak S_{x,\eta}$.
Since $R$ is one-to-one and the linear spaces $ \mathfrak C_{y,\xi}$ and $\mathfrak S_{x,\eta}$  have the same dimension, we see that
\beq\label{R bijective}
R:\mathfrak C_{y,\xi}\to \mathfrak S_{x,\eta}
\eeq
 is a bijection.
Hence, 
when  ${\bf f}\in \I^{n+1}(Y)$ varies so that the linearized conservation law  (\ref{eq: lineariz. conse. law symbols})  for the principal symbols 
is satisfied, the principal symbol  $\tilde a(x,\eta)$ at $(x,\eta)$ of the
solution $\dot u$ of the linearized Einstein equation
 achieves all values in the $(L+6)$ dimensional space $\mathfrak S_{x,\eta}$.

}

 Below, we denote ${\bf f}\in  \mathcal I_{C}^{n+1}
(Y(x_0,\zeta_0;t_0,s_0))$ when
 the principal symbols of ${\bf f}$  
 satisfies the linearized conservation law 
 for principal symbols, that is, equation (\ref{eq: lineariz. conse. law symbols}).

 \subsection{Microlocal analysis of the non-linear interaction of waves}

 \MATTITEXT{ 
 Next we consider the interaction of four $C^k$-smooth 
 waves  having conormal singularities, where $k\in \Z_+$ is sufficiently large. 
 Interaction of such waves produces a ``corner point'' in the spacetime. On related microlocal tools to consider scattering by corners, see \cite{Vasy1,Vasy2}.
 Earlier considerations of interaction of three waves 
 has been done by Melrose and Ritter \cite{MR1,MR2}
 and Rauch and Reed,  \cite{R-R}  for non-linear hyperbolic equations in $\R^{1+2}$
 where the non-linearity appears in the lower order terms. \MTEXT{Recently, the
 interaction of two strongly singular waves has been studied by 
 Luc and Rodnianski \cite{L-Rodnianski}.} 
   
\subsubsection{Interaction of non-linear waves on a general manifold}

Next, we introduce a vector of four $\e$ variables
denoted by $\vec \e=(\e_1,\e_2,\e_3,\e_4)\in \R^4$. 
Let $s_0,t_0>0$ and consider  $u_{\vec \e}=(g_{\vec\e}-\hat g,\phi _{\vec\e}-\hat \phi)$
where $v_{\vec \e}=(g_{\vec \e},\phi_{\vec \e})$ solve the equations (\ref{eq: adaptive model with no source})
with $\F={\bf f}_{\vec \e}$ where 
\beq\label{eq: f vec e sources}
 {\bf f}_{\vec \e}:=\sum_{j=1}^4\e_j {\bf f}_{j},\quad
 {\bf f}_{j}\in \mathcal I_{C}^{n+1} (Y(x_j,\zeta_j;t_0,s_0)),
 \eeq
and $(x_j,\zeta_j)$ are light-like
vectors with $x_j\in U_{\hat g}.$
Moreover,  we assume that for some $0<r_2<r_1$  and $s_-+r_2<s^\prime<s_+$ the sources satisfy
 \beq\label{eq: source causality condition}
& & \supp({\bf f}_{j})\cap J ^+_{\hat g}(\supp({\bf f}_{k}))=\emptyset,\quad\hbox{for all }j\not=k,\\
\nonumber
& &\supp({\bf f}_{j})\subset  I_{\hat g}( \mu_{\hat g}(s^\prime-r_2),\mu_{\hat g}(s^\prime)),\quad \hbox{for all }j=1,2,3,4,
%
 \eeq
 where $r_1$ is the parameter used
 to define $W_{\hat g}=W_{\hat g}(r_1)$, see (\ref{observer neighborhood with hat}).
  The first condition implies that the supports of the sources are causally independent.
}


The sources   $ {\bf f}_{j}$ give raise to  $\B^L$-section valued solutions
of the linearized wave equations, which we denote by
\ba
u_j:=u^{(1)}_{j}=\p_{\e_j}u_{\vec \e}|_{\vec \e=0}={\bf Q} \,{\bf f}_{j}\in \I( \Lambda(x_j,\zeta_j;t_0,s_0)),
\ea
where   ${\bf Q}={\bf Q}_{\hat g}$.
In the following we use the notations
$\p_{\vec \e}^1u_{\vec \e}|_{\vec \e=0}:=\p_{\e_1}u_{\vec \e}|_{\vec \e=0}$,
$\p_{\vec \e}^2 u_{\vec \e}|_{\vec \e=0}:=\p_{\e_1}\p_{\e_2}u_{\vec \e}|_{\vec \e=0},$
$\p_{\vec \e}^3 u_{\vec \e}|_{{\vec \e}=0}:=\p_{\e_1}\p_{\e_2}\p_{\e_3} u_{\vec \e}|_{{\vec \e}=0}$,
and
\ba
\p_{\vec \e}^4 u_{\vec \e}|_{\vec \e=0}:=\p_{\e_1}\p_{\e_2}\p_{\e_3}\p_{\e_4} u_{\vec \e}|_{{\vec \e}=0}.
\ea

Next we denote the waves produced by the $\ell$-th order interaction by
\beq\label{measurement eq.}
\M^{(\ell)}:=\p_{\vec \e}^\ell u_{\vec \e}|_{\vec \e=0},\quad \ell\in \{1,2,3,4\}.
\eeq

}
Below,  we use the notations $\B^\beta _j$, $j=1,2,3,4$ and $S^\beta_n$, $n=1,2$ to denote operators of the form
\beq\label{extra notation}
& &\hbox{$\B^\beta _j:(v_{p})_{p=1}^{10+L}\mapsto (b^{r,(j,\beta )}_{p}(x)\p_x^{ k(\beta,j)} v_{r}(x))_{p=1}^{10+L}$, and}\\ \nonumber
& &\hbox{$S^\beta _n={\bf Q}$ or $S^\beta _n=I$, or $S^\beta _n=[{\bf Q},a^\beta_n(x)D^\a]$},\quad \a=\a(\beta,n),\  |\a|\leq 4
\eeq
\extension{Here,  $|k(\beta)|=\sum_{j=1}^4k(\beta,j)$ 
and $|\a(\beta)|=|\a(\beta,1)|+|\a(\beta,2)|$ satisfy the bounds described below.}
and the coefficients $b^{r,(j,\beta )}_{p}(x)$ and $a^\beta_j(x)$ 
depend on the derivatives of $\hat g$.
Here, $k(\beta,j)=k_j^\beta=ord(\B^\beta_j)$ are the orders of $\B^\beta_j$ and
$\beta$ is just an index running over
a  finite set $J_4\subset \Z_+$.


 Computing the $\e_j$ derivatives of the equations (\ref{w1-w3}), (\ref{w4}), and (\ref{eq: mixed source terms})
with the sources ${\bf f}_{\vec \e}$, and taking into account
 the condition (\ref{eq: source causality condition}), we obtain: 
\noextension{
  \beq\label{w1-3 solutions}
& &\hspace{-2cm}\M^{(4)}=
{\bf Q}\F^{(4)},\quad  \F^{(4)}=
\sum_{\sigma\in \Sigma(4)}\sum_{\beta\in J_4}
(\mathcal G^{(4),\beta}_\sigma+\tilde {\mathcal G}^{(4),\beta}_\sigma),
\eeq 
where  
$\Sigma(4)$  is the set of permutations, that is, the set of the bijections
 $\sigma:\{1,2,3,4\}\to \{1,2,3,4\}$. 
}
\extension{
 \beq\label{w1-3 solutions A}
 \M^{(1)}&=&u_1,\\
\nonumber
 \M^{(2)}&=&\sum_{\sigma\in \Sigma(2)}{\bf Q}(A[ u_{\sigma(2)},u_{\sigma(1)}]),\\
\nonumber
 \M^{(3)}&=&\sum_{\sigma\in\Sigma(3)}\bigg(
2{\bf Q}(A[ u_{\sigma(3)},{\bf Q}(A[ u_{\sigma(2)}, u_{\sigma(1)}])])+
\\ \nonumber & &\quad +{\bf Q}(B[  u_{\sigma(3)},
 u_{\sigma(2)}, u_{\sigma(1)}])\bigg),
 \eeq
 where $A$ is a bi-linear operator discussed below and $B$  is a trilinear
 operator,
 or 
  \beq\label{w1-3 solutions}
& & \hspace{-2cm}\M^{(1)}=u_1,\\ \nonumber
  & &\hspace{-2cm}\M^{(2)}=\sum_{\sigma\in \Sigma(2)}\sum_{\beta\in J_2 }{\bf Q}(\B^{\beta}_2u_{\sigma(2)}\,\cdotp \B^{\beta}_1u_{\sigma(1)}),\hspace{-2cm}\\ \nonumber
 & &\hspace{-2cm}\M^{(3)}=\sum_{\sigma\in \Sigma(3)}\sum_{\beta\in J_3}{\bf Q}(\B_3^{\beta}u_{\sigma(3)}\,\cdotp  S^{\beta}_1(\B^{\beta}_2u_{\sigma(2)}\,\cdotp \B^{\beta}_1u_{\sigma(1)})),\hspace{-2cm} \nonumber
 \\ \nonumber
& &\hspace{-2cm}\M^{(4)}=
{\bf Q}\F^{(4)},\quad  \F^{(4)}=
\sum_{\sigma\in \Sigma(4)}\sum_{\beta\in J_4}
(\mathcal G^{(4),\beta}_\sigma+\tilde {\mathcal G}^{(4),\beta}_\sigma),
\eeq 
where  
$\Sigma(k)$  is the set of permutations, that is, bijections
 $\sigma:\{1,2,\dots,k\}\to \{1,2,\dots,k\}$, and $J_k\subset \Z_+$ are  finite sets. }
%
%
Then, we have
\beq\label{M4 terms}
& & {\mathcal G}^{(4),\beta}_\sigma=\B_4^{\beta}u_{\sigma(4)} \,\cdotp S_2^{\beta}(\B_3^{\beta}u_{\sigma(3)}\,\cdotp S^{\beta}_1(\B^{\beta}_2u_{\sigma(2)}\,\cdotp \B^{\beta}_1u_{\sigma(1)}))\hspace{-2cm}
\eeq
\MTEXT{where  the orders of the differential operators
satisfy $k_4^\beta+k_3^\beta+k_2^\beta+k_1^\beta+|\a(\beta,1)|+|\a(\beta,2)|\leq 6$ and
$k_4^\beta+k_3^\beta+|\a(\beta,2)|\leq 2$
and 
\beq\label{tilde M4 terms}
& &\tilde  {\mathcal G}^{(4),\beta}_\sigma=\mathcal S_2^{\beta} (\B_4^{\beta}u_{\sigma(4)} \,\cdotp \B_3^{\beta}u_{\sigma(3)})\,\cdotp \mathcal S^{\beta}_1(\B^{\beta}_2u_{\sigma(2)}\,\cdotp \B^{\beta}_1u_{\sigma(1)}),\hspace{-2cm}
\eeq
where  $k_4^\beta+k_3^\beta+k_2^\beta+k_1^\beta+|\a(\beta,1)|+|\a(\beta,2)|\leq 6$,
$k_4^\beta+k_3^\beta+|\a(\beta,2)|\leq 4$, and $k_1^\beta+k_2^\beta+|\a(\beta,1)|\leq 4$.}

We denote below $\vec S_\beta=(S^\beta _1, S^\beta _2)$ and
  $  \M^{(4),\beta}_\sigma={\bf Q}_{\hat g}\mathcal G^{(4),\beta}_\sigma$ and 
$\tilde  \M^{(4),\beta}_\sigma={\bf Q}_{\hat g} \tilde { \mathcal G}^{(4),\beta}_\sigma$.
\motivation{Let us explain how the  terms above appear in the Einstein equations: 
\MTEXT{By taking $\p_{\vec \e}$ derivatives of the wave $u_{\vec\e}$
we obtain terms similar to (\ref{w1-w3})  and (\ref{w4}). 
In particular, we have that   $  \G^{(4),\beta}_\sigma$ and 
$\tilde { \mathcal G}^{(4),\beta}_\sigma$ can be written in the form}
\beq\label{eq: tilde M1}
& & \mathcal G^{(4),\beta}_\sigma=A^\beta_3[  u_{\sigma(4)},S^\beta _2(A^\beta_2[ u_{\sigma(3)},S^\beta _1(A^\beta_1[ u_{\sigma(2)}, u_{\sigma(1)}])],
\\ \label{eq: tilde M2}
& &\tilde  {\mathcal G}^{(4),\beta}_\sigma =A^\beta_3[ S^\beta _2(A^\beta_2[ u_{\sigma(4)}, u_{\sigma(3)}]),
S^\beta _1(A^\beta_1[ u_{\sigma(2)}, u_{\sigma(1)}])],
\eeq
where 
$A^\beta_j[V,W]$ are   2nd  order multilinear operators of the form
\beq\label{A form}
A[V,W]=\sum_{|\a|+|\gamma|\leq 2}a_{\a\gamma}(x)(\p_x^\a V(x))\,\cdotp(\p_x^\gamma W(x)).
\eeq 
By commuting derivatives and the operator ${\bf Q}_{\hat g}$  we obtain 
(\ref{M4 terms}) and (\ref{tilde M4 terms}).
Below, \MTEXT{the terms of particular importance are} the bilinear form that is given
for $V=(v_{jk},\phi)$
and  $W=( w^{jk},\phi^{\prime})$
by
\beq\label{A alpha decomposition2}
& & A_1[V,W]=-\hat g^{jb} w_{ab}\hat g^{ak}\p_j\p_k  v_{pq}.
\eeq}

%


%
%
%
%



%


%

\subsubsection{On the singular support of the non-linear interaction  of three waves}
\label{subsection: sing supp interactions}


Let us next consider four light-like future pointing directions $(x_j,\xi_j)$, $j=1,2,3,4$, and
use below the notations, see (\ref{eq: x(h) notation}),
\ba
(\vec x,\vec\xi)=((x_j,\xi_j))_{j=1}^4,\quad 
(\vec x(h),\vec\xi(h))=((x_j(h),\xi_j(h)))_{j=1}^4.
\ea
We will consider the case when we send distorted plane  waves
propagating on surfaces $K_j=K(x_j,\xi_j;t_0,s_0)$,  $t_0,s_0>0$, cf.\ (\ref{associated submanifold}), and
these waves interact. 

Next we consider the 3-interactions of the waves.
Let $\X((\vec x,\vec\xi);t_0,s_0)$ be set of all light-like vectors $(x,\xi)\in L^+M_0$
that are in  the normal bundles   $N^*(K_{j_1}\cap K_{j_2}\cap K_{j_3})$
with some $1\leq j_1<j_2<j_3\leq 4$.
%

Moreover, we define ${\mathcal Y}((\vec x,\vec\xi);t_0,s_0)$ to be the set of  all  $y\in M_0$
such that there are $(z,\zeta)\in \X ((\vec x,\vec\xi);t_0,s_0)$, and
 $t\geq 0$ such that $\gamma_{z,\zeta}(t)=y$. Finally,
 let 
 \beq\label{Y-set}
 {\mathcal Y}((\vec x,\vec\xi);t_0)=\bigcap_{s_0>0}{\mathcal Y}((\vec x,\vec\xi);t_0,s_0).
 \eeq  
The three wave interaction happens then on $\pi ( \X ((\vec x,\vec\xi);t_0,s_0))$ and, roughly speaking, this interaction sends singularities
to  $\Y ((\vec x,\vec\xi);t_0,s_0)$.

%
%

\motivation{For instance in  Minkowski space, when three plane waves (whose singular supports 
are hyperplanes) collide,
the intersections of the hyperplanes is a 1-dimensional space-like line $K_{123}=K_1\cap K_2\cap K_3$ in the 4-dimensional space-time.
This corresponds to a point  moving continuously in time. Roughly speaking, the point 
seems to move at a higher speed than light (i.e.\ it appears like a tachyonic, moving point source)
and  produces a (conic) shock wave type of singularity (see Fig.\ 2 where the interaction time is only finite).
In this paper we do not analyze carefully the singularities produced by the  three wave interaction 
near $ \Y ((\vec x,\vec\xi);t_0,s_0)$. Our goal is to 
consider the singularities produced by the four wave interaction in the domain
$M_0\setminus \Y ((\vec x,\vec\xi);t_0,s_0)$.}

\extension{

\begin{center}

\psfrag{1}{$x^0$}
\psfrag{2}{$x^1$}
\psfrag{3}{$x^2$}
\includegraphics[width=7.5cm]{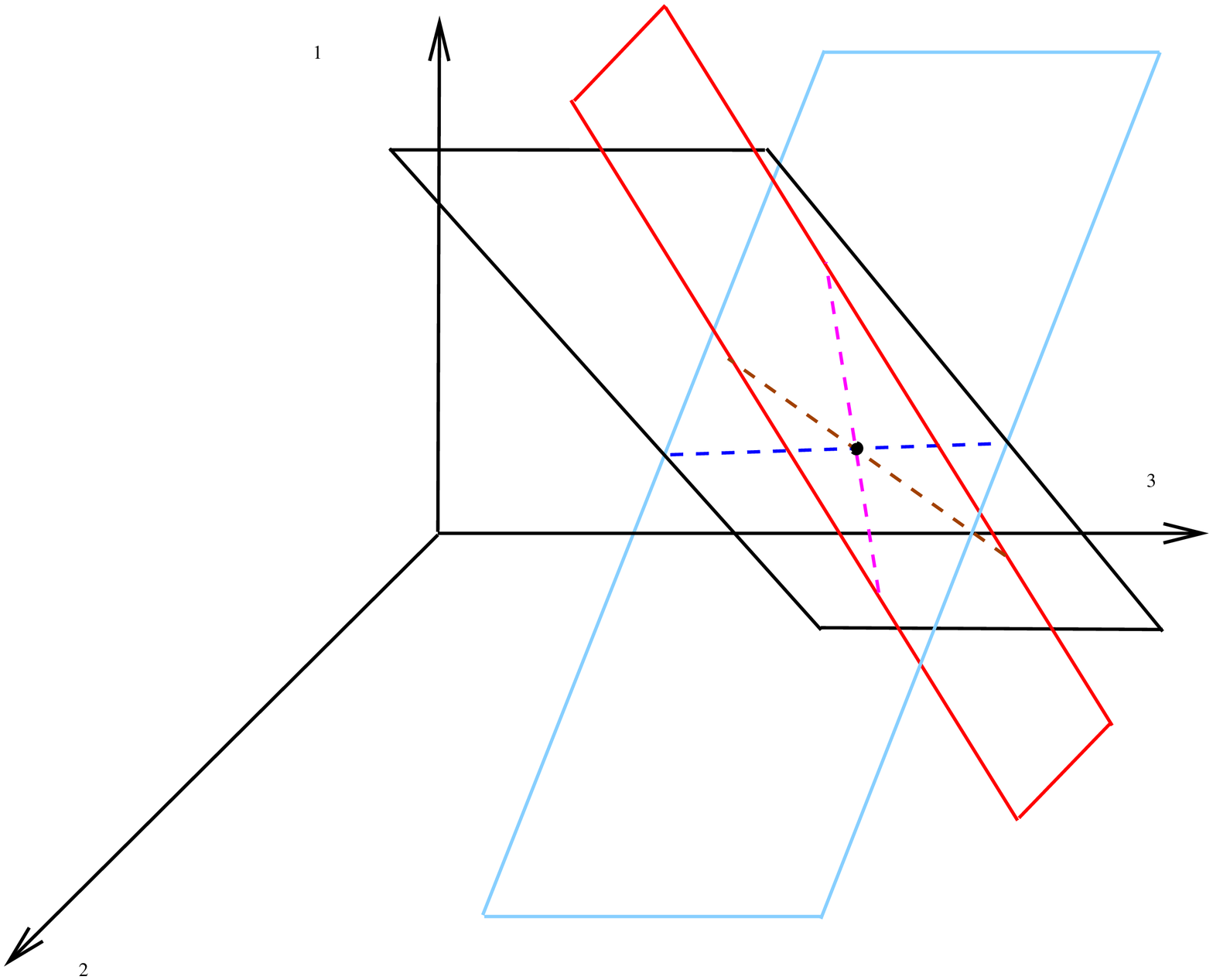}
\end{center}
{\it FIGURE A1: A schematic figure where the space-time is represented as the  3-dimensional set $\R^{2+1}$. 
In the figure 3 pieces of plane waves  have singularities 
on strips of hyperplanes (in fact planes) $K_1,K_2,K_3$, colored
by light blue, red, and black. These planes  have intersections, and
in the figure the sets $K_{12}=K_1\cap K_2$, $K_{23}=K_2\cap K_3$,
and $K_{13}=K_1\cap K_3$ are shown as dashed lines with
dark blue, magenta, and brown colors. These dashed lines  intersect
at a point $\{q\}=K_{123}=K_1\cap K_2\cap K_3$.}

\begin{center}

\psfrag{1}{$\pi(\mathcal U_{z_0,\eta_0})$}
\psfrag{2}{$U_{\hat g}$}
\psfrag{3}{$J_{\hat g}({p^-},{p^+})$}
\psfrag{4}{$W_{\hat g}$}
\psfrag{5}{$\mu_{z,\eta}$}
\includegraphics[width=7.5cm]{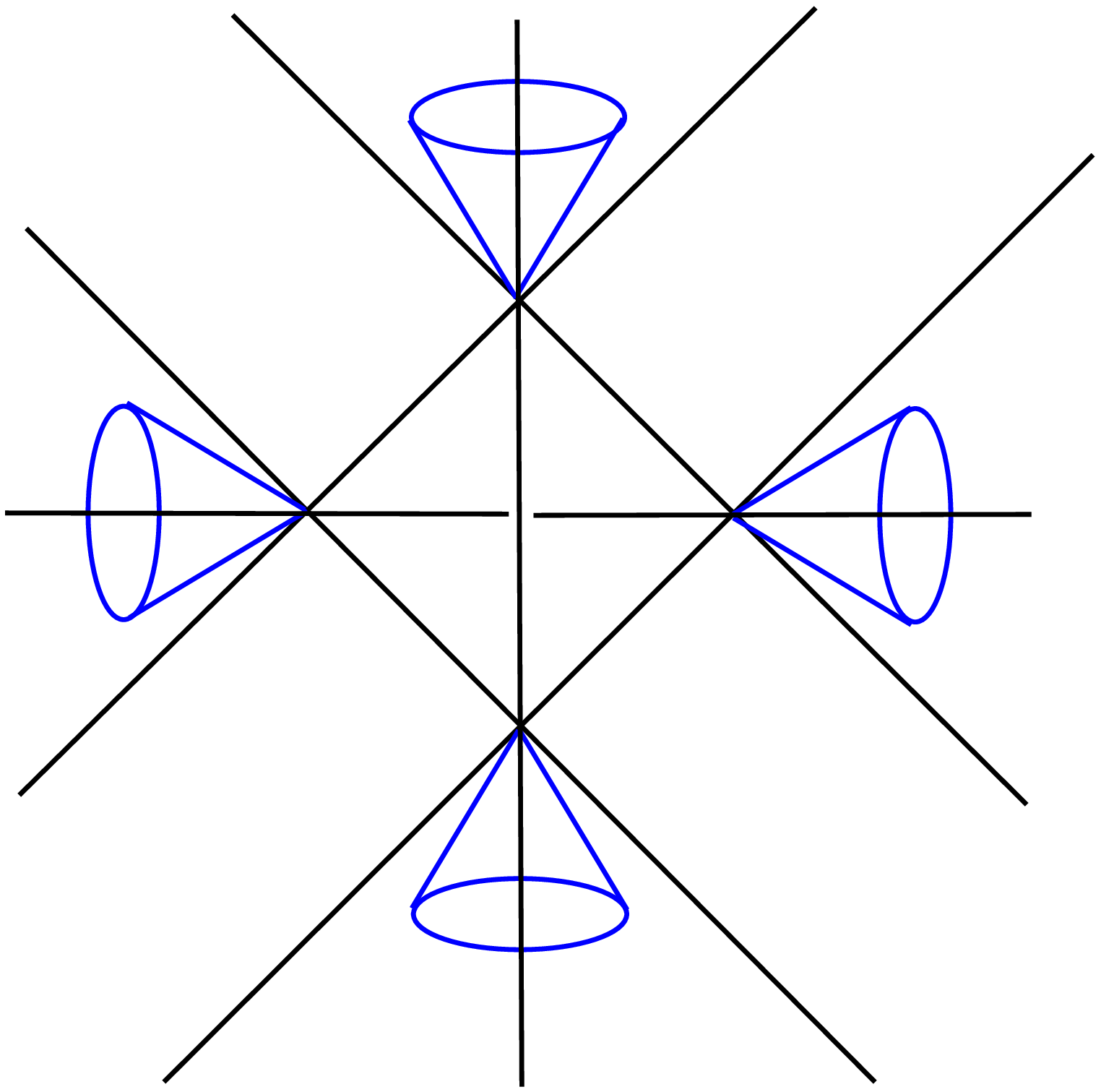}
\end{center}
{\it FIGURE A2: In section \ref{subsection: sing supp interactions} we consider four colliding distorted plane  waves.
In the figure we consider the Minkowski space $\R^{1+3}$ and the figure corresponds 
to a ``time-slice'' $\{T_1\}\times \R^3$.
We assume that the distorted plane  waves are sent 
so far away that the waves look like pieces of plane waves. The plane waves $u_j\in \mathcal I(K_j)$, $j=1,2,3,4$,
are conormal distributions that are  solutions of the linear wave equation and their singular supports are the sets $K_j$,
that are pieces of 3-dimensional planes in the space-time. The sets $K_j$ are not shown in the figure.
The 2-wave interaction wave $\mathcal M^{(2)}$ is singular on the set $\cup_{j\not =k}K_j\cap K_j$.
There are 6 intersection sets $K_j\cap K_j$ that are shown as black line segments. Note that these
lines have 4 intersection points, that is, the vertical and the horizontal black
lines do not intersect  in $\{T_1\}\times \R^3$.
The four intersection points of the black lines are the sets $(\{T_1\}\times \R^3)\cap (K_j\cap K_j\cap K_n)$.
These points correspond to  points moving in time (i.e., they are curves in the space-time) that produce singularities of the 3-interaction wave 
$\mathcal M^{(3)}$. The points
 seem to move faster than the speed of the light (similarly, as  a shadow of a
far away  object may seem to  move faster than the speed of the light).  Such
point sources produce ``shock waves'', and due to this, $\mathcal M^{(3)}$ is
singular on the sets $\mathcal Y((\vec x,\vec \xi),t_0,s_0)$ defined in 
formulas (65)-(66). The set $(\{T_1\}\times \R^3)\cap \mathcal Y((\vec x,\vec \xi),t_0,s_0)$
is the union of the four blue cones shown in the figure. 
}

\begin{center}

\psfrag{1}{$y_{-1}$}
\includegraphics[width=7.5cm]{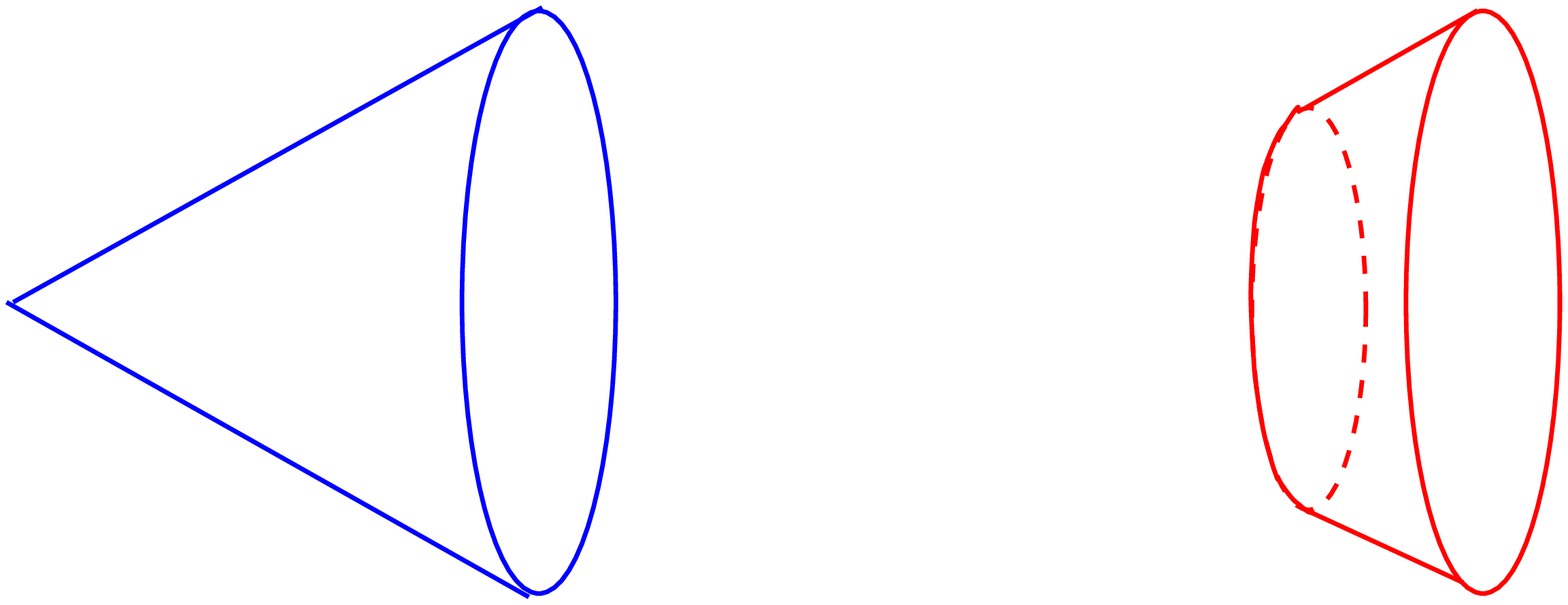}
\end{center}
{\it FIGURE A3: A schematic figure in  3-dimensional Euclidean space  $\R^{3}$
where
we describe the geometry of the wave produced by the interaction of three
waves in the Minkowski space. The figure shows such a 3-interaction wave in a time slice  $\{T_1\}\times \R^{3}$
with two different values of the parameter $s_0$.  On the left, $s_0$ is large
and the 3-interaction wave has singular support on a cone. 
 On the right, $s_0$ is small and the  3-interaction wave is a frustum (a truncated cone, or,  a cone with its apex cut off by a plane), similar to set 
$\{(x^1,x^2,x^3)\in \R^3;\ (x^1)^2+(x^2)^2=c(x^3)^2$, $a<x^3<a+h\}$, with some
$a,c>0$ and a thickness $h=c_1s_0$.
}

\begin{center}

\psfrag{1}{$\pi(\mathcal U_{z_0,\eta_0})$}
\psfrag{2}{$U_{\hat g}$}
\psfrag{3}{$J_{\hat g}({p^-},{p^+})$}
\psfrag{4}{$W_{\hat g}$}
\psfrag{5}{$\mu_{z,\eta}$}
\includegraphics[width=7.5cm]{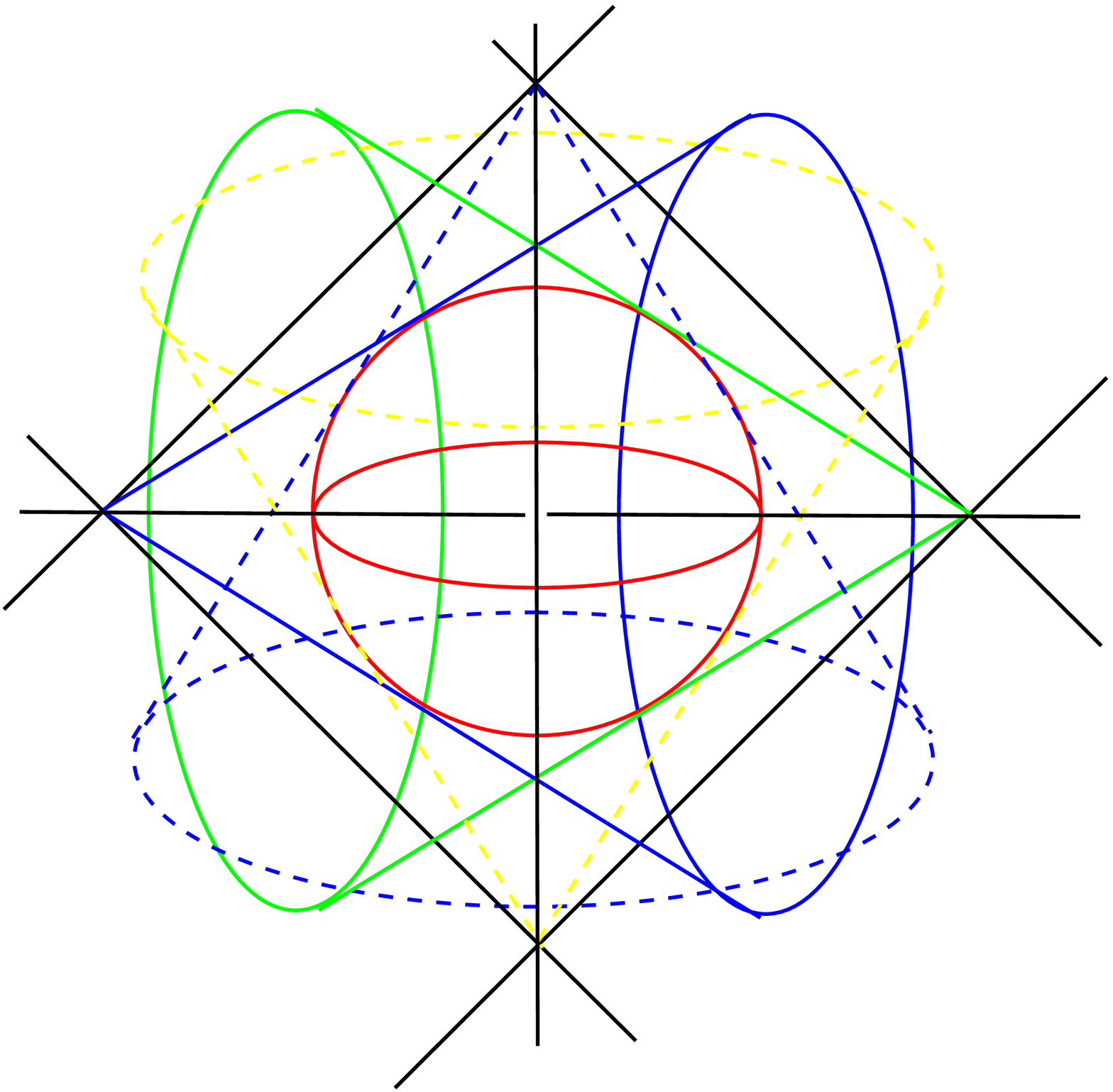}
\end{center}

\noindent
{\it FIGURE A4: The same situation that was described in Fig. A2 
is shown
at a later time, that is, the figure shows the time-slice $\{T_2\}\times \R^3$ with
$T_2>T_1$, when the parameter $s_0$ is quite large. The four  distorted plane waves have now collided and
produced a point source in the space-time at a point $q\in K_1\cap K_2\cap K_3\cap K_4$.
This gives  raise to the singularities of the 4-interaction wave $\mathcal M^{(4)}$.
In the figure $T_1<t<T_2$, where $q$ has the time coordinate $t$.
The four cones in the figure, shown with solid blue and green curves
and dashed blue and yellow curves are   the intersection of  the time-slice $\{T_2\}\times \R^3$ 
and the set $\mathcal Y((\vec x,\vec \xi),t_0,s_0)$. Inside the cones the red sphere
is the set $\mathcal L^+(q)\cap (\{T_2\}\times \R^3)$ that
corresponds to the spherical plane  wave produced by the point source at $q$.
}

\begin{center}

\psfrag{1}{$\pi(\mathcal U_{z_0,\eta_0})$}
\psfrag{2}{$U_{\hat g}$}
\psfrag{3}{$J_{\hat g}({p^-},{p^+})$}
\psfrag{4}{$W_{\hat g}$}
\psfrag{5}{$\mu_{z,\eta}$}
\includegraphics[width=5.5cm]{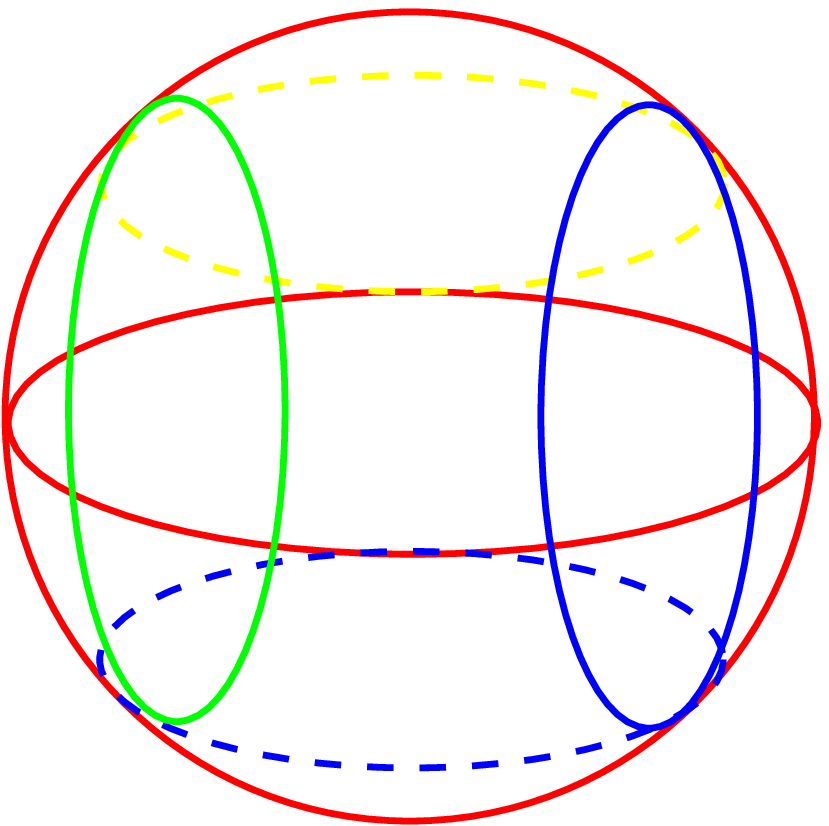}
\end{center}
{\it FIGURE A5: The same situation that was described in Fig.\ A4. 
That is,  the figure shows in 
the time-slice $\{T_2\}\times \R^3$ the singularities produced by four colliding
distorted plane  waves, when the parameter $s_0$ is very small. In this case the truncated cones
degenerate to circles.
Indeed,
 the figure corresponds to the case when the parameter $s_0$ is close
to zero, that is, the pieces of the distorted plane  waves are concentrated near a single
geodesic in the space-time.
As  $s_0$ is small, the distorted plane  waves interact only during a short time.
Due to this, the sets $(K_j\cap K_n)\cap (\{T_2\}\times \R^3)$,
that were black lines in the previous figures, are now empty, and do not
appear at the figure.
In the figure the  red sphere is the set $\mathcal L^+(q)\cap (\{T_2\}\times \R^3)$ corresponding
to  the spherical  wave
 produced by a point source at $q$.
The solid blue and green circles on the sphere
and the dashed blue and yellow  circles on the sphere are
 the intersection of  the time-slice $\{T_2\}\times \R^3$ 
and the set $\mathcal Y((\vec x,\vec \xi),t_0,s_0)$. These circles
are actually not circles but surfaces of frustums (a truncated cone, or,  a cone with its apex cut off by a plane), similar to set 
$\{(x^1,x^2,x^3)\in \R^3;\ (x^1)^2+(x^2)^2=c(x^3)^2$, $a<x^3<a+h\}$, with some
$a,c>0$ and a very small thickness $h>0$. Note that the 2-Hausdorff measure of the frustums
go to zero as $s_0\to 0$ and the 3-Hausdorff measure of the set
$\mathcal Y((\vec x,\vec \xi),t_0,s_0)$ goes to zero as $s_0\to 0$. Thus, when all the distorted plane  waves
intersect at a point $q$, 
the singularities of $\mathcal M^{(4)}$  are observed on a set whose 3-Hausdorff measure in
the space-time is independent of $s_0$ but the singularities of $\mathcal M^{(3)}$ 
are observed on a set which 3-hHausdorff measure in
the space-time is smaller than $Cs_0$. 
}

}

\subsubsection{Gaussian beams}

 \MATTITEXT{ 
Our aim is to consider interactions of  4 waves that produce  a new source, and
to this end we use test sources that produce gaussian beams.

Let $y\in U_{\hat g}$ and  $\eta\in T_{y}M$ be a future pointing light-like vector. 
We choose a complex function $p\in C^\infty(\hattuM _0)$
such that
 $\im  p(x) \geq 0$ and $\im  p(x)$ vanishes
only at $y$, 
$ p(y)=0$, $d(\re p)|_{y}=\eta^\sharp$, 
$d(\im  p)|_y=0$ and
the Hessian of $\im  p$ at $y$ is positive definite.
To simplify notations, we use below  complex sources and waves.
 The physical linearized waves can be obtained by taking the real part of these complex waves. We use a large parameter $\tau\in \R_+$ and define a  test source  
\beq\label{Ftau source}
F_\tau(x)=F_\tau(x;y,\eta),\ \hbox{where}\quad
F_\tau(x;y,\eta)=\tau^{-1}
\exp (  i\tau  p(x))
 h(x) \hspace{-1cm}
\eeq
where $h$ is section on $\mathcal B^L$ supported
in a small neighborhood $W$ of $y$. 
\noextension{The construction of $p(x)$ and $F_\tau$ is discussed
in detail in \cite{preprint}.}
\extension{The construction of $p(x)$ and $F_\tau$ is discussed
in detail below.}


\motivation{
We consider both the usual causal
solutions and the  solutions of the adjoint wave equation
for which time is reversed, that is, we use the anti-causal parametrix
${\bf Q}^*={\bf Q}^*_{\hat g}$ instead of the usual causal   parametrix ${\bf Q}=(\square_{\hat g}+V(x,D))^{-1}$. 

For a function $v:M_0\to \R$ 
the large $\tau$ asymptotics of the $L^2$-inner products $\bra F_\tau,v\cet_{L^2(M_0)}$
can be used to verify if the point $(y,\eta^\sharp)\in T^*M_0$ belongs in the wave front
set $\hbox{WF}(v)$ of $v$, see e.g.\ \cite{Ralston}. Below, we will use
the fact that when $v$ is of the form $v={\bf Q}_{\hat g}b$, we have 
$\bra F_\tau,v\cet_{L^2(M_0)}=\bra {\bf Q}_{\hat g}^*F_\tau,b\cet_{L^2(M_0)}$
where ${\bf Q}_{\hat g}^*F_\tau$ is a Gaussian beam solution to the adjoint
wave equation.
We will use this to analyze the singularities of $\M^{(4)}$ at $y\in U_{\hat g}$.}

The wave   $u_\tau={\bf Q}^*F_\tau$ satisfies by \cite{Ralston}, see also \cite{KKL},
\beq\label{gaussian beam 1}
\|u_\tau-u^N_{\tau}\|_{C^k(J(p^-,p^+))}\leq C_N\tau^{-n_{N,k}}
\eeq
where $n_{N,k}\to \infty$ as $N\to \infty$ and $u^N_\tau$ is a formal Gaussian beam  of order $N$ having
the form
\beq\label{gaussian beam 2}
u^N_{\tau}(x ) = 
e^{i\tau \varphi (x) }a(x,\tau),\quad 
a(x,\tau)=
\sum _{n=0}^N U_n(x) \tau^{-n}.
\eeq
Here
$
\varphi (x)=A(x)+iB(x)
$
and $A(x)$ and $B(x)$ are real-valued functions, $B(x)\geq 0$, and $B(x)$ vanishes
only on $\gamma_{y,\eta} (\R)$, and
for $z= \gamma_{y,\eta}(t)$, and
$\zeta=\dot \gamma_{y,\eta}^\flat (t)$, $t<0$ we have
$dA|_{z}=\zeta$,  $dB|_{z}=0$. Moreover,
the Hessian of $B$ at $z$  restricted to the orthocomplement 
of $\zeta$ 
is positive definite.
The functions $h$ and $U_n$ defined above can be chosen to be smooth so that $h$ is supported in
any neighborhood $V$ of $y$ and $U_N$ are supported in
any neighborhoods of  
$\gamma_{y,\eta} ((-\infty,0])$.  Below, $a(x,\tau)$ is the symbol and $U_0(y)$ is the principal symbol of the gaussian beam $u_\tau(x)$.

\extension{
The source function $F_\tau$ can
be constructed in local coordinates using the asymptotic representation of the 
gaussian beam, namely, by considering first 
\ba
F^{(series)}_\tau(x)=c\tau^{-1/2}\int_\R e^{-\tau s^2}\phi(s)f_\tau(x;s)ds,
\ea
where
$f_\tau(x;s)=\square_{\hat g}(H(s-x^0)u_\tau(x))$,. Here, $\phi\in C^\infty_0(\R)$ 
has value 1 in some neighborhood of zero, and $H(s)$ is the Heaviside function,
and  we have 
\ba
e^{-i\tau p(x)}F^{(series)}_\tau(x)=e^{-i\tau p(x)}F^{(1)}(x)\tau^{-1}+e^{-i\tau p(x)}F^{(2)}(x)\tau^{-2}+O(\tau^{-3}).
\ea Then,
$F_\tau$  can be defined by $F_\tau=F^{(1)}(x)\tau^{-1}.$ 
Indeed
we see that
\ba
\square_{\hat g}^{-1}F^{(series)}_\tau&=&
\tau^{-1}\Phi_\tau(x)\,u_\tau(x),\quad \hbox{with}\\
\Phi_\tau(x)&=&
c\tau^{1/2}\int_\R e^{-\tau s^2}\phi(s)H(s-x^0) ds\\
&=&
\left\{\begin{array}{cl}
O(\tau^{-N}),&\hbox{for }x_0<0,\\
1+O(\tau^{-N}),&\hbox{for }x_0>0.
\end{array}\right.
\ea
Moreover, by writing $\square_{\hat g}=\hat g^{jk}\p_j\p_k+\hat \Gamma^j\p_j$,
$\hat \Gamma^j=\hat  g^{pq}\hat \Gamma^j_{pq}$ we have  that
\ba
f_\tau(x;s)&=&
\hat g^{00}u_\tau(x)\delta^\prime(s-x^0)\\
& &-2\sum_{j=1}^3\hat g^{j0}(\p_ju_\tau(x))\,\delta(s-x^0)
-\hat \Gamma^0u_\tau(x)\,\delta(x^0-s).
\ea
Let us use  normal coordinates centered at the point $y$ so that
$\Gamma^j(0)=0$ and $\hat g^{jk}(0)$ is the Minkowski metric.
Then we define $p(x)=(x^0)^2+\varphi(x)$. The leading
order term of  $F^{(series)}_\tau$ is given by
\ba
F_\tau= c\tau^{-1/2}e^{-\tau p(x)}\left(
2x^0\phi(x^0)U_0(x)+2i(\p_0\varphi(x))U_0(x)\phi(x^0)
\right),
\ea
where $\p_0\varphi|_y$ does not vanish. }}


\subsubsection{Indicator function for singularities produced by interaction of four waves}

Let $y\in U_{\hat g}$ and $\eta\in L^+_{x_0}(M_0,\hat g)$ be  a future pointing light-like vector.
We will next check whether   
$(y,\eta^\flat)\in\WF(\M^{(4)})$.

Using the function $\M^{(4)}$ defined
in (\ref{w1-3 solutions}) 
with  four sources
$ {\bf f}_{j}\in \mathcal I_{C}^{n+1} (Y(x_j,\zeta_j;t_0,s_0))$, $j\leq 4$,  that produce
the  pieces of plane waves 
$u_j\in \I^n(\Lambda(x_{j},\xi_{j};t_0,s_0))$ in $M_0\setminus Y(x_j,\zeta_j;t_0,s_0)$, 
and 
 the source $F_\tau$ in (\ref{Ftau source}),
 we define the indicator functions
\beq\label{test sing}
\Theta_\tau^{(4)}=\bra F_{\tau},\M^{(4)}\cet_{L^2(U)}=\sum_{\beta\in J_4}\sum_{\sigma\in \Sigma(4)}(T_{\tau,\sigma}^{(4),\beta}+\tilde T_{\tau,\sigma}^{(4),\beta}).
\eeq
\motivation{Here, the terms  $T_{\tau,\sigma}^{(\ell),\beta}$ and $\tilde T_{\tau,\sigma}^{(\ell),\beta}$ correspond, for $\sigma$ and $\beta$,  
 to different types of interactions the four waves $u_{(j)}$ can have.}
To define the terms
 $T_{\tau,\sigma}^{(\ell),\beta}$ and $\tilde T_{\tau,\sigma}^{(\ell),\beta}$ 
appearing above,
we use  generic notations where we drop
the index $\beta$, that is, we denote $\B_j=\B_j^\beta$ and $S_j=S_j^\beta.$
With these notations, using 
the decompositions (\ref{w1-3 solutions}), (\ref{M4 terms}), and (\ref{tilde M4 terms}) and
the fact that $u_\tau={\bf Q}^* F_\tau$, we define
\noextension{
\beq\nonumber
T^{(4),\beta}_{\tau,\sigma}&=&\bra F_\tau,{\bf Q}(\B_4u_{\sigma(4)} \,\cdotp  S_2(\B_3u_{\sigma(3)}\,\cdotp   S_1(\B_2u_{\sigma(2)}\,\cdotp \B_1u_{\sigma(1)})))\cet_{L^2(\hattuM _0)}\\ 
\label{T-type source}
&=&\bra (\B_3u_{\sigma(3)})\,\cdotp S_2^*(  (\B_4u_{\sigma(4)}) \,\cdotp u_\tau), S_1(\B_2u_{\sigma(2)}\,\cdotp \B_1u_{\sigma(1)})\cet_{L^2(\hattuM _0)},\\
\nonumber
\tilde T^{(4),\beta}_{\tau,\sigma}&=&\bra F_\tau,{\bf Q}( S_2(\B_4u_{\sigma(4)} \,\cdotp \B_3u_{\sigma(3)})\,\cdotp  S_1(\B_2u_{\sigma(2)}\,\cdotp \B_1u_{\sigma(1)}))\cet_{L^2(\hattuM _0)}\\ \label{tilde T-type source}
&=&\bra  S_2(\B_4u_{\sigma(4)} \,\cdotp  \B_3u_{\sigma(3)}) \,\cdotp u_\tau,  S_1(\B_2u_{\sigma(2)}\,\cdotp \B_1u_{\sigma(1)})\cet_{L^2(\hattuM _0)}.
\eeq
}
\extension{
\beq\nonumber
T^{(4),\beta}_{\tau,\sigma}&=&\bra F_\tau,{\bf Q}(\B_4u_{\sigma(4)} \,\cdotp  S_2(\B_3u_{\sigma(3)}\,\cdotp   S_1(\B_2u_{\sigma(2)}\,\cdotp \B_1u_{\sigma(1)})))\cet_{L^2(\hattuM _0)}\\ 
\label{T-type source}
&=&\bra {\bf Q}^* F_\tau,\B_4u_{\sigma(4)} \,\cdotp   S_2(\B_3u_{\sigma(3)}\,\cdotp  S_1(\B_2u_{\sigma(2)}\,\cdotp \B_1u_{\sigma(1)}))\cet_{L^2(\hattuM _0)}\\ \nonumber
&=&\bra (\B_4u_{\sigma(4)}) \,\cdotp u_{\tau}, S_2(\B_3u_{\sigma(3)}\,\cdotp  S_1(\B_2u_{\sigma(2)}\,\cdotp \B_1u_{\sigma(1)}))\cet_{L^2(\hattuM _0)}
\\ \nonumber
&=&\bra (\B_3u_{\sigma(3)})\,\cdotp S_2^*(  (\B_4u_{\sigma(4)}) \,\cdotp u_\tau), S_1(\B_2u_{\sigma(2)}\,\cdotp \B_1u_{\sigma(1)})\cet_{L^2(\hattuM _0)},
\eeq
and
\beq
\nonumber
\tilde T^{(4),\beta}_{\tau,\sigma}&=&\bra F_\tau,{\bf Q}( S_2(\B_4u_{\sigma(4)} \,\cdotp \B_3u_{\sigma(3)})\,\cdotp  S_1(\B_2u_{\sigma(2)}\,\cdotp \B_1u_{\sigma(1)}))\cet_{L^2(\hattuM _0)}\\ \label{tilde T-type source}
&=&\bra {\bf Q}^* F_\tau, S_2(\B_4u_{\sigma(4)} \,\cdotp  \B_3u_{\sigma(3)})\,\cdotp  S_1(\B_2u_{\sigma(2)}\,\cdotp \B_1u_{\sigma(1)})\cet_{L^2(\hattuM _0)}\\ \nonumber
&=&\bra  S_2(\B_4u_{\sigma(4)} \,\cdotp  \B_3u_{\sigma(3)}) \,\cdotp u_\tau,  S_1(\B_2u_{\sigma(2)}\,\cdotp \B_1u_{\sigma(1)})\cet_{L^2(\hattuM _0)}.
\eeq
}

When  $\sigma$ is the identity, we will omit it in our notations
and denote $T^{(4),\beta}_{\tau}=T^{(4),\beta}_{\tau,id}$, etc.
It turns out later,  in Prop.\ \ref{SL:order}, that the term
$\bra F_\tau,{\bf Q}(\B_4 u_4\,{\bf Q}(\B_3u_3,{\bf Q}(\B_2u_2\,\B_1u_1)]]\cet$,
where the orders $k_j$ of $\B_j$ are $k_1=6$,
 $k_2=k_3=k_4=0$ is the term that is crucial for our considerations. We  enumerate this term with $\beta=\beta_1:=1$, i.e.,
  \beq\label{beta0 index}
\hbox{$\vec S_{\beta_1}=({\bf Q},{\bf Q})$ and $k^{\beta_1}_1=6$, 
 $k^{\beta_1}_2=k^{\beta_1}_3=k^{\beta_1}_2=0$.}
 \eeq
\MTEXT{ 
This term arises from the term $\bra F_\tau,\G^{(4),\beta}_\sigma\cet$,
such that $\G^{(4),\beta}_\sigma$ is of the form
(\ref{eq: tilde M1}), where $\sigma=id$ and all quadratic forms 
$A^{\beta}_3$, $A^\beta_2$, and $A^\beta_1$ are of the from (\ref{A alpha decomposition2}).
More precisely, 
$\bra F_\tau,\G^{(4),\beta}_\sigma\cet$  is the sum of
$T^{(4),\beta_1}_{\tau,id}$  and terms analogous to that (including
terms with commutators of ${\bf Q}$ and differential operators),
where
\beq\label{T-example}
T^{(4),\beta_1}_{\tau,id}=-\bra {\bf Q}^*((F_\tau)_{nm}),
u_{4}^{rs} \,\cdotp  {\bf Q}(u^{ab}_{3}\,\cdotp   {\bf Q}(u_2^{ik}\,\cdotp 
\p_r\p_s\p_a\p_b\p_i\p_k u^{nm}_{1})))\cet.\hspace{-1.6cm}
\eeq
Here, $u^{ik}_{j}=\hat g^{in}(x)\hat g^{km}(x)\,(u_{j}(x))_{nm}$, $j=1,2,3,4$ 
and $(u_{j}(x))_{nm}$
\MTEXT{is the metric  part of the wave $u_j(x)=(((u_{j}(x))_{nm})_{n,m=1}^4,
((u_{j}(x))_{\ell})_{\ell=1}^L)$.}
}


%
%

\subsubsection{Properties of the indicator functions on a general manifold}
Next we analyze
the indicator function for sources
$ {\bf f}_{j}\in \mathcal I_{C}^{n+1} (Y(x_j,\zeta_j;t_0,s_0))$, $j\leq 4$, 
and the source $F_\tau$ related to
 $(y,\eta)\in L^+M_0$,  $y\in U_{\hat g}$.
We denote in the following 
$(x_5,\xi_5)=(y,\eta)$.
%

\begin{definition}\label{def:  4-intersection of rays} Let $t_0>0$.
We say that the geodesics corresponding to vectors $(\vec x,\vec \xi)=((x_j,\xi_j))_{j=1}^4$
intersect  and  the intersection takes place at the point $q\in \hattuM _0$ 
if there are $t_j\in (0,{\bf t}_j)$, ${\bf t}_j=\rho(x_j(t_0),\xi_j(t_0)$ such that 
 $q=\gamma_{x_j,\xi_j}(t_j)$  
 for all $j=1,2,3,4$. 
 \end{definition}
 
  Let  $\Lambda_q^+$ be the Lagrangian manifold 
  \ba
  \Lambda_q^+=\{(x,\xi)\in T^*M_0;\ x=\gamma_{q,\zeta}(t),\ 
  \xi^\sharp=\dot\gamma_{q,\zeta}(t),\ \zeta\in L^+_qM_0,\ t>0\}
  \ea
  Note that the projection of $ \Lambda_q^+$
 on $\hattuM _0$ is  the light cone $\L^+_{\hat g}(q)$. 

Next we consider $x_j\in U_{\hat g}$ and $ \xi_j\in L^+_{x_j}\hattuM _0$,
 and $\vartheta_1,t_0>0$ such that
  $(\vec x,\vec \xi)=((x_j,\xi_j))_{j=1}^4$
  satisfy, see Fig.\ 4(Right),
\beq\label{eq: summary of assumptions 1}& &\hspace{-1.0cm}(i)\ \gamma_{x_j,\xi_j}([0,t_0])\subset W_{\hat g},\ 
 x_j(t_0)\not \in J^+_{\hat g}(x_k(t_0)),
  \hbox{ for $j,k\leq 4$, $j\not =k$,}  \hspace{-5mm}\\ & &\nonumber
 \hspace{-1.0cm}(ii)\ 
\hbox{For all $j,k\leq 4$,
 $d_{\hat g^+}((x_j,\xi_j),(x_k,\xi_k))<\vartheta_1$,} \hspace{-5mm}\\
 & &\nonumber
 \hspace{-1.0cm}(iii)\ 
\hbox{There is $\hat x\in  \hat \mu$ such that \ 
for all $j\leq 4$,
 $d_{\hat g^+}(\hat x,x_j)<\vartheta_1$,} \hspace{-5mm}
 \eeq
 Above, 
$(x_j(h),\xi_j(h))$ are defined in (\ref{eq: x(h) notation}).
We denote
 \beq\label{eq: summary of assumptions 2}
& & \V((\vec x,\vec \xi),t_0)= M_0\setminus \bigcup_{j=1}^4
J^+_{\hat g}(\gamma_{x_j(t_0),\xi_j(t_0)}({\bf t}_j)),
 \eeq
where ${\bf t}_j:=\rho(x_j(t_0),\xi_j(t_0)).$}
  Note that
 two geodesics $\gamma_{x_j(t_0),\xi_j(t_0)}([0,\infty))$ can intersect only once in $\V((\vec x,\vec \xi),t_0)$.
\motivation{We will analyze the 4-wave interaction of waves in the set $\V((\vec x,\vec \xi),t_0)$
where all observed singularities are produced before the 
geodesics $\gamma_{x_j(t_0),\xi_j(t_0)}([0,\infty))$ have conjugate points, that is,
before the waves $u_{(j)}$ have caustics.}
Below we use $\sim$ for  the terms that
have the same  asymptotics up to an error $O(\tau^{-N})$ for all $N>0$, as $\tau \to \infty$.

}

%

  \begin{proposition}\label{lem:analytic limits A}    Let $(\vec x,\vec \xi)=((x_j,\xi_j))_{j=1}^4$  be future pointing light-like vectors 
   satisfying (\ref{eq: summary of assumptions 1})
   and $t_0>0$.
  Let $x_5\in \V((\vec x,\vec \xi),t_0)\cap U_{\hat g}$ and $(x_5,\xi_5)$ be a future pointing light-like vector
such that
 $x_5\not\in {\mathcal Y}((\vec x,\vec \xi);t_0)\cup\bigcup_{j=1}^4\gamma_{x_j,\xi_j}(\R)$,
see (\ref{Y-set}) and (\ref{eq: summary of assumptions 2}).
When  $n\in \Z_+$ is large enough and $s_0>0$ is small enough, 
the function $\Theta^{(4)}_\tau$, see (\ref{test sing}),
 corresponding to  
 $ {\bf f}_{j}\in \mathcal I_{C}^{n+1} (Y(x_j,\zeta_j;t_0,s_0))$, $j\leq 4$, 
and 
 the source $F_\tau$, see (\ref{Ftau source}),  
corresponding to 
$(x_5,\xi_5)$
satisfy

%
%
 \medskip

 \noindent
 (i) If the geodesics corresponding to $(\vec x,\vec\xi)$ either
 do not intersect or intersect at $q$ and $(x_5,\xi_5)\not \in \Lambda^+_q$, then
   $|\Theta^{(4)}_\tau|\leq C_N\tau^{-N}$
 for all $N>0$.

\noindent
(ii) If the geodesics corresponding to $(\vec x,\vec\xi)$  intersect at $q$,
 $q=\gamma_{x_j,\xi_j}(t_j)$ and  the vectors $\dot \gamma_{x_j,\xi_j}(t_j)$, $j=1,2,3,4$ are linearly independent,
{and there are  $t_5<0$ and $\xi_5\in L^+_{x_5}M_0$ such that $q=\gamma_{x_5,\xi_5}(t_5)$,} 
then, with $m=-4n+2$, 
\beq\label{indicator}
\Theta^{(4)}_\tau\sim  
\sum_{k=m}^\infty s_{k}\tau^{-k},\quad\hbox{as $\tau\to \infty$.}
\eeq



Moreover,
let $b_j=(\dot\gamma_{x_j,\xi_j}(t_j))^\flat$, $j=1,2,3,4,5$, and
 $\bsequence=(b_{j})_{j=1}^5\in (T^*_q\hattuM _0)^5$. Let
 $w_j$ be the principal symbols of the waves $u_j={\bf Q}f_j$ at $(q,b_j)$ for $j\leq 4$. 
{Also, let $w_5$ be the principal symbol of $u_\tau={\bf Q}F_\tau$
 at $(q,b_5)$,} and ${\bf w}=(w_j)_{j=1}^5$. 
 Then there is 
  a real-analytic function $\mathcal G(\bsequence,{\bf w})$ such that
 the leading order term in (\ref{indicator})  satisfies 
 \beq\label{definition of G}
s_{m}=
\mathcal  G(\bsequence,{\bf w}).
\eeq

(iii) Under the assumptions in (ii), the point $x_5$ has a neighborhood $V$ such that 
$\M^4$ in $V$ satisfies $\M^4|_V\in  \I(\Lambda^+_q)$.
\end{proposition}

\extension{

\noindent
{\bf Proof.}  
Below, to simplify notations,
we denote $K_j=K(x_{j},\xi_{j};t_0,s_0)$ and 
$K_{123}=K_1\cap K_2\cap K_3$ and $K_{124}=K_1\cap K_2\cap K_4$, etc.
We will denote $ \Lambda_j=\Lambda(x_j,\xi_j;t_0,s_0)$
to consider also the singularities of $K_j$ related to conjugate points.

We will consider below separately the case when the
following linear independency condition is satisfied:
\medskip

(LI) Assume if that  $J\subset \{1,2,3,4\}$ and $y\in J^-(x_6)$
are such that for all $j\in J$  we have $\gamma_{x_j,\xi_j}(t^\prime_j)=y$
with some $t^\prime_j\geq 0$, then  the vectors
$\dot\gamma_{x_j,\xi_j}(t_j^\prime )$, $j\in J$ are linearly independent.
\medskip

Let us first  consider  the case when (LI) is valid.

By the definition of ${\bf t}_j$,  if
the intersection $\gamma_{x_5,\xi_5}(\R_-) \cap(\cap_{j=1}^4\gamma_{x_j,\xi_j}((0,{\bf t}_j)))
$ is non-empty, it can
contain only one point. In the case that such a
point exists, we denote it by $q$. When $q$ exists,
the intersection of $K_j$ at this point are transversal and
we see that when $s_0$ is small enough, the set
$\cap_{j=1}^4 K_j$ consists only of the point $q$.
We assume below that $s_0$ is such that  this is true.

We define    two local coordinates $Z:W_0\to \R^4$
and  $Y:W_1\to \R^4$ such that $W_0,W_1\subset \V((\vec x,\vec \xi),t_0)$,
 see definition after (\ref{eq: summary of assumptions 2}).
 We assume that these local coordinates are such that
$K_j\cap W_0=\{x\in W_0;\ Z^j(x)=0\}$ and $K_j\cap W_1=\{x\in W_1;\ Y^j(x)=0\}$ for $j=1,2,3,4$.
In the Fig.\ 5, $W_0$ is a neighborhood of $z$ and
$W_1$ is a neighborhood of $y$. 
We note that the origin $0=(0,0,0,0)\in \R^4$ does not
necessarily belong to the set $Z(W_0)$ or the set $Y(W_1)$,
for instance in the case when the four geodesics $\gamma_{x_j,\xi_j}$ do not intersect.
However, this is the case when the  four geodesics intersect at the point $q$, we need
to consider the case when $W_0$ and $W_1$ are neighborhoods of $0$. 
 %
Note that
$W_0$ and $W_1$ do not contain any cut points of the geodesic
 $\gamma_{x_j,\xi_j}([t_0,\infty)$.

We will denote $z^j=Z^j(x)$.
We assume that  $Y:W_1\to \R^4$ are similar coordinates
and denote  $y^j=Y^j(x)$. We also denote
 below $dy^j=dY^j$ and  $dz^j=dZ^j$.
Let $\Phi_0\in C^\infty_0(W_0)$ and 
 $\Phi_1\in C^\infty_0(W_1)$.  
 
 Let us next considering the map ${\bf Q}^*:C^\infty_0(W_1)\to
 C^\infty(W_0).$  By \cite{MU1}, ${\bf Q}^*\in I(W_1\times W_0;
   \Delta_{T\hattuM _0}^{\prime}, \Lambda_{\hat g})$ is an  operator 
   with a classical symbol and its canonical
   relation 
$\Lambda_{{\bf Q}^*}^{\prime}$ is associated to a union of two  intersecting Lagrangian manifolds,
$\Lambda_{{\bf Q}^*}^{\prime}= \Lambda_{\hat g}^{\prime}\cup \Delta_{T\hattuM _0}$,
intersecting cleanly \cite{MU1}.
Let $\e_2>\e_1>0$ and $B_{\e_1,\e_2}$ be a pseudodifferential
operator on $\hattuM _0$ which is  a microlocally smoothing
operator (i.e., the full  symbol vanishes  in local coordinates)
outside of the $\e_2$-neighborhood $\V_2\subset T^*\hattuM _0$ of the set of the light like covectors
$L^*\hattuM _0$ and for which 
$(I-B_{\e_1,\e_2})$ is microlocally smoothing
operator 
 in the $\e_1$-neighborhood $\V_2$ of $L^*\hattuM _0$.
 The neighborhoods here are defined 
 with respect to the Sasaki metric of $(T^*\hattuM _0,\hat g^+)$ and
$\e_2,\e_1$ are chosen later in the proof. 
 Let us decompose the operator ${\bf Q}^*={\bf Q}^*_1+{\bf Q}^*_2$
 where  ${\bf Q}^*_1={\bf Q}^*(I-B_{\e_1,\e_2})$ and   ${\bf Q}^*_2={\bf Q}^*B_{\e_1,\e_2}$.
 As  $\Lambda_{{\bf Q}^*}= \Lambda_{\hat g}\cup \Delta_{T\hattuM _0}^{\prime}$,
we see that then there is 
a neighborhood $ \W_2=\W_2(\e_2)$ of $L^*\hattuM _0\times L^*\hattuM _0\subset (T^*\hattuM _0)^2$
such that the Schwartz
kernel ${\bf Q}^*_2(r,y)$ of the operator ${\bf Q}^*_2$ satisfies
\beq\label{eq: W_2 neighborhood}
\hbox{WF}({\bf Q}^*_2)\subset \W_2.
\eeq
Moreover,  $\Lambda_{{\bf Q}^*_1}\subset \Delta_{T\hattuM _0}^{\prime}$ implying
that ${\bf Q}^*_1$ is a pseudodifferential operator with a classical symbol,
 ${\bf Q}^*_1\in I(W_1\times W_0;
   \Delta_{T\hattuM _0}^{\prime})$,
 and
 ${\bf Q}^*_2 
 \in I(W_1\times W_0;
   \Delta_{T\hattuM _0}^{\prime}, \Lambda_{\hat g})$
   is a Fourier integral operator (FIO) associated
   to two cleanly intersecting Lagrangian manifolds, similarly to ${\bf Q}^*$. 
In the case when  $p=1$ we can write ${\bf Q}^*_p$ as 
\beq\label{eq: representation of Q* p}
({\bf Q}^*_1v)(z)=\int_{\R^{4+4}}e^{i(y-z)\cdotp \xi}q_1(z,y,\xi)v(y)\,dyd\xi,
\eeq
where 
{
a classical symbol $q_1 (z,y,\xi)\in S^{-2}(W_1\times W_0; \R^4)$ with a real
valued principal symbol
\ba
\hat q_1 (z,y,\xi)=\frac {\chi(z,\xi)}{\hat g^{jk}(z)\xi_j\xi_k}
\ea
where $\chi(z,\xi)\in \C^\infty$ is a cut-off function vanishing
in a neighborhood of  the set where $g(\xi,\xi)=0$. Note that
then $Q_1-Q_1^*\in \Psi^{-3}(W_1\times W_0)$.}

Furthermore, let us decompose 
${\bf Q}^*_1={\bf Q}^*_{1,1}+{\bf Q}^*_{1,2}$ corresponding to the decomposition
$q_1(z,y,\xi)=q_{1,1}(z,y,\xi)+q_{1,2}(z,y,\xi)$
 of the symbol,
where
 \beq\label{eq: q_{1,1} and q_{1,2}} \hspace{-.3cm}
 & &q_{1,1}(z,y,\xi)=q_{1}(z,y,\xi)\psi_R(\xi),\\
 & & \nonumber q_{1,2}(z,y,\xi)=q_{1}(z,y,\xi)(1-\psi_R(\xi)),
 \hspace{-2cm}
 \eeq
and $\psi_R\in C^\infty_0(\R^4)$
is a cut-off function that is equal to one in  a ball $B(R)$ of radius $R$ specified below.

Next we start to consider the terms $T_{\tau}^{(4),\beta}$ and $\tilde T_{\tau}^{(4),\beta}$ of the type
 (\ref{T-type source}) 
and
(\ref{tilde T-type source}). In these terms,
we can represent the gaussian beam $u_{\tau}(z)$ in $W_1$ in the form 
\beq\label{eq:representation of the wave}
u_\tau(y)=e^{i\tau \varphi (y)}a_5(y,\tau)
\eeq
where the function $ \varphi$ is a complex phase function having non-negative imaginary
part such that $\im   \varphi$, defined on $W_1$, vanishes exactly on the geodesic
 $\gamma_{x_5,\xi_5}\cap W_1$. Note that 
  $\gamma_{x_5,\xi_5}\cap W_1$ may be empty. Moreover,
 $a_5\in
S^{0}_{clas}(W_1;\R)$ is a classical symbol.

 Also, for $y=\gamma_{x_5,\xi_5}(t)\in W_1$
we have that $d\varphi(y)=c\dot \gamma_{x_5,\xi_5}(t)^\flat$,  where $c\in \R\setminus \{0\}$, is light-like.



We consider first the asymptotics of 
 terms $T_\tau^{(4),\beta}$ and $\tilde T^{(4),\beta}_\tau$ of the type
 (\ref{T-type source}) 
and
(\ref{tilde T-type source})
   where $S_1=S_2={\bf Q}$ and, the symbols $a_j(z,\theta_j)$, $j\leq 3$ and  
    $a_j(y,\theta_j)$, $j\in \{4,5\}$ are scalar valued
   symbols written in the $Z$ and $Y$ coordinates.
$\B_1,\B_2,\B_3$ are multiplication operators with $\Phi_0$,
and $\B_4$ is a  multiplication operators with $\Phi_1$
 and consider section-valued symbols and general operators later.

Let us consider functions
$U_j\in \I^{p_j}(K_j)$, $j=1,2,3,$ supported in $W_0$ and 
$U_4\in \I(K_4)$,  supported in $W_1$,  having classical symbols. They have the form
 \beq\label{U1term}
& &U_j(x)=\int_{\R}e^{i\theta_jx^j}a_j(x,\theta_j)d\theta_j,\quad a_j\in
S^{p_j}_{clas}(W_{\kappa(j)};\R), 
\eeq
for all $j=1,2,3,4$ (Note that here the phase function is  $\theta_jx^j=\theta_1x^1$ for $j=1$ etc,
that is, there is no summing over index $j$). We may assume that $a_j(x,\theta_j)$ vanish near $\theta_j=0$.

Since  $x_5\in \V((\vec x,\vec \xi),t_0)\cap U_{\hat g}$ and 
$W_0,W_1\subset \V((\vec x,\vec \xi),t_0)$, we see that
an example of functions (\ref{U1term}) are $U_j(z)=\Phi_{\kappa(j)}(z)u_j(z)$, $j=1,2,3,4$,
where $u_j(z)$ are the distorted plane  waves. Here and below, 
 $\kappa(j)=0$ for $j=1,2,3$  and $\kappa(4)=1$ and we
also denote $\kappa(5)=1$.
Note that $p_j=n$ correspond to the case when $U_j\in \I^n(K_j)=\I^{n-1/2}(N^*K_j)$.

Denote $\Lambda_j=N^*K_j$ and $\Lambda_{jk}=N^*(K_j\cap K_k)$. By \cite[Lem.\ 1.2 and 1.3]{GU1}, the  
pointwise product 
$U_2\,\cdotp U_1\in \I(\Lambda_1,\Lambda_{12})+ \I(\Lambda_2,\Lambda_{12})$  and thus
by \cite[Prop.\ 2.2]{GU1},
${\bf Q}(U_2\,\cdotp U_1)\in \I(\Lambda_1,\Lambda_{12})+ \I(\Lambda_2,\Lambda_{12})$ 
and  it
can be written as 
\beq\label{Q U1U2term}
{\bf Q}(U_2\,\cdotp U_1)=\int_{\R^2}e^{i(\theta_1z^1+\theta_2z^2)}c_1(z,\theta_1,\theta_2)d\theta_1d\theta_2.
\eeq
Note that here $c_1(z,\theta_1,\theta_2)$ is sum of product type symbols, see (\ref{product symbols}).
As  $ N^*(K_1\cap K_2)
\setminus (N^*K_1\cup N^*K_2)$ consists of vectors which
are non-characteristic for 
the wave operator, that is, the wave operator is elliptic in
a neighborhood of this subset of the cotangent bundle,   the principal symbol 
$\hat c_1$ of ${\bf Q}(U_2\,\cdotp U_1)$ on  $N^*(K_1\cap K_2)
\setminus (N^*K_1\cup N^*K_2)$ is  given by
\beq\label{product type symbols}
& &\hat  c_1(z,\theta_1,\theta_2)\sim
 s(z,\theta_1,\theta_2)a_1(z,\theta_1)a_2(z,\theta_2),
 \\ \nonumber
 & &\hspace{-1cm}s(z,\theta_1,\theta_2)
=1/{ \hat g(\theta_1b^{(1)}+\theta_2b^{(2)},\theta_1b^{(1)}+\theta_2b^{(2)})}
 =1/({2 \hat g(\theta_1b^{(1)},\theta_2b^{(2)})}). \hspace{-1cm}
 \eeq
Note that $s(z,\theta_1,\theta_2)$ is a smooth function  
on $N^*(K_1\cap K_2)
\setminus (N^*K_1\cup N^*K_2)$ and homogeneous of order $(-2)$ in $\theta=(\theta_1,\theta_2)$.
Here, we use $\sim$ to denote that the symbols have the same principal symbol.
%
%
%
%
Let us next
make computations in the case when 
 $a_j(z,\theta_j)\in C^\infty(\R^4\times \R)$ is positively homogeneous for $|\theta_j|>1$, that is, we have 
 $ a_j(z,s)=a^{\prime}_j(z)s^{p_j}$, where $p_j\in \N$  and $|s|>1$.
 We consider $T_\tau^{(4),\beta}=\sum_{p=1}^2T_{\tau,p}^{(4),\beta}$
 where $T_{\tau,p}^{(4),\beta}$ is defined as 
 $T_{\tau}^{(4),\beta}$ by replacing the term ${\bf Q}^*(U_4\,\cdotp u_{\tau})$
  by ${\bf Q}^*_p(U_4\,\cdotp u_{\tau})$.

 Let us now consider the case $p=2$ and choose the parameters $\e_1$ and $\e_2$
that determine the
decomposition ${\bf Q}^*={\bf Q}^*_1+{\bf Q}^*_2$. 
First, we observe that for $p=2$ we can write using $Z$ and $Y$ coordinates
 \beq\label{eq; T tau asympt, p=2}
T_{\tau,2}^{(4),\beta}
&=&
\tau^{4}
\int_{\R^{12}}e^{i\tau \Psi_2(z,y,\theta)}
c_1(z,\tau \theta_1,\tau \theta_2)\cdotp\\ \nonumber
& &\hspace{-1.3cm}
\cdotp  a_3(z,\tau \theta_3)
{\bf Q}^*_2(z,y)
a_4(y,\tau \theta_4) a_5(y,\tau )\,d\theta_1d\theta_2d\theta_3d\theta_4 dydz,\\
\nonumber \Psi_2(z,y,\theta)
&=&\theta_1z^1+\theta_2z^2+\theta_3z^3+\theta_4y^4+\varphi(y).
\eeq

Denote $\theta=(\theta_1,\theta_2,\theta_3,\theta_4)\in \R^{4}$.
Consider the case when $(z,y,\theta)$ is a critical point of $\Psi_2$
satisfying $\im \varphi(y)=0$. Then we have
 $\theta^{\prime}=(\theta_1,\theta_2,\theta_3)=0$
and $z^{\prime}=(z^1,z^2,z^3)=0$, 
$y^4=0$,
  $d_y\varphi(y)=(0,0,0,-\theta_4)$,
implying that $y\in K_4$ and 
$(y,d_y\varphi(y))\in N^*K_4$.
As  $\im \varphi(y)=0$, we have that 
  $y=\gamma_{x_5,\xi_5}(\tsparameter)$ with some $\tsparameter\in\R_-$.
 As we have $\dot \gamma_{x_5,\xi_5}(\tsparameter)^\flat=d_y\varphi(y)\in N^*_yK_4$,
we obtain $\gamma_{x_5,\xi_5}([\tsparameter,0])\subset K_4$. However, 
this 
is not possible by our assumption 
$x_5\not \in \cup_{j=1}^4 K(x_{j},\xi_{j};s_0)$ when $s_0$ is small enough.
Thus  the phase function 
$ \Psi_2(z,y,\theta)$ has no critical points satisfying $\im \varphi(y)=0$.

%

When the orders $p_j$ of the symbols $a_j$ are small enough, the integrals
in the $\theta$ variable in  (\ref{eq; T tau asympt, p=2}) are convergent in the classical sense.
We use  now properties of the wave front set to 
compute the asymptotics of  oscillatory integrals
and to this end we introduce the function
\beq
\tilde {\bf Q}^*_2(z,y,\theta)={\bf Q}^*_2(z,y),
\eeq
that is, consider ${\bf Q}^*_2(z,y)$ as a constant function
in $\theta$. Below, denote $\psi_4(y,\theta_4)=\theta_4y^4$ and $r=d\varphi(y)$. 
Note that then $d_{\theta_4}\psi_4=y^4$ and $d_y\psi_4=(0,0,0,\theta_4)$.
Then  in $W_1\times W_0\times \R^4$
\ba
d_{z,y,\theta}\Psi_2&=&(\theta_1,\theta_2,\theta_3,0;
r+d_y\psi_4(y,\theta_4),
z^1,z^2,z^3,d_{ \theta_4}\psi_4(y,\theta_4))
\\
&=&(\theta_1,\theta_2,\theta_3,0;
d\varphi(y)+(0,0,0,\theta_4),
z^1,z^2,z^3,y^4)
\ea
and we see that if $((z,y,\theta),d_{z,y,\theta}\Psi_2)\in \hbox{WF}(\tilde {\bf Q}_2^*)$
and  $\im \varphi(y)=0$, we have $(z^1,z^2,z^3)=0$, $y^4=d_{\theta_4}\psi_4(y,\theta_4)=0$
and $y\in  \gamma_{x_5,\xi_5}$. Thus 
$z\in K_{123}$ and $y\in  \gamma_{x_5,\xi_5}\cap K_4$.

Let  us use the following notations
\beq\label{new notations}
& &z\in K_{123},\quad \omega_\theta:=(\theta_1,\theta_2,
 \theta_2,0)=\sum_{j=1}^3\theta_j dz^j\in T^*_z\hattuM _0,\\
 \nonumber
& &y\in K_4\cap \gamma_{x_5,\xi_5},\quad (y,w):=(y,d_y\psi_4(y,\theta_4))\in N^*K_4,
\\
 \nonumber
& & r=d\varphi(y)=r_jdy^j\in T^*_y\hattuM _0,\quad \kappa:=r+w.
 \eeq
Then, $y$ and $ \theta_4$ satisfy $y^4=d_{ \theta_4}\psi_4(y,\theta_4)=0$ and
$w=(0,0,0,\theta_4)$.
 
 Note that  by definition of the $Y$  coordinates $w$ is a light-like covector.
 By definition of the  $Z$ coordinates,
 $\omega_\theta\in 
 N^*K_1+N^*K_2+N^*K_3=N^*K_{123}$.  
 
 Let us first consider what happens if $\kappa=r+w=d\varphi(y)+(0,0,0,\theta_4)$ is light-like.
In this case, all vectors $\kappa$, $w$, and $r$ are light-like
and satisfy $\kappa =r+w$. This is possible
only if  $r\parallel w$, i.e., $r$ and $w$ are parallel, see \cite[Cor.\ 1.1.5]{SW}.  Thus  $r+w$ is light-like if and only
if $r$ and $w$ are parallel.

Consider next the case when  $(x,y,\theta)\in W_1\times W_0\times \R^4$ 
is such that $((x,y,\theta),d_{z,y,\theta}\Psi_2)\in \hbox{WF}(\tilde {\bf Q}_2^*)$
and   $\im \varphi(y)=0$.
Using the above notations (\ref{new notations}), we obtain
$
d_{z,y,\theta}\Psi_2=(\omega_\theta,
r+w;(0,0,0,d_{ \theta_4}\psi_4(y, \theta_4)))=(\omega_\theta,d\varphi(y)+(0,0,0,\theta_4);(0,0,0,y^4))$,
 where $y^4=d_{ \theta_4}\psi_4(y, \theta_4)=0$,
and thus we have
\ba
((z,\omega_\theta),(y,
r+w))\in \hbox{WF}( {\bf Q}_2^*)=\Lambda_{{\bf Q}^*_2}.
\ea

\begin{center}

\psfrag{1}{$(x_1,\xi_1)$}
\psfrag{2}{\hspace{-5mm}$(x_2,\xi_2)$}
\psfrag{3}{\hspace{-5mm}$(x_3,\xi_3)$}
\psfrag{4}{\hspace{-5mm}$(x_4,\xi_4)$}
\psfrag{5}{\hspace{-5mm}$(x_5,\xi_5)$}
\psfrag{6}{$z$}
\psfrag{7}{$y$}
\psfrag{8}{$r$}
\psfrag{9}{$w$}
\psfrag{0}{$\omega_\theta$}
\includegraphics[height=5.5cm]{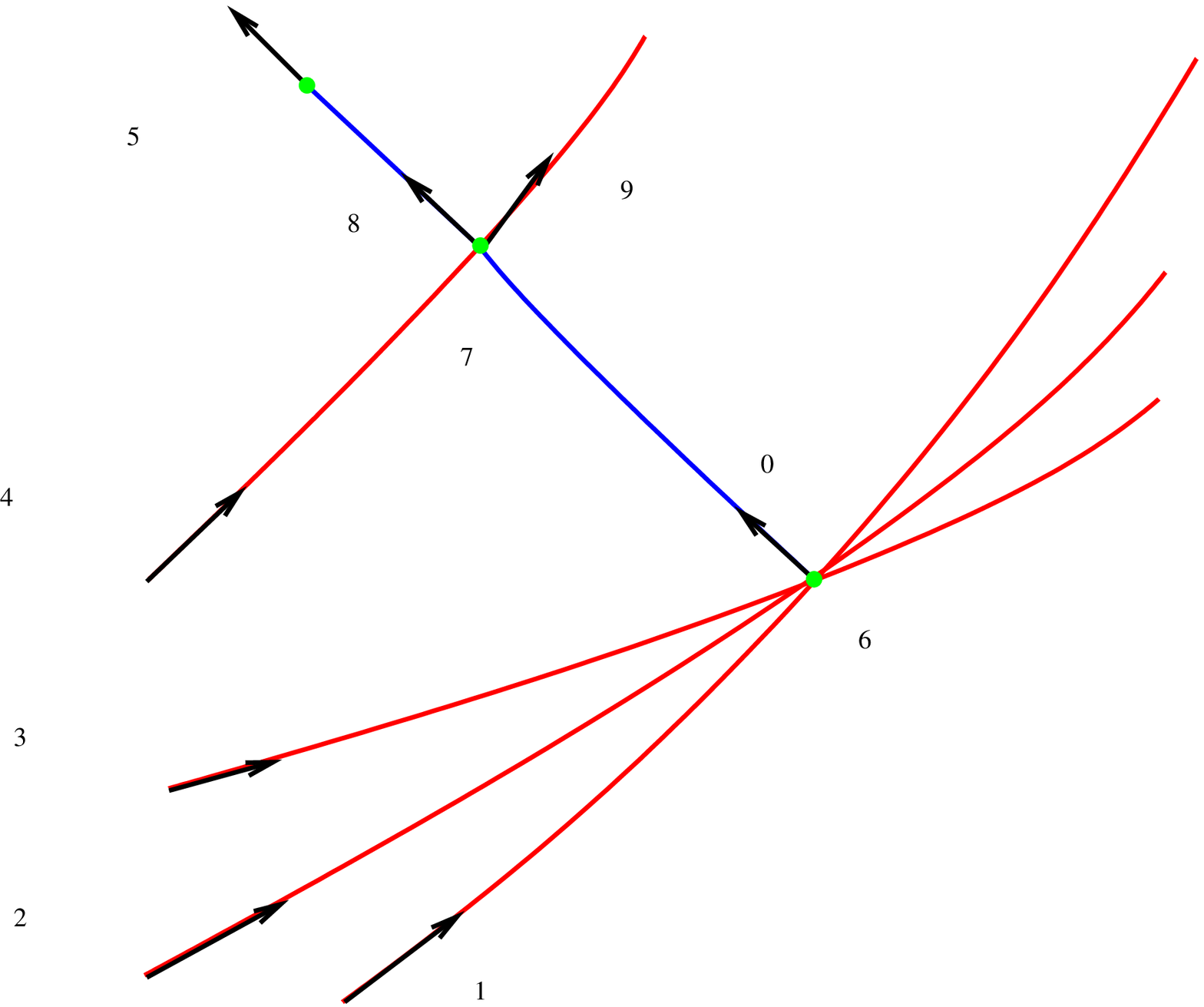}
\psfrag{1}{${\hat x}^\prime$}  
\psfrag{2}{$p_2$}  
\psfrag{3}{$J({p^-},{p^+})$}  
\psfrag{4}{${\hat x}$}  
\psfrag{5}{}  
\psfrag{6}{${\hat x}_1$}  
\psfrag{6}{${\hat x}_1$}  
\psfrag{7}{$\tilde p$}  
\psfrag{8}{$$}  
\includegraphics[width=5.5cm]{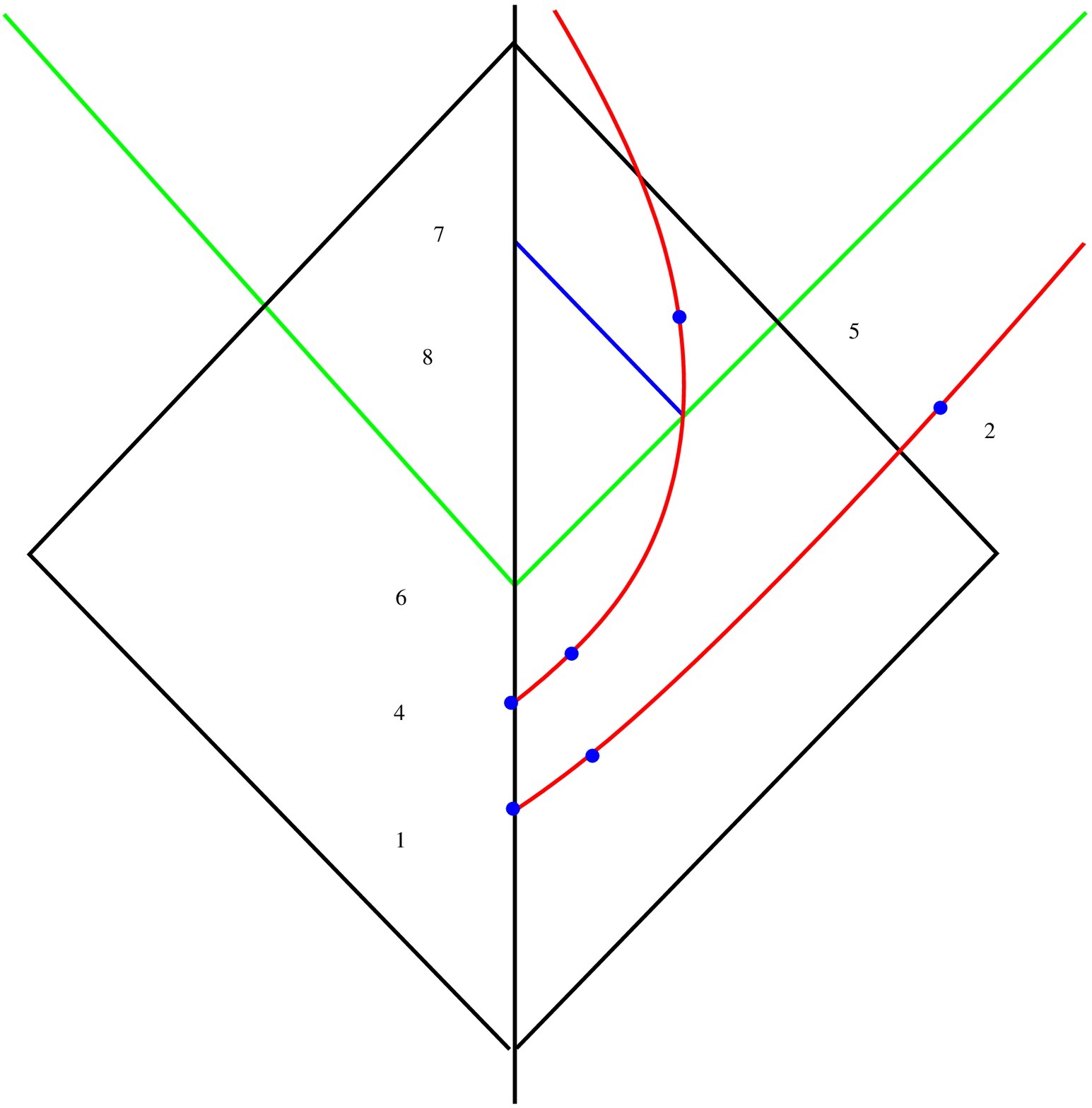}  

\end{center}

\noindent
{\it FIGURE 5. {\bf Left:} In the figure we consider the case A1 where three geodesics intersect at $z$ and 
the waves propagating near these geodesics interact and create a wave that hits
 the fourth geodesic at the point $y$, the produced singularities are detected by the gaussian beam
source at the point $x_5$. Note that $z$ and $y$ can be conjugate points on the geodesic connecting
them. In the case A2 the points $y$ and $z$ are the same.
{\bf Right:} Condition I  is valid.
}
\medskip

Since $\Lambda_{{\bf Q}^*_2}\subset \Lambda_{\hat g}\cup \Delta_{T\hattuM _0}^{\prime}$,
this implies that one of the following conditions are valid:
\ba
& &(A1)\ ((z,\omega_\theta),(y,
r+w))\in \Lambda_{\hat g}, \hspace{5cm} \
\\ & &\hspace{-0.5cm}\hbox{or}\\ 
& &(A2)\ ((z,\omega_\theta),(y,
r+w))\in \Delta_{T\hattuM _0}^{\prime}.
\ea
Let
$\gamma_0$ 
be the geodesic with
$
\gamma_0(0)=z,\ \dot \gamma_0(0) =\omega_\theta^\sharp.
$
Then (A1) and (A2) are equivalent to the following conditions:
\ba
& &(A1)\ \hbox{There is $\tsparameter\in \R$ such
that $(\gamma_0(\tsparameter), \dot \gamma_0(\tsparameter)^\flat) =(y,\kappa)$ and},\\
& &\quad\quad \ \ \hbox{
 the vector $\kappa$ is  light-like,}
\\ & &\hspace{-0.5cm}\hbox{or}\\ 
& &(A2)\ z=y\hbox{ and }
\kappa=-\omega_\theta.
\ea

Consider next the case when (A1)  is valid. 
As  $\kappa$ is light-like,  $r$ and $w$ are parallel. 
Then, since $(\gamma_{x_5,\xi_5}(t_1),\dot \gamma_{x_5,\xi_5}(t_1)^\flat)=(y,r)$ we see that $\gamma_0$ is a continuation of the geodesic $ \gamma_{x_5,\xi_5}$,
that is, for some $t_2$ we
have $(\gamma_{x_5,\xi_5}(t_2),\dot \gamma_{x_5,\xi_5}(t_2))=(z,\omega_\theta)
\in N^*K_{123}$. This implies that $x_5\in {\mathcal Y}$  that is not possible
by our assumptions. Hence (A1) is not possible.

Consider next the case when (A2)  is valid.Then we would also have 
that $r\parallel w$  then $r$ is parallel to $\kappa=-\omega_\theta\in N^*K_{123}$,
 and since $(\gamma_{x_5,\xi_5}(t_1),\dot \gamma_{x_5,\xi_5}(t_1)^\flat)=(y,r)$ we
 would have  that $x_5\in {\mathcal Y}$. As this is not possible by our assumptions,
 we see that $r$ and $w$ are not parallel. This implies that
 $\omega_\theta=-\kappa$ is not light-like.


For any given  $(\vec x,\vec \xi)$  and $(x_5,\xi_5)$ there exists
 $\e_2>0$ so that $(\{(y, \omega_\theta)\}\times T^*M)\cap \W_2(\e_2)=\emptyset$,
 see (\ref{eq: W_2 neighborhood}), and thus
 \beq\label{eq: y and omega}
 (\{(y, \omega_\theta)\}\times T^*M)\cap \hbox{WF}({\bf Q}_2^*)=\emptyset.
 \eeq
 Next we assume that $\e_2>0$ and also $\e_1\in (0,\e_2)$ are chosen
 so that  (\ref{eq: y and omega}) is valid.
 Then there are no $(z,y,\theta)$ such that 
  $((z,y,\theta),d\Psi_2(z,y,\theta))\in  \hbox{WF}(\tilde {\bf Q}_2^*)$ 
 and  $\im \varphi(y)=0$. Thus by 
  Corollary 1.4 in \cite{Delort} or
  \cite[Lem.\ 4.1]{Ralston}
 yields  $T_{\tau,2}^{(4),\beta}=O(\tau^{-N})$ for all $N>0$.
Alternatively, one can use 
 the complex version of \cite[Prop. 1.3.2]{Duistermaat}, 
 obtained using combining the proof of  \cite[Prop. 1.3.2]{Duistermaat} and 
 the method of stationary phase with a complex phase, see \cite[Thm.\ 7.7.17]{H1}.
 
 Thus to analyze the asymptotics of  $T_{\tau}^{(4),\beta}$ 
 we need to consider only $T_{\tau,1}^{(4),\beta}$.
Next, we analyze the case when $U_4$ is a conormal
distribution and has the form (\ref{U1term}).

%
%

 Let us thus consider the case $p=1$. Now
 \beq
U_4(y)\,\cdotp u_{\tau}(y)&=&
\int_{\R^{1}}e^{i\theta_4y^4+i\tau \varphi(y)}
a_4(y,\theta_4)a_5(y,\tau )
\,d\theta_4. \label{QU4Ut term PRE}
\eeq
 We obtain
by (\ref{eq: representation of Q* p}) 
\beq\nonumber
({\bf Q}^*_1(U_4\,\cdotp u_{\tau}))(z)&=&
\int_{\R^{9}}e^{i((y-z)\cdotp \xi+\theta_4y^4+\tau \varphi(y))}q_1(z,y,\xi)\cdotp\\
& &\quad
\cdotp a_4(y,\theta_4)a_5(y,\tau )
\,d\theta_4\,dy d\xi. \label{QU4Ut term}
\eeq
Then
 $T_{\tau,1}^{(4),\beta}=T_{\tau,1,1}^{(4),\beta}+T_{\tau,1,2}^{(4),\beta}$, 
 cf.\ (\ref{eq: q_{1,1} and q_{1,2}}),
 where 
  \beq\label{eq; T tau asympt pre}
& &\hspace{-.4cm}T_{\tau,1,k}^{(4),\beta}
=
\int_{\R^{16}}e^{i(\theta_1z^1+\theta_2z^2+\theta_3z^3+(y-z)\cdotp \xi+\theta_4y^4+\tau \varphi(y))}
c_1(z,\theta_1,\theta_2)\cdotp
\hspace{-1cm}
\\ \nonumber
& &\cdotp  a_3(z,\theta_3)
q_{1,k}(z,y,\xi)a_4(y,\theta_4) a_5(y,\tau )\,d\theta_1d\theta_2d\theta_3d\theta_4dydz d\xi,
\eeq
or
 \beq\label{eq; T tau asympt}
& &\hspace{-.1cm}T_{\tau,1,k}^{(4),\beta}=\tau^{8}
\int_{\R^{16}}e^{i\tau (\theta_1z^1+\theta_2z^2+\theta_3z^3+(y-z)\cdotp \xi+\theta_4y^4+\varphi(y))}
c_1(z,\tau \theta_1,\tau \theta_2)\cdotp\hspace{-1cm}\\ \nonumber
& &\hspace{-0.3cm}
\cdotp  a_3(z,\tau \theta_3)
q_{1,k}(z,y,\tau \xi)a_4(y,\tau \theta_4) a_5(y,\tau )\,d\theta_1d\theta_2d\theta_3d\theta_4dydz d\xi.\hspace{-1cm}
\eeq

%
Let $(\cirt z,\cirt \theta,\cirt y,\cirt \xi)$ be a critical point of
 the phase function
  \beq\label{Psi3 function}
 \Psi_3(z,\theta,y,\xi)=\theta_1z^1+\theta_2z^2+\theta_3z^3+(y-z)\,\cdotp \xi+\theta_4y^4+\varphi(y).\hspace{-1cm}
 \eeq
 Denote $\cirt  w=(0,0,0,\cirt \theta_4)$ and 
$ \cirt r=d\varphi(\cirt y)=\cirt r_jdy^j $. 
Then
 \beq\label{eq: critical points} 
  & &\p_{\theta_j}\Psi_3=0,\ j=1,2,3\quad\hbox{yield}\quad \cirt z\in K_{123},\\
  \nonumber
& &\p_{\theta_4}\Psi_3=0\quad\hbox{yields}\quad \cirt  y\in K_4,\\
  \nonumber
 & &\p_z\Psi_3=0\quad\hbox{yields}\quad\cirt  \xi=\omega_{\cirt  \theta},\\
   \nonumber
& &\p_{\xi}\Psi_3=0\quad\hbox{yields}\quad\cirt  y=\cirt z,\\
  \nonumber
 & &\p_{y}\Psi_3=0\quad\hbox{yields}\quad \cirt \xi=
 -d\varphi(\cirt  y)-\cirt  w.
  \eeq
  The critical points we need to consider for the asymptotics satisfy also
\beq\label{eq: imag.condition} \\ \nonumber
\im \varphi(\cirt y)=0,\quad\hbox{so that }\cirt y\in \gamma_{x_5,\xi_5},\ \im d\varphi(\cirt y)=0,\
\re d\varphi(\cirt y)\in L^{*,+}_{\cirt y}\hattuM _0.\hspace{-1cm}
\eeq


Next\hiddenfootnote{Alternatively, to analyze  the term $T_{\tau,1,2}^{(4),\beta}$
we could use the fact that that  $e^{\Psi_3}=|\xi|^{-2} (\nabla_z-\omega_\theta)^2 e^{\Psi_3}$
where $\omega_\theta=(\theta_1,\theta_1,\theta_1,0)$. This can be used in the 
integration by parts in the $z$ variable. 
Possibly, we could use this to analyze directly
$T_{\tau,1}^{(4),\beta}$ and omit the decomposition ${\bf Q}^*_1={\bf Q}^*_{1,1}+{\bf Q}^*_{1,2}$.  
} we analyze      the terms $T_{\tau,1,k}^{(4),\beta}$ starting with $k=2$.
 Observe that 
 the 3rd and 5th equations in (\ref{eq: critical points})  imply that at the critical points
 $\xi=( \p_{y_1}\varphi(y), \p_{y_2}\varphi(y),$ $ \p_{y_3}\varphi(y),0)$.
 Thus the critical points are bounded in the $\xi$ variable. 
 Let us now fix the parameter $R$ determining $\psi_R(\xi)$ in
 (\ref{eq: q_{1,1} and q_{1,2}}) so that 
 $\xi$-components of the critical
 points are in a ball $B(R)\subset \R^4$.
  Using the identity  $e^{\Psi_3}=|\xi|^{-2} (\nabla_z-\omega_\theta)^2 e^{\Psi_3}$
where $\omega_\theta=(\theta_1,\theta_1,\theta_1,0)$
we can include the operator $|\xi|^{-2}(\nabla_z-\omega_\theta)^2$
in
the integral (\ref{eq; T tau asympt}) with $k=2$
and integrate by
parts. Doing this  two times we can show
that this oscillatory  integral (\ref{eq; T tau asympt}) with $k=2$ becomes an integral of a Lebesgue-integrable
 function. 
 Then, by using method of stationary phase
 and  
%
  the fact that $\psi_R(\xi)$ vanishes at all critical points of $\Psi_3$
 where $\im \Psi_3$ vanishes, we see that   $T_{\tau,1,2}^{(4),\beta}=O(\tau^{-n})$
 for all $n>0$.  
%
%


Above, we have shown that  the term $T_{\tau,1,1}^{(4),\beta}$  has the same
asymptotics as $T_\tau^{(4),\beta}$. Next we analyze  this term.
Let $(\cirt z,\cirt \theta,\cirt y,\cirt \xi)$ be a critical point
of $ \Psi_3(z,\theta,y,\xi)$  such that $y$ satisfies (\ref{eq: imag.condition}).
Let us next use the same notations (\ref{new notations}) which
we used above.
Then (\ref{eq: critical points}) and (\ref{eq: imag.condition}) imply
\beq\label{eq: psi1 consequences}
\cirt z=\cirt y\in \gamma_{x_5,\xi_5}\cap \bigcap_{j=1}^4 K_j,\ \cirt \xi=\omega_{\cirt \theta}=-\cirt r-\cirt w.
\eeq
Note that in this case all the four geodesics $\gamma_{x_j,\xi_j}$ intersect
at the point $q$ and by our assumptions, $\cirt r=d\varphi(\cirt y)$ is such a co-vector
that in the $Y$-coordinates $\cirt r=(\cirt r_j)_{j=1}^4$ with $\cirt r_j\not =0$ for all $j=1,2,3,4$.
In particular, this shows that
 the existence of the critical point  of $ \Psi_3(z,\theta,y,\xi)$
implies that there exists 
an intersection point of $\gamma_{x_5,\xi_5}$ and $ \bigcap_{j=1}^4 K_j$.
Equations (\ref{eq: psi1 consequences}) imply also that 
\ba
\cirt r=\sum_{j=1}^4\cirt r_jdy^j=-\omega_{\cirt \theta}-\cirt w=-\sum_{j=1}^3\cirt \theta_jdz^j-
\cirt \theta_4 dy^4.
\ea

%


To consider the case when $\cirt y=\cirt z$, let us assume for a while that that $W_0=W_1$
and that the $Y$-coordinates and $Z$-coordinates coincide, that is, $Y(x)=Z(x)$.
Then the covectors $dz^j=dZ^j$ and $dy^j=dY^j$ coincide for $j=1,2,3,4$.
Then we have
\beq\label{eq: r is -theta}
\cirt r_j=-\cirt \theta_j,\quad\hbox{i.e.,\ } \cirt \theta:=\cirt \theta_jdz^j=-\cirt r=\cirt r_jdy^j\in T^*_{\cirt y}\hattuM _0.
\eeq


Let us apply
the method of stationary phase to $T^{(4),\beta}_{\tau,1,1}$ as $\tau\to \infty$. 
Note that
as $c_1(z,\theta_1,\theta_2)$ is a product type symbol, we need to use the fact $\theta_1\not=0$ and $\theta_2\not=0$
for the critical points as we have by (\ref{eq: r is -theta}) 
and the fact that $\cirt r=d\varphi(y)|_{\cirt y}\not \in N^*K_{234}\cup N^*K_{134}$
as $x_5\not \in  {\mathcal Y}$ and assuming that $s_0$ is small enough.

{
In local $Y$ and $Z$ coordinates where $\cirt z=\cirt y=(0,0,0,0)$ we
can use the method of  stationary phase,
 similarly to the proofs of \cite[Thm.\ 1.11 and 3.4]{GrS},
 to compute 
the asymptotics of  (\ref{eq; T tau asympt}) with $k=1$.
Let us explain this  computation in detail.
To this end, let us start with some preliminary considerations.

Let $\phi_1(z,y,\theta,\xi)$ be a smooth bounded
 function in $C^\infty(W_0\times W_1\times \R^4\times \R^4)$ that is 
homogeneous of degree zero in the $(\theta,\xi)$ variables in the set  $\{|(\theta,\xi)|>R_0\}$ with some $R_0>0$. 
Assume that $\phi_1(z,y,\theta,\xi)$ is equal to one  in 
a conic neighborhood, with respect to $(\theta,\xi)$, of the points where some of the 
$\theta_j$ or $\xi_k$ variable is zero. Note that in this set 
the positively homogeneous  functions $a_j(z,\theta_j/|\theta_j|)$ may be non-smooth.
Let $\Sigma\subset W_0\times W_1\times \R^4\times \R^4$ 
be a conic 
 neighborhood 
of the critical points of the phase function $\Psi_3(y,z,\theta,\xi)$.
Also, assume that the function  $\phi_1(y,z,\theta,\xi)$ vanishes in 
the intersection of $\Sigma$ 
and the set  $\{|(\theta,\xi)|>R_0\}$.
Let  
\ba
\Psi_{(\tau)}(z,y,\theta,\xi)= 
\theta_1z^1+\theta_2z^2+\theta_3z^3+(y-z)\cdotp \xi+\theta_4y^4+\tau \varphi(y)\ea 
 be the phase function appearing in (\ref{eq; T tau asympt pre}).
 Note that $\Sigma$ contains  the critical points of $\Psi_{(\tau)}$ for all $\tau>0$.
 Let
$$L_\tau=\frac {\phi_1(z,y,\tau^{-1}\theta,\tau^{-1} \xi)}{ |d_{z,y,\theta,\xi}
\Psi_{(\tau)}(z,y,\theta,\xi)|^{2}}\, (d_{z,y,\theta,\xi}
\overline{\Psi_{(\tau)}(z,y,\theta,\xi)})\cdotp d_{z,y,\theta,\xi},$$
so that 
\beq\label{kaav A}
L_\tau \exp(\Psi_{(\tau)})= \phi_1(z,y,\tau^{-1} \theta,\tau^{-1}\xi)\,\exp(\Psi_{(\tau)}).
\eeq
Note that if $\im \varphi(y)=0$, then $y\in \gamma_{x_5,\xi_5}$
and hence $d\varphi(y)$ does not vanish and that when $\tau$ is large enough,  the function $\Psi_{(\tau)}$ has 
no critical points in 
the support of $\phi_1$. Using
these we see that $ \gamma_{x_5,\xi_5}\cap W_1$
has a neighborhood $V_1\subset W_1$ where $|d\varphi(y)|>C_0>0$ and
there are $C_1,C_2,C_3>0$ so that if
  $\tau>C_1$, $y\in V_1$, and $(z,y,\theta,\xi)\in \supp(\phi_1)$
then
\beq\label{str est1}
\,|d_{z,y,\theta,\xi}
\Psi_{(\tau)}(z,y,\theta,\xi)|^{-1}\leq \frac {C_2}{\tau-C_3}.
\eeq
After these preparatory steps, we are ready to  compute the asymptotics of
 $T^{(4),\beta}_{\tau,1,1}$. To this end, we first transform the integrals in  (\ref{eq; T tau asympt})
 to an integral of a Lebesgue integrable function by using integration by parts of $|\xi|^{-2}(\nabla_z-\omega_\theta)^2$
 as explained above. Then
 we decompose  $T^{(4),\beta}_{\tau,1,1}$ into three
 terms $T^{(4),\beta}_{\tau,1,1}=I_1+I_2+I_3$. To obtain the first term $I_1$ we include the factor
 $(1-\phi_1(z,y,\theta,\xi))$ in the  integral  
 (\ref{eq; T tau asympt}) with $k=1$. 
  The integral $I_1$ can then be computed
 using the method of stationary phase similarly to the proof of \cite[Thm.\ 3.4]{GrS}. 
 Let $\chi_1\in C^\infty(W_1)$ be a function that 
 is supported in $V_1$ and vanishes on $\gamma_{x_5,\xi_5}$. 
 The terms $I_2$ and $I_3$
 are obtained by 
 including the factor
 $\phi_1(z,y,\tau^{-1}\theta,\tau^{-1}\xi)\chi_1(y)$ and  $\phi_1(z,y,\tau^{-1}\theta,\tau^{-1}\xi)(1-\chi_1(y))$
  in the  integral 
 (\ref{eq; T tau asympt pre}) with $k=1$, respectively.
(Equivalently, 
 the terms $I_2$ and $I_3$
 are obtained by 
 including the factor
 $\phi_1(z,y,\theta,\xi)\chi_1(y)$ and  $\phi_1(z,y,\theta,\xi)(1-\chi_1(y))$
  in the  integral 
 (\ref{eq; T tau asympt}) with $k=1$, respectively.)
 Using integration by parts in the integral (\ref{eq; T tau asympt pre}) and inequalities  (\ref{kaav A}) and (\ref{str est1}), we see that that $I_2=O(\tau^{-N})$
 for all $N>0$. Moreover, the fact that $\im \varphi(y)>c_1>0$ in  
$W_1\cap V_1$ implies that $I_3=O(\tau^{-N})$
 for all $N>0$.}


 Combining the above we obtain
the asymptotics
 \beq\label{definition of xi parameter1}
& &\quad \quad T^{(4),\beta}_\tau
\sim \tau^{4+4-16/2-2+\rho-2}\sum_{k=0}^\infty c_k
\tau^{-k} =\tau^{-4+\rho}\sum_{k=0}^\infty c_k\tau^{-k}
, \\ \nonumber
& &\hspace{-8mm}c_0=h(\cirt z)c_1(0,-\cirt r_1,-\cirt r_2)\hat a_3(0,-\cirt r_3)
\hat q_1(0,0,-(\cirt r_1,\cirt r_2,\cirt r_3,0))\hat a_4(0,-\cirt r_4)\hat a_5(0,1),\hspace{-1cm}
\eeq
where $\rho=\sum_{j=1}^5 p_j$
and $\hat a_j$ is the principal
symbol of $U_j$ etc.   The factor $h(\cirt z)$ is non-vanishing and is
determined by the determinant of the Hessian of the phase function
$\varphi$ at $q$.
A direct computation shows that $\det(\hbox{Hess}_{z,y,\theta,\xi} 
\Psi_{3}(\cirt z,\cirt y,\cirt \theta,\cirt \xi))=1$. Above,
$(\cirt z,\cirt y,\cirt \theta,\cirt \xi)$ is the critical point satisfying
(\ref{eq: critical points}) and (\ref{eq: imag.condition}), 
where in the local coordinates $(\cirt z,\cirt y)=(0,0)$ and
$h(\cirt z)$ is constant times powers of values of the cut-off functions $\Phi_0$
and $\Phi_1$ at zero. Recall that we considered above the case
when $\mathcal B_j$ are multiplication operators with these cut-off functions.
 The term $c_k$ depends on the 
derivatives of the symbols $a_j$ and $q_1$ of order less
or equal to $2k$ at the critical point.
If  $ \Psi_3(z,\theta,y,\xi)$  has no critical points, that is,
 $q$ is not an intersection point of all five geodesics $\gamma_{x_j,\xi_j}$, $j=1,2,3,4,5$ we
obtain the asymptotics $T^{(4),\beta}_{\tau,1,1}=O(\tau^{-N})$ for all $N>0$.

 For  future reference we note that if we use the method of  stationary
 phase in the last integral of (\ref{eq; T tau asympt}) only in the integrals with respect
 to $z$ and $\xi$, yielding that at the critical point we have  $y=z$ and $\xi=\omega_\beta(\theta)=(\theta_1,\theta_2,\theta_3,0)$, we 
 see that 
 $T_{\tau,1,1}^{(4),\beta}$ can
 be written as
 \beq \label{eq; T tau asympt modified}
& &T_{\tau,1,1}^{(4),\beta}
=c\tau^{4}
\int_{\R^{8}}e^{i(\theta_1y^1+\theta_2y^2+\theta_3y^3+\theta_4y^4)+i\tau \varphi(y)}
c_1(y, \theta_1, \theta_2)\cdotp\\ \nonumber
& &
\cdotp a_3(y, \theta_3)
q_{1,1}(y,y,\omega_\beta(\theta))a_4(y,\theta_4) a_5(y,\tau )\,d\theta_1d\theta_2d\theta_3d\theta_4dy\\ \nonumber
& &=c\tau^{8}
\int_{\R^{8}}e^{i\tau (\theta_1y^1+\theta_2y^2+\theta_3y^3+\theta_4y^4+\varphi(y))}
c_1(y, \tau \theta_1, \tau \theta_2)\cdotp\\ \nonumber
& &
\cdotp a_3(y,\tau  \theta_3)
q_{1,1}(y,y,\tau  \omega_\beta(\theta))a_4(y,\tau \theta_4) a_5(y,\tau )\,d\theta_1d\theta_2d\theta_3d\theta_4dy.
\eeq

Next, we consider the terms $\tilde T^{(4),\beta}_\tau$ of the type
(\ref{tilde T-type source}). Such term
is an integral of the product of $u_\tau$ and two other factors 
${\bf Q}(U_2\,\cdotp U_1)$ and ${\bf Q}(U_4\,\cdotp U_3)$.
As  the last two factors can be written in the form (\ref{Q U1U2term}),
one can see using the method of stationary phase that  $\tilde T^{(4),\beta}_\tau$ 
has similar asymptotics to  $T^{(4),\beta}_\tau$ as $\tau\to \infty$,
with the leading order coefficient $\tilde c_0=\tilde h(\cirt z)\hat c_1(0,-\cirt r_1,-\cirt r_2)
\hat c_2(0,-\cirt r_3,-\cirt r_4)$ where $\hat c_2$ is given as in (\ref{product type symbols})
with symbols $\hat a_3$ and $\hat a_4$, and moreover,
$\tilde h(\cirt z)$ is a constant times powers of values of the cut-off functions $\Phi_0$
and $\Phi_1$ at zero.

This proves the claim in the special case where
$u_j$ are conormal distributions supported in the coordinate neighborhoods $W_{\kappa(j)}$,
$a_j$ are positively homogeneous scalar valued symbols, $S_j={\bf Q}$, and $\B_j$ are 
multiplication functions with smooth cut-off functions.

By using a suitable partition of unity and 
summing the results of the  above computations, 
similar results to the above follows when  
$a_j$ are general classical symbols that are $\mathcal B$-valued and 
the waves $u_j$ are supported on
$J^+_{\hat g}(\supp({\bf f}_j))$.
Also, $S_j$ can be replaced by operators of type (\ref{extra notation}) and
$\B_j$ can replaced by differential operator without other essential changes
expect that the highest order power of $\tau$ changes.
Then, in the asymptotics of terms $T^{(4),\beta}_\tau$  the function
$h(\cirt z)$ in (\ref{definition of xi parameter1}) is a section in  dual bundle $(\B_L)^4$. The coefficients
of $h(\cirt z)$ in local coordinates are polynomials of $\hat g^{jk}$, $\hat g_{jk}$, 
$\hat \phi_{\ell}$, and 
their derivatives at $\cirt z$. Similar representation is obtained for the asymptotics
 of the terms $\tilde T^{(4),\beta}_\tau$.

As we integrated by parts two times the operator
$(\nabla_z-\omega_\theta)^2$ and the total order of 
of $\B_j$ is less or equal to 6, we see that it is
enough to assume above that the symbols $a_j(z,\theta_j)$
are of  order $(-12)$ or less.
The leading order asymptotics come from the term
where the sum of orders of $\B_j$ is 6 and $p_j=n$ for $j=1,2,3,4$,
$p_5=0$ so that $m=-4-\rho+6=4n+2$. 
We also see that the terms containing permutation $\sigma=\sigma_\beta$ of the indexes
of the distorted plane  waves can be analyzed analogously. This proves (\ref{indicator}).

Note that above $\cirt r$  depends only on $\hat g$ and ${\bf b}$ and we write
$\cirt r=\cirt r({\bf b})$, omiting the $\hat g$ depdency.
Making the above computations explicitly, we obtain an explicit formula for
the leading order coefficient $s_{m}$ in (\ref{indicator})  in terms of
$\bsequence$ and ${\bf w}$.
This show that $s_{m}$
coincides with some real-analytic function $G(\bsequence,{\bf w})$.
%
%
This proves the claims (i) and (ii) in the case when the linear independency condition (LI) is valid.

Next, consider the case when the linear independency condition (LI) is not valid.
Again, by the definition of ${\bf t}_j$,  if
the intersection $\gamma_{x_5,\xi_5}(\R_-) \cap(\cap_{j=1}^4\gamma_{x_j,\xi_j}((0,{\bf t}_j)))
$ is non-empty, it can
contain only one point. In the case that such a
point exists, we denote it by $q$.

When (LI) is not valid, we have that the linear space $\hbox{span}(b_j;\ j=1,2,3,4)\subset T_q^*M_0$ has dimension  3 or less.
We use the facts that for $w\in I(\Lambda_1,\Lambda_2)$
we have WF$(w)\subset \Lambda_1\cup \Lambda_2$
and the fact, see \cite[Thm.\ 1.3.6]{Duistermaat}  
 \ba
& &\hbox{WF}(v\,\cdotp w)\subset \\
& &\hbox{WF}(v)\cup \hbox{WF}(w)
\cup\{(x,\xi+\eta);\ (x,\xi)\in \hbox{WF}(v),\ (x,\eta)\in \hbox{WF}(w)\}.
\ea
Let us next consider the terms corresponding to the permutation $\sigma=Id$.
The above facts imply that $\tilde {\mathcal G}^{(4),\beta}$ 
in (\ref{w1-3 solutions}) satisfies
\ba
\hbox{WF}(\tilde {\mathcal G}^{(4),\beta})
\cap T_q^*M_0\subset \mathcal Z_{s_0}:=\X_{s_0}\cup \bigcup _{1\leq j\leq 4} N^*K_j\cup \bigcup_{1\leq j<k\leq 4} N^*K_{jk}, 
\ea 
where $\X_{s_0}=\X((\vec x,\vec\xi);t_0,s_0)$.
Also, for
\ba
w_{123}=\B_3^{\beta}u_{3}\,\cdotp \mathcal C_1^{\beta} S^{\beta}_1(\B^{\beta}_2u_{2}\,\cdotp \B^{\beta}_1u_{1}),
\ea
appearing in (\ref{M4 terms}), we have WF$(w_{123})
\subset \mathcal Z_{s_0}$ and thus using H\"ormander's theorem \cite[Thm.\ 26.1.1]{H4}, we see
that  WF$(S_2^{\beta}(w_{123}))\subset \Lambda^{(3)}$,
where  $\Lambda^{(3)}$ is the flowout of $\mathcal Z_{s_0}$ in the canonical relation 
of ${\bf Q}$. Then
\ba
\pi(\Lambda^{(3)})\subset \Y_{s_0}\cup \bigcup _{1\leq j\leq 4} K_j\cup \bigcup_{1\leq j<k\leq 4} K_{jk}, 
\ea 
where $\Y_{s_0}=\Y((\vec x,\vec\xi);t_0,s_0)$ and $\pi:T^*M_0\to M_0$ is the projection to the
base point.

Observe that $E=\hbox{span}(b_j;\ j=1,2,3,4)\subset T_q^*M_0$
has dimension 3 or less, $\Lambda^{(3)}\cap T_q^*M_0\subset E$
and WF$(u_{4})\cap T_q^*M_0\subset E$. Thus,
${\mathcal G}^{(4),\beta}=\B_4^{\beta}u_{4}\,\cdotp C_2^{\beta}S_2^{\beta}(w_{123})$ 
satisfies WF$({\mathcal G}^{(4),\beta})
\cap T_q^*M_0\subset E$.  Now, $E\subset \mathcal Z_{s_0}.$ By our assumption,
$(q,b_5)\not \in \mathcal X((\vec x,\vec\xi);t_0)$, and thus, cf.\ (\ref{w1-3 solutions}),  we see that 
\ba
\bra u_\tau,{\bf Q}(\sum_{\beta\leq n_1} {\mathcal G}^{(4),\beta}+\tilde {\mathcal G}^{(4),\beta})\cet=O(\tau^{-N})
\ea 
for all $N>0$ when $s_0$ is small enough.  The terms where
the permutation $\sigma$ is not the identity can be analyzed similarly.
This proves the claims (i) and (ii) in the case when the linear independency condition (LI) is not valid. This proves (i) and (ii).

(iii)  Let us fix $(x_j,\xi_j)$, $j\leq 4$, $s_0$, and the
 waves $u_j\in \I(K_j)$, $j\leq 4$.

 First we observe that if the condition (LI) is not valid, we see similarly to the above that
$\M^{(4)}$ is $C^\infty$ smooth in $\V\setminus (\Y\cup \bigcup_{j=1}^4 K_j)$.
Thus to prove the claim of the proposition we can assume that (LI) is valid.

Let us decompose  $\F^{(4)}$, given by
 (\ref{w1-3 solutions}) and (\ref{M4 terms})-(\ref{tilde M4 terms}) 
 as  $\F^{(4)}=\F_{1}^{(4)}+\F_{2}^{(4)}$
 where $\F^{(4)}_{p}$ is defined similarly to
 $\F^{(4)}$ in  (\ref{w1-3 solutions}) and 
 (\ref{M4 terms})-(\ref{tilde M4 terms})
by modifying these formulas so that 
the operator $S_1^\beta$ is replaced by 
$S_{1,p}^\beta$, where $S_{1,p}^\beta={\bf Q}_p$, when
 $S_{1}^\beta={\bf Q}$, and $S_{1,p}^\beta=(2-p)I$, when $S_{1}^\beta=I$.
 Here, the operators ${\bf Q}_p$  are defined as above using the parameters
  $\e_2$ and $\e_1$ defined below.
 
%

 Using formulas (\ref{w1-3 solutions}), (\ref{U1term}), (\ref{Q U1U2term}), and (\ref{QU4Ut term})
 we see that near $q$ in the $Y$ coordinates $\M_1^{(4)}={\bf Q} \F_1^{(4)}$ 
 can be calculated using that
 \beq\label{f4-formula}
 \F_1^{(4)}(y)=\int_{\R^4} e^{iy^j\theta_j}b(y,\theta)\,d\theta,
 \eeq
 where $K_j$ in local coordinates is given by $\{y^j=0\}$ and
  $b(y,\theta)$  is a finite sum of terms  that are products of some of the following
 terms: at most one  product type symbol $c_l(y,\theta_j,\theta_k)\in S(W_0;\R\times (\R\setminus \{0\}))$ (they appear in the terms  (\ref{T-type source})-(\ref{tilde T-type source})  where the $S_j^\beta$ operators are ${\bf Q}$ 
 and do not appear if these operators are the identity),
 and one ore more term which is either
 the symbols $a_j(y,\theta_j)\in S^n(W_0;\R)$,
or the functions $q_1(y,y,\omega_\beta(\theta))$,
  cf.\  (\ref{eq; T tau asympt modified}),
 where $\omega_\beta(\theta)$ is equal to some of the vectors
$(\theta_1,\theta_2,\theta_3,0)$,
 $(\theta_1,\theta_2,0,\theta_4)$, $(\theta_1,0,\theta_3,\theta_4)$, or
 $(0,\theta_1,\theta_3,\theta_4)$, depending on the permutation $\sigma$.

Let us consider next the source $F_\tau$ is determined by the functions $(p,h)$ in (\ref{Ftau source}). 
 Then using the  method of  stationary phase gives  the asymptotics, c.f. (\ref{eq; T tau asympt modified}),
 \ba
 \bra u_\tau,\F_1^{(4)}\cet \hspace{-1mm}\sim \hspace{-1mm}\tau^8\hspace{-1mm}\int_{\R^8} e^{i\tau (\varphi(y)+y^j\theta_j)}(a_5(y,\tau),b(y,\tau \theta))_{\hat G}d\theta dy
\sim \sum_{k=m}^\infty  
 s_k(p,h) \tau^{-k}
 \ea
where $\hat G$ is a Riemannian metric
on the fiber of $\B^L$ at $y$, that is isomorphic to $\R^{10+L}$, and
the critical point of the phase function is $y=0$ and $\theta=-d\varphi(0)$.
  As we saw above, we have that
 when $\e_2>0$ is small enough then for $p=2$ we have
$\bra u_\tau,\F_p^{(4)}\cet=O(\tau^{-N})$ for all $N>0.$
%

Let us choose sufficiently
small $\e_3>0$ and choose a function $\chi(\theta)\in C^\infty(\R^4)$ that 
vanishes in a $\e_3$-neighborhood (in the $\hat g^+$ metric) of
$\mathcal A_q$, 
\beq\label{set Yq}
\mathcal A_q&:=&
N_q^*K_{123}
\cup N_q^*K_{134}\cup N_q^*K_{124}\cup N_q^*K_{234}
\eeq
and is equal to 1 outside the $(2\e_3)$-neighborhood of this set.

Let $\phi\in C^\infty_0(W_1)$ be a function that is one near $q$.
Also, let  \ba\hbox{$b_0(y,\theta)=\phi(y)\chi(\theta)b(y,\theta)$}\ea be a classical symbol,
$p=\sum_{j=1}^4 p_j$,
and let  $
\F^{(4),0}(y)
\in \I^{p-4}(q) 
$ be the
 conormal distribution 
 that is given by the formula (\ref{f4-formula}) with
$b(y,\theta)$ being replaced by 
 $b_0(y,\theta)$.

When $\e_3$ is small enough (depending on the point $x_5$), we have that 
$F_\tau$ is determined by functions $(p,h)$ 
and the corresponding 
gaussian beams $u_\tau$
propagating on the geodesic $\gamma_{x_5,\xi_5}(\R)$ such
that the geodesic
passes through 
$x_5\in V$, we have 
\ba
\bra u_\tau,\F^{(4),0}\cet\sim \sum_{k=m}^\infty s_k(p,h) \tau^{-k},
\ea
that is, we have $\bra u_\tau,\F^{(4),0}\cet-\bra u_\tau,\F^{(4)}\cet=O(\tau^{-N})$
for all $N$.
When $\gamma_{x_5,\xi_5}(\R)$  does not pass through $q$,
we have that  $\bra u_\tau,\F^{(4)}\cet$ and $\bra u_\tau,\F^{(4),0}\cet$ 
are both of order $O(\tau^{-N})$ for all $N>0$.

Let $V\subset \V((\vec x,\vec \xi),t_0)\setminus \bigcup_{j=1}^4\gamma_{x_j,\xi_j}([0,\infty))$, see   (\ref{eq: summary of assumptions 2}), be an open set.
By varying  the source $F_\tau$,  defined in  (\ref{Ftau source}),  
 we see, by multiplying the solution  with a smooth
 cut of function and  using  Corollary 1.4 in \cite{Delort} in local coordinates,
 or \cite{MelinS},
%
 we have
 that the function
$\M^4-{\bf Q}\F^{(4),0}$ has no wave front set in $T^*(V)$ and it is
thus $C^\infty$-smooth function in $V$.

Since by \cite{GU1}, ${\bf Q}:\I^{p-4}(\{q\})\to \I^{p-4-3/2,-1/2}(N^*(\{q\}),\Lambda_q^+)$, the above implies that
\beq\label{eq: conormal}
\M^4|_{V\setminus {\mathcal Y}}\in \I^{p-4-3/2}(V\setminus  {\mathcal Y};\Lambda_q^+),
\eeq
where $  {\mathcal Y}=  {\mathcal Y}((\vec x,\vec \xi),t_0,s_0)$. When $x_5$ if fixed, choosing $s_0$ to
be small enough, we obtain the claim (iii).
\hfill \Box \medskip
}
 
Next we will show that $\mathcal G({\bf b},{\bf w})$ is not  vanishing
identically.

%

\generalizations
{{\bf Remark 3.6.A.} 
Similar considerations to the 
proof of  Prop.\ \ref{lem:analytic limits C}, can be used to show that $\M^3|_{(U\setminus \cup_j K_j)\cap I^-(x_5)}$
is a Lagrangian distribution. As an example of this, let us analyze the term
\ba
m={\bf Q}(\B_3 u_3\,\cdotp v_{12}),\quad \hbox{where }v_{12}={\bf Q}(\B_2 u_1\,\cdotp \B_2 u_3),
\ea
where $u_j\in \I^{p_j}(K_j)$ are solutions
of the linear wave equation.
Let $S=(\bigcup_j N^*K_j)\cup(\bigcup _{j,k}N^*K_{jk})$.
Next, we do computations in local coordinates $X:V\to \R^n$,  $x^j=X^j(x)$, such that $V\subset  I^-(x_6)$ 
and $K_j\cap V =\{x\in V;\ x^j=0\}$. 
As  $v_{12}=v_{12}^1+v_{12}^2\in \I(K_1,K_{12})+\I(K_2,K_{12})$, we can write the function
$f=\B_3 u_3\,\cdotp v_{12}$  in
$V$ as an oscillatory integral
\beq\label{def; f}
f(x)=\int_{\R^3}e^{i(\theta_1x^1+\theta_2x^2+\theta_3x^3)}b(x,\theta_1,\theta_2,\theta_3)d\theta_1 d\theta_2d\theta_3.
\eeq
Here, $b(x,\theta_1,\theta_2,\theta_3)\in S_{cl}^{p_1+p_2+p_3-2}(V;\R^3)$ away from
the set where some  some of the variables $\theta_j$ vanish and $\xi=\sum_{j=1}^3 \theta_jdx^j$
is light-like.
Let $\eta(x,\theta_1,\theta_2,\theta_3)=\phi(x,\sum_{j=1}^3 \theta_jdx^j),$
where $\phi$ is  a cut-off function
that
vanishes  in a conic neighborhood $S$ 
and is one outside a suitable conic neighborhood of $S$.
Then 
 $c(x,\theta_1,\theta_2,\theta_3)=\eta(x,\theta_1,\theta_2,\theta_3)b(x,\theta_1,\theta_2,\theta_3)$ is a classical symbol. Writing 
$f=f_1+f_2$ where $f_1$ and $f_2$ have the representation 
(\ref{def; f}) with the symbol $b$ is replaced by $c$ and $b-c$, correspondingly.
Then $f_1\in \I^{p_1+p_2+p_3-2}(K_{123})= \I^{p_1+p_2+p_3-5/2}(N^*K_{123})$ and 
we see that ${\bf Q}f_1\in
\I^{p_1+p_2+p_3-4}(Y;\tilde \Lambda),$  is a Lagrangian distribution on  $Y=I^-(x_5)\setminus ( \cup_j K_j)$,
where $\tilde \Lambda$ is the flow-out of $N^*K_{123}$,
and ${\bf Q}f_2$ is a $C^\infty$-smooth in  $Y$.
Hence, 
\beq\label{3-singularity order}
\M^3|_{Y}\in \I^{\tilde p-10}(Y;\tilde \Lambda),\quad \tilde p=p_{1}+p_{2}+p_{3}.
\eeq
}

\subsection{Non-vanishing interaction in the Minkowski space}
\subsubsection{WKB computations and the indicator functions in the  Minkowski space}
 \MATTITEXT{ 
To show that the function  $\mathcal  G(\bsequence,{\bf w})$, see (\ref{definition of G}),  is not  vanishing identically, we will next consider 
waves in Minkowski space.

In this section,
$x=(x^0,x^1,x^3,x^4)$ are the standard coordinates
in Minkowski space and
 $\hat g_{jk}=\diag(-1,1,1,1)$ denote the metric
in the standard coordinates of the Minkowski space $\R^4$. Below we call
the principal symbols of the linearized waves the {\it polarizations} to emphasize their physical meaning.
\MTEXT{We denote ${\bf w}=(w_{(j)})_{j=1}^5$.
Then for $j\leq 4$, the polarizations $w_j=(v_{(j)}^{met},v_{(j)}^{scal})$,
represented as a pair of the metric part of the polarization $v_{(j)}^{met}\in\hbox{sym}(\R^{4\times 4})\equiv \R^{10}$ and the scalar field part 
of the polarization $v_{(j)}^{scal}\in \R^L$, are such that for the metric 
 part  $v_{(j)}^{met}$ of the polarization
has to satisfy 4 linear conditions (\ref{harmonicity condition for symbol}) with the Minkowski metric (that follow
from the linearized harmonicity condition\noextension{ (\ref{harmonicity condition A})).} 
We study the special case when all polarizations of the scalar fields  $\phi_\ell$,
$\ell=1,2,\dots,L$, vanish,
that is, $v^{scal}_{(j)}=0$ for all $j$. 
To simplify the notations,  we denote below $v_{(j)}^{met}=v_{(j)}$
 and  ${\bf v}=(v_{(j)})_{j=1}^5$ so that ${\bf w}=({\bf v},0)$.}
In this case, in Minkowski space the function $\mathcal G({\bf v},0,\bsequence)$
can
be analyzed by assuming that there are no matter fields, which we do next.
Later we return to the case  of general  polarizations.
%


We assume that  the waves $u_j(x)$, $j=1,2,3,4$, solving the linear wave equation
in the Minkowski space,  are of the form
\beq\label{eq. plane waves}
 u_j(x)= v_{(j)}\, \bigg(b^{(j)}_px^p \bigg)^{a}_+,\quad t^a_+=|t|^aH(t),
\eeq
where $b^{(j)}_pdx^p$, $p=1,2,3,4$ are  {four linearly independent} light-like co-vectors of $\R^4$,
$a>0$ and $v_{(j)}$ 
are constant $4\times 4$ matrices. 
{We also assume that $b^{(5)}$ is not in the linear span
of any three vectors $b^{(j)}$, $j=1,2,3,4$.} 
In the following,
we denote 
$b^{(j)}\,\cdotp x:=b^{(j)}_p x^p$
and $\bsequence=(b^{(j)})_{j=1}^5$. 
Let us next consider the wave produced by interaction of two plane wave solutions
in the  Minkowski space.\hiddenfootnote{We start with the observation
 that
\ba
& &(\p_0^2-\p_1^2-\p_2^2-\p_3^2)\bigg(  (x_0-x_1)^a_\pm\,\cdotp (x_0-x_2)^c_\pm   \bigg)\\
& &=2ac\,
(x_0-x_1)^{a-1}_\pm\,\cdotp (x_0-x_2)^{c-1}_\pm ,  \\
& &(\p_0^2-\p_1^2-\p_2^2-\p_3^2)\bigg(  (x_0-x_1)^a\,\cdotp (x_0-x_2)^c_\pm   \bigg)\\
& &=2ac\,
(x_0-x_1)^{a-1}_\pm\,\cdotp (x_0-x_2)^{c-1}.
\ea
More generally,  for $\tilde w=(0,\tilde w^1,\tilde w^2,\tilde w^3)\in \R^4$, $\tilde w\,\cdotp \tilde w=1$ and  $w=(0,w^1,w^2,w^3)\in \R^4$, $w\,\cdotp w=1$,
we have
\ba
& &(\p_0^2-\p_1^2-\p_2^2-\p_3^2)\bigg(  (x_0-x_1)^a_+\,\cdotp (x_0-\tilde w\,\cdotp x)^c_+   \bigg)=
\\&=&2ac\,
(x_0-x_1)^{a-1}_+\,\cdotp (x_0-\tilde w\,\cdotp x^{\prime})^{c-1}_+ + \\ 
& &-2ac\,
(x_0-x_1)^{a-1}_+\,\cdotp \tilde w_1(x_0-\tilde w\,\cdotp x^{\prime})^{c-1}_+
\ea
and
\ba
& &(\p_0^2-\p_1^2-\p_2^2-\p_3^2)\bigg(  (x_0-w\,\cdotp x^{\prime})^a_+\,\cdotp (x_0-\tilde w\,\cdotp x^{\prime})^c_+   \bigg)=
\\&=&2ac\,
(x_0-w\,\cdotp x^{\prime})^{a-1}_+\,\cdotp (x_0-\tilde w\,\cdotp x^{\prime})^{c-1}_+ + \\ 
& &+\sum_{j=1}^3 2ac\tilde w_jw_j
(x_0-w\,\cdotp x^{\prime})^{a-1}_+\,\cdotp \tilde w_1(x_0-\tilde w\,\cdotp x^{\prime})^{c-1}_+
\\&=&2ac(1-\tilde w\,\cdotp w)\,
(x_0-w\,\cdotp x^{\prime})^{a-1}_+\,\cdotp (x_0-\tilde w\,\cdotp x^{\prime})^{c-1}_+.
\ea}

Let $b^{(1)}$ and $b^{(2)}$ be light like co-vectors.
We use the notations
\ba
& &u^{a_1,a_2}(x;{b^{(1)},b^{(2)}})=(b^{(1)}\,\cdotp x)^{a_1}_+\,\cdotp (b^{(2)}\,\cdotp x)^{a_2}_+
\ea
for a product of two
plane waves. 
We define the formal parametrix ${\bf Q} _0$,
\beq\label{B3 eq}
{\bf Q}_0(u^{a_1,a_2}(x;{b^{(1)},b^{(2)}}))=\frac {u^{a_1+1,a_2+1}(x;{b^{(1)},b^{(2)}})}{2(a_1+1)(a_2+1)\,\hat g( b^{(1)},b^{(2)})} .
\eeq
Then $\square_{\hat g}({\bf Q}_0(u^{a_1,a_2}(x;{b^{(1)},b^{(2)}})))=u^{a_1,a_2}(x;{b^{(1)},b^{(2)}})$.
{
%
%
%
Also, let
\ba
& &u^{a}_{\tau} (x;b^{(4)},b^{(5)})= u_4(x)\, u_\tau(x),
\quad 
u_4(x)= (b^{(4)}\,\cdotp x)^a_+,\quad  u_\tau(x)=e^{i\tau\,b^{(5)}\cdotp x},
\ea
so that
$\square_{\hat g} ( u^{a}_{\tau} (x;b^{(4)},b^{(5)}))=2a\,\hat g(b^{(4)},b^{(5)})i\tau \, u^{a-1,0}_{\tau} (x;b^{(4)},b^{(5)}).
$
Let 
\beq\label{Q0 def}{\bf Q}_0( u^{a}_{\tau} (x;b^{(4)},b^{(5)}))=\frac 1{2i(a+1)\,\hat g(b^{(4)},b^{(5)})\tau}u^{a+1}_{\tau} (x;b^{(4)},b^{(5)}).\hspace{-1cm}
\eeq
}
%
%


 
We will prove that the indicator function
$\mathcal G({\bf v},0,\bsequence)$ in (\ref{definition of G}) does not vanish identically
  by showing that it coincides with  the {\it formal}
  indicator function $\mathcal G^{({\bf m})}({\bf v},\bsequence)$, {defined below,}
 which is a real-analytic function that does not vanish identically.
We define  the (Minkowski) indicator function (c.f.\ (\ref{test sing}) and (\ref{definition of G}))
 by \ba
\mathcal G^{({\bf m})}({\bf v},{\bf b})=
\lim_{\tau\to\infty} \tau^{m}(\sum_{\beta\leq n_1}
\sum_{\sigma\in \Sigma(4)} 
T^{({\bf m}),\beta}_{\tau,\sigma}+\tilde T^{({\bf m}),\beta}_{\tau,\sigma}),
\ea 
where the super-index $({\bf m})$ refers
to the word ``Minkowski''.
Above, 
 $\sigma$  runs over all permutations of the set $\{1,2,3,4\}$. The functions
$T^{({\bf m}),\beta}_{\tau,\sigma}$ and $\tilde T^{({\bf m}),\beta}_{\tau,\sigma}$ are counterparts
of the functions $T^{(4),\beta}_{\tau}$ and $\tilde T^{(4),\beta}_{\tau}$,
see (\ref{T-type source})-(\ref{tilde T-type source}),
obtained by replacing the distorted plane waves and
the gaussian beam by the plane waves. We also replace 
the parametrices ${\bf Q}$ and ${\bf Q}^*$ by a formal parametrix ${\bf Q}_0$.
Also, we include 
 a smooth cut off function   $h\in C^\infty_0(M)$ which is one near 
 the intersection point $q$ of  the $K_j$. Thus,
\beq
\label{Term type 1}
& &\hspace{- .5cm} T^{({\bf m}),\beta}_{\tau,\sigma}=\bra S^0_2(u_\tau \,\cdotp \B_{4}u_{\sigma(4)}), h\,\cdotp \B_{3}u_{\sigma(3)}\,\cdotp 
S_1^0(\B_{2}u_{\sigma(2)}\,\cdotp \B_{1}u_{\sigma(1)})\cet_{L^2(\R^4)},\hspace{-1cm} \\
\label{Term type 2}
& &\hspace{- .5cm} \tilde T^{({\bf m}),\beta}_{\tau,\sigma}=\bra u_\tau,h\,\cdotp S_2^0(\B_4u_{\sigma(4)} \,\cdotp \B_{3}u_{\sigma(3)})\,\cdotp 
S^0_1(\B_2u_{\sigma(2)}\,\cdotp \B_1u_{\sigma(1)})\cet_{L^2(\R^4)},\hspace{-1cm} 
\eeq
where $u_j$ are given by (\ref{eq. plane waves}) with $a=-n-1$, $j=1,2,3,4$.
Here, the differential operators $\B_j=\B^\beta_{j}$ are in Minkowski space
 constant coefficient operators and finally,
 $S_j^0=S^0_{j,\beta}\in \{{\bf Q}_0,I\}$.

}

\extension{

Let us now consider the orders of the differential  operators appearing above.
The orders $k_j=ord(\B^\beta_j)$ 
of the differential operators $\B^\beta_j$, 
defined in (\ref{extra notation}), 
 depend on 
$\vec S_\beta^0=(S_{1,\beta}^0,S_{2,\beta}^{0})$ 
as follows: 
When $\beta$ is such that $\vec S_\beta^0=(  {\bf Q}_0,  {\bf Q}_0)$,  for the terms $ T^{({\bf m}),\b}_{\tau,\sigma}$
 we have 
\beq\label{k: s for (Q_0,Q_0) and T}
& &k_1+k_2+k_3+k_4\leq 6,\quad  k_3+k_4\leq 4,\quad
k_4\leq 2
\eeq
and for the terms $\tilde T^{({\bf m}),\b}_{\tau,\sigma}$ we have 
\beq\label{k: s for (Q_0,Q_0) and tilde T} 
k_1+k_2+k_3+k_4 \leq 6,\quad
k_1+k_2\leq 4,\quad
k_3+k_4\leq 4.\hspace{-1.5cm}
\eeq

When  $\beta$ is such that $\vec S_\beta^0=(I,  Q_0)$
 we have  for the terms $T^{({\bf m}),\b}_{\tau,\sigma}$
\beq\label{k: s for (I,Q_0) and T} 
k_1+k_2+k_3+k_4 \leq 4,\quad k_4\leq 2,
\eeq
 and for terms $\tilde T^{({\bf m}),\b}_{\tau,\sigma}$ we have 
\beq\label{k: s for (I,Q_0) and tilde T}  k_1+k_2+k_3+k_4 \leq 4,\quad
k_1+k_2\leq 2.
\eeq

When  $\beta$ is such that $\vec S_\beta^0=(  {\bf Q}_0,I)$, 
both for the terms $T^{({\bf m}),\b}_{\tau,\sigma}$ 
and $\tilde  T^{({\bf m}),\b}_{\tau,\sigma}$
we have 
\beq\label{k: s for (Q_0,I) and tilde T} 
& &k_1+k_2+k_3+k_4\leq 4,\quad
k_3+k_4\leq 2.
\eeq
Finally, when  $\beta$ is such that $\vec S_\beta^0=(I,I)$,
 for the terms   $T^{({\bf m}),\b}_{\tau,\sigma}$ and $\tilde T^{({\bf m}),\b}_{\tau,\sigma}$ we have
$k_1+k_2+k_3+k_4 \leq 2$. 

Let us next summarize these formulas in different notations:
}

\noextension{Let us now consider the orders of the differential  operators appearing above.}
Recall that $k_j=ord(\B^\beta_j)$  are  the orders 
of the differential operators $\B^\beta_j$, 
defined in (\ref{extra notation}). 
For $j=1,2$, we  define $K_{\beta,j}=1$ when $S_{\beta,j}^0={\bf Q}_0$
and $K_{\beta,j}=0$ when $S_{\beta,j}^0=I$.
%
%
Then the allowed values of $\vec k=(k_1,k_2,k_3,k_4)$ 
 depend on $K_{\beta,1}$ and $K_{\beta,2}$
as follows: 
We require that  
 \beq\label{k: collected}
& &k_1+k_2+k_3+k_4\leq 2K_{\beta,1}+2K_{\beta,2}+2,\quad  
k_3+k_4\leq 2K_{\beta,2}+2,\\ \nonumber
 & &\hbox{ $k_4\leq 2$,
 for all terms  $T^{({\bf m}),\b}_{\tau,\sigma}$, and}\\
 \nonumber
 & &\hbox{ $ k_1+k_2\leq 2K_{\beta,1}+2$,
 for all  terms $\tilde T^{({\bf m}),\b}_{\tau,\sigma}$}
\eeq
cf. (\ref{M4 terms}) and (\ref{tilde M4 terms}).
\extension{
\medskip

\noindent
{\bf Lemma 3.A.1.}  {\it When
 $b^{(j)}$, $j=1,2,3,4$ are   linearly independent light-like co-vectors and light-like
 co-vector $b^{(5)}$ is not in the linear span
of any three vectors $b^{(j)}$, $j=1,2,3,4$ we have
$\mathcal G({\bf w},{\bf b})=\mathcal G^{({\bf m})}({\bf v},{\bf b})$
when  $w_{(j)}=(v_{(j)},0)\in \R^{10}\times \R^L$.}
\medskip

\noindent{\bf Proof.}
Let us start by considering the relation of ${\bf Q}_0$ with the causal inverse ${\bf Q}$
in (\ref{Q0 def}).
Let
\ba
\nonumber
w_{\tau,0}&=&{\bf Q}_0(u^{a,0}_{\tau} (\,\cdotp;b^{(4)},b^{(5)}))\\
&=&\int_{\R}e^{i\theta_4z^4+ i\tau P\cdotp z}
a_4(z,\theta_4)d\theta_4,\\
w_\tau&=&{\bf Q}^*(J ),\\
J&=&
 u_4 \,\cdotp (\chi\,\cdotp u^\tau).\ea
Here 
\ba
J(z)&=&\chi(X^0(z))\, \square w_{\tau,0}\\
&=&\int_{\R}e^{i\theta_4z^4+ i\tau P\cdotp z}
(\tau b_1(z,\theta_4)+
b_2(z,\theta_4))a_5(z,\tau)d\theta_4,
\ea
where $a_5(z,\tau)=1$.
Then 
\ba
& &\square(w_{\tau}-\chi w_{\tau,0})=J_1,\quad \hbox{for } x\in \R^4,\\
& &J_1=[\square,\chi]w_{\tau,0}\\
&&\quad\ =\int_{\R}e^{i\theta_4z^4+ i\tau P\cdotp z}
(\tau b_3(z,\theta_4)+
b_4(z,\theta_4))d\theta_4,
\ea
where $b_3(z,\theta_4)$ and $b_4(z,\theta_4)$ are
supported in the domain $ T_0< X^0(z)<T_0+1$ and
$w_{\tau}-\chi w_{\tau,0}$ is supported in the domain $X^0(z)<T_0+1$ .
Thus 
\ba
w_{\tau}=\chi w_{\tau,0}+{\bf Q}^*J_1.
\ea
Here, 
we can write 
\beq\label{modified J}
& &J_1(z)= u_4^{(1)}(z) u^{\tau,(1)}(z)+
u_4^{(2)}(z) u^{\tau,(2)}(z)
,\quad\hbox{where}\\ \nonumber
& &u_4^{(1)}(z)= \int_{\R}e^{i\theta_4z^4}
b_3(z,\theta_4)d\theta_4,\\ \nonumber
& & u^{\tau,(1)}(z)=\tau u^\tau(z),\\ \nonumber
& &u_4^{(2)}(z)= \int_{\R}e^{i\theta_4z^4}
b_4(z,\theta_4)d\theta_4,\\ \nonumber
& & u^{\tau,(2)}(z)=u^\tau(z).
\eeq

Let us now substitute this in to the above  microlocal computations
done in the proof of Prop.\ \ref{lem:analytic limits A}. 
%
%

Recall that $b^{(j)}$, $j=1,2,3,4$, are four linearly independent co-vectors.
This means that a condition analogous to (LI) in the proof of Prop.\ \ref{lem:analytic limits A}
is satisfied, and that $b^{(5)}$ is not in the space spanned by any of three  of the
co-vectors $b^{(j)}$, $j=1,2,3,4$.
Also, observe that the hyperplanes $K_j=\{x\in \R^4;\ b^{(j)}\,\cdotp x=0\}$
intersect at origin of $\R^4$. Thus, we see that
the arguments in  the proof of Prop.\ \ref{lem:analytic limits A}
are valid mutatis mutandis if the phase function of the gaussian beam 
$\varphi(x)$ is replaced
by the phase function of the plane wave, $b^{(5)}\,\cdotp x$, 
and  the geodesic $\gamma_{x_5,\xi_5}$, on
which the gaussian beam propagates, is replaced by the whole space $\R^4$.
In particular, as $b^{(5)}$ is not in the space spanned by any of three  of those
co-vectors, the case (A1) in the proof of Prop.\ \ref{lem:analytic limits A} cannot occur.
In particular, we have that
 the  leading order asymptotics of the terms
$T^{\b}_{\tau,\sigma}$ and $\tilde T^{\b}_{\tau,\sigma}$ do not change
as these asymptotics are  obtained using
the method of stationary phase for the integral 
(\ref{eq; T tau asympt modified}) and the other analogous integrals
at the critical point $z=0$.
In other words, we can replace the gaussian beam by a plane wave in our considerations
similar to those in the proof of Prop.\ \ref{lem:analytic limits A}.

Using (\ref{modified J})  and the fact that $b_3(z,\theta_4)$ and $b_4(z,\theta_4)$ 
vanish near $z=0$, we see that if $u_4$ and $u^{\tau}$ are replaced by
 $u_4^{(j)}$ and $u^{\tau,(j)}$, respectively, where $j\in \{1,2\}$ and
we can do similar computations based on the method of stationary phase  as are done in the 
proof of Proposition \ref{lem:analytic limits A}. Then both terms
$T^{\b}_{\tau,\sigma}$ and $\tilde T^{\b}_{\tau,\sigma}$
have asymptotics  $O(\tau^{-N})$ for all $N>0$ as $\tau\to \infty$.
In other words, in the proof of Prop.\ \ref{lem:analytic limits A} the term
$w_{\tau}={\bf Q}^*( u_4u^\tau )$ can be replaced by $\chi w_{\tau,0}$
without changing the leading order asymptotics. 
This shows
that $\mathcal G({\bf w},{\bf b})=\mathcal G^{({\bf m})}({\bf v},{\bf b})$,
where $w_{(j)}=(v_{(j)},0)\in \R^{10}\times \R^L$.  
\hfill \Box \medskip

}

Summarizing the above: 
Let
 $b^{(j)}$, $j=1,2,3,4$ be  linearly independent light-like co-vectors and 
 $b^{(5)}$  be  a light-like  co-vector that is not in the linear span
of any three vectors $b^{(j)}$, $j=1,2,3,4$. Then, \extension{using
the above lemma and}
by analyzing the   microlocal computations
done in the proof of Prop.\ \ref{lem:analytic limits A}, we see
that
$\mathcal G({\bf w},{\bf b})=\mathcal G^{({\bf m})}({\bf v},{\bf b})$
when  $w_{(j)}=(v_{(j)},0)\in \R^{10}\times \R^L$,  ${\bf w}=(w_{(j)})_{j=1}^5,$ and ${\bf v}=(v_{(j)})_{j=1}^5.$


\begin{proposition}\label{singularities in Minkowski space}
Let 
${\mathbb X}$  be the set of $({\bf b},v_{(2)},v_{(3)},v_{(4)})$,
where ${\bf b}$ is a  5-tuple of
light-like covectors ${\bf b}=(b^{(1)},
b^{(2)},b^{(3)},b^{(4)},b^{(5)})$ and 
$v_{(j)}\in \R^{10}$, $j=2,3,4$ are the polarizations that satisfy the equation (\ref{harmonicity condition for symbol})
with respect to $b^{(j)}$,
i.e.,  the harmonicity condition for the principal symbols.
\MTEXT{For $\hat b^{(5)}\in \R^4$, let 
 ${\mathbb X}(\hat b^{(5)})$  be the set elements in ${\mathbb X}$ where $b^{(5)}=\hat b^{(5)}$.}
 Then  for any light-like $\hat b^{(5)}$ there is 
a  generic (i.e. open and dense)
subset ${\mathbb X}^\prime(\hat b^{(5)})$ of ${\mathbb X}(\hat b^{(5)})$ such that for all
$({\bf b},v_{(2)},v_{(3)},v_{(4)})\in {\mathbb X}^\prime({\hat b^{(5)}})$ 
 there exist linearly independent vectors  $v_{(5)}^q$, $q=1,2,3,4,5,6,$
with the following property:
\smallskip

 If
 $v_{(5)}\in \hbox{span}(\{v_{(5)}^q;\ q=1,2,3,4,5,6\})$  is non-zero, then 
 there exists 
 a vector $v_{(1)}$ for which 
the pair 
$(b^{(1)},v_{(1)})$
satisfies the equation (\ref{harmonicity condition for symbol})
 and
 $\mathcal G^{({\bf m})}({\bf v},{\bf b})\not =0$
 with ${\bf v}=(v_{(1)},v_{(2)},v_{(3)},v_{(4)},v_{(5)})$.
%
%
%
\end{proposition}

\extension{

\noindent{\bf Proof.} 
In the proof below, let $a\in \Z_+$ be  large enough.
To show that the coefficient $\mathcal G^{({\bf m})}({\bf v},{\bf b})$  of the 
leading order term in the asymptotics  
is non-zero, we consider a special case when the direction vectors of the
intersecting plane waves  in the  Minkowski space are 
the linearly independent light-like vectors of the form
\ba
b^{(5)}=(1,1,0,0),\quad b^{(j)}=(1,1-\frac 12\rhoepsilon_j^2,\rhoepsilon_j+O(\rhoepsilon_j^3),\rhoepsilon_j^3),\quad j=1,2,3,4,
\ea
where $\rhoepsilon_j>0$ are small parameters for which
\beq\label{b distances}
& &\|b^{(5)}-b^{(j)}\|_{(\R^4,\hat g^+)}= \rhoepsilon_{j}(1+o(\rhoepsilon_{j})),\quad j=1,2,3,4.
\eeq
With an appropriate choice of $O(\rhoepsilon_k^3)$ above, the vectors $b^{(k)}$, $k\leq 5$
are  light-like and
\ba
\,\hat g(b^{(5)},b^{(j)})&=&-1+(1-\frac 12 \rhoepsilon_j^2)=-\frac 12 \rhoepsilon_j^2,\\
\,\hat g( b^{(k)},b^{(j)})&=&
-\frac 12 \rhoepsilon_k^2-\frac 12 \rhoepsilon_j^2+O(\rhoepsilon_k\rhoepsilon_j).
\ea
Below, we denote $\omega_{kj}=\,\hat g( b^{(k)},b^{(j)})$. 
We consider the case when the orders of $\rhoepsilon_j$
are   
\beq\label{eq: ordering of epsilons}
\rhoepsilon_4=\rhoepsilon_2^{100},\ \rhoepsilon_2=\rhoepsilon_3^{100},\hbox{ and }\rhoepsilon_3=\rhoepsilon_1^{100},  
\eeq
{
so that $\rho_4<\rho_2<\rho_3<\rho_1$. \MTEXT{When $\rho_1$ is small enough,
$b_j,$ $j\leq 4$ are linearly independent.}
Note that when $\rho_1$ is small enough, $b^{(5)}$ is not a linear combination
of any three vectors $b^{(j)}$, $j=1,2,3,4$.}

The coefficient $\mathcal G^{({\bf m})}$ of the leading order asymptotics
is computed by analyzing the leading order terms of 
all 4th order interaction terms, similar to those 
given in (\ref{Term type 1}) and (\ref{Term type 2}).
We will start by analyzing the most important terms $ T^{({\bf m}),\beta}_\tau$
of the type (\ref{Term type 1})
%
%
when $\beta$ is such that $\vec S_\b=({\bf Q}_0,{\bf Q}_0)$.
%
When $k_j=k_j^\b$ is the order of $\B_j$, 
and we denote $\vec k_\b=(k_1^\b,k_1^\b,k_3^\b,k_4^\b)$,
we see that 
\beq\label{first asymptotical computation}
 & &T^{({\bf m}),\beta}_\tau
=\bra {\bf Q}_0(\B_4 u_4\,\cdotp u^\tau), h\,\cdotp \B_3u_3\,\cdotp {\bf Q}_0(\B_2u_2\,\cdotp \B_1 u_1)\cet\\
\nonumber
& &\hspace{-1cm}=C\frac {\P_\beta} {\omega_{45}
\tau} \frac 1{\omega_{12}}
\bra u^{a-k_4+1,0}_{\tau} (\,\cdotp ;b^{(4)},b^{(5)}), h\,\cdotp u_3\,\cdotp u^{a-k_2+1,a-k_1+1} (\,\cdotp ;b^{(2)},b^{(1)})\cet
\\
\nonumber&&\hspace{-1cm} =C\frac {\P_\beta}{\omega_{45}\tau} \frac 1{\omega_{12}} 
\int_{\R^4}
(b^{(4)}\,\cdotp x)^{a-k_4+1}_+e^{i\tau(b^{(5)}\,\cdotp x)} 
h(x)(b^{(3)}\,\cdotp x)^{a-k_3}_+\cdotp\\
\nonumber& &\quad\quad\quad\cdotp
(b^{(2)}\,\cdotp x)^{a-k_2+1}_+
(b^{(1)}\,\cdotp x)^{a-k_1+1}_+\,dx,
\eeq
where $\P=\P_\beta$ is a polarization factor involving the coefficients of $\B_j$, 
the directions $b^{(j)}$, and the polarization $v_{(j)}$. Moreover, $C=C_a$ is a generic constant
 depending on $a$ and $\beta$ but not on  $b^{(j)}$ or $v_{(j)}$.

We will analyze the polarization factors later, but as a sidetrack, 
let us already explain now the nature of the polarization
term 
when $\beta=\beta_1$,
see (\ref{beta0 index}).  Observe that this term appear only when  we analyze
 the term
 $\bra F_\tau,{\bf Q}(A[u_4,{\bf Q}(A[u_3,{\bf Q}(A[u_2,u_1])])])\cet$ where
 all operators $A[v,w]$ are of the type $A_2[v,w]=\hat g^{np}\hat g^{mq}v_{nm}\p_p\p_q w_{jk}$,
 cf.\ (\ref{eq: tilde M1}) and (\ref{eq: tilde M2}).
 Due to this, we have the polarization factor 
 \beq\label{pre-eq: def of D}
\P_{\beta_1} =(v_{(4)}^{rs}b^{(1)}_rb^{(1)}_s)(v_{(3)}^{pq}b^{(1)}_pb^{(1)}_q)(v_{(2)}^{nm}b^{(1)}_nb^{(1)}_m)\D,
\eeq
 where $v_{(\ell)}^{nm}=\hat g^{nj}\hat g^{mk}v^{(\ell)}_{jk}$ and
\beq\label{eq: def of D}
 \D =\hat g_{nj}\hat g_{mk}v_{(5)}^{nm}v_{(1)}^{jk}.
 \eeq
We will postpone the analysis of  the polarization factors
$\P_{\beta}$ in  $ T^{({\bf m}),\beta}_\tau$  with $\beta\not=\beta_1$ later.

Let us now return back to the computation (\ref{first asymptotical computation})\hiddenfootnote{
In the computations below, we analyze the following terms:

  1. The $T^{{\bf Q}, {\bf Q}}$ term in formula (\ref{first asymptotical computation}) comes from the 2nd term in (\ref{w4 solutions})
 
  2. $\tilde T^{{\bf Q}, {\bf Q}}$ term on p. 55 that comes from the 1st term in (\ref{w4 solutions}) inside
  the large brackets,
  
  3. $T^{I, {\bf Q}}$ term on p. 56  that comes from the 3rd term in (\ref{w4 solutions}) inside
  the large brackets,

  4. $\tilde T^{I, {\bf Q}}$ term on p. 56  that comes from the 4th term in (\ref{w4 solutions}) inside
  the large brackets,

  5. $T^{{\bf Q}, I}$ term on p. 56  that comes also from the 4th term in (\ref{w4 solutions}) inside
  the large brackets
  if we use
another permutation of $\{1, 2, 3, 4\}$,

 6. $T^{{\bf Q}, I}$ term on p. 56  that comes from the 5th term in (\ref{w4 solutions}) inside
  the large brackets.
}.
We next use in $\R^4$ the coordinates $y=(y^1,y^2,y^3,y^4)^t$
where $y^j=b^{(j)}_kx^k$, i.e., and let $A\in \R^{4\times 4}$
be the matrix for which $y=A^{-1}x$.
Let ${\bf p}=(A^{-1})^tb^{(5)}$.
 {In the $y$-coordinates,  $b^{(j)}=dy^j$ for $j\leq 4$ and
 $b^{(5)}=\sum_{j=1}^4{\bf p}_jdy^j$ and
 \ba
 {\bf p}_j
 =\hat g(b^{(5)},dy^j)=\hat g(b^{(5)},b^{(j)})= \omega_{j5}=-\frac 12 \rhoepsilon_j^2.
 \ea
Then   $b^{(5)}\,\cdotp x={\bf p}\,\cdotp y$.
 We use the  notation
$
{\bf p}_j=\omega_{j5}=-\frac 12 \rhoepsilon_j^2,
$ 
that is, we denote the same object with several symbols,
to clarify the steps we do in the computations.

%
%
%
 
 }

 Then  $\det (A)=2\rhoepsilon_1^{-3}\rhoepsilon_2^{-2}\rhoepsilon_3^{-1}(1+O(\rhoepsilon_1))$ and
\ba
 T^{({\bf m}),\beta}_\tau
=\frac {C\P_\beta}{ \omega_{45}\tau} \frac {\det (A)}{ \omega_{12}} \int_{(\R_+)^4}
e^{i\tau {\bf p}\cdotp y} 
h(Ay)y_4^{a-k_4+1}y_3^{a-k_3}y_2^{a-k_2+1}
y_1^{a-k_1+1}\,dy.
\ea
Using repeated integration by parts 
we see that
%
%
%

\beq\label{t 4 formula}\\ \nonumber
 T^{({\bf m}),\beta}_\tau&=& {C\det (A)\,{\P_\beta}}
\frac {(i\tau)^{-(12+4a-|\vec k_\b|)}(1+O(\tau^{-1}))}{ \rhoepsilon_4^{2(a-k_4+1+2)}\rhoepsilon_3^{2(a-k_3+1)}\rhoepsilon_2^{2(a-k_2+2)}\rhoepsilon_1^{2(a-k_1+1+2)}}\,.
\eeq
{\newtext Note that here and below $O(\tau^{-1})$ may depend also on $\rhoepsilon_j$, that is,
we have  $|O(\tau^{-1})|\leq C(\rhoepsilon_1,\rhoepsilon_2,\rhoepsilon_3,\rhoepsilon_4)\tau^{-1}.$}

To show that $ \mathcal G^{({\bf m})}({\bf v},{\bf b})$ is  non-vanishing 
we need to estimate $\P_{\beta_1}$ from below. In doing this we encounter 
the difficulty that $\P_{\beta_1}$ can go to zero, and moreover,
simple computations show that as  the pairs
$(b^{(j)},v^{(j)})$
satisfies the harmonicity condtion (\ref{harmonicity condition for symbol})
we have $ v_{(r)}^{ns} b^{(j)}_nb^{(j)}_s=O(\rhoepsilon_r+\rhoepsilon_j).$
%
%
%
%
However, to show that $ \mathcal G^{({\bf m})}({\bf v},{\bf b})$ is  non-vanishing 
 for some  ${\bf v}$ we consider a 
  particular choice of polarizations $v^{(r)}$,
 namely  
\beq\label{chosen polarization}
 v^{(r)}_{mk}= b^{(r)}_mb^{(r)}_k,\quad \hbox{for $r=2,3,4,$ but not for $r=1,5$}
\eeq
so that for $r=2,3,4$, we have
\ba
\hat g^{nm}b^{(r)}_nv^{(r)}_{mk}=0,\quad
\hat g^{mk}v^{(r)}_{mk}=0,\quad
\hat g^{nm}b^{(r)}_nv^{(r)}_{mk}-\frac 12 (\hat g^{mk}v^{(r)}_{mk})b^{(r)}_k=0.
\ea
Note that for this choice of $ v^{(r)}$ the linearized harmonicity conditions hold.
Moreover, for this choice of $ v^{(r)}$ we see that for $\rhoepsilon_j\leq \rhoepsilon_r^{100}$
\beq\label{vbb formula}
v_{(r)}^{ns} b^{(j)}_nb^{(j)}_s
 = \hat g(b^{(r)},b^{(j)})\, \hat g(b^{(r)},b^{(j)})
 = \rhoepsilon_r^4+O (\rhoepsilon_r^5). 
\eeq


In particular, when $\beta=\beta_1$, so that
 $k_{\b_1}=(6,0,0,0)$ and  
 the  polarizations are given by (\ref{chosen polarization}), we have
\ba
\P_{\beta_1}= (\D +O(\rhoepsilon_1))\rhoepsilon_1^4\cdotp\rhoepsilon_1^4\cdotp\rhoepsilon_1^4,
\ea
where $\D$ is the inner product of $v^{(1)}$ and $v^{(5)}$ given in (\ref{eq: def of D}). 
Then\hiddenfootnote{Indeed, when $\beta=\beta_1$ we have
\ba
T^{\beta_1,0}_\tau
\ba
T^{(1)}_\tau
&=&\frac {c_1^{\prime}}{4 } \cdotp {(i\tau)^{-(12 +4a-|\vec k_\b|)}(1+O(\tau^{-1}))}\vec \e^{-2\vec a}
 (\omega_{45}\omega_{12})^{-1} \rhoepsilon_4^{2(k_4-1+1)}\rhoepsilon_2^{2(k_2-1+1)}\rhoepsilon_3^{2(k_3+1)}\rhoepsilon_1^{2(k_1-1+1)}{\P_\beta}\\
 &=&{c_1^{\prime}}
 {(i\tau)^{-(12 +4a-|\vec k_\b|)}(1+O(\tau^{-1}))}\vec \e^{-\vec a-\vec 1}
 (\omega_{45}\omega_{12})^{-1}\rhoepsilon_4^{-2}\rhoepsilon_2^{-2}\rhoepsilon_3^{0}\rhoepsilon_1^{10}{\P_\beta}
 \\
 &=&{c_1^{\prime}}
 {(i\tau)^{-(12 +4a-|\vec k_\b|)}(1+O(\tau^{-1}))}\vec \e^{-2(\vec a+\vec 1)}
\rhoepsilon_4^{-4}\rhoepsilon_2^{-2}\rhoepsilon_3^{0}\rhoepsilon_1^{20}\D.
\ea
} 
the term $T^{({\bf m}),\beta_1}_\tau$, which later turns out to
have the strongest asymptotics in our considerations, has the asymptotics
\beq\label{leading term}
& &\hspace{-1cm}T^{({\bf m}),\beta_1}_\tau=\L_\tau,\quad\hbox{where}\\
\nonumber
& &\hspace{-1cm}\L_\tau=
C\det (A)
{(i\tau)^{-(6 +4a)}(1+O(\tau^{-1}))}\vec\rhoepsilon^{\,\vec d}
\rhoepsilon_4^{-4}\rhoepsilon_2^{-2}\rhoepsilon_3^{0}\rhoepsilon_1^{20}\D,\hspace{-1.5cm}
\eeq 
where $\vec\rhoepsilon=(\rho_1,\rho_2,\rho_3,\rho_4)$,
$\vec a=(a,a,a,a)$, and $\vec 1=(1,1,1,1)$.
To compare different terms,
we express
$\rhoepsilon_j$ in powers of $\rhoepsilon_1$ as explained in formula (\ref{eq: ordering of epsilons}), 
that is, we write $\rhoepsilon_4^{n_4}\rhoepsilon_2^{n_2}\rhoepsilon_3^{n_3}\rhoepsilon_1^{n_1}=\rhoepsilon_1^m$
with $m=100^3n_4+100^2n_2+100n_3+n_1$.
In particular, we will below write
\ba
& &\L_\tau=C_{\beta_1}(\vec\rhoepsilon)\,\tau^{n_0}
(1+O(\tau^{-1}))\quad\hbox{as $\tau\to \infty$ for each fixed $\vec \e$, and}\\ 
& &C_{\beta_1}(\vec\e)=c^{\prime}_{\beta_1}\,\rhoepsilon_1^{m_0}(1+o(\rhoepsilon_1))
\quad\hbox{as  } \rhoepsilon_1\to 0. 
\ea
 Below we will show that $c^{\prime}_{\beta_1}$
does not vanish for generic $(\vec x,\vec \xi)$ and $(x_5,\xi_5)$  and polarizations ${\bf v}$.
We will  consider below $ \beta\not= \beta_1$
and show that also  these terms have the asymptotics 
\ba
& &T^{({\bf m}),\beta}_\tau=C_{\beta}(\vec\rhoepsilon)\,\tau^{n}
(1+O(\tau^{-1}))\quad\hbox{as $\tau\to \infty$ for each fixed $\vec \e$, and}\\ 
& &C_{\beta}(\vec\e)=c^{\prime}_{\beta}\,\rhoepsilon_1^{m}(1+o(\rhoepsilon_1))
\quad\hbox{as  } \rhoepsilon_1\to 0. 
\ea
When\hiddenfootnote{
Alternatively, the ordering could be formulated as follows: 
We say that
\beq \label{Slava 1}
C' \tau^{-c'_\tau} \rhoepsilon_4^{c'_4} \rhoepsilon_3^{ c'_3} \rhoepsilon_2^{ c'_2} \rhoepsilon_1^{ c'_1}
\prec C \tau^{-c_\tau} \rhoepsilon_4^{c_4} \rhoepsilon_3^{ c_3} \rhoepsilon_2^{ c_2} \rhoepsilon_1^{ c_1} 
\eeq
if $C \neq 0$ and some of the following conditions is valid:
\begin{itemize}
\item [(i)] if $c_\tau < c'_\tau$;
\item [(ii)] if $c_\tau = c'_\tau$ but $c_4 <c_4'$;
\item [(iii)] if $c_\tau = c'_\tau$ and $c_4 = c_4'$ but $c_2 <c_2'$;
\item [(iv)] if $c_\tau = c'_\tau$ and $c_4 = c_4',\,c_2 =c_2' $ but $c_3 <c_3'$;
\item [(v)] if $c_\tau = c'_\tau$ and $c_4 = c_4',\,c_3 =c_3',\, c_2 =c_2' $ but $c_1 <c_1'$.
\end{itemize}
Note that this is an ordering not just a partial ordering.} we have that 
either $n\leq n_0$ and $m<m_0$,
or $n<n_0$, we say that 
$T^{({\bf m}),\beta}_\tau$ has weaker asymptotics than
$T^{({\bf m}),\beta_1}_\tau$ and denote
$T^{({\bf m}),\beta}_\tau \prec \L_\tau.$

{

As we consider here the asymptotic of five
small parameters $\tau^{-1}$ and $\e_i,$ $i=1,2,3,4,$ and compare in which order
we make them tend to $0$, let us explain the above ordering
in detail. Above, we have chosen the order: first $\tau^{-1}$, then $\e_4$, $\e_2$, $\e_3$ and finally,
 $\e_1$. 
 In correspondence with this choice we can introduce an ordering on all monomials $c \tau^{-c_\tau} \e_4^{n_4} \e_3^{n_3} \e_2^{n_2} \e_1^{n_1}$. Namely, we
say that
\beq \label{SL1}
C^\prime  \tau^{-n^\prime _\tau} \e_4^{n^\prime _4} \e_3^{ n^\prime _3} \e_2^{ n^\prime _2} \e_1^{ n^\prime _1}\prec C \tau^{-n_\tau} \e_4^{n_4} \e_3^{ n_3} \e_2^{ n_2} \e_1^{ n_1} 
\eeq
if $C \neq 0$ and one of the following holds
\begin{itemize}
\item [(i)] if $n_\tau < n^\prime _\tau$;
\item [(ii)] if $n_\tau = n^\prime _\tau$ but $n_4 <n_4^\prime $;
\item [(iii)] if $n_\tau = n^\prime _\tau$ and $n_4 = n_4^\prime $ but $n_2 <n_2^\prime $;
\item [(iv)] if $n_\tau = n^\prime _\tau$ and $n_4 = n_4^\prime ,\,n_2 =n_2^\prime  $ but $n_3 <n_3^\prime $;
\item [(v)] if $n_\tau = n^\prime _\tau$ and $n_4 = n_4^\prime ,\,n_3 =n_3^\prime ,\, n_2 =n_2^\prime  $ but $n_1 <n_1^\prime $.
\end{itemize}

Then, we can analyze terms $T^{({\bf m}),\beta}_{\tau,\sigma}$ and $
\tilde T^{({\bf m}),\beta}_{\tau,\sigma}$ in the formula for 
\beq \label{SL2}
\Theta_\tau^{(4)}=
\Theta_{\tau,\vec\e}^{(4)}=\sum_{\beta\in J_\ell}\sum_{\sigma\in \cell}
\left(T^{({\bf m}),\beta}_{\tau,\sigma} + 
\tilde T^{({\bf m}),\beta}_{\tau,\sigma} \right).
\eeq
Note that here the terms, in which the permutation $\sigma$ is either the identical permutation $id$ or
the permutation  $\sigma_0=(2,1,3,4)$, are the same.

\medskip

{\bf Remark 3.3.} We can find the leading order asymptotics of the strongest terms in the decomposition
(\ref{SL2}) using the  following algorithm.
 First, let us multiply $\Theta_{\tau,\vec\e}^{(4)}$ by $\tau^{\hat n_\tau}$,
where $\hat n_\tau= \min_\beta n_\tau(\beta)$. Taking then $\tau \to \infty$ will give non-zero contribution from
only those terms $T^{({\bf m}),\beta}_{\tau,\sigma}$ and 
$\tilde T^{({\bf m}),\beta}_{\tau,\sigma} $ where $n_\tau(\beta)=\hat n_\tau$.
This corresponds to step $(i)$ above.
Multiplying next by $\e_4^{-\hat n_4}$, where $\hat n_4= \min_\beta n_4(\beta)$ under
the condition that $n_\tau(\beta)=\hat n_\tau$ and taking $\e_4 \to 0$ corresponds 
to selecting terms $T^{({\bf m}),\beta}_{\tau,\sigma} $ with $n_\tau(\beta)=\hat n_\tau$ and $n_4(\beta)=\hat n_4$
and terms $ \tilde T^{({\bf m}),\beta}_{\tau,\sigma} $ with  $ \, \tilde n_\tau(\beta)=
\hat n_\tau$ and $\tilde n_4(\beta)=\hat n_4$.
This corresponds to step $(ii)$.
Continuing this process we  obtain a scalar value that gives 
the leading order asymptotics of the strongest terms in the decomposition
(\ref{SL2}).\medskip

The next results  tells  what are the strongest terms in (\ref{SL2}).

\begin{proposition} \label{SL:order} 
Assume that $v^{(r)}$, $r=2,3,4$ are given by (\ref{chosen polarization}). In 
(\ref{SL2}), the strongest term are $T^{({\bf m}),\beta_1}_{\tau}=T^{({\bf m}),\beta_1}_{\tau,id}$
and $T^{({\bf m}),\beta_1}_{\tau,\sigma_0}$
in the sense that for all $(\beta,\sigma)\not \not\in \{ (\beta_1,id),(\beta_1,\sigma_0)\}$ we have  $T^{({\bf m}),\beta}_{\tau,\sigma}\prec 
T^{({\bf m}),\beta_1}_{\tau,id}$. 
\end{proposition}

{\bf Proof.}} When  $\vec S_\b=(S_1,S_2)=({\bf Q}_0,{\bf Q}_0)$,
similar computations to the above ones yield
\ba
\tilde  T^{({\bf m}),\beta}_\tau
&=&\bra  u^\tau,h\,\cdotp  {\bf Q}_0(\B_4u_4\,\cdotp \B_3u_3)\,\cdotp {\bf Q}_0(\B_2u_2\,\cdotp \B_1u_1)\cet\\
&=&C\det (A) {{\P_\beta}} \frac {(i\tau)^{-(12+4a-|\vec k_\b|)}(1+O(\tau^{-1}))}{ \rhoepsilon_4^{2(a-k_4+2)}\rhoepsilon_3^{2(a-k_3+1+2)}\rhoepsilon_2^{2(a-k_2+2)}\rhoepsilon_1^{2(a-k_1+1+2)}}\,.
\ea

Let us next consider the case when  $\vec S_\b=(S_1,S_2)=(I,{\bf Q}_0)$.
Again, the computations similar to the above ones show that
\ba
 T^{({\bf m}),\beta}_\tau
&=&\bra {\bf Q}_0(u^\tau\,\cdotp \B_4 u_4 ), h\,\cdotp \B_3u_3\,\cdotp I(\B_2u_2\,\cdotp \B_1u_1)\cet\\
&=& {i\,C{\P_\beta} \det (A)}\frac {(i\tau)^{-(10+4a-|\vec k_\b|)}(1+O(\tau^{-1}))}{
 \rhoepsilon_4^{2(a-k_4+1+2)}\rhoepsilon_3^{2(a-k_3+1)}\rhoepsilon_2^{2(a-k_2+1)}\rhoepsilon_1^{2(a-k_1+1)}}\,
\ea
and\hiddenfootnote{In fact, when $(S_1,S_2)=(I,{\bf Q}_0)$, the term $\tilde  T^{({\bf m}),\beta}_\tau$ is similar to the 
corresponding term in the above case $(S_1,S_2)=({\bf Q}_0,I)$
up to a permutation of indexes, and thus there is
 no need to analyze these terms separately.} 
\ba
& &\hspace{-1cm}\tilde  T^{({\bf m}),\beta}_\tau
=\bra  u^\tau, h\,\cdotp {\bf Q}_0(\B_4u_4\,\cdotp \B_3\,u_3)\,\cdotp I(\B_2\,u_2\,\cdotp \B_1\, u_1)\cet\\
&=&{{\P_\beta} C\det (A)} 
\frac {(i\tau)^{-(10+4a-|\vec k_\b|)}(1+O(\tau^{-1}))}
{\rhoepsilon_4^{2(a-k_4+2)}\rhoepsilon_3^{2(a-k_3+1+2)}\rhoepsilon_2^{2(a-k_2+1)}\rhoepsilon_1^{2(a-k_1+1)}}.
\ea
When $\vec S_\b=(S_1,S_2)=({\bf Q}_0,I)$ we have $\tilde  T^{({\bf m}),\beta}_\tau
=T^{({\bf m}),\beta}_\tau
$ and
\ba
 T^{({\bf m}),\beta}_\tau
&=& \bra I(u^\tau \,\cdotp  \B_4 u_4), h\,\cdotp  \B_3 u_3\,\cdotp {\bf Q}_0(\B_2 u_2\,\cdotp \B_1 u_1)\cet,\\
&=& {i{\P_\beta}\,\det (A)C}
\frac {(i\tau)^{-(10+4a-|\vec k_\b|)} (1-O(\tau^{-1}))}{\rhoepsilon_4^{2(a-k_4+1)}\rhoepsilon_3^{2(a-k_3+1)}\rhoepsilon_2^{2(a-k_2+2)}\rhoepsilon_1^{2(a-k_1+1+2)}}\, ,
\ea
and finally when $\vec S_\b=(S_1,S_2)=(I,I)$ 
\ba
\tilde  T^{({\bf m}),\beta}_\tau
&=&\bra  u^\tau, h\,\cdotp I(\B_4u_4\,\cdotp \B_3u_3)\,\cdotp I(\B_2u_2\,\cdotp \B_1 u_1)\cet\\
&=&{{\P_\beta} C_a\det (A)}
\frac {(i\tau)^{-(8+4a-|\vec k_\b|)}(1+O(\tau^{-1}))}{\rhoepsilon_4^{2(a-k_4+1)}\rhoepsilon_3^{2(a-k_3+1)}\rhoepsilon_2^{2(a-k_2+1)}\rhoepsilon_1^{2(a-k_1+1)}}.
\ea

We consider all $\beta$ such that  $\vec S^\beta=({\bf Q}_0,{\bf Q}_0)$
but
$\b\not=\beta_1$. Then
\ba
 \tilde T^{({\bf m}),\beta}_\tau
=C\det (A)
(i\tau)^{-(6 +4a)}(1+O(\tau^{-1}))
\vec\rhoepsilon^{\,\vec d+2\vec k_\b}\rhoepsilon_4^{-2}\rhoepsilon_2^{-2}\rhoepsilon_3^{-4}\rhoepsilon_1^{-4}\cdotp{\P_{\beta}}
\ea where
$\vec k_\b$ is as in (\ref{k: collected}).  Note that for $\b=\beta_1$ we have 
$\P_{\beta_1}= (\D +O(\rhoepsilon_1))\rhoepsilon_1^4\cdotp\rhoepsilon_1^4\cdotp\rhoepsilon_1^4$
while $\b\not=\beta_1$ we just use an estimate
 ${\P_\beta}=O(1)$.
Then we see that
$\tilde  T^{({\bf m}),\beta}_\tau\prec \L_\tau$.

When  $\beta$ is such that $\vec S^\beta=({\bf Q}_0,I)$, we see that 
\ba
 T^{({\bf m}),\beta}_\tau
=C\det (A)
(i\tau)^{-(10 +4a-|\vec k_\b|)}(1+O(\tau^{-1}))
\vec\rhoepsilon^{\,-\vec d+2\vec k_\b}
\rhoepsilon_4^{0}\rhoepsilon_2^{-2}\rhoepsilon_3^{0}\rhoepsilon_1^{-4}{\P_\beta},\\
\tilde  T^{({\bf m}),\beta}_\tau
=C\det (A)
(i\tau)^{-(10 +4a-|\vec k_\b|)}(1+O(\tau^{-1}))
\vec\rhoepsilon^{\,-\vec d+2\vec k_\b}
\rhoepsilon_4^{-2}\rhoepsilon_2^{0}\rhoepsilon_3^{-4}\rhoepsilon_1^{0}{\P_\beta}\
\ea where ${\P_\beta}=O(1)$ and 
$\vec k_\b$ is as in (\ref{k: collected}) 
and hence $ T^{({\bf m}),\beta}_\tau\prec \L_\tau$ and 
$\tilde  T^{({\bf m}),\beta}_\tau\prec \L_\tau$.

When  $\beta$ is such that $\vec S^\beta=(I,{\bf Q}_0)$ or 
$\vec S^\beta=(I,I)$, using inequalities  of the type (\ref{k: collected}) and
 (\ref{k: collected}) in appropriate cases,
 we see that $T^{({\bf m}),\beta}_\tau\prec \L_\tau$ 
 and $\tilde  T^{({\bf m}),\beta}_\tau\prec \L_\tau$.

%
%
%
%
%


The above shows that all terms 
$T^{({\bf m}),\beta}_\tau$ and $\tilde  T^{({\bf m}),\beta}_\tau$ with
maximal allowed $k$.s
have asymptotics with the same power of $\tau$ but their $\vec \rhoepsilon$ asymptotics vary, and when the
asymptotic orders of $\rhoepsilon_j$
are given as explained in after (\ref{eq: ordering of epsilons}), there
is only one term, namely $\L_\tau=T^{({\bf m}),\beta_1}_\tau$, that has the strongest 
order asymptotics given in (\ref{leading term}).

Next we analyze the effect of the permutation $\sigma:
\{1,2,3,4\}\to \{1,2,3,4\}$  of the indexes $j$ of the waves $u_{j}$. We assume below that the permutation 
 $\sigma$ is not the identity map.
%

 Recall that in the computation (\ref{first asymptotical computation})
 there appears a term $\omega_{45}^{-1}\sim \rhoepsilon_4^{-2}$. 
Since this term does not appear 
  in the computations of the terms $\tilde  T^{({\bf m}),\beta}_{\tau,\sigma}$,
 we see that  $\tilde  T^{({\bf m}),\beta}_{\tau,\sigma}\prec \L_\tau$. 
Similarly, if
  $\sigma$  is such that $\sigma(4)\not =4$, the
 term $\omega_{45}^{-1}$ does not appear in the computation of
$ T^{({\bf m}),\beta}_{\tau,\sigma}$
 and hence
$ T^{({\bf m}),\beta_1}_{\tau,\sigma}\prec  \L_\tau$. 
Next we consider the permutations
 for which  $\sigma(4) =4$.

Next we consider $\sigma$ that is either
$\sigma=(3,2,1,4)$ or
$\sigma=(2,3,1,4)$.
 These terms are very similar
and thus we analyze the case when $\sigma=(3,2,1,4)$.
First we consider the case when $\beta=\beta_2$ is such that $\vec S^{\beta_2}=({\bf Q}_0,{\bf Q}_0)$,
 $\vec k_{\beta_2}=(2,0,4,0)$.  
This term appears in the analysis of
the term  $A^{(1)}[  u_{\sigma(4)},{\bf Q}(A^{(2)}[ u_{\sigma(3)},{\bf Q}(A^{(3)}[ u_{\sigma(2)}, u_{\sigma(1)}])])]$ when $(A^{(1)},A^{(2)},A^{(3)})=(A_2,A_1,A_2)$, see (\ref{A alpha decomposition2}). 
By a permutation of the indexes in (\ref{first asymptotical computation}) we obtain 
the formula 
\beq\label{Formula sic!}
T^{({\bf m}),{\beta_2}}_{\tau,\sigma}&=&
c_1^{\prime}\det (A){(i\tau)^{-(6 +4a)}(1+O(\tau^{-1}))}\vec \rhoepsilon^{\,\vec d}
\\ \nonumber
& & \cdotp (\omega_{45}\omega_{32})^{-1}  \rhoepsilon_4^{2(k_4-1)}\rhoepsilon_1^{2k_3}\rhoepsilon_2^{2(k_2-1)}\rhoepsilon_3^{2(k_1-1)}
 \P_{\beta_2},\\
\nonumber
\tilde \P_{{\beta_2}}&=&(v^{pq}_{(4)}b^{(1)}_pb^{(1)}_q)
(v^{rs}_{(3)}b^{(1)}_rb^{(1)}_s)
(v^{nm}_{(2)}b^{(3)}_nb^{(3)}_m)\D. 
\hspace{-1.5cm}
\eeq
Hence, in the case when we use the polarizations (\ref{chosen polarization}), we obtain
\ba
  T^{({\bf m}),\beta_2}_{\tau,\sigma}=
c_1{(i\tau)^{-(6 +4a)}(1+O(\tau^{-1}))}\vec \rhoepsilon^{\,\vec d}
 \rhoepsilon_4^{-4}\rhoepsilon_2^{-2}\rhoepsilon_3^{0+4}\rhoepsilon_1^{6+8}\D.
 \ea
Comparing the power of $\rhoepsilon_3$ in the above expression,
we see that in this case $ T^{({\bf m}),\beta_2}_{\tau,\sigma}\prec \L_\tau$. 
When $\sigma=(3,2,1,4)$, we see in  a straightforward way  also for 
 other $\beta$  for which  $\vec S^{\beta}=({\bf Q}_0,{\bf Q}_0)$, that  $ T^{({\bf m}),\beta}_{\tau,\sigma}\prec \L_\tau$.


When $\sigma=(1,3,2,4)$,
we see that for all $\beta$ with  $|\vec k_\beta|=6$,
\ba
  T^{({\bf m}),\beta}_{\tau,\sigma}\hspace{-1mm}&=&\hspace{-1mm}
C\det (A){(i\tau)^{-(6 +4a)}(1+O(\frac 1\tau))}\vec \rhoepsilon^{\vec d+2\vec k)}\omega_{45}^{-1}\omega_{13}^{-1}
 \rhoepsilon_4^{-2}\rhoepsilon_2^{0}\rhoepsilon_3^{-2}\rhoepsilon_1^{-2}\P_{\beta}\\
 \hspace{-1mm}
 &=&\hspace{-1mm}
C\det (A){(i\tau)^{-(6 +4a)}(1+O(\frac 1\tau))}\vec \rhoepsilon^{\,\vec d+2\vec k}
 \rhoepsilon_4^{-4}\rhoepsilon_2^{0}\rhoepsilon_3^{-2}\rhoepsilon_1^{-4}\P_{\beta}. \ea
 Here,  $\P_{\beta}=O(1)$.
Thus when $\sigma=(1,3,2,4)$, by comparing the powers of $\rhoepsilon_2$ we see that $ T^{({\bf m}),\beta}_{\tau,\sigma}\prec \L_\tau$. The same holds in the case when
$|\vec k_\beta|<6$.
The case when  $\sigma=(3,1,2,4)$ is similar to $\sigma=(1,3,2,4)$. 
This proves Proposition \ref{SL:order}\hfill\Box\medskip

\hiddenfootnote{
THIS WAS
Next, in the  case  when $\sigma=(3,1,2,4),$ we need to analyze the term
\ba
 & &\bra u^\tau ,A_2[u_4,{\bf Q}(A_2[u_2,{\bf Q}(A_1[u_1,u_3])])]\cet\\
&=& \bra u^\tau, u^{pq}_4\,{\bf Q} ((\p_p\p_q\p_r\p_s \p_j\p_k u_1)\,{\bf Q}(u^{rs}_2\, u^{jk}_3))\cet
+\hbox{similar terms}.
\ea
The first term on the right hand side, that we index by  $\beta=\beta_2$,
corresponds to $k_{\beta_2}=(0,0,6,0)$. Then we see that $\P_{\beta_2}=\P_{\beta_0}$ and
\ba
  T^{({\bf m}),\beta_2}_{\tau,\sigma}\hspace{-1mm}&=&\hspace{-1mm}
C\det (A){(i\tau)^{-(6 +4a)}(1+O(\frac 1\tau))}\vec \rhoepsilon^{\vec d}\omega_{45}^{-1}\omega_{32}^{-1}
 \rhoepsilon_4^{-2}\rhoepsilon_2^{0}\rhoepsilon_3^{-2}\rhoepsilon_1^{4}\P_{\beta_2}\\
 \hspace{-1mm}&=&\hspace{-1mm}
C\det (A){(i\tau)^{-(6 +4a)}(1+O(\frac 1\tau))}\vec \rhoepsilon^{\,\vec d}
 \rhoepsilon_4^{-4}\rhoepsilon_2^{0}\rhoepsilon_3^{-4}\rhoepsilon_1^{4}\P_{\beta_2}. \ea
Comparing the powers of $\rhoepsilon_2$ we see that $ T^{({\bf m}),\beta_2}_{\tau,\sigma}\prec \L_\tau$.
When $\sigma=(3,1,2,4)$ we see in a straightforward manner for also values of 
$\beta\not=\beta_2$ 
 that   $ T^{({\bf m}),\beta}_{\tau,\sigma}\prec \L_\tau$.}

}

Summarizing; we have analyzed the terms $ T^{({\bf m}),\beta}_{\tau,\sigma}$ 
corresponding to any $\b$ and all $\sigma$
except $\sigma=\sigma_0=(2,1,3,4)$.
Clearly, the sum  $ \sum_\b T^{({\bf m}),\beta}_{\tau,\sigma_0}$ 
is equal to the sum $ \sum_\b T^{({\bf m}),\beta}_{\tau,id}$.
Thus, when
the asymptotic orders of $\rhoepsilon_j$
are given  in (\ref{eq: ordering of epsilons})
and
the polarizations satisfy (\ref{chosen polarization}), we have
\beq\nonumber
\mathcal G^{({\bf m})}({\bf v},{\bf b})&=&\lim_{\tau\to \infty}
\sum_{\b,\sigma} \frac {T^{({\bf m}),\beta}_{\tau,\sigma}}{(i\tau)^{(6 +4a)}}
=\lim_{\tau\to \infty}
 \frac {2 T^{({\bf m}),\beta_1}_{\tau,id}
(1+O(\rhoepsilon_1))} {(i\tau)^{(6 +4a)}}\\
\label{eq: summary}
&=&2c_1\det (A)
(1+O(\rhoepsilon_1))\,\vec \rhoepsilon^{\,\vec d}
\rhoepsilon_4^{-4}\rhoepsilon_2^{-2}\rhoepsilon_3^{0}\rhoepsilon_1^{20}\D.\hspace{-1.5cm}
\eeq

Next we consider general polarizations.
Let $Y=\hbox{sym}(\R^{4\times 4})$ 
and consider the non-degenerate, symmetric, bi-linear form $B:(v,w)\mapsto  \hat g_{nj}\hat g_{mk} v^{nm}w^{jk}$ in $Y$.
Then
$ \D =B(v^{(5)},v^{(1)})$. 

 
Let $L(b^{(j)})$ denote the subspace of dimension $6$ of the symmetric
 matrices $v\in Y$ that satisfy 
 equation (\ref{harmonicity condition for symbol}) with covector $\xi=b^{(j)}$.
%
%

%

%
%


Let $\W({b^{(5)}})$ be the real analytic \MTEXT{submanifold 
of $(\R^4)^5\times (\R^{4\times 4})^3\times  (\R^{4\times 4})^6\times (\R^{4\times 4})^6$
consisting  of elements  $\eta=({\bf b}, {\underline v},V^{(1)},V^{(5)})$,
where ${\bf b}=(b^{(1)},b^{(2)},b^{(3)} ,b^{(4)},
b^{(5)})$} is a sequence of  light-like vectors {with the given vector $b^{(5)}$}, 
and ${\underline v} =(v^{(2)},v^{(3)},v^{(4)})$ 
satisfy $v^{(j)}\in L(b^{(j)})$ for all $j=2,3,4$, 
and $V^{(1)}=(v_p^{(1)})_{p=1}^6\in  (\R^{4\times 4})^6$ is a basis of $L(b^{(1)})$ and 
 $V^{(5)}=(v_p^{(5)})_{p=1}^6\in  (\R^{4\times 4})^6$ is sequence of vectors in $Y$ such 
 that $B(v_p^{(5)},v_q^{(1)})=\delta_{pq}$ for $p\leq q$.
 \MTEXT{Note that $\W({b^{(5)}})$ has two components where the orientation of
 the basis $V^{(1)}$ is different.}
 
By (\ref{Term type 1}) and (\ref{Term type 1}),  $\mathcal G^{({\bf m})}({\bf v},{\bf b})$ is linear in each $v^{(j)}$. For $\eta\in \W({b^{(5)}})$, we define
\ba
\kappa(\eta):=\det\bigg(\mathcal G^{({\bf m})}({\bf v}_{(p,q)},{\bf b})\bigg) _{p,q=1}^6,\ \hbox{where}\ {\bf v}_{(p,q)}=(v^{(1)}_p,v^{(2)},v^{(3)},v^{(4)},v_q^{(5)}).
\ea
Then $\kappa(\eta)$   \MTEXT{can be 
written as $\kappa(\eta)=A_1(\eta)/A_2(\eta)$
where $A_1(\eta)$ and $A_2(\eta)$
are real-analytic functions on $\W({b^{(5)}})$.
 In fact,  $A_2(\eta)$ is a product of terms
 $\hat g(b^{(j)},b^{(k)})^p$ with some positive integer $p$,
cf.\ \noextension{(\ref{BBB eq}) and} (\ref{t 4 formula}).}   
 
 Let us next  consider  the case when the sequence of the light-like vectors,
 ${\bf b}= ( b^{(1)}, b^{(2)},b^{(3)} , b^{(4)}, b^{(5)})$ given in
 (\ref{b distances}) with $\vec \rhoepsilon$ 
 given in (\ref{eq: ordering of epsilons})
 with some small $\rhoepsilon_1>0$ 
 and let the polarizations $ {\underline v}=( v^{(2)}, v^{(3)}, v^{(4)})$ be  such that $ v^{(j)}\in  L(b^{(j)})$, $j=2,3,4$, are those given
 by (\ref{chosen polarization}),
  and
 $ V^{(1)}=( v_p^{(1)})_{p=1}^6$  be a basis of $ L(b^{(1)})$.
 Let $ V^{(5)}=( v_p^{(5)})_{p=1}^6$ be vectors in $Y$ such 
 that $B( v_p^{(5)}, v_q^{(1)})=\delta_{pq}$ for $p\leq q$.
When $\rhoepsilon_1>0$ is small enough,  formula (\ref{eq: summary}) 
 yields that $\kappa( \eta)\not =0$ for   $ \eta=( {\bf b}, {\underline v}, V^{(1)}, V^{(5)})$.
 Since $\kappa(\eta)= \kappa(\eta)=A_1(\eta)/A_2(\eta)$ where
  \MTEXT{$A_1(\eta)$  and $A_2(\eta)$ are real-analytic} 
on $\L({b^{(5)}})$,
 we have that $\kappa(\eta)$ is non-vanishing and finite on an open and dense
 subset of the component of $\W({b^{(5)}})$ containing $ \eta$.
\MTEXT{The fact that  $\kappa(\eta)$ is  non-vanishing on a generic
 subset of the other component of $\W({b^{(5)}})$ can be seen by
 changing the orientation of $V^{(1)}$. This yields that claim.} 
\hfill \Box \medskip

{

\section{Observations in normal coordinates}

\label{sec: normal coordinates}
We have considered above the singularities of the metric $g$ in the wave gauge
coordinates.
 As the wave gauge coordinates may  also be non-smooth,
we do not know if the observed singularities are caused by the metric or
the coordinates. Because of this, we consider next the metric in normal coordinates.

Let  $v^{\vec \e}=(g^{\vec \e},\phi^{\vec \e})$ be the solution 
 of the  $\hat g$-reduced Einstein equations  (\ref{eq: adaptive model with no source}) with
 the source ${\bf f}_{\vec \e}$ given in (\ref{eq: f vec e sources}). 
 We emphasize that  $g^{\vec \e}$ is the metric in the
$(g,\hat g)$-wave gauge coordinates.


Let $(z,\eta)\in \U_{(z_0,\eta_0)}(\hat h)$, $\mu_{\vec \e}=\mu_{g^{\vec \e},z,\eta}$ 
and $(Z_{j,\vec \e})_{j=1}^4$ be  a frame of vectors obtained by
$g^{\vec \e}$-parallel continuation of some $\vec \e$-independent frame along the geodesic  
$\mu_{\vec \e}([-1,1])$ from $\mu_{\vec \e}(-1)$ to a point $p_{\vec \e}=\mu_{\vec \e}(\tilde r)$. Recall that $g$ and $\hat g$  coincide in the set $(-\infty,0)\times N$ that
contains the point $\mu_{\vec \e}(-1)$.
Let  
$\Psi_{\vec\e}: W_{\vec \e}\to  \Psi_{\vec \e}(W_{\vec \e})\subset
\R^4$ denote  normal coordinates
of $(\hattuM _0,g^{\vec \e})$ defined using the center $p_{\vec \e}$
and the frame $Z_{j,\vec \e}$.
We say that $\mu_{\vec \e}([-1,1])$ are observation geodesics,
and that  $\Psi_{\vec \e}$  are the 
 normal coordinates  associated to $\mu_{\vec \e}$ and $p_{\vec\e}=\mu_{\vec \e}(\tilde r)$.
%
 Denote $\Psi_{0}= \Psi_{\vec \e}|_{\vec \e=0}$, and  $W=W_0$. 
%
%
%
%
%
%
We also denote  $g_{\vec \e}= g^{\vec \e}$ and
$U_{\vec \e}=U_{g^{\vec \e}}$. 

%


\begin{lemma}\label{lem: sing detection in wave gauge coordinates 1}
Let $v^{\vec \e}=(g^{\vec \e},\phi^{\vec \e})$  and $\Psi_{\vec \e}$ be as above.
   Let
 $S\subset U_{\hat g}$ be a smooth 3-dimensional surface such that $p_0=p_{\vec \e}|_{\vec \e=0}\in S$ 
and
\beq\label{eq. measured field}
& & g^{(\alpha)}=\p_{{\vec \e}}^\alpha g^{\vec \e}|_{\vec \e=0},\quad
\phi_\ell^{(\alpha)}=\p_{{\vec \e}}^\alpha\phi_\ell^{\vec \e}|_{\vec \e=0},\quad\hbox{for } |\alpha|\leq 4,\ \alpha\in \{0,1\}^4,
\eeq
and assume that $ g^{(\alpha)}$ and $\phi_l^{(\alpha)}$ are in $C^\infty(W)$ for $|\a|\leq 3$ and
 $g^{(\alpha_0)}_{pq}|_W\in\I^{m_0}(W\cap S)$ and $ \phi_l^{(\alpha_0)}|_W\in\I^{m_0}(W\cap S)$
for $\a_0=(1,1,1,1)$. 
\smallskip

(i) Assume that   $S\cap W$ is empty. Then the tensors
$\p_{{\vec \e}}^{\alpha_0} ((\Psi_{\vec\e})_*g_{\vec \e})|_{\vec \e=0}$ and
 $\p_{{\vec \e}}^{\alpha_0} ((\Psi_{\vec\e})_*\phi_\ell ^{\vec \e})|_{\vec \e=0}$
 are  $C^\infty$-smooth 
in $\Psi_0(W)$.

\smallskip

(ii)
Assume that  
 $ \mu_{0}([-1,1])$ intersets $S$
 transversally  at $p_0$.
 Consider the conditions

 \smallskip

 \noindent (a) There is a 2-contravariant tensor field $v$ that is a smooth section of $TW\otimes TW$  such that
$v(x)\in T_xS\otimes T_xS$ for $x\in S$  
and the principal symbol of
$ \bra v,g^{(\alpha_0)}\cet |_W\in \I^{m_0}(W\cap S)$
is non-vanishing at $p_0$.
 \smallskip

 \noindent (b) The principal symbol of $ \phi_\ell^{(\alpha_0)}|_W\in\I^{m_0}(W\cap S)$ is non-vanishing
at $p_0$ for some $\ell=1,2,\dots,L$.
 \smallskip

If (a) or (b) holds, then either  $\p_{{\vec \e}}^{\alpha_0} ((\Psi_{\vec\e})_*g_{\vec \e})|_{\vec \e=0}$ or
 $\p_{{\vec \e}}^{\alpha_0} ((\Psi_{\vec\e})_*\phi_\ell ^{\vec \e})|_{\vec \e=0}$ is not $C^\infty$-smooth 
in $\Psi_0(W)$.

\end{lemma}


\noindent
{\bf Proof.}  (i)  is obvious.

 (ii) 
%
%
%
%
Denote $\gamma_{\vec \e}(t)=\mu_{\vec \e}(t+\tilde r)$.
Let $X:W_0\to V_0\subset \R^4$, $X(y)=(X^j(y))_{j=1}^4$ be  local coordinates in $W_0$ such
that $X(p_0)=0$ and $X(S\cap W_0)=\{(x^1,x^2,x^3,x^4)\in V_0;\ x^1=0\}$  
and $y(t)=X( \gamma_0(t))=(t,0,0,0)$. 
 Note that the   coordinates $X$ are
independent of ${\vec \e}$.
We assume that the vector fields $Z_{\vec \e,j}$ defining the normal
coordinates are such that 
  $Z_{0,j}(p_0)= \p/\p X^j$.
To do computations
in local coordinates, let us denote
\ba
\tilde g^{(\alpha)}=\p_{{\vec \e}}^\alpha (X_*g^{\vec \e})|_{\vec \e=0},\quad\tilde \phi_\ell^{(\alpha)}=\p_{{\vec \e}}^\alpha (X_*\phi_\ell^{\vec \e})|_{\vec \e=0},
\quad\hbox{for } |\alpha|\leq 4,\ \alpha\in \{0,1\}^4.
\ea
Let $v$ be a tensor field given in (a) such that 
in the $X$ coordinates $v(x)=v^{pq}(x)\frac\p{\p x^p}\frac\p{\p x^q}$ 
so that $v^{pq}(0)=0$ if $(p,q)\not \in\{2,3,4\}^2$ 
at the point $0=X(p_0)$ \MTEXT{and the functions $v^{pq}(x)$ do not depend on $x$, that is, $v^{pq}(x)=\hat v^{mq}\in \R^{4\times 4}$.}
Let    $R^{\vec \e}$  be the curvature tensor  of $g^{\vec \e}$
and  define the functions
\ba
h^{\vec\e}_{mk}(t)=g^{\vec \e}(R^{\vec \e}(\dot \gamma_{\vec \e}(t),Z^{\vec \e}_m(t))\dot \gamma_{\vec \e}(t),Z^{\vec \e}_k(t)),\quad
J_v(t)=
\p_{{\vec \e}}^{\alpha_0}( \hat v^{mq}\,
h^{\vec \e}_{mq}(t))|_{{\vec \e}=0}.
\ea 

\MTEXT{The function $J_v(t)$ is an invariantly defined function on the curve
$\gamma_0(t)$ and thus it can be computed in any coordinates.}
If $\p_{{\vec \e}}^{\alpha_0} ((\Psi_{\vec\e})_*g_{\vec \e})|_{\vec \e=0}$ would
be smooth near $0\in \R^4$, then the
function $J_v(t)$ would be smooth near $t=0$. To show that the $\p_{\vec\e}^{\alpha_0}$-derivatives of the metric tensor  in the normal
coordinates are not smooth, 
we need to show that 
$J_v(t)$ is non-smooth  at $t=0$ for some values of $\hat v^{mq}$.
%
%
%

We will work in the $X$ coordinates and 
denote $\tilde R(x)=\p_{\vec \e}^{\alpha_0} X_*(R_{\vec \e})(x)|_{{\vec \e}=0}$.
Moreover, 
$\tilde \gamma^j$ are  the analogous
4th order $\e$-derivatives and we denote $\tilde g=\tilde g^{(\alpha_0)}$. 
  For simplicity we also denote
   $X_*\hat g$ and $ X_*\hat \phi_l$  by $\hat g$ and $\hat \phi_l$, respectively.

We analyze the functions of $t\in I=(-t_1,t_1)$, e.g., $a(t)$, where  $t_1>0$ is small. We say
that $a(t)$ is of order $n$ if $a(\,\cdotp)\in \I^{n}(\{0\}).$
By the assumptions of the theorem, $\big(\p_x^\beta \tilde g_{jk}\big)(\gamma(t))\in \I^{m_0+\beta_1}(\{0\})$
when $\beta=(\beta_1,\beta_2,\beta_3,\beta_4)$.
 \MTEXT{It follows from (\ref{eq. measured field}) and the linearized equtions
 for the parallel transport that}  $(\tilde \gamma(t),
 \p_t \tilde \gamma(t))$ and $\tilde Z_k$ 
 are in $ \I^{m_0}(\{0\})$. 
The above analysis shows that $\tilde R|_{\gamma_0(I)}\in  \I^{m_0+2}(\{0\})$. Thus in the $X$ coordinates
$\p_{\vec \e }^{\alpha_0}(h^{\vec \e}_{mk}(t))|_{{\vec \e}=0}\in \I^{m_0+2}(\{0\})$ can be written as
\ba
& &\hspace{-4mm}\p_{\vec \e }^{\alpha_0}(h^{\vec \e}_{mk}(t))|_{{\vec \e}=0}
= \hat g( \tilde  R(\dot  \gamma_0(t),\hat Z_m(t))
\dot \gamma_0(t),\hat Z_k(t))+\hbox{s.t.}\\
&\hspace{-4mm}=&\hspace{-4mm}\frac 12
\bigg(\frac \p{\p x^1}(\frac  {\p \tilde  g_{km}}{\p x^1}+\frac  {\p \tilde  g_{k1}}{\p x^m}-\frac  {\p \tilde  g_{1m}}{\p x^k})-\frac \p{\p x^m}(\frac  {\p \tilde  g_{k1}}{\p x^1}+\frac  {\p \tilde  g_{k1}}{\p x^1}-\frac  {\p \tilde  g_{11}}{\p x^k})\bigg)+\hbox{s.t.},
\ea
where all $s.t.=$
 ``smoother terms''
 are in  $\I^{m_0+1}(\{0\})$.
%

Consider next the case (a).
Assume that   for given $(k,m)\in\{2,3,4\}^2$, the principal symbol of $ \tilde g^{(\alpha_0)}_{km}$
is non-vanishing at $0=X(p_0)$. Let $v$ be such a tensor field that $v^{mk}(0)= v^{km}(0)\not=0$
and $ v^{in}(0)=0$ when $(i,n)\not\in \{(k,m),(m,k)\}$.
Then
the above yields (in the formula below, we do not sum over $k,m$)
\ba
J_v(t)=v^{ij}\p_{\vec \e}^{\alpha_0}(h^{\vec \e}_{ij}(t))|_{\vec \e=0}
=\frac {e(k,m)}2\hat v^{km} 
\bigg(\frac \p{\p x^1}\frac  {\p \tilde  g_{km}}{\p x^1}\bigg)
+\hbox{s.t.},
\ea
where $e(k,m)= 2-\delta_{km}$. 
Thus the principal symbol of $J_v(t)$  in $\I^{m_0+2}(\{0\})$ is non-vanishing and
 $J_v(t)$ is not a smooth function. Thus in this case
$\p_{{\vec \e}}^{\alpha_0} ((\Psi_{\vec\e})_*g_{\vec \e})|_{\vec \e=0}$  is not
smooth.


Next, we consider   the case (b).
Assume that there is $\ell$ such that  the principal
symbol of the field 
$\tilde \phi_\ell^{(\alpha_0)}$
 is non-vanishing. 
As  $\p_t\tilde \gamma(t)\in \I^{m_0}(\{0\})$, we see that $\tilde \gamma(t)\in \I^{m_0-1}(\{0\})$.
Then as $\phi_\ell^{\vec \e}$ are scalar fields, 
\ba
& & j_\ell(t)=
\p_{{\vec \e}}^{\alpha_0} \bigg(
 \phi^{\vec \e}_\ell(\gamma_{\vec \e}(t)))\bigg)
\bigg|_{\vec \e=0}
=\tilde  \phi_\ell(\gamma(t))
+
\hbox{s.t.},
\ea
where $ j_\ell \in  \I^{m_0}(\{0\})$ and the smoother terms $(s.t.)$ are in $ \I^{m_0-1}(\{0\})$.
Thus in the case (b),
%
$j_\ell(t)$ is not smooth and hence both $\p_{{\vec \e}}^{\alpha_0} ((\Psi_{\vec\e})_*g_{\vec \e})|_{\vec \e=0}$ and
 $\p_{{\vec \e}}^{\alpha_0} ((\Psi_{\vec\e})_*\phi_\ell ^{\vec \e})|_{\vec \e=0}$ cannot be smooth.
%
\hfill \Box \medskip

\extension{

\begin{center}

\psfrag{1}{$q$}
\psfrag{2}{\hspace{-2mm}$y$}
\includegraphics[width=5.5cm]{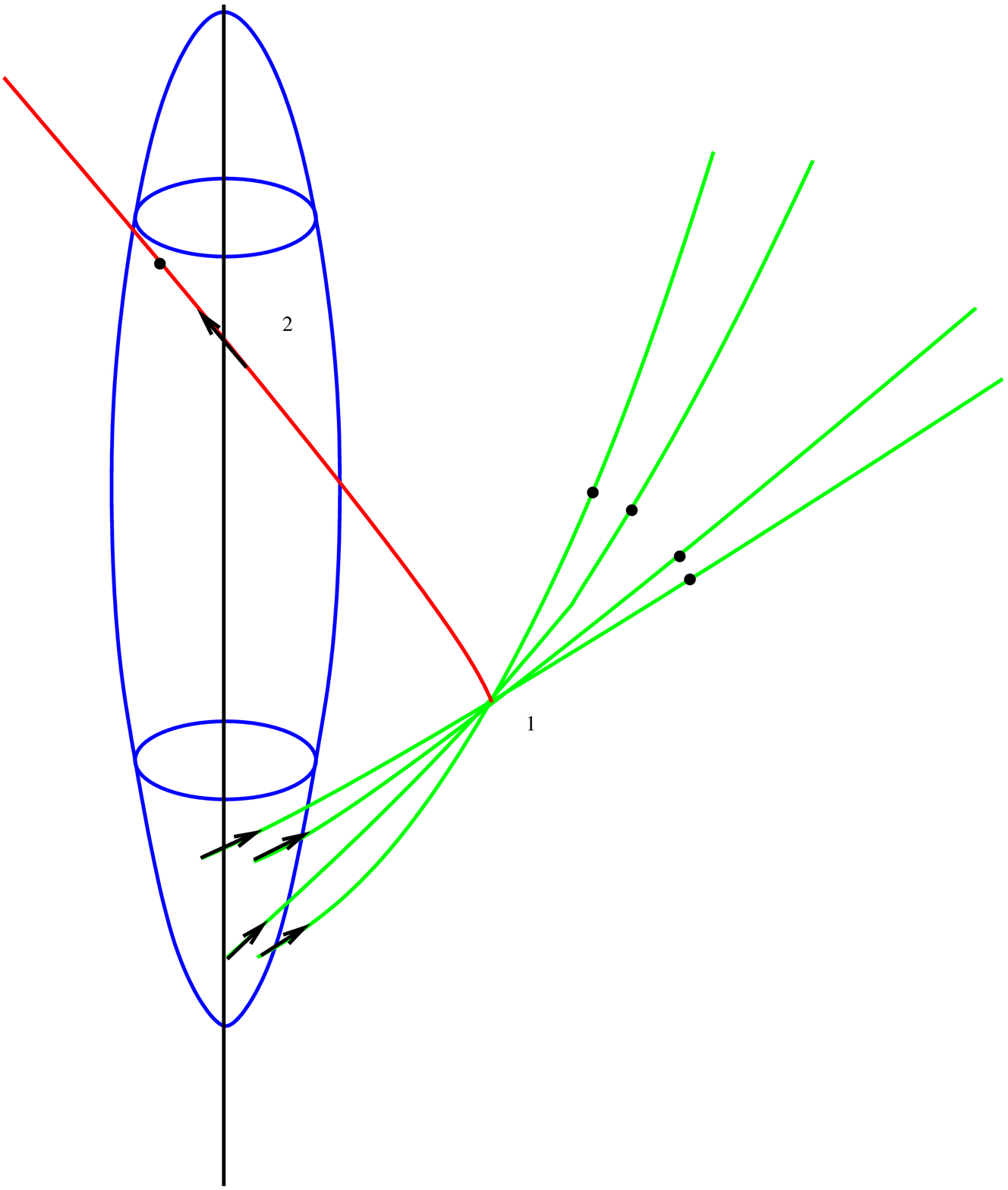}
\end{center}
 {\it FIGURE A6: A schematic figure where the space-time is represented as the  3-dimensional set $\R^{1+2}$. The light-like geodesic emanating from the point $q$  is  shown
as a red curve. The point $q$ is the intersection of
light-like geodesics corresponding to the starting points
and directions $(\vec x,\vec \xi)=((x_j,\xi_j))_{j=1}^4$. A light like geodesic
starting from $q$ passes through the  point $y$ and has the direction $\eta$ at $y$. 
The black points are the first conjugate points on the geodesics
$\gamma_{x_j(t_0),\xi_j(t_0)}([0,\infty))$, $j=1,2,3,4,$ and $\gamma_{q,\zeta}([0,\infty))$.
The figure shows the case when the interaction condition (I) is satisfied for $y\in U_{\hat g}$ with light-like
vectors $(\vec x,\vec \xi)$.
}
\medskip

}

We use now the results above to detect singularities
in normal coordinates.
 We say
that the {\it interaction condition} (I) is satisfied for $y\in U_{\hat g}$ with light-like
vectors $(\vec x,\vec \xi)=((x_j,\xi_j))_{j=1}^4$ and $t_0\geq0$,
with parameters  $(q,\zeta,t)$,
 if
 \medskip
}

{\bf (I)} There exist
 $q\in \bigcap_{j=1}^4\gamma_{x_j(t_0),\xi_j(t_0)}((0,{\bf t}_j))$, ${\bf t}_j=\rho(x_j(t_0),\xi_j(t_0))$,
 $\zeta\in L^+_q(\hattuM _0,\hat g)$ and $t\geq 0$
 such that $y=\gamma_{q,\zeta}(t)$.

 \medskip

 \noindent

{For  ${\bf f}_1\in \mathcal I_{C}^{n+1}(Y((x_1,\xi_1);t_0,s_0))$
the wave ${\bf Q}{\bf f}_1={\mathcal M}^{(1)}$ is a solution of the
linear wave equation. Thus, when ${\bf f}_1$ runs through the set $\mathcal I_{C}^{n+1}(Y((x_1,\xi_1);t_0,s_0))$ the union of the sets 
$\hbox{WF}\,({\bf Q}{\bf f}_1)$ is the manifold 
$\Lambda((x_1,\xi_1);t_0,s_0)$. Thus, 
${\cal D}(\hat g,\hat \phi,\e)$
determines  $\Lambda((x_1,\xi_1);t_0,s_0))\cap T^*U_{\hat g}$.
In particular, using  these data we can determine the geodesic segments 
$\gamma_{x_1,\xi_1}(\R_+)\cap U_{\hat g}$
 for all $x_1\in U_{\hat g}$, $\xi_1\in L^+M_0$.

Below, in  $T\hattuM _0$ 
we use the Sasaki metric corresponding to $\hat g^+$.
Moreover, let
 $s^\prime\in (s_-+r_2,s_+)$, $0<r_2<r_1$ and
 $B_j\subset U_{\hat g}$ be open sets
such that, 
cf.\  (\ref{observer neighborhood with hat}) and (\ref{eq: source causality condition}), 
for some $r_0^\prime\in (0,r_0)$  we have 
\beq\label{eq: admissible Bj}
& &\hbox{
$B_j\subset \subset I_{\hat g}( \mu_{\hat g}(s^\prime-r_2),\mu_{\hat g}(s^\prime)),$ and }\\
\nonumber
& &\hbox{
   for all $g^\prime\in \V(r_0^\prime)$, 
$B_j\cap J ^+_{g^\prime}(B_k)=\emptyset$ for all 
$j\not =k$.}
\eeq

} 

We say that $y\in U_{\hat g}$ 
  satisfies  the singularity {\it detection condition} ({D}) with light-like directions $(\vec x,\vec \xi)$
  and $t_0,\hat s>0$ 
 if 
 \medskip  

{\bf (D)} \MTEXT{For any $s,s_0\in (0,\hat s)$ and   $j=1,2,3,4$
there are $(x_j^{\prime},\xi_j^{\prime})$
 in the
$s$-neighborhood of $(x_j,\xi_j)$, 
$(2s)$-neighborhoods  $B_j$ of  $x_j$, in $(U_{\hat g},\hat g^+)$, satisfying (\ref{eq: admissible Bj}), and 
 sources ${\bf f}_j\in \mathcal I_{C}^{n+1}(Y((x_j^{\prime},\xi_j^{\prime});t_0,s_0))$
 having the LS property (\ref{eq: LSp}) in $C^{s_1}(M_0)$ with a family $\F_j(\e)$. This family is supported in $B_j$ and satisfies $\p_\e \F_j(\e)|_{\e=0}={\bf f_j}$.
Moreover, let 
$u_{\vec \e}$ be the 
 solution    of (\ref{eq: adaptive model with no source}) 
with the source 
$\F_{\vec \e}=\sum_{j=1}^4 \F_j(\e_j)$ and
$\mu_{\vec \e}([-1,1])$  be
observation geodesics with $y=\mu_{0}(\tilde r)$,
and   $\Psi_{\vec \e}$ be the normal coordinates
associated to $\mu_{\vec\e}$ at $\mu_{\vec\e}(\tilde r)$}.
 Then
  $\p_{{\vec \e}}^{\alpha_0} ((\Psi_{\vec\e})_*g_{\vec \e})|_{\vec \e=0}$ or
 $\p_{{\vec \e}}^{\alpha_0} ((\Psi_{\vec\e})_*\phi_\ell ^{\vec \e})|_{\vec \e=0}$ is not $C^\infty$-smooth near $0=\Psi_0(y)$.

\begin{lemma}\label{lem: sing detection in in normal coordinates 2}
Let $(\vec x,\vec \xi)$,   and ${\bf t}_j$ with $j=1,2,3,4$ and $t_0>0$
satisfy (\ref{eq: summary of assumptions 1})-(\ref{eq: summary of assumptions 2}).
Let $t_0,\hat s>0$ be sufficiently small and assume  that   $y\in  \V((\vec x,\vec \xi),t_0)\cap U_{\hat g}$ satisfies $y\not \in 
 {\mathcal Y}((\vec x,\vec \xi);t_0,\hat s)\cup\bigcup_{j=1}^4
 \gamma_{x_j,\xi_j}(\R) $. 
Then

(i) If $y$ does not satisfy condition (I)  with  $(\vec x,\vec \xi)$ and  $t_0$, then 
$y$ does not satisfy
condition ({D}) with  $(\vec x,\vec \xi)$ and  $t_0,\hat s>0$.

(ii) Assume  
$y\in U_{\hat g}$ satisfies condition (I) with  $(\vec x,\vec \xi)$ and  $t_0$ and parameters $q,\zeta$, and  $0<t<\rho(q,\zeta)$.
Then $y$ satisfies condition ({D}) with  $(\vec x,\vec \xi)$, $t_0$, and  any sufficiently
small $\hat s>0$.

(iii) Using the 
data set ${\cal D}(\hat g,\hat \phi,\e)$ we can determine
 \MTEXT{whether} the condition ({D}) is valid  for the given point $y\in W_{\hat g}$  with the 
parameters $(\vec x,\vec \xi)$, $\hat s$, and $t_0$  or not.

   \end{lemma}

   \noindent   {\bf Proof.}  
      (i)   
   If $y\not \in 
{\mathcal Y}((\vec x,\vec \xi);t_0,\hat s)\cup\bigcup_{j=1}^4
 \gamma_{x_j,\xi_j}(\R) $, the same condition holds also for $(\vec x^\prime,\vec \xi^\prime)$ close 
 to ($ \vec x,\vec \xi)$.
   Thus
Prop.\ \ref{lem:analytic limits A} 
and Lemma \ref{lem: sing detection in wave gauge coordinates 1} imply that (i) holds.

(ii) Let  $\mu_{\vec \e}([-1,1])$ be  observation geodesics
with $y=\mu_{0}(\tilde r)$ and $\Psi_{\vec \e}$ be the normal coordinates associated 
to $\mu_{\vec \e}([-1,1])$ \MTEXT{at $\mu_{\vec\e}(\tilde r)$.}
  Our aim is to show that there is a source ${\bf f}_{\vec \e}$ described in (D)
such that 
$\p_{\vec \e}^4((\Psi_{\vec \e})_*g_{\vec \e})|_{\vec \e=0}$  or $\p_{\vec \e}^4((\Psi_{\vec \e})_*\phi_{\vec \e})|_{\vec \e=0}$  is not
$C^\infty$-smooth at $y=\gamma_{q,\zeta}(t)$.

Let $\eta=\p_t \gamma_{q,\zeta}(t)$ and 
 denote
$(y,\eta)=(x_5,\xi_5)$.
\MTEXT{Let $t_j>0$ be such that  $\gamma_{x_j,\xi_j}(t_j)=q$
and denote $b_j=(\p_t \gamma_{x_j,\xi_j}(t_j))^\flat$, $j=1,2,3,4$. 
Also, let us denote  $b_5=\zeta^\flat$ and $t_5=-t$, so that
$q=\gamma_{x_5,\xi_5}(t_5)$ and $b_5=(\dot\gamma_{x_5,\xi_5}(t_5))^\flat$.}


  By
assuming that $V\subset U_{\hat g}$ is a sufficiently small neighborhood of $y$,
we have that $S:=\L_{\hat g}^+(q)\cap V$ is a smooth 3-submanifold,
{as $t<\rho(q,\zeta)$}.

\MTEXT{Let $u_\tau={\bf Q}_{\hat g}^* F_\tau$ be a gaussian beam, produced by a source 
$F_\tau(x)=
F_\tau(x;x_5,\xi_5)$  and function $h(x)$
 given in (\ref{Ftau source}).}
Then    we can use the techniques of \cite{KKL,Ralston},
   see also \cite{Babich,Katchalov}, 
to obtain a result analogous to Lemma
\ref{lem: Lagrangian 1} for the propagation of singularities  along the geodesic
$\gamma_{q,\zeta}([0,t])$: 
We have that when $h(x_5)=\Hsymbol$ and  $w$ 
 is  the principal symbol of $u_\tau$ at $(q,b_5)$,
then
$w=(R_{(5)})^*\Hsymbol$, where  $R_{(5)}$ is 
a bijective linear map.
\extension{
The map $(R_{(5)})^*$  is similar to the map $R_{(1)}(q,b_5;x_5,\xi_5)$ considered in the formula (\ref{eq: R propagation}) and it is obtained by
solving a system of linear ordinary differential equations along the geodesic
connecting $q$ to $x_5$.} 

Let $s\in (0,\hat s)$ be sufficiently small {and denote $b_5^{\prime}=b_5$.}
Using Propositions\ \ref{lem:analytic limits A} and  \ref{singularities in Minkowski space},  
we see that there are
$(b_j^{\prime})^\sharp\in L^+_q\hattuM _0$, $j=1,2,3,4$,
in the  $s$-neighborhoods of  $b_j^\sharp\in L^+_q\hattuM _0$,
vectors $v^{(j)}\in \R^{10+L}$, $j\in \{2,3,4\}$, 
linearly independent  vectors $v^{(5)}_p\in \R^{10+L}$, $p=1,2,3,4,5,6,$
and linearly independent vectors $v^{(1)}_r\in \R^{10+L}$, $r=1,2,3,4,5,6$,
that have the following properties:
\medskip

(a) All  $v^{(j)}$, $j=2,3,4$, and $v^{(1)}_r$, $r=1,2,\dots,6$  satisfy the harmonicity 
conditions for the symbols  (\ref{harmonicity condition for symbol})
with the covector $\xi$ being $b_j^{\prime}$ and $b_1^{\prime}$,  respectively.
\smallskip

(b)
Let  $X_1=\hbox{span}(\{v^{(1)}_r;\ r=1,2,3,\dots,6\})$ and $X_5=\hbox{span}(\{v^{(5)}_p;\ p=1,2,3,\dots,6\})$.
If
 $v^{(5)}\in X_5\setminus\{0\}$ then 
 there exists
 a vector $v^{(1)}\in X_1$
such that for  ${\bf v}=(v^{(1)},v^{(2)},v^{(3)},v^{(4)},v^{(5)})$
 and ${\bf b}^{\prime}=(b_j^{\prime})_{j=1}^5$
 we have
 $\mathcal G({\bf v},{\bf b}^{\prime})\not =0$.
 \medskip

Let $Y_5:=\hbox{sym}(T_yS\otimes T_yS)\times \R^L$.
Since the codimension  of $((R_{(5)})^*)^{-1}
Y_5$ in $\hbox{sym}(T_qM\otimes T_qM)\times \R^L $ is 4 and the dimension
of $X_5$ is 6, we see that dimension of the intersection $Z_5=X_5\cap 
((R_{(5)})^*)^{-1}Y_5$  is at least 2.
 \MTEXT{Thus there exist 
$v^{(j)}\not =0$, $j=1,2,3,4,5$  that satisfy
the above conditions (a) and $(b)$ and $v^{(5)}\in Z_5$.
Let $\Hsymbol^{(5)}=((R_{(5)})^*)^{-1}v^{(5)}\in 
((R_{(5)})^*)^{-1}Z_5 $.}

  Let $s_0\in (0,\hat s)$ and 
$x^{\prime}_j= \gamma_{q,(b^{\prime}_j)^\sharp}(-t_j)$, and
$\xi^{\prime}_j=\p_t \gamma_{q,(b^{\prime}_j)^\sharp}(-t_j)$, $j=1,2,3,4$.
 We denote $(\vec x^{\prime},\vec \xi^{\prime})=((x_j^{\prime},\xi_j^{\prime}))_{j=1}^4$.



%
%
%

Moreover,  let $B_j$ be 
neighborhoods of  $x_j^{\prime}$ satisfying (\ref{eq: admissible Bj}).
Then,
by condition $\mu$-LS  there are 
 sources ${\bf f}_j\in \mathcal I_{C}^{n+1}(Y((x_j^{\prime},\xi_j^{\prime});t_0,s_0))$
 having the LS property (\ref{eq: LSp}) in $C^{s_1}(M_0)$ with some family $\F_j(\e)$
 supported in $B_j$ such that the principal symbols of ${\bf f}_j$  at 
 $(x_j^{\prime}(t_0),(\xi_j^{\prime}(t_0)^\flat)$
are equal to $w^{(j)}=R_j^{-1}v^{(j)}
$,
where $R_j=R_{(1)}(q,b_j^\prime ;x_j^{\prime}(t_0),(\xi_j^{\prime}(t_0))^\flat)$ are defined by formula (\ref{eq: R propagation}). 
Then the principal symbols of
${\bf Q}_{\hat g}{\bf f}_j$
at $(q,b_j^{\prime})$ are equal to $v^{(j)}$. 
Let $u_{\vec \e}=(g_{\vec \e}-\hat g,\phi_{\vec \e}-\hat \phi)$ be the solution of (\ref{eq: notation for hyperbolic system 1})
\MTEXT{corresponding to  $\F_{\vec \e}=\sum_{j=1}^4 \F_j(\e_j)$ with
$\p_{ \e_j} \F_j(\e_j)|_{\e_j=0}= {\bf f}_j $.}
When $s_0$ is small enough,
$\M^{(4)}=\p_{\vec \e}^4u_{\vec \e}|_{\vec \e=0}$  is a conormal
distribution in the neighborhood $V$ of $y$ and 
$\M^{(4)}|_V\in \I(S)$. 
By Propositions\ \ref{lem:analytic limits A} and  \ref{singularities in Minkowski space}, the inner product $\bra F_\tau,\M^{(4)}\cet_{L^2(U_{\hat g})}$ is not of order $O(\tau^{-N})$
for all $N>0$, so  that $\M^{(4)}$  is not smooth
near $y$, see \cite{GrS}.  \MTEXT{Let $h(x)\in C_0^\infty(U_{\hat g})$ be
such that $h(x_5)=
\Hsymbol^{(5)}$.}
The above implies that the principal symbol of 
the function $x\mapsto \bra h(x),\M^{(4)}(x)\cet_{\B^L}$ is
not vanishing at $(y,\eta^\flat)$. Thus
the principal symbols of the functions (\ref{eq. measured field}) are not vanishing
\MTEXT{and as $h(x_5)\in 
((R_{(5)})^*)^{-1}Z_5$, 
%
%
conditions (ii) in Lemma \ref{lem: sing detection in wave gauge coordinates 1} 
are satisfied. Thus
either
$\p_{\vec \e}^4((\Psi_{\vec \e})_*g_{\vec \e})|_{\vec \e=0}$ or $\p_{\vec \e}^4((\Psi_{\vec \e})_*\phi_{\vec \e})|_{\vec \e=0}$
is not $C^\infty$-smooth at} $0=\Psi_0(y)$.
Thus,  condition 
({D}) is valid for $y$. This proves (ii). 
(iii) Below, we will assume that $ {\cal D}(\hat g,\hat \phi,\hat \e)$ is given with $\hat \e>0$.

To verify (D), we need to consider the solution $v_{\vec \e}=(g_{\vec \e},\phi_{\vec \e})$ in  the wave gauge coordinates. For a general element
 $[(U_g,g|_{U_g},\phi|_{U_g},F|_{U_g})]\in {\cal D}(\hat g,\hat \phi,\e)$
  we encounter the difficulty that we do not know the wave gauge coordinates in the set $(U_g,g)$. However,   we construct the wave map (i.e.\ the wave gauge) 
  coordinates \MTEXT{for sources $F$  of a special form.} Below, we give the proof in several steps.
  
  {\bf Step 1.} Let $s_-\leq s^\prime \leq s_0$ and $r_2\in (0,r_1)$ be so small that  
   $I_{\hat g}( \mu_{\hat g}(s^\prime-r_2),\mu_{\hat g}(s^\prime))\subset W_{\hat g}$.
  Assume that $\hat \e$  is small enough and that we are given an element
$[(U_g,g|_{U_g},\phi|_{U_g},F|_{U_g})]\in {\cal D}(\hat g,\hat \phi,\hat \e)$
such that  $F$  
is supported in 
 $I_{g}( \mu_{ g}(s^\prime-r_2),\mu_{ g}(s^\prime))\subset W_g$.
As  we know $(U_g,g|_{U_g})$, assuming that $\hat \e$ is small enough,
we can find the wave map $\Psi:
 I_{g}( \mu_{ g}(s^\prime-r_2),\mu_{ g}(s^\prime))\to U_{\hat g}$ by solving (\ref{C-problem 1 pre}) in
 $ I_{g}^-(\mu_{ g}(s^\prime))\cap U_g$. \MTEXT{Then $\Psi$ is the restriction of the wave map $f$
 solving (\ref{C-problem 1 pre})-(\ref{C-problem 2 pre}).
 Then we can determine in the wave gauge coordinates  $\Psi$ the source
$\Psi_*F$ and the
solution $(\Psi_*g,\Psi_*\phi)$
in $\Psi(  I_{g}^-(\mu_{ g}(s^\prime))\cap U_g)$.
Observe that the function $\Psi_*F$ vanishes on $U_{\hat g}$ outside the set $\Psi(  I_{g}^-(\mu_{ g}(s^\prime))\cap U_g)$.
Thus the source $f_*F$  is determined in the wave gauge coordinates $f$ 
in the whole set $U_{\hat g}$, where $f$ solves (\ref{C-problem 1 pre})-(\ref{C-problem 2 pre}).
This construction can be done for all equivalence
classes $[(U_g,g|_{U_g},\phi|_{U_g},F|_{U_g})]$ in ${\cal D}(\hat g,\hat \phi,\hat \e)$ such that
$F$ 
is supported in 
 $I_{g}( \mu_{ g}(s^\prime-r_2),\mu_{ g}(s^\prime))\subset W_g$. 
 Next, we consider sources on the set $U_{\hat g}$ endowed with
 the background metric $\hat g$.
 Let 
 $\F$ be a source function on the set $U_{\hat g}$ such that $J^-_{\hat g}(\supp(\F))\cap J^+_{\hat g}(\supp(\F))\subset 
I_{\hat g}( \mu_{ \hat g}(s^\prime-r_2),\mu_{\hat g}(s^\prime))$
and assume that  $\F$  is sufficiently small in the $C^{s_1}$-norm.
Then, using the above considerations,
 we can determine
the unique equivalence class $[(U_g,g|_{U_g},\phi|_{U_g},F|_{U_g})]$ in  ${\cal D}(\hat g,\hat \phi,\hat \e)$ for which $f_*F$ is equal to $\F$.} 
In this case we say that the equivalence class corresponds to $\F$ in
the wave gauge coordinates and that $\F$ is an admissible source.

  {\bf Step 2.}
Let   $k_1\geq 8$, $s_1\geq k_1+5$, and $n\in \Z_+$ be large enough. 
Let $(x_j^{\prime},\xi_j^{\prime})$ be covectors in an $s$-neighborhood
of $(x_j,\xi_j)$ and
  $B_j$ be
neighborhoods of  $x_j^{\prime}$ satisfying (\ref{eq: admissible Bj}).
Consider then
  ${\bf f}_j\in \mathcal I_{C}^{n+1}(Y((x_j^{\prime},\xi_j^{\prime});t_0,s_0))$,
  and a family
$\F_j(\e_j)\in  C^{s_1}(M_0)$, $\e_j\in [0,\e_0)$ 
 of functions supported in $B_j$ that depend smoothly
on $\e_j$. Moreover, assume that
 $\p_{\e_j} \F_j(\e_j)|_{\e_j=0}={\bf f}_j$.
 Then, we can use step 1 to test
if all sources ${\bf f}_j$ and $\F_j(\e_j)$, $\e_j\in [0,\e_0)$ are admissible.

  {\bf Step 3.} 
  Next, assume that ${\bf f}_j$ and $\F_j(\e_j)$, $\e_j\in [0,\e_0)$, $j=1,2,3,4$ are admissible
  and are compactly supported in 
neighborhoods $B_j$ of  $x_j^{\prime}$ satisfying (\ref{eq: admissible Bj}).

  \MTEXT{Let us next consider  $\vec a=(a_1,a_2,a_3,a_4)\in (-1,1)^4$ and define
 $\F(\tilde \e,a)=\sum_{j=1}^4 \F_j(a_j\tilde \e)$.}
 Let $(g_{\tilde \e,\vec a},\phi_{\tilde \e,\vec a})$ be the solution of 
  (\ref{eq: adaptive model with no source}) with  the source  $\F(\tilde \e,\vec a)$.
By  (\ref{eq: admissible Bj}), $B_j\cap J ^+_{g_{\tilde \e,\vec a}}(B_k)=\emptyset$ for $j\not=k$ and $\tilde \e$ small enough.
Hence
we have that for sufficiently small $\tilde \e$ the conservation law (\ref{conservation law0})
is satisfied for $(g_{\tilde \e,\vec a},\phi_{\tilde \e,\vec a})$
in the set $Q=(M_0\setminus J^+_{g_{\tilde \e,a}}(p^-))\cup \bigcup_{j=1}^4J^-_{g_{\tilde \e,a}}(B_j)$, see Fig. 6(Left). Solving the
Einstein equations   (\ref{eq: adaptive model with no source})  in
a neighborhood of the closure of 
 the complement of the set $Q$, where  the source  $\F(\tilde \e,a)$ vanishes,
we see that   the conservation law (\ref{conservation law0}) 
is satisfied in the whole set $M_0$, see \cite[Sec.\ III.6.4.1]{ChBook}.
%
%
Hence
  $\sum_{j=1}^4 a_j{\bf f}_j$  has the LS property (\ref{eq: LSp})    in $C^{s_1}(M_0)$ with the family
  $\F(\tilde \e,\vec a)$. When $\tilde \e>0$ is small
  enough, the sources   $\F(\tilde \e,\vec a)$ are admissible.


  {\bf Step 4.}
 Let $v_{\tilde\e,\vec a}=(
g_{\tilde \e,\vec a},\phi_{\tilde \e,\vec a})$ be the solutions of 
the
Einstein equations   (\ref{eq: adaptive model with no source}) with the source $\F(\tilde \e,\vec a)$.

Using Step 1, we can find 
for
all $\vec a\in (-1,1)^4$ and sufficiently small $\tilde \e$   the equivalence classes
$[(U_{\tilde \e,\vec a},g_{\tilde \e,\vec a}|_{U_{\tilde \e,\vec a}},\phi_{\tilde \e,\vec a}|_{U_{\tilde \e,\vec a}},F(\tilde \e,\vec a)|_{U_{\tilde \e,\vec a}})]$
correspond to  $\F(\tilde \e,\vec a)$ in the wave guide coordinates. Then
we can determine, 
using the normal coordinates $\Psi_{\tilde \e,\vec a}$ associated to the observation geodesics $\mu_{\tilde \e\vec\a}$ and 
the solution $v_{\tilde\e,\vec a}$, the function $
(\Psi_{\tilde \e,\vec a})_*v_{\tilde \e,\vec a}$. Also, we can compute the derivatives of this function
with respect to $\tilde\e$ and $a_j$.
%
%

Observe that $\p_{\vec\e }^4f(\e_1,\e_2,\e_3,\e_4)|_{\vec \e=0}=
\p_{\vec a }^4(\p_{\tilde \e }^4f(a_1\tilde \e,a_2\tilde \e,a_3\tilde \e,a_4\tilde \e)|_{\tilde \e=0})|_{\vec a=0}$ so that
 $\M^{(4)}=\p_{a_1}\p_{a_2}\p_{a_3}\p_{a_4}\p^4_{\tilde \e}v_{\tilde \e,\vec a}|_{\tilde \e=0,\vec a=0}$, where  $\M^{(4)}$ given in 
(\ref {measurement eq.}). 
Thus, using the solutions $v_{\tilde\e,a}$ we  determine if the 
function $\p_{\vec\e}^4((\Psi_{\vec \e})_*v_{\vec \e})|_{\vec\e=0}$, corresponding to the source   $f_{\vec \e}=\sum_{j=1}^4 \e_j{\bf f}_j$ 
is singular, where $\Psi_{\vec\e}$ are the
the normal coordinates associated to any observation geodesic.
%
Thus we can verify if the condition (D) holds.
\hfill \Box \medskip

\section{Determination of earliest light observation sets }\label{subsection combining}

Below, we use only the metric $\hat g$ and often denote $\hat g=g$, $U=U_{\hat g}$.  
  
Our next  
aim is to consider the global problem of constructing the set of the earliest light observations of all points $q\in J(p^-,p^+)$. To this end, we need to handle the technical problem that  
in the set $ {\mathcal Y}((\vec x,\vec\xi))$ we  
have not analyzed if we observe singularities or not.  \MTEXT{Also, we have not  
analyzed the waves in the set 
where singularities caused by caustics or their interactions may appear, see (\ref{eq: summary of assumptions 2}).
To avoid these  difficulties,} 
we define next the  
 sets $\S_{reg} ((\vec x,\vec \xi),t_0)$ of points near which we observe   
singularities in a 3-dimensional set.

\begin{definition}\label{def: Sreg}  
Let $(\vec x,\vec \xi)=((x_j,\xi_j))_{j=1}^4$ be a collection of light-like  
vectors with $x_j\in U_{\hat g}$ and $t_0>0$. We define 
 $\S((\vec x,\vec \xi),t_0)$   
be the set of those $y\in U_{\hat g}$    
that satisfies   
the property  ({D}) with $(\vec x,\vec \xi)$ and $t_0$ and some $\hat s>0$.  
Moreover, let   
 $\S_{reg} ((\vec x,\vec \xi),t_0)$ be the set  
of the points  $y_0\in \S((\vec x,\vec \xi),t_0)$ having a neighborhood $W\subset U_{\hat g}$  
such that the intersection $W\cap \S((\vec x,\vec \xi),t_0)$  
is a non-empty smooth 3-dimensional submanifold.  
\MTEXT{We denote   (see (\ref{Earliest element sets}) and Def.\ \ref{def. O_U}) 
\beq  \label{EE1}
& &{\Sclo}(\vec x,\vec \xi),t_0)\hspace{-0.5mm}=\hspace{-0.5mm}\hbox{cl}\,(\S_{reg} ((\vec x,\vec \xi),t_0))\cap U_{\hat g},\\ \label{EE2}
& &{\Scle}((\vec x,\vec \xi),t_0)\hspace{-0.5mm}=
\bigcup
 _{(z,\eta)\in \U_{z_0,\eta_0}}\pointear_{z,\eta}({\Sclo}((\vec x,\vec \xi),t_0)).
\eeq
}\end{definition}

The  data  ${\cal D}(\hat g,\hat \phi,\e)$ determines the set 
${\Scle}((\vec x,\vec \xi),t_0)$.  
Below, we fix $t_0$ to  be ${t_0}=4\kappa_1$, cf.\ Lemma \ref{lem: detect conjugate 0}.  

  \extension{
  
Our goal is to show that ${\Sclo}((\vec x,\vec \xi),t_0)$ coincides with the intersection   
of the light cone $\L^+_{\hat g}(q)$ and $U_{\hat g}$ where $q$ is the intersection  
point of the geodesics corresponding to $(\vec x,\vec \xi)$, see Fig.\ 3(Left).

Let us next motivate the analysis we do below:   
We will consider how to create an artificial point source using interaction  
of distorted plane  waves propagating along light-like geodesics $\gamma_{x_j,\xxi_j}(\R_+)$  
where $(\vec x,\vec \xxi)=((x_j,\xxi_j))_{j=1}^4$ are perturbations  
of a light-like $(y,\zeta)$, $y\in \hat \mu$. 
We will use the fact that for all $q\in J(p^-,p^+)\setminus \hat \mu$ there  
is a light-like geodesic $\gamma_{y,\zeta}([0,t])$ from $y=\hat \mu(f^-_{\hat \mu}(q))$  
to $q$ with $t\leq \rho(y,\zeta)$.  
We will next  
show that   
when we choose $(x_j,\xxi_j)$ to be suitable   
perturbations of $\p_t \gamma_{y,\zeta}(t_0)$, $t_0>0$,  
it is possible that all geodesics $\gamma_{x_j,\xxi_j}(\R_+)$ intersect at $q$ before   
their first cut points,  
that is, $\gamma_{x_j,\xxi_j}(r_j)=q$, $r_j<\rho(x_j,\xxi_j)$.   
We note that we cannot analyze the interaction of the waves  
if the geodesics intersect after the cut points as then  
 the distorted plane   
waves can have caustics. Such interactions of wave caustics can, in principle,  
cause propagating singularities. Thus the sets ${\Sclo}((\vec x,\vec \xi),t_0)$  
contain singularities propagating along the light cone $\L^+_{\hat g}(q)$  
and in addition that they may contain singularities produced by caustics that we do not know how to analyze (that could  
be called ``messy waves''). Fortunately, near an open and dense set of geodesics  
$\mu_{z,\eta}$ the nice singularities propagating along the light cone $\L^+_{\hat g}(q)$  
arrive before the ``messy waves''. This is the reason why we consider below the  
first observed singularities on geodesics  
$\mu_{z,\eta}$.  
 Let us now return to the rigorous analysis.  
 
Below in this section we fix $t_0$ to have the value  
 ${t_0}=4\kappa_1$, cf.\ Lemma \ref{lem: detect conjugate 0}. 
 Recall the notation that  
\ba  
& &(x({t_0}),\xxi({t_0}))=(\gamma_{x,\xxi}({t_0}),\dot \gamma_{x,\xxi}({t_0})),\\  
& &(\vec x({t_0}),\vec \xxi({t_0}))=((x_j({t_0}),\xxi_j({t_0})))_{j=1}^4.  
\ea  
\medskip

\noindent
{\bf Lemma 5.1.A}  
{\it
Let $\vartheta>0$ be arbitrary, $q\in J(p^-,p^+)\setminus {\hat \mu}$ and   
let $y={\hat \mu}(f_{\hat \mu}^-(q))$  
and $\zzeta\in L^+_y\hattuM _0$, $\|\zeta\|_{g^+}=1$ be such that $\gamma_{y,\zzeta}([0,r_1])$, $r_1>t_0=4\kappa_1$ is   
a longest causal (in fact, light-like) geodesic connecting $y$ to $q$.   
Then there exists  a set $\mathcal G$ of 4-tuples   of light-like vectors $(\vec x,\vec \xxi)=((x_j,\xxi_j))_{j=1}^4$   
 such that the points  $x_j$  
 and the directions $\xxi_j$ and the points $x_j(t_0)=\gamma_{x_j,\xxi_j}(t_0)$
  have the following properties:  
\begin{itemize}  
\item[(i)] $x_j(t_0)\in U_{\hat g}$, $x_j(t_0)\not \in J^+(x_k(t_0))$ for $j\not =k$,  
  
\item[(ii)] $d_{g^+}  
((x_l,\xxi_l),(y,\zeta))<\vartheta$ for $l \leq 4$,  
  
\item[(iii)] $q=\gamma_{x_j,\xxi_j}(r_j)$ and $\rho(x_j(t_0),\xxi_j(t_0))+t_0>r_j$,  
  
\item[(iv)] when $(\vec x,\vec \xxi)$ run through the set  $\mathcal G$,  
the directions  
 $ (\dot \gamma_{x_j,\xxi_j}(r_j))_{j=1}^4$ form an open  
 set in $(L^+_qM_0)^4$.  
 
 \end{itemize}  
   
In addition, $\mathcal G$ contains elements $(\vec x,\vec \xxi)$ for which $(x_1,\xxi_1)=  
(y,\zzeta)$.  
}
 \medskip
 
 \noindent 
{\bf Proof.}   
Let  $\eta=\dot \gamma_{y,\zzeta}(r_0)\in L^+_q\hattuM _0$.  
Let us choose light-like directions  $\eta_j\in T_q\hattuM _0$, $j=1,2,3,4,$   
close to $\eta$ so that $\eta_j$ and $\eta_k$ are not parallel for $j\not  =k$.  
In particular, it is possible (but not necessary) that   
$\eta_1=\eta$. Let ${\bf t}:M\to \R$ be a time-function on $M$ that can be  
used to identify $M$ and $\R\times N$. Moreover, let  us choose  
 $T\in ((\gamma_{y,\zeta}(r_0-t_0)),{\bf t}(\gamma_{y,\zeta}(r_0-2t_0)))$ and  for $j=1,2,3,4$, let $s_j>0$   
be such that   
${\bf t}(\gamma_{q,\eta_j}(-s_j))=T$. Choosing  first $T$ to be sufficiently close ${\bf t}(\gamma_{y,\zeta}(r_0-t_0))$  
and then all $\eta_j$, $j=1,2,3,4$ to be sufficiently  
close to $\eta$     
 and defining  
$x_j=\gamma_{q,\eta_j}(-s_j-t_0)$ and $\xxi_j=\dot\gamma_{q,\eta_j}(-s_j-t_0)$  
we obtain the pairs $(x_j,\xxi_j)$ satisfying the properties stated in the claim.  
Indeed, this follows from the fact that $\rho(x,\xi)$ is lower semicontinuous function,
$\rho(q,\eta)\geq r_0$
and the geodesics $\gamma_{q,\eta_j}$ can not intersect before their first
cut points.
As vectors $\eta_j$ can be varied in sufficiently small open sets   
so that the properties stated in the claim stay valid, we obtain   
the  claim  concerning the open set of light-like directions.  
  
The last claim follows from the fact that $\eta_1$ may be equal to $\eta$ and $T={\bf t}(y)$.  
 \hfill \Box \medskip

Next we  analyze the   
set ${\Sclo}((\vec x,\vec \xi),t_0)$.

  }

 If the set 
 $\cap_{j=1}^4 \gamma_{x_j,\xi_j}([t_0,\infty))$ is non-empty we denote  its earliest point by  
$Q((\vec x,\vec \xi),t_0)$.  
If such intersection point does not exists, we define $Q((\vec x,\vec \xi),t_0)$  
to be the empty set.  \MTEXT{Next we consider the relation of ${\Scle}((\vec x,\vec \xi),t_0)$ and $\be_U(q)$, $q=Q((\vec x,\vec \xi),t_0)$, see
Def.\ \ref{def. O_U}.}

\begin{lemma}\label{lem: sing detection in in normal coordinates 3}  
Let $(\vec x,\vec \xi)$, $j=1,2,3,4$ and  $t_0>0$ 
satisfy (\ref{eq: summary of assumptions 1})-(\ref{eq: summary of assumptions 2})  
and assume that $\vartheta_1$ in  (\ref{eq: summary of assumptions 1}) and   
Lemma \ref{lem: detect conjugate 0}   is so small  
that  for all $j\leq 4$,   
$x_j\in  I( \hat \mu(s_1),\hat \mu(s_2))\subset W_{\hat g}$ with some $s_1,s_2\in (s_-,s_+)$.

Let $\V= \V((\vec x,\vec \xi),t_0)$ be the set defined in  (\ref{eq: summary of assumptions 2}).
Then 
 \smallskip  
  
{(i)  
Assume that $y\in \V\cap U_{\hat g}$ satisfies the condition (I) with $(\vec x,\vec \xi)$ and $t_0$ and parameters   
$q,\zeta$, and $t$ such that $0\leq t\leq \rho(q,\zeta)$. 
 Then   
$y\in  {\Sclo}((\vec x,\vec \xi),t_0)$. 
 \smallskip   
  
(ii) Assume  $y\in \V\cap U_{\hat g}$ does not satisfy condition (I) with $(\vec x,\vec \xi)$ and $t_0$.  
Then $y\not \in {\Sclo}((\vec x,\vec \xi),t_0)$.   
\smallskip   
  
(iii) 
If $Q((\vec x,\vec \xi),t_0)\not =\emptyset$ and $q=Q((\vec x,\vec \xi),t_0)\in   \V$,
we have \ba
\quad{\Scle}((\vec x,\vec \xi),t_0)=\be_U(q)\subset \V.  
\ea
Otherwise, if $Q((\vec x,\vec \xi),t_0)\cap  \V=\emptyset$,
then ${\Scle}((\vec x,\vec \xi),t_0)\cap  \V=\emptyset$.
}
   \end{lemma}

   \noindent  
   {\bf Proof.}  (i) Assume first that   $y$ is not in ${\mathcal Y}((\vec x,\vec\xi))$
   and $t<\rho(q,\zeta)$. Then the
   assumptions in (i) and Lemma
   \ref{lem: sing detection in in normal coordinates 2} (ii) and  (iii)
   imply that $y\in {\Sclo}((\vec x,\vec \xi),t_0)$. 
   Consider next a general point $y$ satisfying  the
   assumptions in (i) and let $q=Q((\vec x,\vec \xi),t_0).$ 
   Then $y\in \be_U(q)$.
      Since $\rho(x,\xi)$ is
   lower semi-continuous and 
   the set $\be_U(q)\setminus {\mathcal Y}((\vec x,\vec\xi))$ is  dense in $\be_U(q)$,
   and  $y$ is a limit point of points $y_n\not \in {\mathcal Y}((\vec x,\vec\xi))$ 
   that satisfy the assumptions in (i). Hence $y\in  {\Sclo}((\vec x,\vec \xi),t_0)$.

   The claim (ii) follows from  Lemma
   \ref{lem: sing detection in in normal coordinates 2} (i).

  (iii)   
Suppose
    $q=Q((\vec x,\vec \xi),t_0)\in  \V$ and
     $y\in \be_U(q)\setminus
  \bigcup_{j=1}^4\gamma_{x_j,\xi_j}([t_0,\infty))$. Let $\gamma_{q,\eta}([0,l])$ be a light-like
  geodesic that is the longest causal geodesic from $q$ to $y$ {with $l\leq \rho(q,\eta)$}, and let 
  $p_j=\gamma_{x_j,\xi_j}(t_0+{\bf t}_j)$, ${\bf t}_j=\rho(x_j(t_0),\xi_j(t_0))$, be the first cut point on the geodesic $\gamma_{x_j,\xi_j}([t_0,\infty))$.
 To show that 
   $y$ is in $\V$, we assume the opposite, $y\not \in \V$. Then  for some $j$ there is
  a causal geodesic $\gamma_{p_j,\theta_j}([0,l_j])$ from $p_j$ to $y$.  Now we can use
  a short-cut argument: Let $q=\gamma_{x_j,\xi_j}(t^\prime)$. 
  As  $q\in \V$, we have $t^\prime<t_0+{\bf t}_j$. 
  Moreover, as 
  $y\not \in \gamma_{x_j,\xi_j}([t_0,\infty))$, the union of the geodesic
  $\gamma_{x_j,\xi_j}([t^\prime,t_0+ {\bf t}_j])$ from $q$ to $p_j$ and 
  $\gamma_{p_j,\theta_j}([0,l_j])$ from $p_j$ to $y$ does not form a light-like geodesic
  and thus $\tau(q,y)>0$. As     $y\in \be_U(q)$, this is not possible. Hence
  $y\in \V$.   Thus by (i), $y\in  {\Sclo}((\vec x,\vec \xi),t_0)$
  and hence $\be_U(q)\setminus ( \bigcup_{j=1}^4\gamma_{x_j,\xi_j}([t_0,\infty)))\subset {\Sclo}((\vec x,\vec \xi),t_0)$. 
  Since  the set $\be_U(q)\setminus (
  \bigcup_{j=1}^4\gamma_{x_j,\xi_j}([t_0,\infty)))$ is dense in the closed set
  $\be_U(q)$, the above shows that  $\be_U(q)\subset {\Sclo}((\vec x,\vec \xi),t_0)$.
  Also
  by (ii),  $ {\Sclo}((\vec x,\vec \xi),t_0)\subset \L^+(q).$ 
 {Using Definition  \ref{def. O_U} and  (\ref{EE2}), we see
 that $ {\Scle}((\vec x,\vec \xi),t_0)= \be_U(q).$}

  On the other hand, if $Q((\vec x,\vec \xi),t_0)\cap \V =\emptyset$, we can apply
  (ii) for all $y\in \V\cap U$ and see that  $ {\Sclo}((\vec x,\vec \xi),t_0)\cap \V=\emptyset$.
  This proves (iii). 
\hfill \Box \medskip

   \MTEXT{
 Let 
 \beq\label{set Kt0}
 \K_{t_0}=\{x \in U_{\hat g}&;&
x= \gamma_{{\hat x},\xi}(r),\
{\hat x}=\hat \mu(s), \ s\in [s^-,s^+],\\
\nonumber
& & 
 \xi\in L^+_{\hat x}M_0,\ \|\xi\|_{\hat g^+}=1,\ r\in [0,2{t_0})\}.
\eeq}
\extension{ 
Recall that $\U_{z_0,\eta_0}=\U_{z_0,\eta_0}(\hat h)$  
 was defined using the parameter $\hat h$. We see 
 using compactness of $\hat \mu ([s_-,s_+])$, the continuity of $\tau(x,y)$, and 
the existence of convex neighborhoods \cite[Prop.\ 5.7]{ONeill},
  c.f.\ Lemma 2.A.1,
 that if $t_0^\prime >0$  
 and  $\hat h^\prime\in (0, \hat h)$ are small enough and $(z,\eta)\in \U_{z_0,\eta_0}(\hat h^\prime)$,  
 then the  
 longest geodesic from $x\in \overline \K_{t_0^\prime}$ to the point  
 $\be_{z,\eta}(x)$ is contained in $\U_g$
 and hence we can determine the point $\be_{z,\eta}(x)$ for such $x$ and $(z,\eta)$.  
 Let us replace the parameters $\hat h$ and $t_0$ by $\hat h^\prime$  
 and $t_0^\prime$, correspondingly in our considerations below. Then we  
 may  assume  
 that in addition to the data given in the original formulation of the problem,  
 we are given also the set $\be_U(\K_{t_0})$. Next we do this.
}
\noextension{\MTEXT{Recall that we are given the set $U_{\hat g}$ that is determined via the parameter $\hat h>0$.   
By using compactness of $\hat \mu ([s_-,s_+])$, the continuity of $\tau(x,y)$, and 
the existence of convex neighborhoods \cite[Prop.\ 5.7]{ONeill},
we can determine the sets $\be_U(\K_{t_0^\prime}\cap J^+(\hat \mu(s)))$ for some $t_0^\prime>0$ and 
all $s\in [s_-,s_+]$. Here, recall that  $\be_U(V)=\{\be_U(q);\ q\in V\}\subset 2^{U}$.
 Thus,
 by making $\hat h$ and $t_0=4\kappa_1$ smaller, we may assume  below that we are given the 
sets $\be_U(\K_{t_0}\cap J^+(\hat \mu(s)))$ for $s\in [s_-,s_+]$.}}
%
Below, we may assume that $\vartheta_1$ is so \MTEXT{small that
$\gamma_{y ,\zeta}([0,t_0]) \cap J(p^-,p^+)\subset \K_{t_0}$ when 
$y\in J(p^-,p^+)$, $d_{g^+}(y,\hat \mu)<\vartheta_1$ and
$\zeta\in L^+_yM$, $\|\zeta\|_{g^+}\leq 1+\vartheta_1$.}

Let $\kappa_1,\kappa_2$ be constants  given as in Lemma \ref{lem: detect conjugate 0}.
{Let  $s_0\in [s_-,s_+]$   
be so close to $s_+$  
 that $J^+(\hat \mu(s_0))\cap J^-({p^+})\subset \K_{t_0}$. Then the given data ${\cal D}(\hat g,\hat \phi,\e)$ 
  determines $\be_{U}(J^+(\hat \mu(s_0))\cap J^-({p^+}))$.  

Next we use \MTEXT{a step-by-step construction}: We    
consider  
 $s_1\in (s_-,s_+)$ and assume that  we are given  $\be_{U}(J^+({\hat x}_1)\cap J^-({p^+}))$ with ${\hat x}_1={\hat \mu}(s_1)$.} Then, let  $s_2\in (s_1-\kappa_2,s_1)$.  
Our next aim is to find the earliest light observation sets $\be_{U}(J^+({\hat x}_2)\cap J^-({p^+}))$ with ${\hat x}_2={\hat \mu}(s_2)$. To this end   
we need to make the following definitions  (see Fig.\ 6). 
\medskip

\begin{center}  
  
 \psfrag{1}{\hspace{-2mm}$J^+(B_1)$}  
\psfrag{2}{$J^+(B_2)$}  
\psfrag{3}{$B_1$}  
\psfrag{4}{$B_2$}  
\psfrag{5}{$p^-$}

  \includegraphics[width=5.5cm]{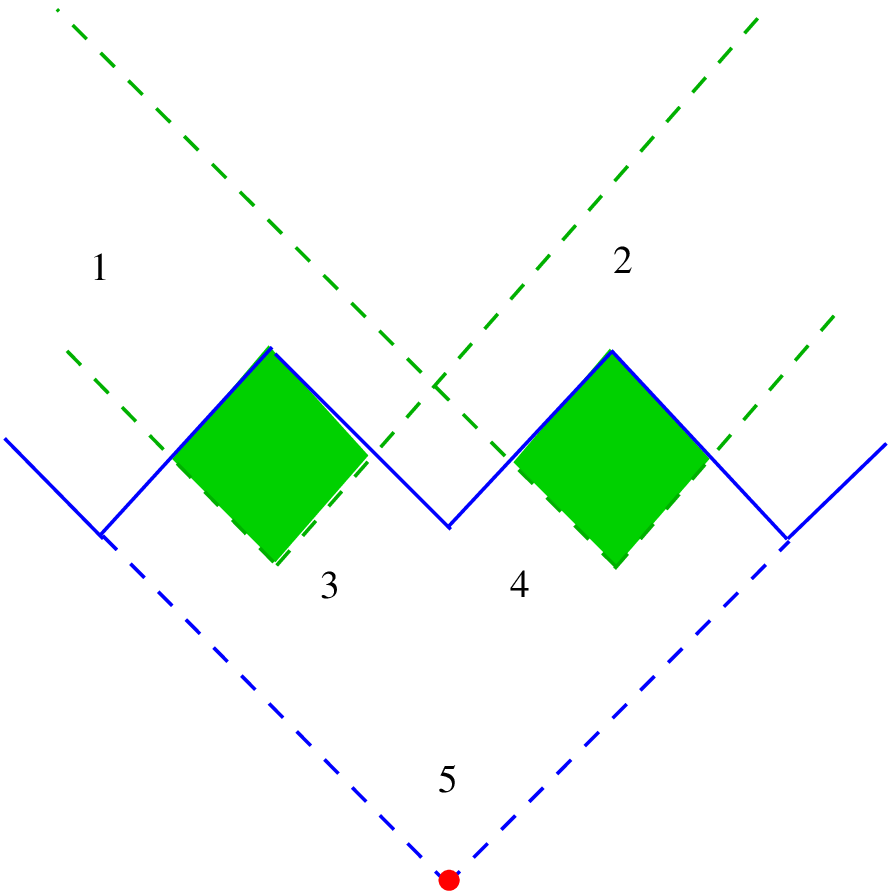}  
\psfrag{1}{$\hspace{-2mm}{\hat x}_2$}  
\psfrag{3}{\hspace{-1.5mm}$p_2$}  
\psfrag{2}{$p_2^\prime$}  
\psfrag{4}{\hspace{-2mm}${\hat x}$}  
\psfrag{5}{y}  
\psfrag{6}{${\hat x}_1$}  
\psfrag{7}{$z$}  
\psfrag{8}{$y^\prime$}  
\psfrag{9}{$\hspace{-2mm}p^-$}  
\psfrag{0}{$p^+$}  
\includegraphics[width=5.5cm]{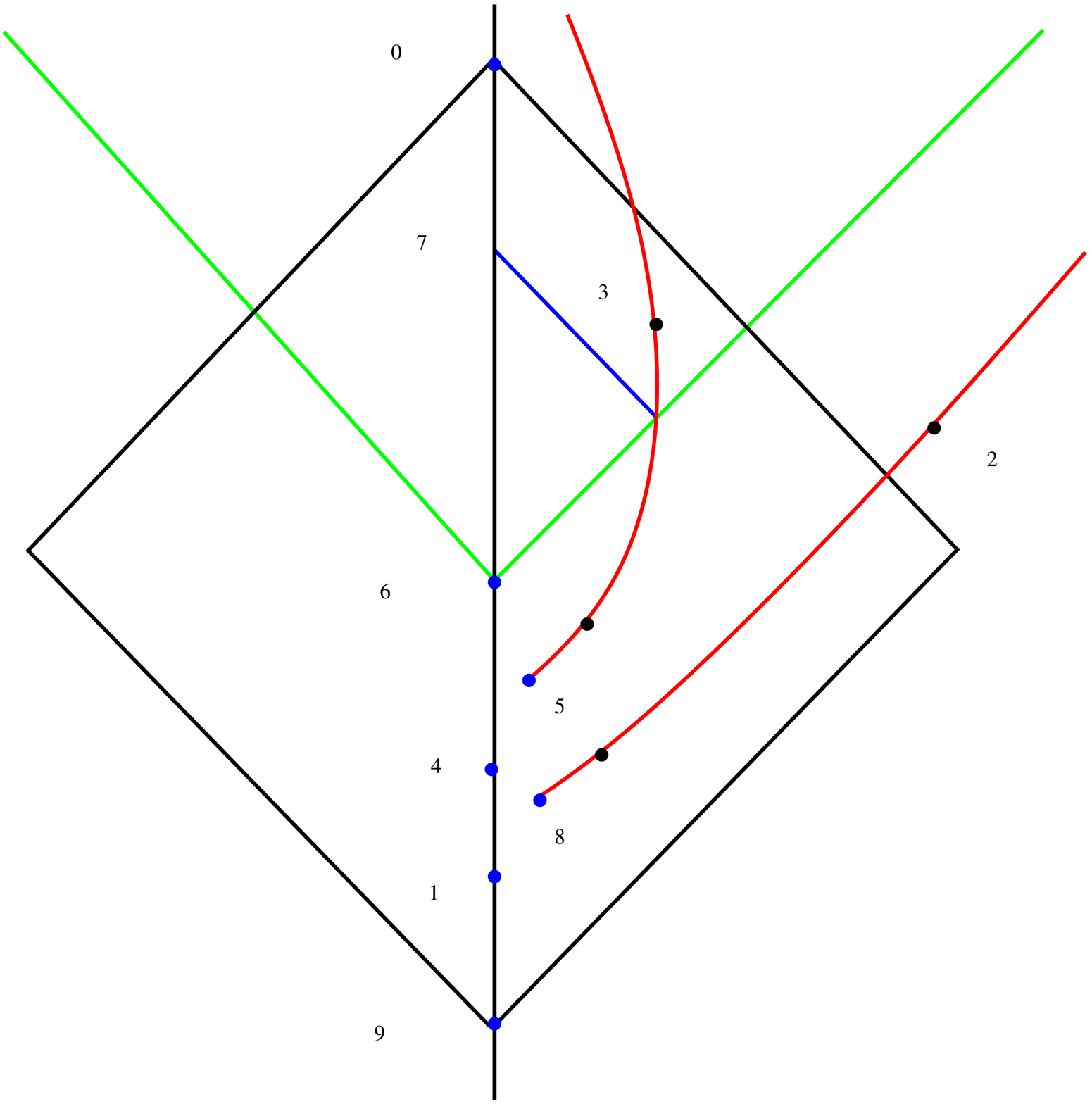}  
\end{center}  
{\it FIGURE 6. {\bf Left:} The setting in the proof of Lemma \ref{lem: sing detection in in normal coordinates 2}.
{\bf   :} The blue points on $\hat \mu$ are ${\hat x}_1=\hat \mu(s_1)$, ${\hat x}_2=\hat \mu(s_2)$, and ${\hat x}=\hat \mu(s)$.
The blue points $y$ and $y^\prime$ are close to $\hat x$.  The set with the green boundary is $J^+({\hat x}_1)$.
 We consider the 
 geodesics $\gamma_{y,\zeta}([0,\infty))$ and $\gamma_{y^\prime,\zeta^\prime}([0,\infty))$.   
 These geodesics  corresponding to the cases when
the geodesic  $\gamma_{y,\zeta}([0,\infty))$ enters in $J^-({p^+})\cap J^+({\hat x}_1)$,
 and the case
when the geodesic  $\gamma_{y^\prime ,\zeta^\prime}([0,\infty))$ does not
enter this set. The point $p_2$ is the cut point of  $\gamma_{y,\zeta}([t_0,\infty))$
and  $p_2^\prime$ is the cut point of  $\gamma_{y^\prime ,\zeta^\prime}([t_0,\infty))$.
At the point $z=\hat \mu(\mathbb{S} (y,\zeta,s_1))$  we observe for the first time  on the geodesic $\hat \mu$  that  
the geodesic $\gamma_{y,\zeta}([0,\infty))$ has entered   
  $J^+({\hat x}_1)$.  
}  

\extension{

\vfill

\goodbreak 
\begin{center}  
  
\psfrag{1}{$\tilde y^\prime$}  
\psfrag{2}{$p_2$}  
\psfrag{3}{$J_{\hat g}(\hat p^-,\hat p^+)$}  
\psfrag{4}{$\tilde y$}  
\psfrag{5}{}  
\psfrag{6}{$y_1$}  
\psfrag{7}{}  
\psfrag{8}{}
\psfrag{A}{\hspace{-3mm}$\mu_{z,\eta}$}  
\psfrag{B}{$\mu_{z_0,\eta_0}$}  
\includegraphics[width=7cm]{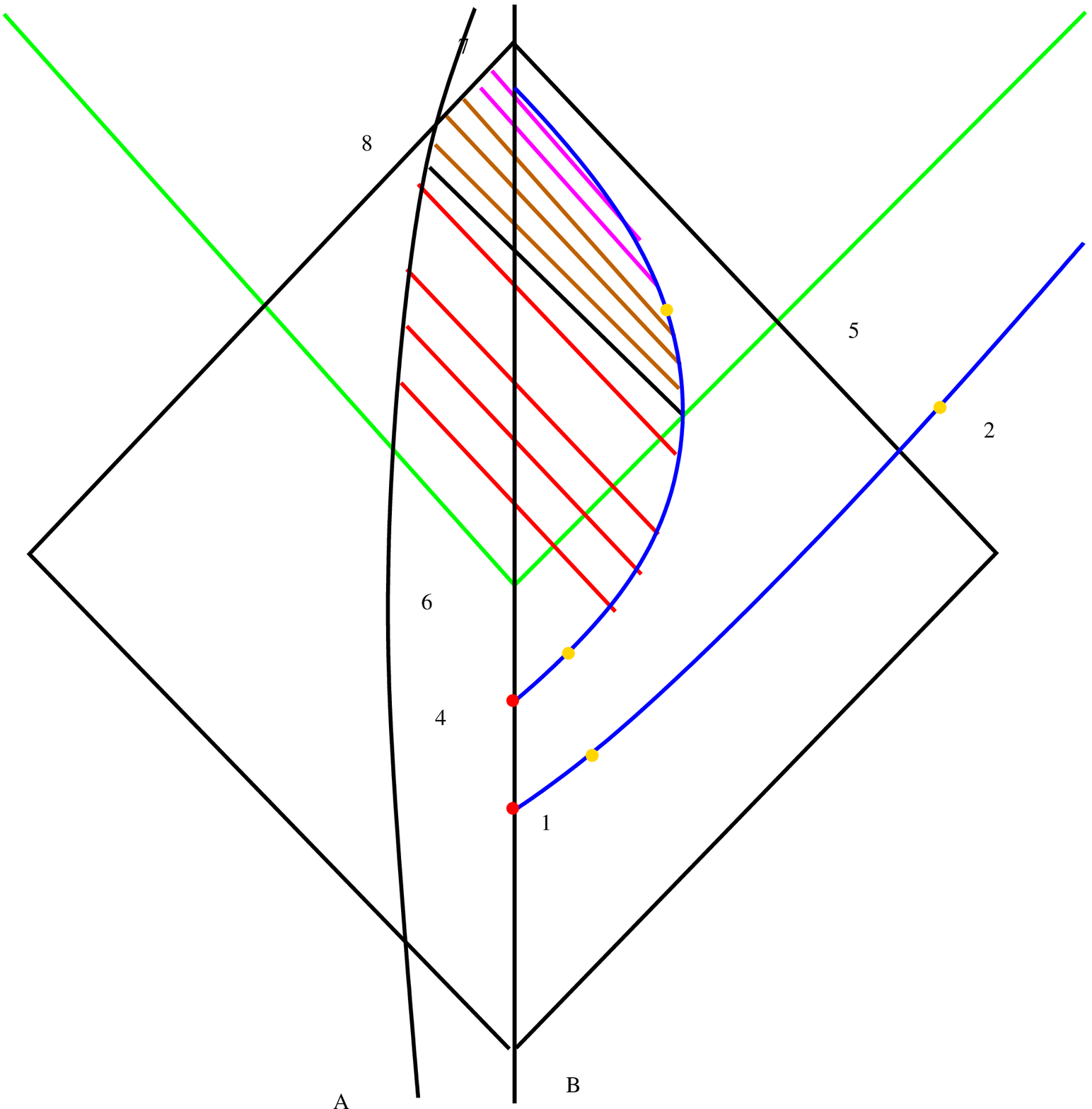}  
\end{center}  
  {\it FIGURE A8: A schematic figure where the space-time is represented as the  2-dimensional set $\R^{1+1}$.   
We consider  
the geodesics emanating from a point $x(r)=\gamma_{\tilde y,\tilde \zeta}(r)$,  
When $r$ is smallest value for which $x(r)\in J^+(y_1)$, a light-like geodesic  (black line segment)  
emanating from $x(r)$ is observed at the point $\tilde p\in \hat \mu $.  
Then $\tilde p= \hat\mu(\mathbb{S} (\tilde y,\tilde \zeta,s_1))$.  
When $r$ is small enough so that $x(r)\not \in J^+(y_1)$,  
the light-geodesics (red line segments) can be observed on $\hat \mu$ in the set  $J^-(\tilde p)$.  
Moreover, when  $r$  is such that $x(r)\in J^+(y_1)$, the light-geodesics can be observed at $\hat \mu$ in the set  $J^+(\tilde p)$. The golden point is the cut point  on $\gamma_{\tilde y,\tilde \zeta}([t_0,\infty))$ and the singularities on the light-like geodesics   
starting before this point (brown line segments) can be   
analyzed, but after the cut point the singularities on the light-like geodesics (magenta line segments)  
are not analyzed in this paper.  
}

\vfill

\goodbreak 

\psfrag{1}{$x_1$}  
\psfrag{2}{$x_2$}  
\psfrag{3}{$q$}  
\psfrag{4}{$p$}  
\psfrag{5}{}  
\psfrag{6}{$y_1$}  
\psfrag{7}{$x$}  
\psfrag{8}{$z$}  
\begin{center}  
\includegraphics[width=6.5cm]{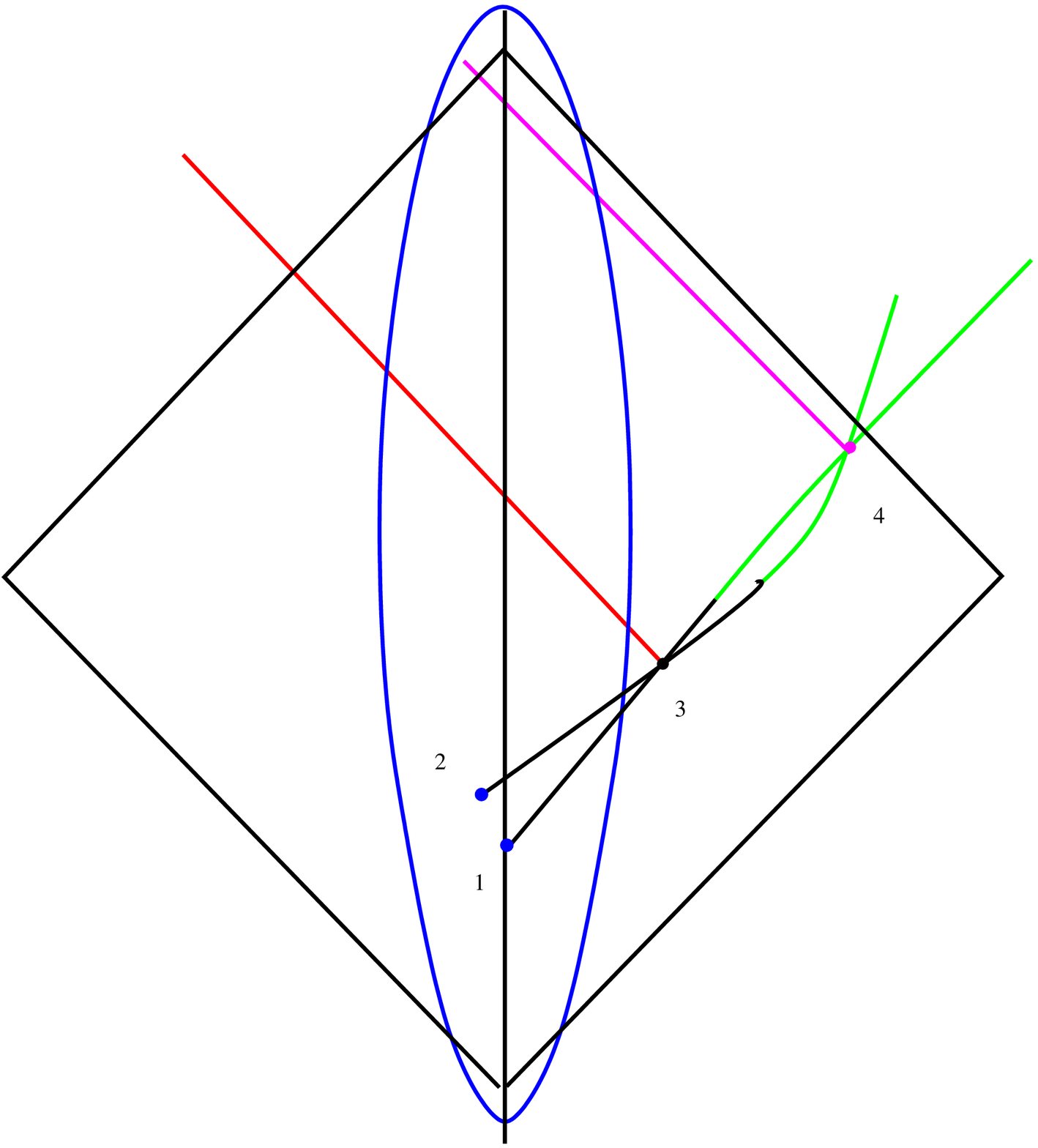}  
\end{center}  
{\it FIGURE A9:  
A schematic figure where the space-time is represented as the  2-dimensional set $\R^{1+1}$.   
In section 5 we consider geodesics $\gamma_{x_j,\xi_j}([0,\infty))$, $j=1,2,3,4,$ that all intersect for the first time at a point  
$q$. 
When we consider geodesics $\gamma_{x_j,\xi_j}([t_0,\infty))$, $j=1,2$ with $t_0>0$, they may  
have cut points at $x^{cut}_j=\gamma_{x_j,\xi_j}({\bf t}_j)$.  
In the figure $\gamma_{x_j,\xi_j}([0,{\bf t}_j))$ are colored by black and the geodesics   
 $\gamma_{x_j,\xi_j}([{\bf t}_j,\infty))$ are colored by green. In the figure  
 the geodesics intersect at the point $q$ before the cut point and for  
the second time at $p$  
 after the cut point. We can analyze the singularities caused by the distorted plane  waves   
 that interact at $q$ but not the interaction of waves after the cut points of the geodesics.  
 It may be that e.g. the intersection at $p$ causes new singularities to appear and we observe  
 those in $U_{\hat g}$. As we cannot analyze these singularities, we consider these singularities  
 as ``messy waves''. However, the ``nice'' singularities caused by the interaction at $q$ propagate along  
 the future light cone of the point $q$ and in
$U_{\hat g}\setminus \cup_{j=1}^4\gamma_{x_j,\xi_j} $ these ``nice'' singularities are observed before the ``messy waves''. Due to this the first singularities we observe near $\hat \mu$ 
come from the point $q$.  
}  
}

\begin{definition}\label{def: hat s 1}  
Let $s_-\leq s_2\leq s<s_1\leq s_+$ satisfy  $s_1<s_2+\kappa_2$, 
 ${\hat x}_j={\hat \mu}(s_j)$, $j=1,2$, and  ${\hat x}={\hat \mu}(s)$,
 ${\hat \zeta} \in L^+_{{\hat x}}U$, 
    $\|{\hat \zeta}\|_{g^+}=1$.   
   \MTEXT{Let $(y,\zeta) \in L^+U$ be in $\vartheta_1$-neighborhood of
    $({\hat x},{\hat \zeta})$ such that
     $y\in J^+({\hat x}_2)$ and 
  the geodesic $\gamma_{y,\zeta}(\R_+)$ does not intersect $\hat \mu$.  

Let $\rkaksi(y,\zeta,s_1)=\inf\{r>0;\ \gamma_{y,\zeta}(r)\in J^+(\hat\mu (s_1))\}$.
 Also, define $\ryksi(y,\zeta)=\inf\{r>0;\ \gamma_{y,\zeta}(r)\in M\setminus I^-(\hat\mu(s_{+2}))\}$
  and $r_0(y,\zeta,s_1)=\min(\ryksi(y,\zeta),\rkaksi(y,\zeta,s_1))$.}
 
When  
 $\gamma_{y ,\zeta }(\R_+)$ intersects $J^+(\hat \mu(s_1))\cap J^-({p^+})$  
we define  
\beq\label{special S formula}  
\mathbb{S} (y ,\zeta ,s_1)=f^+_{\hat \mu}(q_0),
\eeq  
where $ q_0=\gamma_{y ,\zeta }(r_0)$ and   
$r_0=r_0(y,\zeta,s_1).$
In the case when  
$\gamma_{y ,\zeta }(\R_+)$ does not intersect $J^+(\hat \mu(s_1))\cap  {J}^-({p^+})$,  
we define $\mathbb{S} (y ,\zeta ,s_1)=s^+$.

\end{definition}  
  
{We note that above 
  $\ryksi(y,\zeta)$ is finite by \cite[Lem.\ 14.13]{ONeill}.}
  Below we use the notations used in Def.\ \ref{def: hat s 1}.
We saw in Sec.\ \ref{sec: normal coordinates} that  using solutions of the linearized Einstein equations
  we can find $\gamma_{x,\xi}\cap U_{\hat g}$ for any $(x,\xi)\in L^+U$. Thus we 
  can check for given $(x,\xi)$ if $\gamma_{x,\xi}\cap\hat \mu=\emptyset$.
%
%
%

\begin{definition}\label{def: hat s}  
Let $0<\vartheta<\vartheta_1$ and 
$\qP_{\vartheta} (y ,\zeta )$  be the set of $(\vec x,\vec \xxi)=((x_j,\xxi_j))_{j=1}^4$ 
that satisfy (i) and (ii) in formula (\ref{eq: summary of assumptions 1}) with $\vartheta_1$ replaced by $\vartheta$
{and $(x_1,\xi_1)=(y,\zeta)$.}
  We say that the set $S\subset U_{\hat g}$ is a genuine observation
 associated to the geodesic $\gamma_{y,\zeta}$ if  there is $\hat \vartheta>0$
 such that for all $\vartheta\in (0,\hat \vartheta)$
  there are  $(\vec x,\vec \xxi)\in \qP_{\vartheta} (y ,\zeta )$ such that
  $S={\Scle}((\vec x,\vec \xi),t_0)$.
\end{definition}  
  

%

{\begin{lemma}\label{conjugatepoints are nice}   \MTEXT{Suppose  $\max(s_-,s_1-\kappa_2)\leq s<s_1<s^+$
and let
${\hat x}=\hat \mu(s)$, ${\hat x}_1=\hat \mu(s_1)$, and ${\hat \zeta}\in L^+_{{\hat x}}M$, $\|{\hat \zeta}\|_{g^+}=1$.
Moreover, let
  $(y,\zeta)$ be in  a $\vartheta_1$-neighborhood of $({\hat x},{\hat \zeta})$. Assume that the geodesic $\gamma_{y,\zeta}(\R_+)$ does not intersect $\hat \mu$.}

 \noindent (A) Then the cut point  $p_0=\gamma_{y(t_0),\zeta(t_0)}({\bf t}_*)$, ${\bf t}_*= \rho(y(t_0),\zeta(t_0))$ of   
 the geodesics $ \gamma_{y(t_0),\zeta(t_0)}([0,\infty))$, if it exists, satisfies either  

  \smallskip  
   
(i)  $p_0\not \in J^-(\hat \mu(s_{+2}))$,  
 \smallskip  
  
\noindent

or  
 \smallskip

(ii)  $r_0=r_0(y,\zeta,s_1)<\ryksi(y,\zeta)$ 
and $p_0\in I^+({\hat x}_1)$.
 \smallskip  
 
 \noindent
(B)
 There is $\vartheta_2(y,\zeta,s_1)\in (0,\vartheta_1)$   
   such that if $0<\vartheta<\vartheta_2(y,\zeta,s_1)$, 
 $(\vec x,\vec \xi)\in \qP_\vartheta(y,\zeta)$,
  \MTEXT{ and the geodesics $ \gamma_{x_j(t_0),\xi_j(t_0)}([0,\infty))$,
   $j\in \{1,2,3,4\}$,
 has a cut point   
 $p_j=\gamma_{x_j(t_0),\xi_j(t_0)}({\bf t}_j)$,
 then the following holds: \smallskip  
   
If either the point $p_0$ does not exist or it exists and (i) holds
then 
  $p_j\not \in J^-({p^+})$.
On the other hand, if $p_0$ exists and
 (ii) holds, then $f^+_{\hat \mu}(p_j)
>f^+_{\hat \mu}(q_0)$, where $q_0=\gamma_{y,\zeta}(r_0(y,\zeta,s_1))$.}

Note that
$f^+_{\hat \mu}(q_0)=\mathbb{S} (y ,\zeta ,s_1)$.

  \end{lemma}

\noindent  {\bf Proof.}  
(A) Assume that (i) does not hold, that is, $p_0
=\gamma_{y,\zeta}(t_0+{\bf t}_*) \in J^-(\hat \mu(s_{+2}))$. 
By  Lemma \ref{lem: detect conjugate 0}  (ii) we have
$f^-_{\hat \mu}(p_0)>s+2\kappa_2\geq s_1$
that yields $p_0
\in I^+({\hat x}_1)$.
Thus, the geodesic
$\gamma_{y(t_0),\zeta(t_0)}([0,\rho(y(t_0),\zeta(t_0)))$ intersects
 $J^+({\hat x}_1)\cap J^-(\hat \mu(s_{+2}))$. Hence the alternative (ii)
 holds {with $0<r_0
<\ryksi(y,\zeta)$ and moreover, $r_0< t_0+\rho(y(t_0),\zeta(t_0))$.}

(B) If (i) holds, the claim follows since the function $(x,\xi)\mapsto \rho(x,\xi)$ is lower semi-continuous and $(x,\xi,t)\mapsto \gamma_{x,\xi}(t)$ is continuous.  

In the case (ii), we saw above that $r_0<t_0+\rho(y(t_0),\zeta(t_0))$.
Let  $q_0=\gamma_{y,\zeta}(r_0)$. Then by using a short cut argument and the 
fact that $\gamma_{y,\zeta}(\R_+)$ does not intersect $\hat \mu$
we see similarly to the above that $f^+_{\hat \mu}(p_0)>f^+_{\hat \mu}(q_0)=\mathbb{S} (y ,\zeta ,s_1)$.
Since the function $(x,\xi,t)\mapsto f^+_{\hat \mu}(\gamma_{x,\xi}(t))$ is continuous
and $t\mapsto f^+_{\hat \mu}(\gamma_{x,\xi}(t))$ is non-decreasing, and  the function $(x,\xi)\mapsto \rho(x,\xi)$ is lower semi-continuous, 
 we have that the function $(x,\xi)\mapsto f^+_{\hat \mu}(\gamma_{x(t_0),\xi(t_0)}(\rho(x(t_0),\xi(t_0))))$
 is lower semi-continuous, and
 the claim follows.
  \hfill \Box \medskip  
}
  
   \begin{definition}  \label{def 56}
\MTEXT{Let $s_-\leq s_2\leq s<s_1\leq s_+$ satisfy  $s_1<s_2+\kappa_2$, 
 ${\hat x}_j={\hat \mu}(s_j)$, $j=1,2$, and  ${\hat x}={\hat \mu}(s)$,
 ${\hat \zeta} \in L^+_{{\hat x}}U$, 
    $\|{\hat \zeta}\|_{g^+}=1$.   
 Also, let $(y,\zeta) \in L^+U$ be in $\vartheta_1$-neighborhood of
    $({\hat x},{\hat \zeta})$ and $\mathcal G(y ,\zeta ,s_1)$  be the set
    of the
     genuine observations  $S\subset U_{\hat g}$
associated to the geodesic $\gamma_{y,\zeta}$ 
such that $S\in \be_{U}(J^+({\hat x}_1)\cap J^-({p^+}))$.
Moreover, define
  $T (y ,\zeta ,s_1)$ to be the infimum of $s^\prime\in [-1, s_+]$ 
such that $ \hat \mu(s^\prime)\in S\cap \hat \mu$ with
some  $S\in \mathcal G(y ,\zeta ,s_1)$.
If no such $s^\prime$ exists,
we define $T (y ,\zeta ,s_1)=s^+$.}
  \end{definition}  

Let us next  consider  $(\vec x,\vec\xi)\in \qP_\vartheta(y,\zeta)$   
where
$0<\vartheta<\vartheta_2(y,\zeta,s_1)$. Here, 
$\vartheta_2(y,\zeta,s_1)$ is defined in  
Lemma \ref{conjugatepoints are nice}.    
  Assume that for some $j=1,2,3,4$ we have that $\rho(x_j({t_0}),\xxi_j({t_0}))<\mathcal T(x_j({t_0}),\xxi_j({t_0}))$ and consider the cut point  
$p_j=\gamma_{x_j({t_0}),\xxi_j({t_0})}(\rho(x_j({t_0}),\xxi_j({t_0})))$.  
Then either the case (i) or (ii) of Lemma \ref{conjugatepoints are nice}, (B) holds.
If (i) holds, $p_j$  
 satisfies  
$p_j\not\in J^-({p^+})$ and thus $f^+_{\hat \mu}(p_j)>s_+\geq \mathbb{S} (y ,\zeta,s_1)$.
 If (ii) holds, there exists $r_0=r_0(y,\zeta,s_1)<\rkaksi(y,\zeta)$ such that  
  $q_0=\gamma_{y,\zeta}(r_0)\in J^+({\hat x}_1)$ and 
$f^+_{\hat \mu}(p_j)
  >f^+_{\hat \mu}(q_0)=\mathbb{S} (y ,\zeta ,s_1).$ Thus both in case (i) and (ii) we have 
 \beq\label{eq: cut points observed after S}
  f^+_{\hat \mu}(p_j)>  \mathbb{S} (y ,\zeta ,s_1).
 \eeq

Next, consider  a point 
$q=\gamma_{y,\zeta}(r)\in J^-({p^+}),$
where $t_0<r\leq r_0=r_0(y,\zeta,s_1)$.
Let $\theta=-\dot\gamma_{y,\zeta}(r)\in T_qM$.
By Lemma \ref{lem: detect conjugate 0} (ii), the geodesic
$\gamma_{y,\zeta}([t_0,r])=\gamma_{q,\theta}([0,r-t_0])$ has no
cut points. 
Denote  $\tilde x=\hat \mu(\mathbb{S} (y ,\zeta,s_1))$.
Note that then $q\in J^-(\tilde x)$.
Then, consider four  geodesics that emanate  from $q$ to the past, 
in the light-like direction $\eta_1=\theta$ and in the light-like directions
 $\eta_j\in  T_qM$,  $j=2,3,4$ that are sufficiently close to the
direction $\theta$. Let  
$\gamma_{q,\eta_j}(r_j)$ be the intersection points  of $\gamma_{q,\eta_j}$ with
the surface $\{x\in M;\ {\bf t}(x)=c_0\}$, on which the time function  ${\bf t}(x)$ 
has the constant value 
$c_0:={\bf t}(\gamma_{y,\zeta}(t_0))$. 
\MTEXT{Note that such $r_j=r_j(\eta_j)$ exists by the inverse function theorem 
when $\eta_j$ is sufficiently close to $\eta_1$.} 
Choosing  $x_j=\gamma_{q,\eta_j}(r_j+t_0)$ and
$\xi_j=-\dot\gamma_{q,\eta_j}(r_j+t_0)$ we see that
when  $\vartheta_3\in (0,\vartheta_2(y,\zeta,s_1))$ is small enough,
 for all $\vartheta\in (0,\vartheta_3)$, there is
$(\vec x,\vec\xi)=((x_j,\xi_j))_{j=1}^4\in \qP_\vartheta(y,\zeta)$   
such
that the geodesics corresponding to $(\vec x,\vec\xi)$ intersect
at $q$. As  the set $(U_{\hat g},\hat g)$  is known, 
that for sufficiently small $\vartheta$
one can verify if given vectors $(\vec x,\vec\xi)$ satisfy 
$(\vec x,\vec\xi)=((x_j,\xi_j))_{j=1}^4\in \qP_\vartheta(y,\zeta)$, {cf.\  \cite[Prop.\ 5.7]{ONeill}.}
Also, note that as then $\vartheta<\vartheta_2(y,\zeta,s_1)$, 
the inequality (\ref{eq: cut points observed after S}) yields 
 $\tilde x\in \V((\vec x,\vec \xi),t_0)$ and thus $q\in 
 J^-(\tilde x)\subset \V((\vec x,\vec \xi),t_0)$.
 Then  Lemma \ref{lem: sing detection in in normal coordinates 3} (iii) yields
 ${\Scle}((\vec x,\vec \xi),t_0)=\be_{U}(q)$. As   $\vartheta\in (0,\vartheta_3)$
 above can  be arbitrarily small, 
  we have that for 
  any $q=\gamma_{y,\zeta}(r)\in J^-({p^+})$
where $t_0<r\leq r_0=r_0(y,\zeta,s_1)$, we obtain
\beq\label{aaaa}
& &\hspace{-1cm}
\hbox{
$S=\be_{U}(q)$ is a genuine 
observation associated to  $\gamma_{y,\zeta}$ and}
\hspace{-1cm}
 \\
\nonumber
& &\hspace{-1cm}
S\cap 
\hat \mu=\{\hat \mu(\hat s)\},
\ \hat s:=f^+_{\hat \mu}(q)\leq 
\mathbb{S} (y ,\zeta,s_1).\hspace{-1cm}
\eeq

  {
\begin{lemma}\label{lem: Determination of hat s}  
Assume that $\gamma_{y,\zeta}(\R_+)$ does not intersect $\hat \mu$. Then
we have  $T (y ,\zeta ,s_1)=\mathbb{S} (y ,\zeta,s_1)$. 
  
\end{lemma}

\noindent{\bf Proof.}   
  Let us first prove that $T (y ,\zeta ,s_1)\geq \mathbb{S} (y ,\zeta,s_1)$.
To this end, let $s^\prime<\mathbb{S} (y ,\zeta,s_1)$
and  $x^\prime =\hat \mu(s^\prime)$.  
Assume  $S\in  \be_{U}(J^+({\hat x}_1)\cap J^-({p^+}))$ is a genuine 
observation associated to the geodesic $\gamma_{y,\zeta}$ 
and $S\cap \hat \mu=\{\hat \mu(s^\prime)\}$.
 Let   $q \in J^+({\hat x}_1)\cap J^-({p^+})$ be such that
$S=\be_U(q)$.

Then for arbitrarily small  $0<\vartheta<\vartheta_2(y,\zeta,s_1)$ 
there is
 $(\vec x,\vec\xi)\in \qP_\vartheta(y ,\zeta)$
 satisfying   
 ${\Scle}((\vec x,\vec \xi),t_0)=S$.  Let $\V=\V((\vec x,\vec \xi),t_0)$.
 Then by (\ref{eq: cut points observed after S}), we have 
 $x^\prime \in \V$.

 If the geodesics  corresponding to $(\vec x,\vec\xi)$
   intersect  
at some point $q^\prime\in J^-(x^\prime)$, then  by  Lemma \ref{lem: sing detection in in normal coordinates 3} (i), (ii) we have
$S=\be_{U}(q^\prime)$. 
Then 
$\hat \mu(s^\prime)=\hat \mu(f^+_{\hat \mu}(q^\prime))$
implying  $ f^+_{\hat \mu}(q^\prime)=s^\prime$. Moreover, we have
 then that $\be_{U}(q)\
=\be_{U}(q^\prime)$ and 
Theorem \ref{main thm paper part 2} (i) yields 
$q^\prime=q$.
Since $(\vec x,\vec\xi)\in \qP_\vartheta(y,\zeta)$ 
implies $(x_1,\xi_1)=(y,\zeta)$, 
we see that $q^\prime\in \gamma_{y,\zeta}([t_0,\infty))$.
As  $q\in J^+({\hat x}_1)$,
we see that $q=q^\prime\in 
\gamma_{y,\zeta}([t_0,\infty)) \cap J^+({\hat x}_1)=
\gamma_{y,\zeta}([r_0(y ,\zeta,s_1),\infty))$.
However, then $f^+_{\hat \mu}(q^\prime)\geq
\mathbb{S} (y ,\zeta ,s_1)  >s^\prime$ and thus $S\cap \hat \mu=\be_{U}(q)\cap \hat \mu$
can not be equal to $\{\hat \mu(s^\prime)\}$.

On the other hand, if the geodesics  corresponding to $(\vec x,\vec\xi)$
do not   intersect  at any point in 
 $ J^-(x^\prime)\subset \V$, \MTEXT{then either they intersect in some
 $q^\prime_1\in (M\setminus J^-(x^\prime))\cap \V$, do not intersect at all, or intersect at 
 $q^\prime_2\in M\setminus \V$. In the first case, $S=\be_U(q^\prime_1)$ do not satisfy
 $S\cap \hat \mu\in \hat\mu((-1,s^\prime))$. In the other cases,}
  Lemma \ref{lem: sing detection in in normal coordinates 3} (iii)
 yields 
 ${\Scle}((\vec x,\vec \xi),t_0)\cap \V=\emptyset$.
As  $x^\prime=\hat \mu(s^\prime)\in \V$, we see that     $S\cap \hat \mu$
can not be equal to $\{\hat \mu(s^\prime)\}$. Since above 
$s^\prime< \mathbb{S} (y ,\zeta,s_1)$ is arbitrary, this shows that
$T (y ,\zeta ,s_1)\geq \mathbb{S} (y ,\zeta,s_1)$.

Let us next show that $T (y ,\zeta ,s_1)\leq \mathbb{S} (y ,\zeta,s_1)$.
\MTEXT{Assume the opposite. Then, if $\mathbb{S} (y ,\zeta,s_1)=s_+$,
 we see by Def.\ \ref{def 56} that 
$T (y ,\zeta ,s_1)=\mathbb{S} (y ,\zeta,s_1)$ which leads to a contradiction. 
However, if 
$\mathbb{S} (y ,\zeta,s_1)<s_+$,  by Def.\ \ref{def: hat s 1},
we have
(\ref{eq: cut points observed after S}).  
This implies the existence of  $q_0=\gamma_{y,\zeta}(r_0)$, 
$r_0=r_0(y,\zeta,s_1)$ such that 
 $q_0\in  J^+({\hat x}_1)
\cap J^-({p^+})$ and
by (\ref{aaaa}), 
$S=\be_{U}(q_0)$ is a genuine 
observation associated to the geodesic $\gamma_{y,\zeta}$.
By Lemma \ref{conjugatepoints are nice} (ii), 
$\mathbb{S} (y ,\zeta,s_1)=f^-_{\hat \mu}(q_0)$ which implies,
by Def.\ \ref{def 56} that $T (y ,\zeta ,s_1)\leq \mathbb{S} (y ,\zeta,s_1)$.
%
%
Thus, $T (y ,\zeta ,s_1)=\mathbb{S} (y ,\zeta,s_1)$.}
 \hfill \Box \medskip  
  }

Next we reconstruct $\be_{U}(q)$ when $q$ runs over a geodesic segment.

\begin{lemma}\label{lem: Determination of be-sets}  
Let $s_-\leq s_2\leq s<s_1\leq s_+$ with  $s_1<s_2+\kappa_2$, let  
 ${\hat x}_j={\hat \mu}(s_j)$, $j=1,2$, and  ${\hat x}={\hat \mu}(s)$,
 ${\hat \zeta} \in L^+_{{\hat x}}U$, 
    $\|{\hat \zeta}\|_{g^+}=1$.   
    Let $(y,\zeta) \in L^+U$ be in the $\vartheta_1$-neighborhood of
    $({\hat x},{\hat \zeta})$ such that $y\in J^+({\hat x}_2)$.  Assume that $\gamma_{y,\zeta}(\R_+)$ does not intersect
    $\hat \mu$.   
Then, if we are given the data set  ${\cal D}(\hat g,\hat \phi,\e)$, we can determine the collection $\{\be_{U}(q) 
;\ q\in G_0(y,\zeta,s_1)\}$, where $G_0(y,\zeta,s_1)=\{q\in \gamma_{y ,\zzeta }([{t_0},\infty))\cap  (I^-({p^+}) \setminus J^+({\hat x}_1))\}$.
\end{lemma}  
  
 \noindent
{\bf Proof.} Let $s^\prime=\mathbb{S} (y ,\zzeta ,s_1)$, $x^\prime=\hat\mu(s^\prime)$, and
$\Sigma$ be the set of all genuine observations $S$ associated
to the geodesic $\gamma_{y,\zeta}$ such that $S$ intersects
$\hat\mu([-1,s^\prime))$.
%

Let $q=\gamma_{y,\zeta}(r)\in G_0(y,\zeta,s_1)$. Since
$\gamma_{y,\zeta}(\R_+)$ does not intersect
    $\hat \mu$,  
using a short cut argument for the geodesics from  $q$ to $q_0=\gamma_{y,\zeta}(r_0(y,\zeta,s_1))$
and from $q_0$ to $x^\prime$,
 we see that $f^-_{\hat \mu}(q)<s^\prime$.
 Then, 
$q\in I^-({p^+}) \setminus J^+({\hat x}_1)$ and $r<r_0(y,\zeta,s_1)$,
and we have using  (\ref{aaaa}) that
$S=\be_{U}(q)$ is a genuine 
observation associated to the geodesic $\gamma_{y,\zeta}$
and $S\cap \hat \mu=\{\hat \mu(f^-_{\hat \mu}(q))\}$ with $
f^-_{\hat \mu}(q)<s^\prime$. Thus
$\be_{U}(q)\in \Sigma$ and we conclude that $\be_{U}(G_0(y,\zeta,s_1))\subset \Sigma.$

%

Next, suppose $S\in \Sigma$. Then
there is $\hat \vartheta\in (0,\vartheta_2(y,\zeta,s_1))$ 
 such that for all $\vartheta\in (0,\hat \vartheta)$ 
 there is $(\vec x,\vec\xi)\in \qP_\vartheta(y,\zeta)$   so that
${\Scle}((\vec x,\vec \xi),t_0)=S$. 
Observe that by (\ref{eq: cut points observed after S}) we have  $\hat\mu([-1,s^\prime))\subset J^-(x^\prime)\subset 
 \V((\vec x,\vec \xi),t_0)$. 

First, consider the case when the geodesics corresponding to $(\vec x,\vec\xi)$ do not intersect
 at any point in $I^-(x^\prime)$. Then
 Lemma \ref{lem: sing detection in in normal coordinates 3} (iii) yields
 that ${\Scle}((\vec x,\vec \xi),t_0)$ is either empty or does not intersect $I^-(x^\prime)$. 
 Thus $S\cap I^-(x^\prime)={\Scle}((\vec x,\vec \xi),t_0)\cap I^-(x^\prime)$ is empty 
 and $S$ does not intersect
  $\hat\mu([-1,s^\prime))$. Hence $S$ cannot be in $\Sigma$.

Second, consider the case when the geodesics corresponding to $(\vec x,\vec\xi)$
intersect at some point $q\in I^-(x^\prime)$. Then, Lemma \ref{lem: sing detection in in normal coordinates 3} (iii)   
yields $S=\be_{U}(q)$. 
 {Since $(\vec x,\vec\xi)\in \qP_\vartheta(y,\zeta)$ 
implies $(x_1,\xi_1)=(y,\zeta)$,
 the intersection point $q$ has
a representation $q=\gamma_{x_1,\xi_1}(r)$. As  $q\in I^-(x^\prime)$, this yields
$q\in G_0(y,\zeta,s_1)$ and $S\in \be_{U}(G_0(y,\zeta,s_1))$.
 Hence $\Sigma\subset \be_{U}(G_0(y,\zeta,s_1))$.
 }

%
%
%
%

Combining the above arguments, we see  that $\Sigma=\be_{U}(G_0(y,\zeta,s_1))$.
As  $\Sigma$ is determined by the data set, the claim follows.
\hfill \Box \medskip

 
 Let $B(s_2,s_1)$ be the set of all $(y,\zeta,t)$ such that  there are
 $ {\hat x}={\hat \mu}(s)$, $s\in [s_2,s_1)$ and 
 ${\hat \zeta} \in L^+_{{\hat x}}U$,      $\|{\hat \zeta}\|_{g^+}=1$ so that
 $(y,\zeta) \in L^+U$ in $\vartheta_1$-neighborhood of    $({\hat x},{\hat \zeta})$, $y\in J^+({\hat x}_2)$,  and $t\in  [t_0,r_0(y,\zeta,s_1)]$.
   Moreover, let   $B_0(s_2,s_1)$ be the set of all $(y,\zeta,t)\in B(s_2,s_1)$
      \MTEXT{such that  $t<r_0(y,\zeta,s_1)$}
     and 
    $\gamma_{y,\zeta}(\R_+)\cap \hat \mu=\emptyset$.
\MTEXT{Lemma \ref{lem: Determination of be-sets} and the fact that 
${\cal D}(\hat g,\hat \phi,\e)$ determines $\gamma_{y,\zeta}\cap U_{\hat g}$}
show that,  when we are  given 
  the data  set ${\cal D}(\hat g,\hat \phi,\e)$, we  can determine 
the collection
$\Sigma_0(s_2,s_1):=\{\be_{U}(q);  
\ q=\gamma_{y ,\zeta }(t),\ (y,\zeta,t)\in B_0(s_2,s_1)\}$. 
We denote also $\Sigma(s_2,s_1):=\{\be_{U}(q);  
\ q=\gamma_{y ,\zeta }(t),\ (y,\zeta,t)\in B(s_2,s_1)\}$.

Note that the sets $\be_{U}(q)\subset {U}$, where $q\in J:=J^-(p^+)\cap J^+(p^-)$, can be identified
with the function, $F_q:\mathcal U_{z_0,\eta_0}\to \R$,
$F_q(z,\eta)=f^+_{\mu(z,\eta)}(q)$, c.f.\ (\ref{kaava D}). 
Let  $\U=\U_{z_0,\eta_0}$. 
When we  endow the set $\R^\U$ of maps $\U\to \R$
with the topology of pointwise convergence, 
Lemma \ref{B: lemma} yields that $F:q\mapsto F_q$ is continuous
map $F:J\to \R^\U$. 
By Theorem \ref{main thm paper part 2}, $F$ is one-to-one, and since $ J$
is compact  and $\R^\U$ is Hausdorff, we have that $F:J\to F(J)$ is homeomorphism.
Next, we identify $\be_{U}(q)$ and $F_q$.
%

{Using standard results of differential topology, we see that any neighborhood
of $(y,\zeta)\in L^+U$ contains  $(y^\prime,\zeta^\prime)\in L^+U$
such that the geodesic $\gamma_{y^\prime,\zeta^\prime}([0,\infty))$
does not intersect $\hat \mu$. 
Since $(y,\zeta)\mapsto r_0(y,\zeta,s_1)$ is lower semicontinuous, this implies that $\Sigma_0(s_2,s_1)$ is dense in $\Sigma(s_2,s_1)$. 
Hence we obtain
the closure $\overline \Sigma(s_2,s_1)$ of $\Sigma(s_2,s_1)$
as  
the limits points of $ \Sigma_0(s_2,s_1)$.

\MTEXT{Then, we obtain the set we $\be_U( J^+({\hat \mu}(s_2 ))\cap J^-({p^+}))$
as
the union  $\overline\Sigma(s_2,s_1)\cup
\be_U(J^+({\hat \mu}(s_1))\cap J^-({p^+}))\cup \be_U(\K_{t_0}\cap J^+({\hat \mu}(s_2 )))$,  
see (\ref{set Kt0}).

Let  
$s_0,\dots,s_K\in [s_-,s_+]$ be such that  
 $s_j>s_{j+1}>s_j-\kappa_2$ and $s_K=s_-$.
Then, by iterating the above construction so that the values of the parameters $s_1$ and $s_2$ are replaced by
$s_j$ and $s_{j+1}$, respectively,
we can construct the set $\be_{U}( J^+(\hat \mu (s_-))\cap J^-(\hat \mu (s_+)))$.}}

Moreover, similarly to the above construction, we can find the   
sets $\be_{U}( J^+(\hat \mu (s^\prime))\cap J^-(\hat \mu (s^{\prime\prime}))$  
for all $s_-<s^\prime<s^{\prime\prime}<s_+$,
and taking their union, we find the set 
$\be_{U}( I(\hat \mu (s_-),\hat \mu (s_+))$.
   By Theorem \ref{main thm paper part 2}  we can reconstruct the manifold  
 $I(\hat \mu (s_-),\hat \mu (s_+))$  and the conformal structure on it.    
This proves   
Theorem \ref{alternative main thm Einstein}. \hfill \Box \medskip  
  

\noindent
{{\bf Remark 5.1.} The proof of Theorem \ref{alternative main thm Einstein}
can be used to analyze approximative reconstruction of the set $I_{\hat g}(p^-,p^+)=I^+_{\hat g}(p^-)\cap I^-_{\hat g}(p^+)$ and its conformal structure with only one
measurement: We choose one suitably constructed source $\F$, supported in $W_{\hat g}$, measure the fields
$(f,\phi)$ produced by this source in $U_{\hat g}$  and aim to construct an approximation
of the conformal class of the metric $g$ in $I_{\hat g}(p^-,p^+)$. 
In the proof above we showed that it is possible
to use the non-linearity to create an artificial point source at 
a point $q\in J_{\hat g}(p^-,p^+)$.
Using the same method we see that it is possible to create with a single source $\F$ an arbitrary number
of artificial point sources. To see this, let  $\delta>0$ and $P,Q\in \Z_+$ and  consider 
 points  $x_p\in U_{\hat g}\cap  J_{\hat g}(p^-,p^+)$, $p=1,2,\dots,P$.
Let $\xi_{p,k}\in \Sigma_p=\{\xi\in T^*_{x_p}M;\ \|\xi\|_{\hat g}=1\}$, $k=1,2,\dots,n_Q$ be a maximal $1/Q$ net
on the set $\Sigma_p$ and $\Sigma_{p,k}(R)$ be an $R$-neighborhood of $\xi_{p,k}$ in
$\Sigma_p$. Let $R_1=1/Q$ and $R_2=2/Q$, and consider real numbers
$a(p,k)\in (-31,-30)$, chosen to be $a(p,k)=-31+1/j(p,k)$,
where $j:\Z_+^2\to \P$ is a bijection from $\Z_+^2$ to the set $\P$ of
the prime numbers.
Let $F_{p,k}\in \I^{a(p,k)}(\Sigma_{p,k}(R_2))$ be Lagrangian distributions
whose principal symbols are non-vanishing on  $\Sigma_{p,k}(R_1)$.

Assume next that $(M,g)$ is in generic manifold (i.e., it is in the intersection
of countably many open and dense sets in a suitable space of smooth manifolds), points $x_p$ have generic positions,
and let $\delta>0$.
Note that the sets $\Sigma_{p,k}(R_1)$ are a covering of the unit sphere 
$\Sigma_{p}$ and the linearized waves $u_{p,k}={\bf Q}_{\hat g}(F_{p,k})$
are singular on a subset of the light-cones $\L^+_{\hat g}(x_p)$.
Let $\e>0$ be small enough and consider a suitable source $F_\e$,
with $\p_\e F_\e|_{\e=0}=\sum_{p,k}F_{p,k}$, 
that produces the perturbation $u_\e(x)=\sum_{n=1}^4 \e^n u_n(x)+O(\e^5)$ for
$(\hat g,\hat \phi)$.
Assume that we measure the singular supports of the waves $u_n,$ $n=1,2,3,4$
produced by $n$-th order interaction of the waves. When $\e$ is small enough, this could be
done e.g.\  by using thresholding of the curvelet coefficients, of a suitable order, of the solution \cite{Hoop1,Hoop2}. 
Then in $I_{\hat g}(p^-,p^+)$, outside the singular support of the wave $u_3$,  the wave $u_4$ is a sum of
a smooth wave and the waves produced by artificial point sources $\F_\ell$ located at $q_\ell$, 
$\ell=1,2,\dots,L$.  Here $q_\ell\in J_{\hat g}(p^-,p^+)$
 are the intersection points of any four light cones $\L_{\hat g}^+(x_{p_1(\ell)})$,
$p_1(\ell),p_2(\ell),p_3(\ell),p_4(\ell)\leq N$. When $Q$ and $P$ are large enough, so that $N$ is large,
the points $q_\ell$ are a $\delta$-dense subset of $J_{\hat g}(p^-,p^+)$.
Moreover,  for the chosen orders $a(p,q)$, the orders of the sources
$\F_\ell$ at points $q_\ell$ are all different and thus the waves
${\bf Q}_{\hat g}(\F_\ell)$ have different orders. As a very rough analogy,
we can produce in $J_{\hat g}(p^-,p^+)$ an arbitrarily dense collection
of artificial points sources having all different colors.
Using this observation one
can show, using similar methods to those in \cite{AKKLT},
that  in a suitable  compact class of Lorentzian manifolds 
having no conjugate points
the measurement with the source
$F_\e$, defined using sufficiently large $P$ and $Q$ and generic values of principal symbols
on $\Sigma_{p,q}$, determines 
a $\delta$-approximation (in a suitable sense) of the conformal class of the
 manifold $(J_{\hat g}(p^-,p^+),\hat g)$. The details of this construction will
 be considered elsewhere.
}


  \observation{   
{\bf Remark 5.1.}   
When the adaptive source functions $\mathcal S_\ell$ are the  
ones given in Appendix $C$ we can consider an  
improvement of   Theorem \ref{main thm Einstein}:  
When we are given not  ${\cal D}(\hat g,\hat \phi,\e)$ but assume  
that we know ${\cal D}_{0}(\hat g,\hat \phi,\e)$,  
we can use Remarks  
3.1, 3.2, 3.3, 4.1, and 4.2 and the proofs of this section  
 to show that  $ I({p^-},{p^+})$  
and the conformal structure on it can be reconstructed.  
This yields similar result to Theorem \ref{main thm Einstein}.  
}\medskip

\generalizations{
 \section{Generalizations and outlook}\label{sec: Generalization}

\HOX{T23. This section needs to be checked and made less messy. - Matti }
{
In this section we present a sketch of the analysis how the above results can
be applied for the single measurement inverse problem, that is, how we can  obtain an approximation
of the space-time with a single measurement. Recall that  
 Theorem \ref{main thm Einstein} concerns the case when all possible
measurements near the freely falling observer $\mu$ are known. In principle this implies that the perfect
spacetime cloaking with a smooth globally hyperbolic metric is not possible, that is, no two
space-times having different conformal structure appear the same in
all measurements. However, one can ask if 
one can make an approximative image of the space-time knowing
 only one measurement.  In general, in many
inverse problems several measurements can be packed together to 
one measurement.
For instance, for the wave equation with a  time-independent simple metric
this is done in \cite{HLO}. Similarly, Theorem \ref{main thm Einstein} and its proof
make it possible to do approximate reconstructions as discuss below.
\HOX{Add a remark that we can produce any number $N\in \Z_+$ of point
sources in the space time.}

\subsection{Pre-compact collection of Lorentzian manifolds}
For $m\in \Z_+$ and $\Lambda_{m}\in \R_+$, let $\M_{m,\Lambda}$
be the set of pointed  globally hyperbolic Lorentzian manifolds $(M,g,x_0)$ such that
the corresponding Riemannian manifold $(M,g^+)$
has the property that in the closure of the ball $B_{g^+}(x_0,m)$ with center $x_0$
and radius $r=m$
the covariant derivatives of the curvature tensor $R=R(g^+)$ of the Riemannian
metric $g^+$ satisfy $\|\nabla^j R\|_{g^+}\leq \Lambda$
for all $j\leq m$, the  injectivity radius $i_0(x)$ of $(M,g^+)$ satisfies
$i_0(x)>\Lambda^{-1}$ for all $x\in \hbox{cl}(B_{g^+}(x_0,m))$. 
Let us choose some numbers $\Lambda_m$ for all $m\in \Z_+$,
and let 
\ba
\M=\bigcap_{m\in \Z_+}\M_{m,\Lambda_m}
\ea
be the set endowed with the initial topology for which all
identity maps $\iota_m:\M\to \M_{m,\Lambda_m}$ are continuous.
Using the Cheeger-Gromov theory and the Cantor diagonalization
procedure, we see that the set $\B_R$ of the balls  $B_{g^+}(x_0,R)$ of radius $R$ of the pointed manifolds  
 $(M,g,x_0)\in \M$ form a precompact set in the  topology given by the Gromov-Hausdorff metric
 and the closure of $\B_R$, denoted $\hbox{cl}_{GH}(\B_R)$, consists of smooth manifolds,
 see \cite{And1,Cheeger}, see  also the applications for inverse problems in  \cite{AKKLT}.
 Indeed, we see that the balls of radius $R$ of the pointed manifolds in $\M_{m,\Lambda_m}$ are
 precompact for all $m$ and thus  using a Cantor diagonalization argument
 we see claimed is precompactness.
  
Let us fix $R>0$ and define  $\mathcal N$ to be the set 
$((M,g),(x_0,\xi_0))$  such that  $(M,g,x_0)$ is a pointed  Lorentzian manifold for which  $B_{g^+}(x_0,R)\in \hbox{cl}_{GH}(\B_R)$,
and a time-ike vector $(x_0,\xi_0)\in TM$ is such that the  
 freely falling observer $\hat \mu$ in $(M,g,x_0)$, such that $\hat \mu(-1)=x_0$ and 
 $\p_s\hat \mu(-1)=\xi_0$ satisfies
$U_{g}\cup J_{g}(\hat \mu(-1-\e_0),\hat \mu(1+\e_0))\subset B_{g^+}(x_0,R)$, with fixed parameters 
$\e_0>0$ and  $\hat h>0$ (that determines $U_g$), 
and finally that $\hat \mu([-1,1])$ is such that no light-like geodesic has cut points
in $J_{g}(\hat \mu(-1),\hat \mu(1))$.

\subsubsection{Geometric preparations for single source measurement}
Let $S^{g^+}(y,r)$ denote the sphere of $T_yM$ of radius $r$ with respect to the metric
$g^+$ and  $S^{g^+}(y)= S^{g^+}(y,1)$.
For  $(y,\xi)\in L^+M$ we define spherical surfaces
\ba
& &\L((y,\xi);\rho)=\{\exp_y(t\theta);\ \theta\in L_y^+M,\ 
 \|\theta-\frac 1{\|\xi\|_{g^+}}\xi\|_{g^+}<\rho,\ t>0\},\\
 & &\L^c((y,\xi);\rho)=\{\exp_y(t\theta);\ \theta\in L_y^+M,\  
 \|\theta-\frac 1{\|\xi\|_{g^+}}\xi\|_{g^+}\geq \rho,\ t>0\}.
  \ea
  We say that $y$ is the center point of these surfaces.

In particular,
$((M,g),(x_0,\xi_0))\in \mathcal N$ implies that
no light-like geodesic starting from $\mu_{z,\eta}$ can intersect $\mu_{z,\eta}$ again.

 Using compactness of $\mathcal N$, \HOX{We need to check (\ref{consequence of compactness})
 and existence of $\rho>0$  and add parameter $r_0$ below. Idea below is first
 to choose $\e$, i.e. the accuracy on which we want to find $(M,g)$, then
 suitable $L$, then perturbation so that light-like geodesics do not intersect
 ${\hat x}_a$, then $r_0$ and then $\rho$.}
for any $r_0>0$ we can choose $\rho=\rho(r_0)>0$ such that 
for all  $((M,g),(x_0,\xi_0))\in \mathcal N$  and any ${\hat x}_j\in J_{g}(\hat \mu(s_-),\hat \mu(s_+)) $, $j=1,2,3,4$ and $\xi_j\in L^+_{{\hat x}_j}M$ such that 
\beq\label{far away}
d_{\hat g^+}({\hat x}_j,
\gamma_{{\hat x}_k,\xi_k}(\R_+))\geq r_0,\quad \hbox{ when $j\not   =k$},
\eeq
the  set
\beq\label{consequence of compactness}
\bigcap_{j=1}^4 \L(({\hat x}_j,\xi_j);\rho)\cap J_{g}(\hat \mu(s_-),\hat \mu(s_+))\hbox{ contains at most one point.}\hspace{-2.1cm}
 \eeq
 Indeed, if there are no such $\rho$ we find a sequence of manifolds $(M_k,g_k)$
 and $(\vec {\hat x}_k,\vec \xi_k)\in (TM_k)^4$ such that the sets
 (\ref{consequence of compactness}) contain at least two points with the parameter $\rho_k$
 and $\rho_k\to 0$ as $k\to \infty$. Considering suitable subsequences of $(M_k,g_k,x_k)$
 that converge to $(M,g,x)$ 
 and $(\vec {\hat x}_k,\vec \xi_k)$ that converge to $(\vec y,\vec \xi)$ we 
 see that some  geodesics in $(M,g,x)$ have cut points in $J_{g}(\hat \mu(s_-),\hat \mu(s_+))$
 that is not possible due to the definition of $\mathcal N$. Hence required parameter $\rho>0$ exists.

Next, for $((M,g),(x_0,\xi_0))\in  \hbox{cl}( \mathcal N)$ we denote 
\ba
J((M,g),(x_0,\xi_0))=J_{g}(\hat \mu(s_-),\hat \mu(s_+))
\ea
and 
\ba
\mathcal J=\{(J((M,g),(x_0,\xi_0)),[g],x_0,\xi_0);\ ((M,g),(x_0,\xi_0))\in \hbox{cl}( \mathcal N)\},
\ea
where $[g]$ denotes the conformal class of $g$. Equivalence classes $[g]$
are closed sets and we use Hausdorff distance for to define
the distance for the equivalence classes $[g_1]$ and $[g_2]$ using representatives
of the equivalence classes for which the determinant of the metric tensor in the Fermi
coordinates are equal to the determinant of the metric tensor of the Euclidean
metric in the Fermi coordinates associated to a line. We call such representative of 
the equivalence class $[g]$ a Fermi-normalized metric and denote it by $g^{(n)}$.
The corresponding Riemannian metric is denoted by $g^{(n),+}$

If $p\in J_{g}(\hat \mu(s_-),\hat \mu(s_+))$, then for given $(z,\eta)$ there
is unique point $y\in \mu_{z,\eta}(s)$
for which there exists a light like geodesic
 $\gamma_{y,\xi}([0,t])$ that connects $y$ to $p$. Moreover, then
 $s=f^-_{z,\eta}(p)$. Let $(z_j,\eta_j)\in \U_{z_0,\eta_0}$, $j=1,2,\dots$,
 be dense set 
 $(z_j,\eta_j)\not= (z_0,\eta_0)$. We consider these points to be defined so that they have fixed 
 representation in the Fermi coordinates. \HOX{Add details on the points in the Fermi coordinates}
 Using compactness of $J_{ g}(\hat \mu(s_-),\hat \mu(s_+))$,
  we see that  there is $m_0=m_0(M,g,x_0,\xi_0)$ such that
 for all $p\in J_{ g}(\hat \mu(s_-),\hat \mu(s_+))$ there are four pairs $(z_{j_k},\eta_{j_k})$
 with $j_k\leq m_0$ such
 that $q\mapsto ( f^-_{z_{j_k},\eta_{j_k}}(q))_{k=1}^4$ forms smooth coordinates
 in a neighborhood of $p$.

Let choose dense set 
 of numbers $s_j^k\in (-1,1)$, $j,k\in \Z_+$ and let
 ${\hat x}_{j,k}=\gamma_{z_j,\eta_j}(s_j^k)$. For each
 ${\hat x}_{j^\prime,k^\prime}$ we choose a dense
 set of directions, $\xi_{j^\prime,k^\prime,n}$, $n=1,2,\dots,$ of $S^{g^+}({\hat x}_{j^\prime,k^\prime})\cap L^+_{{\hat x}_{j^\prime,k^\prime}}M$. 
 After this, let us use the pairs $({\hat x}_{j^\prime,k^\prime},\xi_{j^\prime,k^\prime,n})$ to define
the pairs $({\hat x}_{j^\prime,k^\prime}(t_{j^\prime,k^\prime}),\xi_{j^\prime,k^\prime,n}(t_{j^\prime,k^\prime}))$ with some 
 $t_{j^\prime,k^\prime}>0$, and 
 re-enumerate the obtained pairs 
 $({\hat x}_{j^\prime,k^\prime}(t_{j^\prime,k^\prime}),\xi_{j^\prime,k^\prime,n}(t_{j^\prime,k^\prime}))\in L^+M$,
 $j^\prime,k^\prime,n\in \Z_+$ as 
$({\hat x}_a,\xi_a)$, $a\in \Z_+$. 
 
 Let $\A$ be the set of 4-tuples $\vec a=(a_\ell)_{\ell=1}^4\in \Z_+^4$ that contain
 four pairwise non-equal positive integers.  
 Let
 \ba
 q(\vec a)=\bigcap_{\ell=1}^4 \L(({\hat x}_{a_\ell},\xi_{a_\ell});\rho)\cap J_{ g}(\hat \mu(s_-),\hat \mu(s_+)).
 \ea
 Note that $q(\vec a)$ contains then only one point or is empty
 if (\ref{far away}) is satisfied with some $r_0$ and $\rho<\rho(r_0)$. Also,
 let 
 \ba
 & & \Psi^{(1)}_k(\vec a)=\prod_{\ell,\kappa\in \{1,2,3,4\},\ \ell\not=\kappa } \phi(\frac 1k \dist_{g^{(n),+}}({\hat x}_{a_\kappa},
\gamma_{{\hat x}_{a_\ell},\xi_{a_\ell}}(\R_+))),\\
& & \Psi^{(2)}_k(\vec a)=\prod_{\ell=1}^4 \phi(\frac 1k \dist_{g^{(n),+}}(q(\vec a), 
 \L^c(({\hat x}_{a_\ell},\xi_{a_\ell});\rho)),\\
& & \Psi_k(\vec a)=\Psi^{(1)}_k(\vec a)\,\cdotp \Psi^{(2)}_k(\vec a), 
  \ea
where $\phi\in C(\R\cup\{\infty\})$ is a function with $\phi(t)=0$ for $t\leq 1$ and $\phi(t)=1$
for $t>2$ and $t=\infty$ and $\dist_{g^{(n),+}}$ is the distance
with respect to the metric $g^{(n),+}$ that is conformal to the metric $g^+$
and which determinant in the Fermi coordinates is normalized as above.

Using the existence
 of $m_0(M,g,x_0,\xi_0)$ and local coordinates associated to four pairs $(z_{j_k},\eta_{j_k})$,
we see that $\{q(\vec a);\ \vec a\in \A,\ q(\vec a)\not=\emptyset\}$ is a dense set of $J_{ g}(\hat \mu(s_-),\hat \mu(s_+))$.

 Let us then consider the continuous map $G:\mathcal J\to \R^{\Z_+\times \Z_+}$,
 \ba
G: (J((M,g),(x_0,\xi_0)),[g],x_0,\xi_0)\mapsto (\Psi_k(\vec a)\,f^+_{z_j,\eta_j}(q(\vec a)))_{\vec a\in \A,\ j\in \Z_+,\ k\in \Z_+}.
 \ea
  \HOX{Check that $G$ is continuous}
In the case when the set $q(\vec a)$ is empty, we define 
 $f^+_{z_j,\eta_j}(q(\vec j(n),\vec k(n))$ to be equal to $\sup\{s;\ \mu_{z_j,\eta_j}(s)\in 
 J_{ g}(\hat \mu(s_-),\hat \mu(s_+))\}$. In fact, we can define this value
 to be arbitrary because of the "cut-off" function $\Psi_k(\vec a)$.
 
  Note that the sets $\P_V(q)$ are  closed. We see using the arguments of the main text that the map $G$ is injective in $\mathcal J$. 
Since 
 $\mathcal J$ is compact and $G$ is continuous, we see for $G_N:\mathcal J\to \R^{\Z_+\times \Z_+}$,
 \ba
G_N: (J((M,g),(x_0,\xi_0)),[g],x_0,\xi_0)\mapsto (\Psi_k(\vec a)\,f^+_{z_j,\eta_j}(q(\vec a)))_{|\vec a|\leq N,\ j\leq N,\ k\leq N}
 \ea
that for any $\e>0$ there is  $N$ 
 such that $G_N(g)$ determines the conformal type $(J((M,g),(x_0,\xi_0)),[g])$ up to $\e$-error in the  metric  in $\mathcal J$. 
 Here, for on $ \mathcal J $  we use
 the distance function defined using 
 $((B_1,[g_1]),(B_2,[g_2]))\mapsto
 d_{GH}((B_1,g_1^{(n),+}),(B_2,g_2^{(n),+}))$.
  \HOX{We need to define the topology of $\mathcal J$ carefully}

 \subsection{Single source measurement}
 Let us now consider distorted plane  waves propagating
 on the spherical surfaces  $\L(({\hat x}_{a_\ell},\xi_{a_\ell});\rho)$, $|\ell |\leq L$.
 When $N$ is large, these waves interact and produce 3-wave and 4-wave interactions.
 These waves can be produced by sources ${\bf f}_{a_\ell}$ supported in an arbitrarily
 small neighborhoods  $W_{a_\ell}\subset U_{\hat g}$ of the "conic" points ${\hat x}_{a_\ell}$ of the spherical surfaces  $\L(({\hat x}_{a_\ell},\xi_{a_\ell});\rho)$.

  We see that for generic set  (i.e.\ in an open and dense set) of  manifolds in $\mathcal N$ it
holds that   there are no points ${\hat x}_{a_\ell}$  \HOX{This "generic" statement needs to be checked}
 and ${\hat x}_{a_\ell^\prime}$ with ${a_\ell}\not= {a_\ell^\prime}$ that
 can be connected with a light-like geodesic. This is due to
 the fact that this condition holds in an open set and we see that if for any $L\in \Z_+$
 we can order the points ${\hat x}_{a(\ell)}$, $\ell=1,2,\dots,L$ in causal order, such that ${\hat x}_{a(\ell)}\not \in J^+({\hat x}_{a(\ell+1)})$, and then making a small perturbation to the metric near
each point ${\hat x}_{a(\ell)}$, in causal order, we can  choose such perturbations
that the point ${\hat x}_{a(\ell)}$ is not on the surface $\cup_{j<\ell}\L^+({\hat x}_{a(j)})$. 
{Below we will assume that neighborhoods \HOX{This generic statement needs to be checked}
 $W_{a_\ell}$ of ${\hat x}_{a_\ell}$ of the  points ${\hat x}_{a_\ell}$ are so small
 that $W_{a_\ell}\subset U_{\hat g}$ does not intersect any $\L^+({\hat x}_{a(j)})$, where
$j\not =\ell$ and $j\leq L$.

In the main text we used sources ${\bf f}_{a_\ell}$ that are elements of  the space
$\I(Y_{a_\ell})$  and produce waves $u_{a_\ell}\in \I (N^*Y_{a_\ell},N^*K_{a_\ell})$,
where $K_{a_\ell}\subset \L^+({\hat x}_{a(\ell)})$ are 3-dimensional subsets  and 
 $Y_{a_\ell}= \L^+({\hat x}_{a(\ell)})\cap {\bf t}^{-1}(T_{a(\ell)})$ are 2-dimensional sub-manifolds.
 As the 
products of  distributions associated to two Lagrangian manifolds are difficult to analyze, we modify the construction by 
replacing $K_a$ by surfaces $K_a\subset \L^+_{\hat g}({\hat x}_a)$ and 
  ${\bf f}_{a_\ell}\in \I(Y_{a_\ell})$ by sources $\tilde {\bf f}_{a_\ell}\in \I(K_{a_\ell})$,
supported in the  neighborhood  $W_{a_\ell}$ of ${\hat x}_{a_\ell}$.
As  then $\tilde {\bf f}_{a_\ell}\in \I(N^*K_{a_\ell})$ and $N^*K_{a_\ell}$ is invariant in the future
bicharacteristic flow of $\square_{\hat g}$, the sources $\tilde {\bf f}_{a_\ell}$
 produce by \cite[Prop. 2.1]{GU1} waves  $ u_{a_\ell}\in \I (N^*K_{a_\ell})$,
where the principals symbols are not anymore evaluations of the principal
symbol of the source at some point but integrals of the principal symbol along a bicharacteristic.

The sources $\tilde {\bf f}_{a_\ell}\in \I (K_{a_\ell})$  can be chosen so that there are 
 neighborhoods $V_{a_\ell}$ of the surfaces  $\L(({\hat x}_{a_\ell},\xi_{a_\ell});\rho)$  that satisfy
$W_{a_\ell}\subset V_{a_\ell}$,
the linearized waves $ u_{a_\ell}={\bf Q}_{ \hat g}\tilde {\bf f}_{a_\ell}$ satisfy  singsupp$( u_{a_\ell})\subset V_{a_\ell}$,
and   
\beq\label{eq: V ehto}
\hbox{$V_{a_\kappa}\cap W _{a_\ell}=\emptyset$  for all $\ell,\kappa\leq L$,
 $\ell\not =\kappa$.}
\eeq
We also assume that
$W_{a_\kappa}$  are such that $J^+_{\hat g}(W_{a_\ell})\cap
J^-_{\hat g}( W_{a_\ell})\subset  U_{\hat g}$ and that the wave map coordinates
corresponding to the Euclidean metric are well defined in
this set. Note that the above properties of $W_{a_\kappa}$ can be obtained
by choosing all surfaces $K_{a_\kappa}$, $\kappa\leq L$ with a the same parameter
$s_0$ that is small enough. Below, we assume that $s_0$ is chosen in such a way. 

Consider now the non-linear interactions produced by the source
\beq\label{Fe source}
F_\e(x)=\e\left(\sum_{\ell=1}^N\tilde {\bf f}_{a_\ell}(x)\right),
\eeq
where $\e>0$ is small.  
This source produces the wave
\beq\label{nl wave}
u(x)=u^{(4)}_\e(x)+O(\e^5),\quad u^{(4)}_\e(x)=\sum_{j=1}^4 \e^j\M^{(j)}(x),
\eeq
where we consider $O(\e^5)$ as a term which is so small that it can be considered to be negligible. \HOX{Maybe we could use $\vec \e$ with $\e_j=10^{-100-10j}$ etc and justuse the considerations of main text}
Let us next assume that we can measure  the singularities of $\M^{(j)}$ for all $j\leq 4$.
Let $S_j$ be the intersection of $U_{\hat g}$ and singular support of  $\M^{(j)}$, $j\leq 4$.

\subsubsection{Analysis of 1st, 2nd, and 3rd order waves}
As in the source $F_\e$ we use only one parameter $\e$ instead of four parameters
$\vec\e=(\e_1,\e_2,\e_3,\e_4)$ we need to consider self-interaction terms,
that is, in our computations the permutations $\sigma:\{1,2,3,4\}\to \{1,2,3,4\}$
need to be replaced by general (non-injective) maps $\sigma:\{1,2,3,4\}\to \{1,2,3,4\}$.
This does not cause difficulties in analyzing the terms $\M^{(j)}$ for all $j\leq 3$
as because of condition (\ref{eq: V ehto}) the source terms do not cause
self-interactions with terms corresponding to different Lagrangian manifolds
$\Lambda_{a_\ell}=N^*K_{a_\ell}$ and the interactions of linearized waves
corresponding to two such manifolds produce  waves that are in
$\I(N^*K_{a_\ell},N^*K_{a_\ell,a_\kappa})$.

For generic sources singsupp$(\M^{(1)}) $ is equal to $S_1=\cup_{\ell\leq L} K_{a_\ell}\cap U_{\hat g}$ and 
singsupp$(\M^{(2)})$  is equal to $S_2=\cup_{\ell,\kappa\leq L} K_{a_\ell,a_\kappa}\cap U_{\hat g}
\subset S_1$. 

Let us then consider singularities of $\M^{(3)}$ using similar methods that we used to analyze the fourth
order terms. Again, observe that we need to consider also the terms for 
which the map $\sigma:\{1,2,3\}\to \{1,2,3\}$ is not bijective. For such terms,
e.g., if $\sigma(1)=r_1$, $\sigma(2)=r_2$ and $\sigma(3)=r_3$ with $r_1=r_3$,   
we see using \cite[Lemmas 1.2 and 1.3]{GU1} 
 \cite[Thm.\ 1.3.6]{Duistermaat}  
 that e.g. $v\in u_{r_1}\cdotp {\bf Q}(u_{r_1}u_{r_2})$ has
 wave front set that is subset of $N^*K_{r_1r_2}$, implying
 that wave front set of ${\bf Q}v$ is subset of  $N^*K_{r_1r_2}$.
 Using this we see that outside $S_2$ only  the terms for which map $\sigma:\{1,2,3\}\to \{1,2,3\}$
 is bijective cause singularities. Let us next show that this kind of singularities can be observed.

To this end, let us consider the term
 $\bra F_\tau,{\bf Q}(A[u_3,{\bf Q}(A[u_2,u_1])])\cet$ where
 all operators $A[v,w]$ are of the type $A_2[v,w]=\hat g^{np}\hat g^{mq}v_{nm}\p_p\p_q w_{jk}$,
 cf.\ (\ref{eq: tilde M1}) and (\ref{eq: tilde M2}), we obtain  formulas similar to
  (\ref{eq: def of D}), where the first factor in 
   (\ref{def of D}) is removed. Next we show that using this and the similar analysis we did
   for the fourth order terms with the WKB analysis, we can see that the 
   principal symbol of  $\M^{(3)}$ in the normal coordinates does not vanish
   for generic interacting distorted plane  waves. Indeed, the indicator function
   $\Theta^{(3)}_\tau$ can be obtained as a sum of terms similar to $T^{(3),\beta}_\tau$ given in
    (\ref{eq: 3rd order term}). For these terms we obtain in the Minkowski space formulas
    analogous to (\ref{first asymptotical computation}) and (\ref{t 4 formula}). Let us consider  the case when
     $b^{(j)}$, $j\leq 5$ are  light-like co-vectors such that
    $b^{(5)}\in \hbox{span}(b^{(1)},b^{(2)},b^{(3)})$ and vector $b^{(4)}$ is linearly independent of $b^{(1)},b^{(2)},b^{(3)}$  
 so that
    $y=A^{-1}x$ are coordinates such  that
   ${\bf p}_4=(A^{-1})^tb^{(4)}=0$. The particular example of such co-vectors that we consider, are
   \ba
& &b^{(5)}=(1,1,0,0),\quad b^{(1)}=(1,1-\frac 12\rhoepsilon_1^2+O(\rhoepsilon_1^3),\rhoepsilon_1,\rhoepsilon_1^3),\\
& &b^{(3)}=(1,-1,0,0),\quad b^{(2)}=(1,1-\frac 12\rhoepsilon_2^2+O(\rhoepsilon_2^3),\rhoepsilon_2,\rhoepsilon_2^{201})
\ea
 and $ \rhoepsilon_1=\rhoepsilon^{100}$, $\rhoepsilon_2=\rho$, and $b^{(4)}$ is some fixed vector.
 Then $b^{(5)}=\sum_{j=1}^3 a_j b^{(j)}$,
 where $a_1=O(1)$, $a_2\sim -a_1\rho^{99}$, and $a_3=O( \rho^{2})$ as $\rho\to 0$.   For
 $T^{(3),\beta}_\tau$ given in
    (\ref{eq: 3rd order term}), we have then 
   \HOX{WKB needs to be checked.}
   \ba
T^{(3),\beta}_\tau&=&
\bra  u_\tau, h\,\cdotp \B_3u_3\,\cdotp {\bf Q}_0(\B_2u_2\,\cdotp \B_1 u_1)\cet\\
& &\hspace{-1cm}=C\frac {\P_\beta^\prime} {\omega_{12}}
\bra u_\tau, h\,\cdotp u_3\,\cdotp v^{a-k_2+1,a-k_1+1} (\,\cdotp ;b^{(2)},b^{(1)})\cet
\\
\nonumber&&\hspace{-1cm} =C\frac {\P^\prime_\beta}{\omega_{12}} 
\int_{\R^4}
e^{i\tau(b^{(5)}\cdotp x)} 
h(x)(b^{(3)}\,\cdotp x)^{a-k_3}_+\cdotp\\
\nonumber& &\quad\quad\quad\cdotp
(b^{(2)}\,\cdotp x)^{a-k_2+1}_+
(b^{(1)}\,\cdotp x)^{a-k_1+1}_+\,dx,\\
& &\hspace{-1cm}
=\frac {C\P^\prime_\beta\det (A^\prime)}{ \omega_{12}} \int_{(\R_+)^3\times \R}
e^{i\tau {\bf p}\cdotp y} 
h(Ay){\hat x}_3^{a-k_3}{\hat x}_2^{a-k_2+1}
{\hat x}_1^{a-k_1+1}\,dy
\\
\nonumber&&\hspace{-1cm} =
 {C\det (A^\prime)\,{\P^\prime_\beta}}
\frac {(i\tau)^{-(10+3a-|\vec k_\b|)}(1+O(\tau^{-1}))}{{\bf p}_3^{(a-k_3+1)}{\bf p}_2^{(a-k_2+2)}{\bf p}_1^{(a-k_1+2)}\omega_{12}}\,,
\ea
   where ${\bf p}_j=g(b^{(5)},b^{(j)})=-\frac 12 \rhoepsilon_j^2+O(\rhoepsilon_j^3)$, $j\leq 3$ and $\omega_{12}=g(b^{(1)},b^{(2)})=
   -\frac 12 \rhoepsilon_2^2+O(\rhoepsilon_2^3)$. 
Moreover, with the appropriate choice of polarizations   $v_{(\ell)}^{nm}=\hat g^{nj}\hat g^{mk}v^{(\ell)}_{jk}$ given in  (\ref{chosen polarization}) for $j=2,3,4$, c.f. (\ref{vbb formula}), we have 
\beq\label{pre-eq: def of D prime}
\P_{\beta_1}^\prime =(v_{(3)}^{pq}b^{(1)}_pb^{(1)}_q)(v_{(2)}^{nm}b^{(1)}_nb^{(1)}_m)\D,
\eeq
 and
\ba
 \D =\hat g_{nj}\hat g_{mk}v_{(5)}^{nm}v_{(1)}^{jk},
 \ea 
 where we note that $v_{(5)}^{nm}$ does not need to satisfy any harmonicity type conditions.
We see that the term $T^{(3),\beta_0}_\tau$, where $k_1=k_1^{\beta_0}=4$, dominates the other terms 
$T^{(3),\beta}_\tau$,when 
$\rho\to 0$. Thus, using the real analyticity of the leading order
coefficient in the asymptotic of  $\Theta^{(3)}_\tau$, we see similarly to the main text that the principal symbol of $\M^{(3)}$ in the normal
coordinates does not vanish in the generic case on the manifolds
$\Lambda_{a_1,a_2,a_3}^{(3)}$ that are the flow out of the light-like
directions of $N^*K_{a_1,a_2,a_3}$.
     
   \HOX{We need to check here arguments!!}

   The above analysis implies that for generic manifold
   and for generic sources singsupp$(\M^{(3)}) $ is equal to 
   \ba
   S_3=U_{\hat g}\cap (\bigcup_{(\vec x_j,\vec \xi_j)=({\hat x}_{a_\ell},\zeta_{a_\ell}),\ \ell\leq L}
    \mathcal Y((\vec x,\vec \xi),t_0,s_0)).
   \ea
   In the generic
   case we can thus assume the set $S_3$ has the above form and that the given set
  singsupp$(\M^{(3)}) $ determines it.
  
\subsubsection{Analysis of the 4th order terms}

Let us  next consider singularities of $\M^{(4)}$ outside the set $S_1\cup S_2\cup S_3$.
Outside
these surfaces we look singularities of  $\M^{(4)}$ using the indicator functions  $\Theta^{(4)}_\tau$,  constructed using gaussian
beams.
Assume that $x_5\in U_{\hat g}\setminus (S_1\cup S_2\cup S_3)$.
As pointed out before we need to consider general (non-injective) maps $\sigma:\{1,2,3,4\}\to \{1,2,3,4\}$.
The case when $\sigma$ is a permutation is considered in the main text.
Let us thus assume that $\sigma$ is not a permutation.

We need to consider cases where the linearized waves $u_1,u_2,u_3$ are replaced by
$u_{r_1},u_{r_2},u_{r_3}$, where $r_1,r_2,r_3\in \{1,2,3,4\}$ are arbitrary,
and $u_4$, kept the same (changing it would be just
renumeration of indexes.) Also, we can assume that $\{r_1,r_2,r_3\}$
is not the set  $\{1,2,3\}$, as this case we have already analyzed.

We start by  considering the $\tilde T$-functions.
We need to consider terms similar to $\tilde T$ in the main text containing e.g.\
$\bra u_\tau,{\bf Q}(u_{r_1}u_{r_2})\cdotp {\bf Q}(u_{r_3}u_{r_4})\cet$,
when some $r_j$ are the same. Assume next that
e.g. $r_1$ and $r_3$ are the same.
We use that fact that  $v\in \I(\Lambda_1, \Lambda_2)$,
we have
WF$(v)\subset \Lambda_1\cup \Lambda_2$. Moreover,
by \cite[Thm.\ 1.3.6]{Duistermaat},
\ba
\hbox{WF}(v\,\cdotp w)\subset 
\hbox{WF}(v)\cup \hbox{WF}(w)
\cup\{(x,\xi+\eta);\ (x,\xi)\in \hbox{WF}(v),\ (x,\eta)\in \hbox{WF}(w)\}
\ea
and thus we see that
the wave front set of term ${\bf Q}(u_{r_1}u_{r_2})\cdotp {\bf Q}(u_{r_1}u_{r_4})$
is in the union of the sets $N^*K_{r_j}$,
$N^*(K_{r_1}\cap K_{r_2})$, and 
$N^*(K_{r_1}\cap K_{r_2}\cap K_{r_4})$. The elements in
the first two sets do not propagate in the bicharacteristic flow
and the last one propagates along $\mathcal Y((\vec x,\vec \xi),t_0,s_0)$.
Observe that in the generic case $\mathcal Y((\vec x,\vec \xi),t_0,s_0)\subset S_3$.
The other  $\tilde T$-functions can be analyzed similarly.
Thus in the generic case the  $\tilde T$-functions do not cause observable 
singularities in the complement of $S_3$.

Consider next the $T$-functions.
We need to consider terms similar to $T$ in the main text containing e.g.\
$\bra u_\tau,u_{r_1}\,\cdotp {\bf Q}(u_{r_2}\,\cdotp {\bf Q}(u_{r_3}\,\cdotp u_{r_4}))\cet$,
when some $r_j$ are the same. 
 Again we decompose ${\bf Q}={\bf Q}_1+{\bf Q}_2$ and consider cases $p=1,2$ as in the main text. In the case $p=2$ we need to consider
critical points of modified phase functions 
\ba
\nonumber \Psi^{\vec r}_2(z,y,\theta)
&=&\theta_1z^{r_1}+\theta_2z^{r_2}+\theta_3z^{r_3}+\theta_4y^4+\varphi(y).
\ea
Under the assumption that $x_5\not \in K_4$, we see similarly
as in the main text that there are no critical points. Now
\ba
d_{z,y,\theta}\Psi^{\vec r}_2&=&(\omega^{\vec r};
r+d_y\psi_4(y,\theta_4),
z^{r_1},z^{r_2},z^{r_3},d_{ \theta_4}\psi_4(y,\theta_4))
\ea
where $\omega^{\vec r}=(\omega^{\vec r}_j)_{j=1}^4$,
$\omega^{\vec r}_j=\sum_{n=1}^3\theta_n\delta_{n\in R(j)}$, $R(j)=\{n\in \{1,2,3\}; r_n=j\}.$
Note that $(z,\omega^{\vec r})\in \bigcap_{n\in \{r_1,r_2,r_3\}}N^*K_{n}$ and some of $r_j$
may be the same. 
We see as in the main text that as $x_5\not \in {\mathcal Y},$
alternative (A1) can not hold.
Note that if (A2) holds then $\omega^{\vec r}\in \X$, and we see as in 
the main text that this is not possible.  Thus the term with $p=2$ causes no observable singularities
near $x_5$.

In the case $p=1$, let $(z,\theta,y,\xi)$ be a critical point of
 the phase function
  \ba
 \Psi^{\vec r}_3(z,\theta,y,\xi)=\theta_1z^{r_1}+\theta_2z^{r_2}+\theta_3z^{r_3}+(y-z)\,\cdotp \xi+\theta_4y^4+\varphi(y).
 \ea
Then
 \beq\label{eq: critical points 1}
  & &\p_{\theta_j}\Psi_3=0,\ j=1,2,3\quad\hbox{yield}\quad z\in K_{{r_1}}\cap K_{r_2}\cap K_{r_3},\\
  \nonumber
& &\p_{\theta_4}\Psi_3=0\quad\hbox{yields}\quad y\in K_4,\\
  \nonumber
 & &\p_z\Psi_3=0\quad\hbox{yields}\quad\xi=\omega^{\vec r},\\
   \nonumber
& &\p_{\xi}\Psi_3=0\quad\hbox{yields}\quad y=z,\\
  \nonumber
 & &\p_{y}\Psi_3=0\quad\hbox{yields}\quad \xi=
 -\p_y\varphi(y)-\theta_4 w.
  \eeq
  The critical points we need to consider for the asymptotics satisfy also
\beq\label{eq: imag.condition copu} \\ \nonumber
\im \varphi(y)=0,\quad\hbox{so that }y\in \gamma_{x_5,\xi_5},\ \im d\varphi(y)=0,\
\re d\varphi(y)\in L^{*,+}_y\hattuM _0.\hspace{-1cm}
\eeq
Now, as we assume that $\{r_1,r_2,r_3\}$
is not the set  $\{1,2,3\}$, some of $r_j$ is equal to $4$ or two of these are the same.   For simplicity, assume that
$r_1$ is either 4 or alternatively, $r_1$ coincides with $r_2$ or $r_3$. As  $w\in N^*K_4$, these imply that 
$\p_y\varphi(y)\in N^*K_{r_2}\cap  N^*K_{r_3}\cap N^*K_4\subset \mathcal X((\vec x,\vec\xi);t_0)$.
Hence the assumption that $x_5\not \in  {\mathcal Y}$ implies that
terms with $p=1$ do not cause any singularities near $x_5$. Thus only the singularities of
$T$-type terms for which $\sigma$ is bijection,
that were analyze in the main text, are visible in $U_{\hat g}\setminus (S_1\cup S_2\cup S_3)$.
The singularities could be analyzed  in much more careful way
if we could analyze products of $v\in \I(K_1,K_{13})$ and $w\in \I(K_2,K_{23})$
and show that such product is a sum of conormal distributions.\hiddenfootnote{HERE ARE SOME EXTRA COMPUTATIONS:
By applying method of stationary phase in $\xi$ and $z$ variable we obtain
a integral similar to (\ref{eq; T tau asympt modified}),
 \beq \label{eq; T tau asympt modified2}
& &T_{\tau,1,1}^{(4),\beta}
=c\tau^{4}
\int_{\R^{8}}e^{i(\theta_1y^{r_1}+\theta_2y^{r_1}+\theta_3y^{r_1}+\theta_4y^{r_1})+i\tau \varphi(y)}
c_1(y, \theta_1, \theta_2)\cdotp\\ \nonumber
& &
\cdotp a_3(y, \theta_3)
q_{1,1}(y,y,\omega^{\vec r}_\beta(\theta))a_4(y,\theta_4) a_5(y,\tau )\,d\theta_1d\theta_2d\theta_3d\theta_4dy.
\eeq
where $\omega^{\vec r}_\beta$ is similar to $\omega^{\vec r}$  is defined with 
$R_\beta(j)=\{n\in \{1,2,3\}; r_n=\sigma_\beta(j)\}$ is  similar to $R(j)$.

 Similarly, to the analysis related to the function $\chi(\theta)$ that we vanish
 in a neighborhood of the set $\mathcal A_q$ given in (\ref{set Yq}),
 we can define a fucntion 
 $\tilde \chi(\theta)\in C^\infty(\R^4)$ that 
vanishes in a $\e_3$-neighborhood $\W$ (in the $\hat g^+$ metric) of
$\tilde {\mathcal A_q}$,
\ba
\tilde {\mathcal A_q}&:=&
\bigcup_{1\leq j<k\leq 4}N_q^*K_{jk}.
\ea
The idea of the sets $N_q^*K_{jk}$ is that in the complement of  their  neighborhoods at least three of the coordinates $(\theta_1,\theta_2,\theta_3,\theta_4)$
 are comparable to $|\vec\theta|$, that is, there is $j_0$ such that $|\theta_k|\geq c_j|\vec \theta|$ for $k\not =j_0$.
 Indeed, e.g.\ $(q;(\theta_1,\theta_2,\theta_3,\theta_4))\not \in N_q^*K_{jk}$ implies 
 that $(\theta_j,\theta_k)\not =0$. As this holds for all $(j,k)$, we see
 that $\theta\not \in \W$ and $|\theta|_{\R^4}=1$, only one coordinate of  $\theta$ can vanish. This makes the symbol
$\tilde \chi(\theta)b(y,\theta)$ a sum of a product type symbols $S^{p,l}(M_0;\R^3\times \R\setminus \{0\})$.

Let  \ba\hbox{$\tilde b_0(y,\theta)=\phi(y)\tilde \chi(\theta)b(y,\theta)$}\ea be a  product type symbol determining 
a  distribution 
$
\tilde \F^{(4),0}$ 
 that is given by the formula (\ref{f4-formula}) with
$b(y,\theta)$ being replaced by 
 $b_0(y,\theta)$.
We  can consider
 $\tilde \F^{(4),0}$ as a sum of in three terms: 
 \ba
 \tilde \F^{(4),0}=\tilde \F^{(4),0,1}+\tilde \F^{(4),0,2}+\tilde \F^{(4),0,2}\ea
  where \HOX{Check 3/2 HERE!!}
 \ba
 \tilde \F^{(4),0,1}&\in& \sum_{ \{k_1,k_2,k_3,k_4\}=\{1,2,3,4\}} \I^{p(\vec k)-4+1/2}(K_{k_2,k_3,k_4})\cap  \I(N^*K_{k_2,k_3,k_4}\setminus N^*\{q\})\\
& &\subset 
\sum_{ \{k_1,k_2,k_3,k_4\}=\{1,2,3,4\}} \I^{p(\vec k)-4+1/2-3/2+1}(N^* K_{k_2,k_3,k_4})\
,
\ea
where $p(\vec k)=p_{k_2}+p_{k_3}+p_{k_4}$ and
\ba
\tilde \F^{(4),0,2}&\in& \I^{p-4}(\{q\}\cap  \I(N^*\{q\}\setminus (\cup_{\vec k} N^*K_{k_2,k_3,k_4}))\\
& &\subset \I^{p-4-2+1}(N^*\{q\})
\ea
 and  $\tilde \F^{(4),0,2}$ is microlocally supported in some (arbitrarily) small
 neighborhood $\W$ of  $\cup (N^*K_{k_2,k_3,k_4}\cap N^*\{q\})$. 
  Let $\Lambda^{(3)}$ be the flowout of $\cup_{\vec k} N^*K_{k_2,k_3,k_4}$.
We see that $\M^{(4)}$ is Lagrangian distribution
 with known order in $I^-(x_6)\setminus  {\mathcal Y}((\vec x,\vec \xi),t_0)$, that is outside
 the set $\cap_{s_0>0}{\mathcal Y}((\vec x,\vec \xi),t_0,s_0)$ that has the Hausdorff dimension
 two. Indeed,\HOX{Check 3/2}
\ba
& &{\bf Q}_{\hat g}\tilde \F^{(4),0,1}\in
\sum_{ \{k_1,k_2,k_3,k_4\}=\{1,2,3,4\}} \I^{p(\vec k)-4-3/2}
\I(U_{\hat g}\setminus K_{k_2,k_3,k_4};
\Lambda^{(3)}),\\
& &{\bf Q}_{\hat g}\tilde \F^{(4),0,2}\in 
\I^{p-5-3/2}(U_{\hat g}\setminus\{q\};\Lambda_q),\\
  & &\hbox{WF}({\bf Q}_{\hat g}\tilde \F^{(4),0,3})\subset  \V .\ea
where $\V\subset T^*M_0$ is a neighborhood of $ 
T^*\mathcal Y((\vec x,\vec \xi),t_0)$ that is
the flow out of the set $\W$ in the bicharacteristic flow.
This makes it possible to analyze $
 \bra u_\tau,\F_1^{(4)}\cet$ and see that ${\bf Q}\F_1^{(4)}$ is Lagrangian distribution
 in $U_{\hat g}\setminus {\mathcal Y}((\vec x,\vec\xi),t_0)$.

To study the singularities of $u^{(4)}_\e(x)$, we consider next interaction of four distorted plane  waves, 
corresponding the $\vec a=(a_\ell)_{\ell=1}^4 \in \A$, $a_\ell\leq N$ that are singular on the surfaces
$K_{a_\ell}=\L(({\hat x}_{a_\ell},\xi_{a_\ell});\rho)$.
  Assume now that $\bigcap_{\ell=1}^4 \K_{a_\ell}$ contains a point $q$ and
  let $K_{a_1,a_2}=\bigcap_{\ell=1}^2 \K_{a_\ell}$ and $K_{a_1,a_2,a_3}=\bigcap_{\ell=1}^3 \K_{a_\ell}$ . 
 Then, we see that   $\Lambda_q^+$ and
  the flowout $\Lambda_{a_1,a_2,a_3}$ of $N^*K_{a_1,a_2,a_3}$ in the canonical
  relation of $\square_g$ are in $T^*U_g$ four dimensional Lagrangian
  manifolds that intersect only on a three dimensional set.
  Assume now that the source ${\bf f}_{a_\ell}$ have such orders 
  that $u_{a_\ell}={\bf Q}_g{\bf f}_{a_\ell}\in \I^{p_{a_\ell}}(K_{a_\ell})$ outside support of ${\bf f}_{a_\ell}$.

The linearized term $\M^{(1)}$ in (\ref{nl wave}) has the wave front set that is
a subset of  $S_1=\bigcup_{a\leq N}N^*K_a$,
the term  $\M^{(2)}$  has wave front set is a subset of  $S_1\cup S_2$ where  $S_2=\bigcup_{a_1<a_2\leq N}N^*K_{a_1,a_2}$. Let us note that when analyze the term $\M^4$, we have
to consider the term $\M^{(1)}\,\cdotp\M^{(1)}$ containing terms
$u_{a_1}\,\cdotp u_{a_1}$, where $u_{a_1}\in \I^{p}(K_{a_1})$. 
Denote $b(x^1)=u_{a_1}(x^1,x^\prime)$ and assume that $\hat b(\xi_1)\sim c_0 |\xi_1|^p$ as $\xi_1\to \infty$,
we see that 

\ba
(\hat b*\hat b)(\xi_1)=\int_\R \hat b(\xi_1-t)\hat b(t)dt
\sim 
|\xi_1|^p
\int_\R c_0(1-|\xi_1|^{-1}t)^p\hat b(t)dt
\sim c_0 b(0)|\xi_1|^p,
\ea

****GENERALIZATION***

Let  $B_1\in \I^{p_1}(\{0\})$ and $B_2\in \I^{p_2}(\{0\})$,
Let us write symbol $b_1(x,\xi)$ of $B_1$ as a sum of positive homogeneous term and
a compactly $\xi$-supported term $b_1(x,\xi)=b_1^h(x,\xi)+b_1^c(x,\xi)$ and similar
lower order terms.
Then $B_1^c,B_2^c\in C^\infty$, and
 $B_1^h\,\cdotp B_2^c\in \I^{p_1}(\{0\})$,  $B_1^c\,\cdotp B_2^h\in \I^{p_2}(\{0\})$,
  $B_1^c\,\cdotp B_2^c\in C^\infty$, and symbol of  $B_1^h\,\cdotp B_2^h$ is 
  obtained via substitution $t=|\xi_1|s$,
\ba
(b_1*b_2)(x,\xi_1)&=&\int_\R b_1(x,\xi_1-t) b_2(x,t)dt\\
&=&|\xi_1|^{p_1+p_2+1}
\int_\R b_1(x,1-s) b_2(x,s)ds
\ea
and thus  $B_1^h\,\cdotp B_2^h\in \I^{p_1+p_2+1}(\{0\})$.
Hence,  $B_1\,\cdotp B_2\in \I^{r}(\{0\})$, $r=\max(p_1,p_2)$.

****GENER. 2***

Let  $B_1\in \I^{p_1}(\{0\})$ and $B_2\in \I^{p_2}(\{0\})$,
Let us estimate the symbol $b_1(x,\xi)$ of $B_1$
by
\ba
|b_1(x,\xi)|\leq b_1^h(x,\xi)+b_1^c(x,\xi)
\ea 
where $b_1^c(x,\xi)$ is a compactly $\xi$-supported  symbol and
$b_1^h(x,\xi)$ is $p_1$-homogeneous. Let $B_1^c$ be the conormal distribution
with symbol $b_1^c(x,\xi)$ and  $B_1^h=B_1-B_1^c$, and introduce similar notations for $B_2$.
Then $B_1^c,B_2^c\in C^\infty$, and
 $B_1^h\,\cdotp B_2^c\in \I^{p_1}(\{0\})$,  $B_1^c\,\cdotp B_2^h\in \I^{p_2}(\{0\})$,
  $B_1^c\,\cdotp B_2^c\in C^\infty$, and symbol $c^h$ of  $B_1^h\,\cdotp B_2^h$ can be estimated
using substitution $t=|\xi_1|s$,
\ba
|c^h(x,\xi_1)|&\leq &\int_\R b_1(x,\xi_1-t) b_2(x,t)dt\\
&\leq &|\xi_1|^{p_1+p_2+1}
\int_\R b_1(x,1-s) b_2(x,s)ds
\ea
and thus  $B_1^h\,\cdotp B_2^h\in \I^{p_1+p_2+1}(\{0\})$.
Hence,  $B_1\,\cdotp B_2\in \I^{r}(\{0\})$, $r=\max(p_1,p_2)$.

****GENER. 3***

We need also to analyze the products of 
  $B_1\in \I^{p_1,l_1}(K_1,K_{13})$ and $B_2\in \I^{p_2,l_2}(K_2,K_{13})$.
  When we analyze this in the case when $|\theta_2|>c|\theta|$ and $|\theta_3|>c|\theta|$ 
  and  $|\theta_1|>c|\theta|$, this is equivalent to analyzing product of 
  $ \I(K_{13})$ and  $ \I(K_{12})$.
  The case where only  $|\theta_2|>c|\theta|$ and $|\theta_3|>c|\theta|$ but $\theta_1$ may be
  small  is difficult. Note that it is possible that we do not need to analyze
  this to understand light-like directions on $N^*K_{123}$. Maybe we should
  first find $S_3$ from $\e^3$ terms.

***

Plan: 

1. It seems likely that in generic case we can see $S_1$, $S_2$ and $S_3$ from different $\e$-order of singularities
The difficulty why we need to consider different $\e$ order is that we need
to consider  the product of distributions 
  $v_1\in \I(K_{13})$ and  $ v_1\in\I(K_{12})$ and we can not find the order of ${\bf Q}(v_1v_2)$.

2. Outside the surfaces  $S_1$, $S_2$ and $S_3$ we can detect $S_4$ and can separated different components.
For this we use the product theorem for the wave front sets of products
of   $v_1\in \I(K_{13})$ and  $ v_1\in\I(K_{12})$ to see that the wave front set of ${\bf Q}(v_1v_2)$ is contained
in $\Lambda^{(3)}$. 

***

and thus $u_{a_1}\,\cdotp u_{a_1}\in \I^{p}(K_{a_1})$ has a non-vanishing principal symbol
at $(x,\xi)\in N^*K_{a_1}$ if and only $u_{a_1}$ has a non-vanishing principal symbol at $(x,\xi)$ and 
a non-vanishing value at $x$.

*** ABOVE: We had problem: We can not analyze terms $\M^{(2)}\cdotp \M^{(2)}$
appearing in $\M_4$. The problem can be solved by considering solutions
outside Hausdroff 2-dimensional set $\Z$.
}
\footnote{We could use Piriou,	A.: Calcul	symbolique non	 lineare pour une onde conormale simple. Ann. Inst. Four. 38, 173-188 (1988), or \cite{GU1}, Lemma 1.3. Look also
Bougrini, H.; Piriou, A.; Varenne, J. P. Propagation et interaction des symboles principaux pour les ondes conormales semi-lin'eaires. (French) [Propagation and interaction of principal symbols for semilinear conormal waves] Comm. Partial Differential Equations 23 (1998), no. 1-2, 333--370.}

\subsection{Separation of singularities from different point sources}
Summarizing, above we have seen that in the generic case  
the terms  $\M^{(j)}$, $j\leq 3$   have wave front sets which union is $S_1\cup S_2\cup S_3$ 
and,  $\M^{(4)}$  has wave front set is a subset of   $S_1\cup S_2\cup S_3\cup S_4$ where  $S_4=\bigcup_{a_1<a_2<a_3<a_4\leq N}\Lambda^+_{q(a_1,a_2,a_3,a_4)}$.
Now we notice that for a generic  (i.e.\ in an open and dense set) $((M,g),(x_0,\xi_0))\in \mathcal N$,
we see that when $\vec a$  and $\vec b$ run over any finite subset of  $\A$,
\HOX{Check the the claims on "generic manifolds"}
the sets $\Lambda^+_{q(b_1,b_2,b_3,b_4)}$ intersect
the sets $N^*K_a$, $N^*K_{a_1,a_2}$, $\Lambda_{a_1,a_2,a_3}^{(3)}$, and 
the other sets 
$\Lambda^+_{q(a_1,a_2,a_3,a_4)}$, 
on at most 3-dimensional sets. Indeed,
 by making a small perturbation to the metric in a small neighborhood of
  $q(a_1,a_2,a_3,a_4)$ the sets  $\Lambda^+_{q(a_1,a_2,a_3,a_4)}$ can be perturbed so
  that the above non-intersection condition becomes valid, and the non-intersection condition
is clearly valid in an open set of metric tensors.
Thus the set of manifolds for which the stated condition holds is  open and dense.
Let us then choose the sources $\tilde {\bf f}_{a_\ell}$ to have generic polarizations.
We see that  in a generic set of polarizations the principal symbol
 of the 4th order interaction term $\M^{(4)}$ is non-vanishing on an open and dense set
 of $S_4$. Recall also that   in a generic situation the singular support of $\M^{(3)}$
 determines $S_3$.

 Thus in a generic situation
we can observe an open and dense subset of the surface $S_4$. Let us now consider how
we can decide which points on the union  of the surfaces 
$\Lambda^+_{q(a_1,a_2,a_3,a_4)}$ belong to the same surface.
Recall that the order of singularity (at points that do not belong to the 3-dimensional intersections) \HOX{Check 3/2}
on $\Lambda^+_{q(a_1,a_2,a_3,a_4)}$ is $(p_{a_1}+p_{a_2}+p_{a_3}+p_{a_4})-4-3/2$, see (\ref{eq: conormal}).
Let us now assume that $p_j$, $j\leq N$ are such that all the numbers
\ba
& &(p_{a_1}+p_{a_2}+p_{a_3}+p_{a_4})-1/2-n_1,
\ea
are different,
where  $n_1$ and $n_2$ are arbitrary integers (corresponding to different orders of the
classical symbols of terms) and \HOX{Gunther, do you know if the symbols of of waves $\M_j$ classical, assuming
that the sources have classical symbols, and if so, what would be a reference for this?}
the indexes $a_1,a_2,a_3,a_4\leq N$
satisfy $a_j\not =a_k$ for $j\not=k$. Let $\Sigma$ be the set of points $(x,\xi)\in T^*U_g$, $x\not \in S_3$ on the 
wave front set of  $\M^{(4)}$ such that the point $x$
has a neighborhood $V$ where $\M^{(4)}$  is a conormal
distribution associated smooth surface $S\subset V$ such that
$\M^{(4)}|_V\in \I^p(N^*S)$, where $p=(p_{a_1}+p_{a_2}+p_{a_3}+p_{a_4})-4-3/2$ \HOX{Check 3/2}
but $p$ can not be replaced by any smaller number. Then for
generic set of manifolds  in $\mathcal N$
the closure of the set $\Sigma$ coincides with 
 $\Lambda^+_{q(a_1,a_2,a_3,a_4)}$. \HOX{Check the construction of  $\Lambda^+_{q(a_1,a_2,a_3,a_4)}$}
This implies
that we can determine $\Lambda^+_{q(a_1,a_2,a_3,a_4)}\cap T^*U_g$ for  all
$\vec a=(a_1,a_2,a_3,a_4)$ and thus
we can find the whole smooth surface $\Lambda^+_{q(a_1,a_2,a_3,a_4)}\cap T^*U_g$.
This implies that we can find  
$G_N(J((M,g),(x_0,\xi_0)),[g],x_0,\xi_0)$. This determines the manifold
$(J((M,g),(x_0,\xi_0)),[g])$ up to a small error.

   Summarizing the above  sketch of construction
   means the following:
   \medskip
   
   {\bf Summary:} {\it For any $\e>0$ one can $L$ such that the following holds
   in a generic set of manifolds
    $((M,g),(x_0,\xi_0))\in \mathcal N$.
    Let   $({\hat x}_{a_\ell},\xi_{a_\ell})$, $|\ell |\leq L$ be constructed above and
    $F_\e$ be a source (\ref{Fe source}) that sends distorted plane  waves to directions
       $({\hat x}_{a_\ell},\xi_{a_\ell})$, $|\ell |\leq L $ with generic polarizations  $w_{a_\ell}$  and
       satisfies (\ref{eq: V ehto}). Then,
       by measuring singularities of terms 
   $\M^{(j)}$, $j\leq 4$ of the     
       the wave  $u^{(4)}_\e(x)$ on $U_g$
    produced by the source  $F_\e$,
   we can find the manifold  
    $J_{g}(\hat \mu(s_-),\hat \mu(s_+))$ and the conformal type 
    of the metric $[g]$ on it up to an error $\e$.}
    \medskip
   
 In other words,  when we consider the above measurement
 with single source and approximate the $O(\e^5)$ terms by zero,
 the observations of the singularities of the terms with different $\e^j$ magnitudes $j\leq 4$, 
 make it possible in a generic case to construct
a  discrete approximation of the earliest light observation point
sets, $\{\mathcal P_V(q_k);\ k\leq K\}.$
Such data, under appropriate conditions, could be used to determine a discrete approximation 
for the conformal structure of the space-time in an analogous manner used
in \cite{AKKLT,K2L}. However, making the approach described above to a detailed proof is outside
the scope of this paper  will be considered
elsewhere.}
\medskip
}}

\hiddenfootnote{
\section*{Appendix A: Reduced Einstein equation}

In this section we review known results on Einstein equations and wave maps.

\HOX{The appendices need to be shortened  in the final version of the paper.}

\subsection*{A.1.\ Summary of the used notations}


Let us recall some definitions given in Introduction, in the Subsection
\ref{subsec: Gloabal hyperbolicity}.
Let $(\hattuM ,\hat g)$ be a $C^\infty$-smooth globally hyperbolic Lorentzian
manifold  and  $\tilde g$ be a $C^\infty$-smooth globally hyperbolic
metric on $M$ such that $\hat g<\tilde g$. 

Recall that there is an isometry $\Phi:(M,\tilde g)\to (\R\times N,\tilde  h)$,
where $N$
is a 3-dimensional manifold and the metric $\tilde  h$ can be written as
$\tilde  h=-\beta(t,y) dt ^2+\overline h(t,y)$ where $\beta:\R\times  N\to (0,\infty)$ is a smooth function and 
$\overline h(t,\cdotp)$ is a Riemannian metric on $ N$ depending smoothly
on $t\in \R$. 
As in the main text, we identify these isometric manifolds and denote $\hattuM =\R\times N$.
%
Also, for $t\in \R$, recall that $\hattuM (t)=(-\infty,t)\times N$. We use parameters $t_1>t_0>0$
and denote $\hattuM _j=\hattuM (t_j)$, $j\in \{0,1\}$.
We use  the time-like geodesic    $\hat \mu =\mu_{\hat g}$, $\mu_{\hat g}:[-1,1]\to M_0$ on $(M_0,\hat g)$
and  the set $\K_j:=J^+_{\tilde g}({p^-})\cap M_j$ with ${p^-}=
\hat \mu(s_-)\in (0,t_0)\times N$, $s_-\in (-1,1)$
and recall that  
$J^+_{\tilde g}({p^-})\cap M_j$ is compact and 
there exists $\e_0>0$ such that if $\|g-\hat g\|_{C^0_b(\hattuM _1;\hat g^+)}<\e_0$, then
$g|_{\K_1}<\tilde g|_{\K_1}$,
and in particular, we have  $J^+_{g}(p)\cap \hattuM _1\subset \K_1$
for all $p\in \K_1$.

Let us use local coordinates on $\hattuM _1$ and denote by
$\nabla_k=\nabla_{X_k}$  the covariant derivative with respect to the metric $g$
to the direction $X_k=\frac \p{\p x^p}$
 and by
 $\hat \nabla_k=\hat \nabla_{X_k}$ the covariant derivative with respect to the metric $\hat g$
to the direction $X_k$.

\subsection*{A.2.\ Reduced Ricci and Einstein tensors}
Following \cite{FM} we recall that 
\beq\label{q-formula2copy}
\Ric_{\mu\nu}(g)&=& \Ric_{\mu\nu}^{(h)}(g)
+\frac 12 (g_{\mu q}\frac{\p \Gamma^q}{\p x^{\nu}}+g_{\nu q}\frac{\p \Gamma^q}{\p x^{\mu}})
\eeq
where  $\Gamma^q=g^{mn}\Gamma^q_{mn}$,
\beq\label{q-formula2copyB}
& &\hspace{-1cm}\Ric_{\mu\nu}^{(h)}(g)=
-\frac 12 g^{pq}\frac{\p^2 g_{\mu\nu}}{\p x^p\p x^q}+ P_{\mu\nu},
\\ \nonumber
& &\hspace{-2cm}P_{\mu\nu}=  
g^{ab}g_{ps}\Gamma^p_{\mu b} \Gamma^s_{\nu a}+ 
\frac 12(\frac{\p g_{\mu\nu }}{\p x^a}\Gamma^a  
+ \nonumber
g_{\nu l}  \Gamma^l _{ab}g^{a q}g^{bd}  \frac{\p g_{qd}}{\p x^\mu}+
g_{\mu l} \Gamma^l _{ab}g^{a q}g^{bd}  \frac{\p g_{qd}}{\p x^\nu}).\hspace{-2cm}
\eeq
Note that $P_{\mu\nu}$ is a polynomial of $g_{jk}$ and $g^{jk}$ and first derivatives of $g_{jk}$.
The harmonic  Einstein tensor is  
\beq\label{harmonic Ein}
\Ein^{(h)}_{jk}(g)=
\Ric_{jk}^{(h)}(g)-\frac 12 g^{pq}\Ric_{pq}^{(h)}(g)\, g_{jk}.
\eeq 
The  harmonic  Einstein tensor is extensively used to study
Einstein equations in local coordinates where one can use the Minkowski
space $\R^4$ as the background space. To do global constructions
with a background space $(M,\hat g)$ one uses  the  reduced Einstein tensor.
The $\hat g$-reduced Einstein tensor 
 $\Ein_{\hat g} (g)$ and the $\hat g$-reduced  Ricci tensor
  $\Ric_{\hat g} (g)$  
are given 
by
\beq\label{Reduced Einstein tensor}
& &(\Ein_{\hat g} (g))_{pq}=(\Ric_{\hat g} (g))_{pq}-\frac 12 (g^{jk}(\Ric_{\hat g} g)_{jk})g_{pq},\\& &
\label{Reduced Ric tensor}
(\Ric_{\hat g} (g))_{pq}=\Ric_{pq} g-\frac 12 (g_{pn} \hat \nabla _q\hat F^n+ g_{qn} \hat \nabla _p\hat F^n)
\eeq
where $\hat F^n$ are the harmonicity functions given by
\beq\label{Harmonicity condition BBB}
\hat F^n=\Gamma^n-\hat \Gamma^n,\quad\hbox{where }
\Gamma^n=g^{jk}\Gamma^n_{jk},\quad
\hat \Gamma^n=g^{jk}\hat \Gamma^n_{jk},
\eeq
where $\Gamma^n_{jk}$ and $\hat \Gamma^n_{jk}$ are the Christoffel symbols
for $g$ and $\hat g$, correspondingly.
Note that $\hat \Gamma^n$ depends also on $g^{jk}$.
As $\Gamma^n_{jk}-\hat \Gamma^n_{jk}$ is the difference of two connection
coefficients, it is a tensor. Thus $\hat F^n$ is tensor (actually, a vector field), implying that both
$(\Ric_{\hat g} (g))_{jk}$ and $(\Ein_{\hat g} (g))_{jk}$ are  2-covariant tensors.
Observe that  the $\hat g$-reduced Einstein tensor is sum of the harmonic Einstein tensor 
and a term that is a zeroth order in $g$, 
\beq\label{hat g reduced einstein and reduced einstein}
(\Ein_{\hat g} (g))_{\mu\nu}=\Ein^{(h)}_{\mu\nu} (g)+\frac 12 (g_{\mu q}\frac{\p \hat \Gamma^q}{\p x^{\nu}}+g_{\nu q}\frac{\p \hat \Gamma^q}{\p x^{\mu}}).
\eeq

%

%

%

\subsection*{A.3.\ Wave maps and reduced Einstein equations}
Let us consider the manifold $M_1=(-\infty,t_1)\times N$ with a $C^m $-smooth metric $g^{\prime}$,
$m\geq 8$,
 which
is a perturbation of the metric $\hat g$ and satisfies the Einstein equation
\beq\label{Einstein on M_0}
\Ein(g^{\prime})=T^{\prime}\quad \hbox{on }M_1,
\eeq
or equivalently, 
\ba
\Ric(g^{\prime})=\rho^{\prime},\quad \rho^{\prime}_{jk}=T_{jk}^{\prime}-\frac 12 ((g^{\prime})^{  nm}T^{\prime}_{nm})g^{\prime}_{jk}\quad \hbox{on }M_1.
\ea
Assume also that  $g^{\prime}=\hat g$ in the domain $A$,
where $A=M_1\setminus\K_1$ and  $\|g^{\prime}-\hat g\|_{C^2_b(M_1,\hat g^+)}<\e_0$,
so that $(M_1,g^{\prime})$ is globally hyperbolic. Note that then $T^{\prime}=\hat T$ in the set $A$ 
and that the metric $g^{\prime}$ coincides with $\hat g$
in particular in the set $M^-=\R_-\times N$

We recall next the considerations of \cite{ChBook}.
Let us  consider the Cauchy problem for the wave map
$f:(M_1,g^\prime)\to (M,\hat g)$, namely
\beq
\label{C-problem 1}& &\square_{g^{\prime},\hat g} f=0\quad\hbox{in } M_1,\\
\label{C-problem 2}& &f=Id,\quad \hbox{in  }\R_-\times N,
\eeq
where $M_1=(-\infty,t_1)\times N\subset M$. In (\ref{C-problem 1}), 
$\square_{g^{\prime},\hat g} f=g^\prime\,\cdotp\hat \nabla^2 f$ is the wave map operator, where $\hat \nabla$
is the covariant derivative of a map $(M_1,g^\prime)\to (M,\hat g)$, see \cite[formula (7.32)]{ChBook}.
In  local coordinates
$X:V\to \R^4$ of $V\subset M_1$, denoted
$X(z)=(x^j(z))_{j=1}^4$ and $Y:W\to \R^4$ of $W\subset \hattuM $, denoted
$Y(z)=(y^A(z))_{A=1}^4$, 
the wave map $f:M_1\to \hattuM $ has representation $Y(f(X^{-1}(x)))=(f^A(x))_{A=1}^4$ and
the wave map operator in equation (\ref{C-problem 1}) is given by
\beq\label{wave maps2}
& &(\square_{g^{\prime},\hat g} f)^A(x)=
(g^{\prime})^{  jk}(x)\bigg(\frac \p{\p x^j}\frac \p{\p x^k} f^A(x) -\Gamma^{{\prime} n}_{jk}(x)\frac \p{\p x^n}f^A(x)
\\
& &\quad \quad\quad\quad\quad\quad\quad\quad\nonumber
+\hat \Gamma^A_{BC}(f(x))
\,\frac \p{\p x^j}f^B(x)\,\frac \p{\p x^k}f^C(x)\bigg)\eeq
where $\hat \Gamma^A_{BC}$ denotes the Christoffel symbols of metric $\hat g$
and  $\Gamma^{ {\prime} j}_{kl}$ are the Christoffel symbols of metric $g^{\prime}$.
When (\ref{C-problem 1})
is satisfied, we say that  $f$ is wave map with respect to the pair $(g^{\prime},\hat g)$.
The important property of the wave maps is that, if 
$f$ is wave map with respect to the pair $(g^{\prime},\hat g)$ and
$g=f_*g^{\prime}$ 
then, as follows from (\ref{wave maps2}), the identity map $Id:x\mapsto x$ is a wave map with respect to the pair $(g,\hat g)$
and, the wave map equation for the identity map is equivalent to (cf.\ \cite[p.\ 162]{ChBook})
\beq\label{wave equation id}
\Gamma^n=\hat\Gamma^n,\quad \hbox{where }\Gamma^n=g^{jk}\Gamma^n_{jk},
\quad \hat \Gamma^n=g^{jk}\hat\Gamma^n_{jk}
\eeq
%
where the Christoffel symbols $\hat\Gamma^n_{jk}$ of the metric $\hat g$ are smooth functions.

As   $g=g^{\prime}$ outside a compact set $\K_1\subset (0,t_1)\times N$,
we see that   this Cauchy problem is equivalent to the same equation
restricted to the
set $(-\infty,t_1)\times B_0$, where $B_0\subset N$ is such an open
relatively compact set that $\K_1\subset (0,t_1]\times B_0$
with the boundary condition $f=Id$ on  $(0,t_1]\times \p B_0$.
Then using results of \cite {HKM},  that can be applied for equations on manifold  as is
done  in Appendix B, and combined with the Sobolev embedding theorem, 
we see\footnote{See also:  
Thm.\ 4.2 in App.\ III  of
of \cite{ ChBook}, and its proof for the estimates
for the time on which the solution exists.} that 
 there is $\e_1>0$ such that if $\|g^\prime-\hat g\|_{C^{m}_b(\hattuM _1;\hat g^+)}<\e_1$,  then there  is a map  $f:M_1 \to M$ satisfying
the Cauchy problem (\ref{C-problem 1})-(\ref{C-problem 2})
and the solution depends continuously, 
 in
$C^{m-5}_b([0,t_1]\times N,g^+)$, on the metric $g^{\prime}$.
 Moreover, by the uniqueness of the wave map, we have
$f|_{M_1\setminus \K_1}=id$ so that
 $f(\K_1)\cap M_0\subset \K_0$.

 As
 the inverse function of the wave map $f$ depends continuously,
 in
$C^{m-5}_b([0,t_1]\times N,g^+)$, on the metric $g^{\prime}$
we can also assume that $\e_1$ is  so small that
$\hattuM _0\subset f(M_1)$.

Denote next $g:=f_*g^{\prime}$, $T:=f_*T^{\prime}$, and $\rho:=f_*\rho^{\prime}$ 
and define $\hat \rho=\hat T-\frac 12 (\hbox{Tr}\, \hat T)\hat g.$ 
Then $g$ is $C^{m-6}$-smooth and the equation (\ref{Einstein on M_0}) implies 
\beq\label{Einstein on M1}
\Ein(g)=T\quad \hbox{on }\hattuM _0.
\eeq 
Since $f$ is a wave map, $g$ satisfies (\ref{wave equation id}) and thus  by the
definition of the reduced Einstein tensor,  (\ref{Reduced Einstein tensor}), we have
\ba
\Ein_{pq}(g)=
(\Ein_{\hat g} (g))_{pq}\quad \hbox{on } \hattuM _0.
\ea
This and  (\ref{Einstein on M1}) yield the $\hat g$-reduced Einstein equation
\beq\label{hat g reduced einstein equations}
(\Ein_{\hat g} (g))_{pq}=T_{pq}\quad \hbox{on } \hattuM _0.
\eeq
This equation is useful for our considerations as it is a quasilinear,
hyperbolic equation on $\hattuM _0$. Recall that 
 $g$ coincides with
$\hat g$ in $M_0\setminus \K_0$. The unique solvability of this 
Cauchy problem is studied in e.g.\ \cite[Thm.\ 4.6 and 4.13]{ ChBook} and \cite {HKM} 
and in Appendix B below.


\subsection*{A.4.\ Relation of the reduced Einstein equations and for the original Einstein equation}

The metric $g$ which solves the $\hat g$-reduced Einstein
equation $\Ein_{\hat g} (g)=T$ is a solution of the original
Einstein equations  $\Ein (g)=T$ if the harmonicity
functions $\hat F^n$
vanish identically. Next we recall the result that 
the harmonicity functions vanish on $\hattuM _0$
when 
\beq\label{eq: good system}
& &(\Ein_{\hat g}(g))_{jk}=T_{jk},\quad \hbox{on }M_0,\\
\nonumber & &\nabla_pT^{pq}=0,\quad \hbox{on }M_0,\\
\nonumber & &g=\hat g,\quad \hbox{on }M_0\setminus \K_0.
\eeq
To see this, let
us  denote $\Ein_{jk}(g)=S_{jk}$, $S^{jk}=g^{jn}g^{km}S_{nm}$,
and $T^{jk}=g^{jn}g^{km}T_{nm}$. 
Following the standard arguments,
see \cite{ ChBook}, we see from (\ref{Reduced Einstein tensor}) that in local coordinates
\ba
S_{jk}-(\Ein_{\hat g}(g))_{jk}=\frac12(g_{jn}\hat \nabla_k \hat F^n+g_{kn}\hat \nabla_j\hat F^n-g_{jk}\hat \nabla_n \hat  F^n).
\ea
Using equations (\ref{eq: good system}), the Bianchi identity  $\nabla_pS^{pq}=0$,
and the basic property of Lorentzian connection,
$\nabla_kg^{nm}=0$,
we obtain 
\ba
0&=&
2\nabla_p(S^{pq}-T^{pq})\\
&=&\nabla_p(g^{qk}\hat  \nabla_k F^p+g^{pm}\hat\nabla_m \hat F^q-  g^{pq}\hat  \nabla_n \hat F^n)
\\
&=&g^{pm}\nabla_p \hat \nabla_m \hat F^q+(g^{qk} \nabla_p \hat \nabla_k \hat F^p-  g^{qp} \nabla_p \hat \nabla_n \hat F^n)\\
&=&g^{pm}\nabla_p \hat \nabla_m \hat F^q+W^q(\hat F)
\ea
where $\hat F=(\hat F^q)_{q=1}^4$ and
the operator 
\ba
W:(\hat F^q)_{q=1}^4\mapsto (g^{qk} (\nabla_p \hat \nabla_k \hat F^p- \nabla_k \hat \nabla_p \hat F^p))_{q=1}^4
\ea
is a linear first order differential operator which coefficients are polynomial functions of
$\hat g_{jk}$, $\hat g^{jk}$, $g_{jk}$, $g^{jk}$ and their first derivatives.

Thus the harmonicity functions $\hat F^q$ satisfy on $\hattuM _0$ the hyperbolic initial 
value problem 
\ba
& &g^{pm}\nabla_p \hat \nabla_m \hat F^q+W^q(\hat F)=0,\quad\hbox{on }M_0,\\
& & \hat F^q=0,\quad \hbox{on }M_0\setminus \K_0,
\ea
and as this initial 
Cauchy problem is uniquely solved by \cite[Thm.\ 4.6 and 4.13]{ ChBook}
or \cite {HKM}, we see that $\hat F^q=0$ on $\hattuM _0$. Thus
equations (\ref{eq: good system})
yield that
the Einstein equations   $\Ein(g)=T$ holds on $\hattuM _0$.

}

\hiddenfootnote{
\subsection*{Appendix B: Stability and existence of the  direct problem}

 {

Let us  start by explaining  how we can choose a  $C^\infty$-smooth metric $\tilde g$ 
such that $\hat g<\tilde g$ and  $(M,\tilde g)$ is globally hyperbolic:
When $v(x)$ the eigenvector corresponding to the negative eigenvalue
of $\hat g(x)$, we can choose a smooth, strictly positive function
$\eta:\hattuM \to \R_+$ such that $\tilde g^\prime:=\hat g-\eta v\otimes v<\tilde g$. Then
$(\hattuM ,\tilde g^\prime)$ is globally hyperbolic, $\tilde g^\prime$ is smooth
and $\hat g<\tilde g^\prime$. Thus we can replace $\tilde g$ by the smooth metric
 $\tilde g^\prime$ having the same properties that are required for $\tilde g$.
 
 Let us now return to consider existence and stability of the solutions of the
 Einstein-scalar field equations.
Let $t={\bf t}(x)$ be local time so that there is a diffeomorphism
$\Psi:M\to \R\times N$, $\Psi(x)=({\bf t}(x),Y(x))$, and
 $S(T)=\{x\in \hattuM _0; \ {\bf t}(x)=T\}$, $T\in \R$ are Cauchy surfaces.
 Let $t_0>0$. Next we identify $M$ and $\R\times N$ via the map $\Psi$ and just
 denote $M=\R\times N$.
Let us denote $M(t)=(-\infty,t)\times N$,  and $M_0=M(t_0)$.
By \cite[Cor.\ A.5.4]{BGP} the set $\K=J^+_{\tilde g}({p^-})\cap \hattuM _0$,
where $\hat p^0\in M_0$,
is compact.
Let $N_1,N_2\subset N$ be such
 open relatively compact sets with smooth boundary that
 $N_1\subset N_2$ and  $Y(J_{\tilde g}^+(\hat p^0)\cap \hattuM _0)\subset N_1$.  \HOX{We could improve explanation on $\tilde N$}
 
 To simplify citations to existing literature, 
let us define $\tilde N$ to be a compact manifold without boundary such
that $N_2$ can be considered as a subset of $\tilde N$. Using a construction based on
a suitable partition of unity, the Hopf double of the manifold $N_2$, 
and the Seeley extension of the metric tensor,  we can
endow $\tilde M=(-\infty,t_0)\times \tilde N$ with a smooth Lorentzian
metric $\hat g^e$ (index $e$ is for "extended")
so that $\{t\}\times \tilde N$ are Cauchy surfaces of $\tilde N$
and that $\hat g$ and $\hat g^e$ coincide in  the set $\Psi^{-1}((-\infty,t_0)\times N_1)$ that contains the set $J^+_{\hat g}(\hat p^0)\cap \hattuM _0$.} We extend the metric $\tilde g$ to a (possibly non-smooth) globally
hyperbolic metric $\tilde g^e$ on $\tilde M_0=(-\infty,t_0)\times \tilde N)$
such that $\hat g^e<\tilde g^e$.

 To simplify notations below we denote $\hat g^e=\hat g$ and
 $\tilde g^e=\tilde g$ 
 on the whole
 $\tilde M_0$.
 Our aim is to prove the estimate (\ref{eq: Lip estim}).

Let us denote by $t={\bf t}(x)$ the local time. 
Recall that when $(g,\phi)$ is a solution of 
the scalar field-Einstein equation, we denote $u=(g-\hat g,\phi-\hat \phi)$.
We will consider the equation for $u$, and to emphasize that the metric
depends on $u$, we denote $g=g(u)$ and assume below that  both the metric $g$ is  dominated
by $\tilde g$, that is, $ g<\tilde g$.
We use the pairs ${\bf u}(t)=(u(t,\cdotp),\p_t u(t,\cdotp))\in
H^1(\tilde N)\times L^2(\tilde N)$ and the notations
 ${\bf v}(t)=(v(t),\p_t v(t))$ etc. 
Let us consider a generalization of the system (\ref{eq: notation for hyperbolic system 1}) of the form 
\beq\label{eq: notation for hyperbolic system}
& &\square_{g(u)}u+V(x,D)u+H(u,\p u) =R(x,u,\p u)F+K,\ \ x\in \tilde M_0,\hspace{-1cm}\\
& & \nonumber \supp(u)\subset \K,
\eeq
where $\square_{g(u)}$ is the Lorentzian Laplace operator operating on the sections 
of the bundle $\B^L$ on $M_0$ and 
$\supp(F)\cup \supp(K)\subset \K$.
Note that above  $u=(g-\hat g,\phi-\hat \phi)$ and $g(u)=g$.
Also, $F\mapsto R(x,u,\p u)F$ is a 
 linear first order differential operator which coefficients at $x$ are depending smoothly on 
 $u(x)$, $\p_j u(x)$ and the  derivatives 
of $(\hat g,\hat \phi)$ at $x$ and
\ba V(x,D)=V^j(x)\p_j+V(x)\ea is a linear first order differential operator which coefficients 
at $x$ are depending smoothly on the  derivatives 
of $(\hat g,\hat \phi)$  at $x$, and finally, $H(u,\p u)$ is a polynomial
of $u(x)$ and $\p_j u(x)$ which coefficients 
at $x$ are depending smoothly on the  derivatives 
of $(\hat g,\hat \phi)$   such that
 $\p^\a_v\p^\beta_w H(v,w)|_{v=0,w=0}=0$ for
$|\a|+|\b|\leq 1$. 
By \cite[Lemma 9.7]{Ringstrom}, the equation  (\ref{eq: notation for hyperbolic system})
has at most one solution with given  $C^2$-smooth source functions $F$ and $K$.
Next we consider the existence of $u$ and its dependency on $F$ and $K$.


%

Below we use notations, c.f. (\ref{eq: notation for hyperbolic system 1}) and
(\ref{eq: notation for hyperbolic system 1b}) \ba
{\mathcal R}({\bf u},F)=R(x,u(x),\p u(x))F(x),\quad
{\mathcal H}({\bf u})=H(u(x),\p u(x)).\ea
Note that $u=0$, i.e., $g=\hat g$ and $\phi=\hat\phi$ satisfies (\ref{eq: notation for hyperbolic system}) with $F=0$
and $K=0$.
Let us use the same notations as in  \cite {HKM} cf. also \cite[section 16]{Kato1975}, 
 to consider quasilinear wave equation on $[0,t_0]\times \tilde N$.
Let $\H^{(s)}(\tilde N)=H^s(\tilde N)\times H^{s-1}(\tilde N)$ and 
\ba
Z=\H^{(1)}(\tilde N),\quad  Y=\H^{(k+1)}(\tilde N),\quad X=\H^{(k)}(\tilde N).
\ea
The norms on these space are defined invariantly using the smooth Riemannian 
metric $h=\hat g|_{\{0\}\times \tilde N}$ on $\tilde N$. \HOX{Maybe we should explain the 
bundles, the connection and the wave operator on section in a detailed way.}
Note that $\H^{(s)}(\tilde N)$ are in fact the Sobolev spaces of sections on the bundle $\pi:\B_K\to
\tilde N$, where $\B_K$ denotes also the  pull back bundle of $\B_K$ on $\tilde M$ in the map
$id:\{0\}\times \tilde N\to \tilde M_0$,
or on the bundle $\pi:\B_L\to \tilde N$. Below, $\nabla_h$ denotes the  standard connection of the bundle  $\B_K$
or  $\B_L$ associated to the metric $h$.

Let $k\geq 4$ be an even integer.
By definition of $H$ and $R$ we see that  there are  $0<r_0<1$ and $L_1,L_2>0$,
all depending on $\hat g$, $\hat \phi$, $\K$, and $t_0$, such that if 
$0<r\leq r_0$ and
\beq\label{new basic conditions}
& &\|{\bf v} \|_{C([0,t_0];\H^{(k+1)}(\tilde N))}\leq r,\quad \|{\bf v}^\prime\|_{C([0,t_0];\H^{(k+1)}(\tilde N))}\leq r,\\
& &\|F\|_{C([0,t_0];H^{(k+1)}(\tilde N))}\leq r^2,\quad \|K\|_{C([0,t_0];H^{(k+1)}(\tilde N))}\leq r^2 \nonumber\\
& &\|F^\prime\|_{C([0,t_0];H^{(k+1)}(\tilde N))}\leq r^2,\quad \|K^\prime\|_{C([0,t_0];H^{(k+1)}(\tilde N))}\leq r^2 \nonumber
\eeq
then
\beq\label{new f1f2 conditions}
& &\quad \quad \|g(\cdotp;v)^{-1}\|_{C([0,t_0];H^s(\tilde N))}\leq L_1 ,\\ \nonumber
& &\|{\mathcal H} ({\bf v})\|_{C([0,t_0];H^{s-1}(\tilde N))}\leq L_2r^2,
\quad \|{\mathcal H} ({\bf v}^\prime)\|_{C([0,t_0];H^{s-1}(\tilde N))}\leq L_2r^2,\\ \nonumber
& &\|{\mathcal H} ({\bf v})-{\mathcal H} ({\bf v}^\prime)\|_{C([0,t_0];H^{s-1}(\tilde N))}\leq  L_2r\, \|{\bf v}-{\bf v}^\prime\|_{C([0,t_0];\H^{(s)}(\tilde N))}, 
\\ \nonumber
& &\|\mathcal R ({\bf v}^\prime,F^\prime)\|_{C([0,t_0];H^{s-1}(\tilde N))}\leq L_2r^2,
\quad \|\mathcal R ({\bf v},F)\|_{C([0,t_0];H^{s-1}(\tilde N))}\leq L_2r^2,\hspace{-1cm} \\ \nonumber
& &\|\mathcal R ({\bf v},F)-\mathcal R ({\bf v}^\prime, F^\prime)\|_{C([0,t_0];H^{s-1}(\tilde N))}
\\ \nonumber & &\leq  
L_2r\, \|{\bf v}-{\bf v}^\prime\|_{C([0,t_0];\H^{(s)}(\tilde N))}+L_2\|F-{F}^\prime\|_{C([0,t_0];H^{s+1}(\tilde N))\cap C^1([0,t_0];H^{s}(\tilde N))},\hspace{-2cm}
%
  \eeq
 for all $s\in [1,k+1]$. 
%

%
%

Next we write (\ref{eq: notation for hyperbolic system}) as a first order
system.  To this end,
let $\mathcal A(t,{\bf v}):\H^{(s)}(\tilde N)\to \H^{(s-1)}(\tilde N)$ be the operator 
$\mathcal A(t,{\bf v})=\mathcal A_0(t,{\bf v})+\mathcal A_1(t,{\bf v})$
where in local coordinates and in the local trivialization of the bundle $\B^L$
\ba
\mathcal A_0(t,{\bf v})=-\left(\begin{array}{cc} 0 & I\\
 \frac 1{g^{00}(v)}\sum_{j,k=1}^3 g^{jk}(v)\frac \p{\p x^j}\frac\p{ \p x^k}\quad &
 \frac 1{g^{00}(v)}\sum_{m=1}^3 g^{0m}(v)\frac {\p}{\p x^m} \end{array}\right)
\ea 
with $g^{jk}(v)=g^{jk}(t,\cdotp;v)$ is a function on $\tilde N$ and 
\ba
\mathcal A_1(t,{\bf v})= \frac {-1}{g^{00}(v)}\left(\begin{array}{cc} 0 & 0\\
\sum_{j=1}^3 
 B^j(v)
   \frac \p{\p x^j}\quad &
B^{0}(v)\end{array}\right)
\ea 
where $B^j(v)$ depend on $v(t,x)$ and its first derivatives, and 
the connection coefficients (the Christoffel  symbols) corresponding  to $g(v)$. We denote $\S=(F,K)$ and 
 \ba
& &f_\S(t,{\bf v})=(f_\S^1(t,{\bf v}),f_\S^2(t,{\bf v}))\in \H^{(k)}(\tilde N),
\quad\hbox{where}\\
& & f_\S^1(t,{\bf v})=0,\quad f_\S^2(t,{\bf v})=\mathcal R({\bf v},F)(t,\cdotp)-\mathcal H({\bf v})(t,\cdotp)+K(t,\cdotp).
 \ea
 
 Note that when (\ref{new basic conditions}) are satisfied with $r<r_0$, inequalities (\ref{new f1f2 conditions}) imply that there exists $C_2>0$ so that
 \beq\label{fF bound r2}
 & &\|f_\S(t,{\bf v})\|_Y+\|f_{\S^\prime}(t,{\bf v}^\prime)\|_Y\leq C_2r^2,\\
 & &\|f_\S(t,{\bf v})-f_\S(t,{\bf v}^\prime)\|_Y\leq C_2
 r\,\|{\bf v}-{\bf v}^\prime\|_{C([0,t_0];Y)}. \nonumber
 \eeq

%
%

Let $U^{\bf v}(t,s)$ be the wave propagator corresponding to metric $g(v)$, that is,
$U^{\bf v}(t,s): {\bf h}\mapsto {\bf w}$, where ${\bf w}(t)=(w(t),\p_t w(t))$  solves
\ba
(\square_{g(v)}+V(x,D))w=0\quad\hbox{for } (t,y)\in [s,t_0]\times \tilde N,\quad\hbox{with }
 {\bf w}(s,y)={\bf h}.
 \ea 
 Let $S=(\nabla_{h}^*\nabla_{h}+1)^{k/2}:Y\to Z$ be an isomorphism.
As $k$ is an even integer, we see using multiplication estimates for Sobolev
spaces, see e.g.\
  \cite [Sec.\ 3.2, point (2)]{HKM}, that there exists $c_1>0$ (depending on
  $r_0,L_1,$ and $L^2$) so that  $\A(t,{\bf v})S-S\A(t,{\bf v})=C(t,{\bf v})$,
 where $\|C(t,{\bf v})\|_{Y\to Z}\leq c_1$ for all ${\bf v}$ satisfying
 (\ref{new basic conditions}). This yields that 
 the property (A2) in  \cite {HKM} holds, namely
 that $S\A(t,{\bf v})S^{-1}=\A(t,{\bf v})+B(t,{\bf v})$ where
$B(t,{\bf v})$ extends to a bounded operator in $Z$
 for which $\|B(t,{\bf v})\|_{Z\to Z}\leq c_1$ for all ${\bf v}$ satisfying
 (\ref{new basic conditions}).
Alternatively, to see the mapping properties of $B(t,{\bf v})$ we could use the fact that 
 $B(t,{\bf v})$ is a zeroth order pseudodifferential operator with
$H^{k}$-symbol.\hiddenfootnote{On mapping properties of such
pseudodifferential operators, see J. Marschall, Pseudodifferential operators with coefficients in Sobolev spaces. Trans. Amer. Math. Soc. 307 (1988), no. 1, 335-361.} 
 
 Thus the proof of  \cite [Lemma 2.6]{HKM} shows\hiddenfootnote{ Alternatively, as $U^{\bf v}(t,s)$ are propagators
  for linear wave equations having finite speed of wave propagation, one can prove 
the  estimates (\ref{U bounds}) by considering the wave equation
in local coordinate neighborhoods $W_j\subset W_j^\prime \subset \tilde N$,
$\Phi_j:  W_j^\prime\to \R^3$ and a partition of unity $\psi_j\in C^\infty_0(W_j)$ and
cut-off functions $\psi^\prime_j\in C^\infty_0(W_j^\prime)$.
Then  \cite [Lemma 2.6]{HKM}, used  in local coordinates, shows
that when $t_k-t_{k-1}$ is small enough then 
\ba 
\| \psi^\prime_j U^{\bf v}(t_k,t_{k-1})(\psi_j {\bf w}) \|_{Y}\leq C_{3,jk}^{\prime}e^{C_4(t-s)}\| \psi_j {\bf w} \|_{Y}.
\ea
Using  sufficiently small time steps $t_k-t_{k-1}$ and combining the estimates
for different $j$:s
together, one obtain estimates (\ref{U bounds}).}
 that 
%
there is a constant  $C_3>0$
   so that 
  \beq
  \label{U bounds}\|U^{\bf v}(t,s)\|_{Z\to Z}\leq C_3\quad\hbox{and}\quad
  \|U^{\bf v}(t,s) \|_{Y\to Y}\leq C_3\hspace{-1.5cm}
  \eeq 
  for $0\leq s<t\leq t_0$.
    By interpolation of estimates (\ref{U bounds}), we see also that 
    \beq
  \label{U bounds2}
  \|U^{\bf v}(t,s)\|_{X\to X}\leq C_3,
  \eeq  for $0\leq s<t\leq t_0$.

Let us next modify the reasoning given in \cite{Kato1975}: let  $r_1\in (0,r_0)$
be a parameter which value will be chosen later,  $C_1>0$
and 
$E$ be the space of functions ${\bf u}\in C([0,t_0];X)$ for which
\beq
\label{1 bounded}& &\|{\bf u}(t)\|_Y\leq r_1\quad\hbox{and}\\
\label{2 Lip}& &\|{\bf u}(t_1)-{\bf u}(t_2)\|_X\leq C_1|t_1-t_2|\eeq
for all $t,t_1,t_2\in [0,t_0]$.
The set $E$ is endowed by the metric of  $C([0,t_0];X)$. We 
note that  by \cite[ Lemma 7.3] 
{Kato1975},
a convex $Y$-bounded, $Y$-closed set is closed also in $X$.
Similarly, functions $G:[0,t_0]\to X$  satisfying (\ref{2 Lip}) 
form a closed subspace of $C([0,t_0];X)$. Thus $E\subset X$ is a closed set implying
that $E$ is a complete metric space.

Let
 \ba 
 & &\hspace{-1cm}W=\\ & &\hspace{-1cm}\{(F,K)\in C([0,t_0];H^{k+1}(\tilde N))^2;\ \sup_{t\in [0,t_0]}\|F(t)\|_{H^{k+1}(\tilde N)}+\|K(t)\|_{H^{k+1}(\tilde N)}< r_1\}.\hspace{-1cm}
 \ea

Following \cite[p. 44]{Kato1975}, we see that
the solution of
equation (\ref{eq: notation for hyperbolic system}) with the source $\S\in W$  is found as
a fixed point, if it exists, of the map $\Phi_\S:E\to C([0,t_0];Y)$ where 
$\Phi_\S({\bf v})={\bf u}$ is  given by
$$
{\bf u}(t)=\int_0^t U^{\bf v}(t,\tilde t)f_\S(\tilde t,{\bf v})\,d\tilde t,\quad 0\leq t\leq t_0.
$$ 
Below, we denote  ${\bf u}^{\bf v}=\Phi_\S({\bf v})$.

Since $\Phi_{\S_0}(0)=0$ where
$\S_0=(0,0)$, we see using 
the above and the inequality $\|\,\cdotp\|_X\leq \|\,\cdotp\|_Y$ that the
 function ${\bf u}^{\bf v}$ satifies
  \ba
& &  \|{\bf u}^{\bf v}\|_{C([0,t_0];Y)}\leq C_3C_2t_0 r_1^2,\\
& &  \|{\bf u}^{\bf v}(t_2)-{\bf u}^{\bf v}(t_1)\|_{X}\leq C_3C_2r_1^2|t_2-t_1|,\quad t_1,t_2\in [0,t_0].
  \ea
  When $r_1>0$ is so small that $C_3C_2(1+t_0)<r_1^{-1}$ 
  and $C_3C_2r_1^2<C_1$
  we see that $  \|\Phi_\S({\bf v})\|_{C([0,t_0];Y)}<r_1$
  and $ \|\Phi_\S({\bf v})\|_{C^{0,1}([0,t_0];X)}<C_1$.
  Hence $ \Phi_\S(E)\subset E$
 and we can consider $\Phi_\S$ as a map
 $\Phi_\S:E\to E$.
  
  As  $k>1+\frac 32$, it follows from Sobolev embedding theorem
  that $X=\H^{(k)}(\tilde N)\subset C^1(\tilde N)^2$. This yields that
  by \cite[Thm. 3]{Kato1975}, for the original reference, see Theorems III-IV
  in \cite{Kato1973},
  \ba
 & & \|(U^{\bf v}(t,s)- U^{\bf v^\prime}(t,s)){\bf h}\|_X \\
 & &\leq 
 C_3\bigg(\sup_{t^\prime \in [0,t]}\|\A(t^\prime,{\bf v})-
 \A(t^\prime,{\bf v}^\prime)\|_{Y\to X}  \|U^{\bf v}(t^\prime,0) {\bf h}\|_Y\bigg)\\
 & &\leq 
 C_3^2 \|{\bf v}-{\bf v}^\prime\|_{C([0,t_0];X)}  \|{\bf h}\|_Y.
  \ea
  Thus,
  \ba
&&\|  U^{\bf v}(t,s)f_\S(s,{\bf v})- U^{\bf v^\prime}(t,s)f_{\S}(s,{\bf v}^\prime)\|_X\\
&\leq&\|  (U^{\bf v}(t,s)- U^{\bf v^\prime}(t,s))f_{\S}(s,{\bf v})\|_X+
\|  U^{\bf v^\prime}(t,s)(f_\S(s,{\bf v})- f_{\S}(s,{\bf v}^\prime))\|_X\\
&\leq&
(1+C_3)^2C_2r_1^2 \|{\bf v}-{\bf v}^\prime\|_{C([0,t_0];X)} .
  \ea
This implies that
\ba
& & \| \Phi_\S({\bf v})-\Phi_{\S}({\bf v}^\prime)\|_{C([0,t_0];X)}
\leq t_0(1+C_3)^2C_2r_1^2 \|{\bf v}-{\bf v}^\prime\|_{C([0,t_0];X))}.
\ea

Assume next that  $r_1>0$ is so small that we have also
\ba
t_0(1+C_3)^2C_2r_1^2<\frac 12.
\ea
 cf.\
  Thm.\ I in \cite{HKM} (or
  (9.15)  and (10.3)-(10.5) in \cite{Kato1975}). For $\S\in W$ this 
 implies that $\Phi_\S:E\to E$ is a contraction with a contraction constant 
 $C_L\leq \frac 12$,
  and thus 
 $\Phi_\S$ has a unique fixed point ${\bf u}$  in the space $E\subset C^{0,1}([0,t_0];X)$.
 
 Moreover, elementary considerations related to fixed point
 of the map $\Phi_\S$ show that ${\bf u}$  in  $C([0,t_0];X)$ depends in $E\subset C([0,t_0];X)$
 Lipschitz-continuously on
 $\S\in  W\subset C([0,t_0];H^{k+1}(\tilde N))^2$. Indeed, if $\|\S-\S^\prime\|_{C([0,t_0];H^{k+1})^2}<\e$,
 we see that
  \beq\label{fF bound r2 B}
 \|f_\S(t,{\bf v})-f_{\S^\prime}(t,{\bf v})\|_Y\leq C_2\e,\quad t\in [0,t_0],
 \eeq
 and when  (\ref{fF bound r2}) and (\ref{U bounds}) are satisfied 
 with $r=r_1$,
 we have
  \ba
& & \| \Phi_\S({\bf v})-\Phi_{\S^\prime}({\bf v}^\prime)\|_{C([0,t_0];Y)}
\leq C_3C_2t_0 r_1^2.
\ea
 Hence 
 \ba
 \|\Phi_\S({\bf v})-\Phi_{\S^\prime}({\bf v})\|_{C([0,t_0];Y)}\leq  t_0C_3C_2\e.
 \ea
 This and standard estimates for fixed points, yield that when $\e$ is small enough 
 the fixed point ${\bf u}^\prime$ of the map  $\Phi_{\S^\prime}:E\to E$ 
 corresponding to the source ${\S^\prime}$ and
  the fixed point ${\bf u}$ of the map $\Phi_{\S}:E\to E$ 
  corresponding to the source ${\S}$ satisfy
  \beq\label{stability}
  \|{\bf u}-{\bf u}^\prime\|_{C([0,t_0];X)}\leq  \frac 1{1-C_L}t_0C_3C_2\e.
  \eeq
{
Thus the solution ${\bf u}$ 
  depends in $C([0,t_0];X)$ Lipschitz continuously on
 $\S\in  C([0,t_0];H^{k+1}(\tilde N))^2$ (see also \cite[Sect.\ 16]{Kato1975}, 
 and \cite{Ringstrom})}. 
In fact,  for analogous systems it is possible to show that $u$ is in $C([0,t_0];Y)$,
but one can not obtain  
 Lipschitz or H\"older  stability for $u$ in the $Y$-norm, see  \cite{Kato1975}, Remark 7.2. 
 
 Finally, we note that the fixed point ${\bf u}$ of $\Phi_\S$ can be found as a limit
 $ {\bf u}=\lim_{n\to \infty}{\bf u}_n$ in $C([0,t_0];X)$, where ${\bf u}_0=0$ and  ${\bf u}_n=\Phi_\S({\bf u}_{n-1})$.
Denote ${\bf u}_n=(g_n-\hat g,\phi_n-\hat \phi)$. We see that if
 $\supp({\bf u}_{n-1})\subset J_{\hat g}(\supp(\S))$ 
then also $\supp(g_{n-1}-\hat g)\subset J_{\hat g}(\supp(\S))$.
Hence  for all $x\in M_0\setminus J_{\hat g}(\supp(\S))$ we see that $J_{g_{n-1}}^-(x)
\cap J_{\hat g}(\supp(\S))=\emptyset.$ Then, using the definition of the map $\Phi_\S$ we see that
 $\supp({\bf u}_n)\subset J_{\hat g}(\supp(\S))$. Using induction we see that this holds
for all $n$ and hence we see that the solution ${\bf u}$ satisfies
 \beq\label{eq: support condition for u}
 \supp({\bf u})\subset J_{\hat g}(\supp(\S)).
 \eeq 
 \hiddenfootnote{REMOVE MATERIAL IN THIS FOOTNOTE OR MOVE IT ELSEWHERE:
  Let us next use the above considerations to study system 
   (\ref{eq: notation for hyperbolic system 1}) on $\R\times N$.
   For this end, we denote below the background metric of  $\tilde M=\R\times \tilde N$
   by $\tilde g$ and the background metric of  $\hattuM _0=\R\times  N$
   by $\hat g$.
 Consider a source 
\ba
\S\in {C^1([0,T];H^{s_0-1}(\tilde N))\cap
C^0([0,t_0];H^{s_0}(\tilde N))}
\ea  that is supported in a compact set $\K$ and
small enough in the norm of this space and let $\tilde u$ 
be the solution of the  
system (\ref{eq: notation for hyperbolic system 1}) on $\R\times \tilde N$.
Assume that  $u$ is some $C^2$-solution of  
(\ref{eq: notation for hyperbolic system 1}) on $\R\times N$
with the same source $\S$. Our aim is to show that it
coincides with $\tilde u$ in its support.
%
For this end, let $\tilde g^\prime$ be a metric tensor which
has the same eigenvectors as $\tilde g$ and the same
eigenvalues except that  the unique negative eigenvalue of $\tilde g$
 is multiplied by $(1-h_1)$, $h_1>0$. 
Assume that $h_1>0$ is so small
that $J^+_{\hattuM _0,\tilde g^\prime}(p^+)\subset  \hattuM _0\subset N_0\times (-\infty,t_0]$
and assume that $\e$ is so small that $g(\tilde u)<\tilde g^\prime$ for all 
sources $\S$ with $\|\S\|_{C^1([0,t_0];H^{s_0-1}(\tilde N))\cap
C^0([0,t_0];H^{s_0}(\tilde N))}\leq \e$. Then
 $J^+_{\tilde M,g(\tilde u)}(p^+)\subset  (-\infty,t_0]\times N_0$.
This  proves the estimate (\ref{eq: Lip estim}).

Next us we consider solution $u$ on $\hattuM _0$ corresponding to source $\S$.
Let  $T(h_1,\S)$
be the supremum of all $T^\prime\leq t_0$ for which
$J^+_{g(u)}(p^+)\cap \hattuM _0(T^\prime)\subset N_1\times (-\infty,T^\prime]$.
Assume that $ T(h_1,\S)<t_0$.
Then for $T^\prime= T(h_1,\S)$ we see that
$u$ can be continued by zero from $ (-\infty,T^\prime]\times N_1$
to a function on $(-\infty,T^\prime]\times \tilde N$ that 
satisfies the system (\ref{eq: notation for hyperbolic system 1}). Thus
$u$ coincides in $ (-\infty,T^\prime]\times N_1$
with $\tilde u$ and vanishes outside this set.
Since  $J^+_{\tilde M,g(\tilde u)}(p^+)\subset  (-\infty,t_0]\times N_0$, this implies that
also $u$ satisfies 
$J^+_{g(u)}(p^+)\cap \hattuM _0(T^\prime)\subset (-\infty,T^\prime]\times N_0$.
As $u$ is $C^2$-smooth and $N_0\subset \subset N_1$
we see that there is $T^\prime_1> T(h_1,\S)$ so 
that $J^+_{g(u)}(p^+)\cap \hattuM _0(T^\prime_1)\subset  (-\infty,T^\prime_1]
\times N_1$ that is in contradiction with definition of $ T(h_1,\S)$. Thus
we have to have  $T(h_1,\S)=t_0$.


This shows that when the norm of $\S$ is smaller than the above chosen $\e$
then $J^+_{g(u)}(p^+)\cap \hattuM _0\subset N_0\times (-\infty,t_0]$.
Then $u$ vanishes outside $N_1\times (-\infty,t_0]$ by
the support condition in (\ref{eq: notation for hyperbolic system 1}).
Thus $u$ and $\tilde u$ are both supported on  $N_1\times (-\infty,t_0]$ 
and coincide there. As $\tilde u$ is unique, 
this shows that for a sufficiently small sources $ \S$
 supported on  $\K$   the solution $u$ of
(\ref{eq: notation for hyperbolic system 1})  in $M_0$
exists and is unique.}

\bigskip

\subsection*{Appendix D: An inverse problem for a non-linear wave equation}

In this appendix  we explain how a problem for a scalar wave equation can
be solved with the same techniques that we used for the Einstein equations.

Let $(M_j,g_j)$, $j=1,2$ be two globally hyperbolic  $(1+3)$  
dimensional Lorentzian
manifolds represented using
  global smooth time functions as $M_j=\R\times N_j$, $\mu_j=\mu_j([-1,1])
\subset M_j$ be
a time-like geodesic and $U_j\subset M_j$ be open, relatively compact  
neighborhood
of $\mu_j([s_-,s_+])$, $-1<s_-<s_+<1$. Let $M_j^0=(-\infty,T_0)\times  
N_j$ where $T_0>0$ is such
that $U_j\subset  M_j^0$.
{\mltext Consider the
non-linear wave equation
\beq\nonumber
& &\hspace{-.5cm}  \square_{g_j}u(x)+a_j(x)\,u(x)^2=f(x)
\quad\hbox{on }M_j^0,\hspace{-.5cm}\\
\label{eq: wave-eq general} & &\quad \supp(u)\subset J^+_{g_j}(\supp(f)),
\eeq
where $\supp(f)\subset U_j,$
\ba
\square_gu=
-\sum_{p,q=1}^4\det(-g(x))^{-1/2}\frac  \p{\p x^p}
\left ((-\det(g(x))^{1/2}
g^{pq}(x)\frac \p{\p x^q}u(x)\right),
\ea
$\det(g)=\det((g_{pq}(x))_{p,q=1}^4)$,
  $f\in C^6_0(U_j)$ is a controllable source, and $a_j$ is
a non-vanishing $C^\infty$-smooth  function.}
Our goal is to prove the following result:

\begin{theorem}\label{main thm3}
Let $(M_j,g_j)$, $j=1,2$ be two
open,  smooth, globally hyperbolic    Lorentzian manifolds of  
dimension $(1+3)$.
Let
$p^+_j=\mu_j(s_+), p^-_j=\mu_j(s_-)\in M_j$ the points of a  time-like  
geodesic  $\mu_j=
\mu_j([-1,1])\subset M_j$, $-1<s_-<s_+<1$,
and let  $U_j\subset M_j$ be an open  relatively compact neighborhood
of $\mu_j([s_-,s_+])$  given in (\ref{eq: Def Wg with hat}).  Let  
$a_j:M_j\to \R$, $j=1,2$ be $C^\infty$-smooth functions that are
non-zero on $M_j$.

Let   $L_{U_j}$, $j=1,2$ be measurement operators defined
in an open set $\mathcal W_j\subset C^6_0(U_j)$ containing the zero  
function by setting
\beq\label{measurement operator}
L_{U_j}: f\mapsto u|_{U_j},\quad f\in C^6_0(U_j),
\eeq
where $u$ satisfies the wave equation (\ref{eq: wave-eq general}) on  
$(M_j^0,g_j)$.

Assume that there is a diffeomorphic isometry
$\Phi:U_1\to U_2$ so that $\Phi(p^-_1)=p^-_2$ and
$\Phi(p^+_1)=p^+_2$ and the measurement maps
satisfy
\ba
((\Phi^{-1})^*\circ L_{U_1}\circ \Phi^*) f =L_{U_2}f
\ea
for all $f\in \W$ where $\W$ is  some  neighborhood of the zero function in $C^6_0(U_2)$.

Then there is a diffeomorphism $\Psi:I(p^-_1,p^+_1)\to I(p^-_2,p^+_2)$,
and the metric $\Psi^*g_2$ is conformal to $g_1$ in  
$I(p^-_1,p^+_1)\subset M_1$,
that is, there is $\beta(x)$ such that $g_1(x)=\beta(x)(\Psi^*g_2)(x)$ in  
$I(p^-_1,p^+_1)$.
\end{theorem}

We note that the  smoothness assumptions assumed above on the functions
$a$ and the source $f$ are not optimal.
The proof, presented below, is based on using the interaction
of singular waves.  The techniques used can be modified used to study
different non-linearities, such as the equations
$\square_{\hat g} u+a(x)u^3=f$, $\square_{\hat g} u+a(x)u_t^2=f$,
or $\square_{g(x,u(x))} u=f$, but these considerations are outside the scope
of this paper.

Theorem \ref{main thm3} can be applicable
  for example in the mathematical analysis of  non-destructive
testing or imaging in non-linear medium e.g,
in imaging  the non-linearity  of the acoustic material parameter  inside a
given body  when it
is under large, time-varying, possibly periodic, changes of the external
pressure and at the same time the body is probed with small-amplitude fields.
  Such acoustic measurements are analogous to
  the recently developed Ultrasound Elastography imaging technique where the interaction
  of the elastic shear and pressure
waves is used for medical imaging, see  
e.g.\ \cite{Hoskins,McLaughlin1,McLaughlin2,Ophir}. There, the slowly  
progressing
shear wave is imaged using a pressure wave and the
image of the shear wave inside the body is used to determine
approximately the material parameters. In other words,
the changes which the elastic wave causes in the medium
are imaged using the interaction of the s-wave and p-wave
components of the  elastic wave.

\HOX{Here a corollary is removed due to a problem in proof.}
\medskip

Next we consider the proofs.

\medskip

\noindent
 {\bf Proof.} (of Theorem \ref{main thm3}).
We will explain how the proof of
Theorem \ref{main thm Einstein} for the Einstein equations  needs to
be modified to obtain the similar result for the non-linear wave equation.

Let  $(M,\hat g)$ be a  smooth
globally hyperbolic Lorentzian manifold that we represent
using a global smooth time function as $M=(-\infty,\infty)\times N$,
and consider   $M^0=(-\infty,T)\times N\subset M$.
Assume that
the set $U$, where the sources are supported and where
we observe the waves, satisfies
$U\subset [0,T]\times N$.

The results of  section \ref{subsec: Direct problem} concerning the  
direct problem
for Einstein equations can be modified for the wave equation
\beq\label{PABC eq}
& &\square_{\hat g}u+au^2=f,\quad\hbox{in }M^0=(-\infty,T)\times N,\\
& &u|_{(-\infty,0)\times N}=0,\nonumber
\eeq
where
$a=a(x)$ is a smooth, non-vanishing function. Here we denote
the metric by $\hat g$ to emphasize the fact that it is independent
on the solution $u$.  Below,
let  $Q$ be the causal inverse  operator of $\square_{\hat g}$.

When $f$ in $C_0([0,t_0];H^6_0(B))\cap
C_0^1([0,t_0];H^5_0(B))$
is small enough,
we see by using  \cite[Prop.\ 9.17]{Ringstrom} and   \cite [Thm.\ III]{HKM},
see also (\ref{stability}) in Appendix B, that  the
equation (\ref{PABC eq}) has a unique solution
$u\in  C_0([0,t_0];H^{5}(N))\cap C_0^1([0,t_0];H^{4}(N))$.
Moreover, 
we  can consider the case when $f=\e f_0$ where $\e>0$
is small.
Then, we can write
\ba
u=\e w_1+\e^2 w_2+\e^3 w_3+\e^4 w_4+E_\e
\ea
where $w_j$ and the reminder term $E_\e$ satisfy
\ba
w_1&=&Qf,\\
w_2&=&-Q(a\,w_1\,w_1),\\
w_3&=&-2Q(a\, w_1\,w_2)\\
&=&2Q(a\, w_1\,Q(a\,w_1\,w_1))
,\\
w_4
&=&-Q(a\, w_2\,w_2)-2Q(a\, w_1\,w_3)\\
&=&-Q(a\, Q(a\,w_1\,w_1)\,Q(a\,w_1\,w_1))\\
& &+4Q(a\, w_1\,Q(a\, w_1\,w_2))
\\
&=&-Q(a\, Q(a\,w_1\,w_1)\,Q(a\,w_1\,w_1))\\
& &-4Q(a\, w_1\,Q(a\, w_1\,Q(a\,w_1\,w_1))),\\
& &\hspace{-1.5cm}\|E_\e\|_{C([0,t_0];H^{4}_0(N))\cap  
C^1([0,t_0];H^{3}_0(N))}\leq C\e^5.
\ea

If we consider sources  $f_{\vec\e}(x)=\sum_{j=1}^4\e_j f_{(j)}(x)$,
$\vec\e=(\e_1,\e_2,\e_3,\e_4),$
and the corresponding solution $u_{\vec \e}$ of (\ref{PABC eq}), we see that
\beq \nonumber
\M^{(4)}&=&\p_{\vec \e}^4u_{\vec \e}|_{\vec\e=0}\\
&=& \nonumber
\p_{\e_1}\p_{\e_2}\p_{\e_3}\p_{\e_4}u_{\vec \e}|_{\vec\e=0}\\
\label{4th interaction for wave eq B}
&=&-\sum_{\sigma\in \Sigma(4)}\bigg(Q(a\,  
Q(a\,u_{(\sigma(1))}\,u_{(\sigma(2))})\,Q(a\,u_{(\sigma(3))}\,u_{(\sigma(4))}))\\  
\nonumber
& &\quad+ 4Q(a\, u_{(\sigma(1))}\,Q(a\,  
u_{(\sigma(2))}\,Q(a\,u_{(\sigma(3))}\,u_{(\sigma(4))})))\bigg),
\eeq
where $u_{(j)}=Qf_{(j)}$ and $\cell$ is the set
of permutations of the set $\{1,2,3,\dots,\ell\}$.

The results of Lemma \ref{lem: Lagrangian 1} can
be replaced by the results of \cite[Prop.\ 2.1]{GU1} as follows.
Using the same notations as in  Lemma \ref{lem: Lagrangian 1}, let
$Y=Y(x_0,\zeta_0;t_0,s_0)$, $K=K(x_0,\zeta_0;t_0,s_0)$, and  
$\Lambda_1=\Lambda(x_0,\zeta_0;t_0,s_0)$, and consider a source $f\in  
\I^{n+1}(Y)$.  Then $u=Qf$ satisfies
$u|_{M_0\setminus Y}\in
  \I^{n-1/2} ( M_0\setminus Y;\Lambda_1)$. Assume that  
$(x,\xi),(y,\eta)\in L^+M$  are on the same bicharacteristics of  
$\square_{\hat g}$,
  and $x<y$, that is, $((x,\xi),(y,\eta))\in \Lambda_{\hat g}^\prime$.  
Moreover, assume
  that $(x,\xi)\in N^*Y$.
Let   $\tilde b(x,\xi)$ be the principal
   symbol of $f$ at $(x,\xi)$ and
    $\tilde a(y,\eta)$ be the principal
   symbol of $u$ at $(y,\eta)$. Then $\tilde a(y,\eta)$
   depends linearly on $ \tilde f(x,\xi)$ and
    $\tilde a(y,\eta)$ vanishes if and only if
   $ \tilde f(x,\xi)$ vanishes.

Analogously to the Einstein equations,
we consider the indicator function
\beq\label{test sing 2}
\Theta_\tau^{(4)}=\bra F_{\tau},\M^{(4)}\cet_{L^2(U)},
\eeq
where
$\M^{(4)}$ is given by (\ref{4th interaction for wave eq B})
with  $u_{(j)}=Qf_{(j)}$, $j=1,2,3,4$, where  $f_{(j)}\in  
\I^{n+1}(Y(x_j,\xi_j;t_0,s_0))$, $n\leq -n_1$,
  and $F_\tau$ is the source producing a gaussian beam $Q^*F_\tau$
that propagates to the past along the geodesic $\gamma_{x_5,\xi_5}(\R_-)$,
see (\ref{Ftau source}).

Similar results to the ones given in Proposition \ref{lem:analytic limits A}
are
valid. Let us consider next the case when $(x_5,\xi_5)$ comes from the  
4-intersection
of  rays corresponding to $(\vec x,\vec \xi)=((x_j,\xi_j))_{j=1}^4$ and
$q$ is the corresponding intersection point, that is, $q=\gamma_{x_j,\xi_j}(t_j)$ for
all $j=1,2,3,4,5$.
Then 
\beq\label{indicator 2}
\Theta^{(4)}_\tau\sim 
\sum_{k=m}^\infty s_{k}\tau^{-k}
  \eeq
as $\tau\to \infty$ where   $m=-4n+4$.
Moreover,
let $b_j=(\dot\gamma_{x_j,\xi_j}(t_j))^\flat$ and
  $\bsequence=(b_{j})_{j=1}^5\in (T^*_q\hattuM _0)^5$,
  $w_j$ be the principal symbols of the waves $u_{(j)}$
  at $(q,b_j)$, and ${\bf w}=(w_j)_{j=1}^5$.
Then we see as in Proposition \ref{lem:analytic limits A}  that there is
  a real-analytic function $\mathcal G(\bsequence,{\bf w})$ such that
  the leading order term in (\ref{indicator})  satisfies
  \beq\label{definition of G 2}
s_{m}=
\mathcal  G(\bsequence,{\bf w}).
\eeq

The proof of Prop.\ \ref{singularities in Minkowski space} dealing
with the Einstein equations needs significant changes and we need to prove  
the following:

\begin{proposition}\label{singularities in Minkowski space for wave equation}
The  function $\ \mathcal  G(\bsequence,{\bf w})
$ given in (\ref{definition of G 2}) for the non-linear wave equation
is a non-identically vanishing real-analytic function.
\end{proposition}

\noindent{\bf Proof.}
Let us use the notations introduced in Prop.\ \ref{singularities in  
Minkowski space}.

As for the Einstein equations, we consider light-like vectors
\ba
b_5=(1,1,0,0),\quad b_j=(1,1-\frac  
12\rhoepsilon_j^2,\rhoepsilon_j+O(\rhoepsilon_j^3),\rhoepsilon_j^3),\quad  
j=1,2,3,4,
\ea
  in the Minkowski space $\R^{1+3}$, endowed with the standard metric $g=\diag(-1,1,1,1)$, where the terms
  $O(\rhoepsilon_k^3)$ are such that the vectors $b_j$, $j\leq 5$,
are  light-like. Then
\ba
g(b_5,b_j)= -\frac 12 \rhoepsilon_j^2,\quad
g(b_k,b_j)=-\frac 12 \rhoepsilon_k^2-\frac 12  
\rhoepsilon_j^2+O(\rhoepsilon_k\rhoepsilon_j).
\ea
Below, we denote $\omega_{kj}=g(b_k,b_j)$.
Note that if $\rhoepsilon_j<\rhoepsilon_k^4$, we have
$\omega_{kj}=-\frac 12 \rhoepsilon_k^2+O(\rhoepsilon_k^3).$

{For the wave equation,
we use different parameters $\rhoepsilon_j$ than for the Einstein  
equations, and
define (so, we use here the "unordered" numbering 4-2-1-3)
\beq\label{eq: ordering of epsilons wave equation}
\rhoepsilon_4=\rhoepsilon_2^{100},\  
\rhoepsilon_2=\rhoepsilon_1^{100},\hbox{ and  
}\rhoepsilon_1=\rhoepsilon_3^{100}.
\eeq
Below in this proof, we denote $\vec\rhoepsilon\to 0$ when
$\rhoepsilon_3\to 0$ and
$\rhoepsilon_4,$ $\rhoepsilon_2,$ and $\rhoepsilon_1$ are defined using
$\rhoepsilon_3$ as in (\ref{eq: ordering of epsilons wave equation}).

Let us next consider in Minkowski space
the coordinates $(x^j)_{j=1}^4$ such that
$K_j=\{x^j=0\}$ are light-like hyperplanes and  the waves $u_j=u_{(j)}$ that
satisfy in the Minkowski space $\square u_j=0$ and can be written
as
\ba
u_j(x)=\int_{\R}e^{i x^j \theta }a_j(x,\theta)\,d\theta,
\quad
a_j(x,\theta^{\prime})\in S^{n}(\R^4;\R\setminus 0),\quad j\leq 4,
\ea
and
\ba
u_\tau(x) =\chi(x^0)w_{(5)} \exp(i\tau b^{(5)}\,\cdotp x).
\ea
Note that the singular supports of the waves $u_j$, $j=1,2,3,4,$  
intersect then at the point
$\cap_{j=1}^4 K_j=\{0\}$.
Analogously to the definition  (\ref{definition of G}) we considered
  for  the Einstein equations, we
   define  the (Minkowski) indicator function
\ba
\mathcal G^{({\bf m})}(v,{\bf b})=
\lim_{\tau\to\infty} \tau^{m}(\sum_{\b\leq n_1}
\sum_{\sigma\in \Sigma(4)}
T^{({\bf m}),\b}_{\tau,\sigma}+\tilde T^{({\bf m}),\b}_{\tau,\sigma}),
\ea
where
\ba
T^{({\bf m}),\beta}_{\tau,\sigma}
&=&\bra Q_0(u_\tau\, \cdotp \a u_{\sigma(4)}), h\,\cdotp  \a  
u_{\sigma(3)}\,\cdotp Q_0(\a u_{\sigma(2)}\,\cdotp u_{\sigma(1)})\cet,\\
\tilde T^{({\bf m}),\beta}_{\tau,\sigma}
&=&\bra u_\tau,h\a\,Q_0(\a u_{\sigma(4)}\,\cdotp  \,u_{\sigma(3)})\,\cdotp
  Q_0(\a u_{\sigma(2)}\,\cdotp \, u_{\sigma(1)})\cet.
\ea

As for the Einstein equations, we
see that when $\alpha$ is equal to the value of the function $a(t,y)$

at the intersection point $q=0$ of the waves,
we have
  $\mathcal G^{({\bf m})}(v,{\bf b})=\mathcal G(v,{\bf b})$.

Similarly to the Lemma \ref{lem:analytic limits A}  we analyze next  
the functions
\ba
\Theta_\tau^{({\bf m})}=\sum_{\beta\in J_\ell}\sum_{\sigma\in  
\Sigma(4)}(T_{\tau,\sigma}^{({\bf m}),,\beta}+\tilde  
T_{\tau,\sigma}^{({\bf m}),\beta}).
\ea
Here $({\bf m})$ refers to "Minkowski".
We denote $T_{\tau}^{({\bf m}),\beta}=T_{\tau,id}^{({\bf m}),\beta}$
and  $\tilde T_{\tau}^{({\bf m}),\beta}=\tilde T_{\tau,id}^{({\bf m}),\beta}$.

Let us first consider the case when the permutation
$\sigma=id$. Then, as in the proof of Prop.\ \ref{singularities in  
Minkowski space},
in the case  when $\vec S^\beta =(Q_0,Q_0)$, we have
\ba
T^{({\bf m}),\beta}_\tau\\
&& \hspace{-2cm}=C_1 \det(A) \cdotp
(i\tau)^{m}(1+O(\frac 1\tau))
\vec\rhoepsilon^{\,2\vec n}
(\omega_{45}\omega_{12})^{-1}\rhoepsilon_4^{-4}\rhoepsilon_2^{-4}\rhoepsilon_1^{-4}\rhoepsilon_3^{2}\cdotp\P\\
&& \hspace{-2cm}=C_2\det(A) \cdotp
(i\tau)^{m}(1+O(\frac 1\tau))
\vec\rhoepsilon^{\,2\vec n}
\rhoepsilon_4^{-4-2}\rhoepsilon_2^{-4}\rhoepsilon_1^{-4-2}\rhoepsilon_3^{-2}\cdotp\P
\ea
where $\P$ is the product of the principal symbols of the waves $u_j$
at zero, $\vec\rhoepsilon^{\,2\vec n}=
\rhoepsilon_1^{2n}\rhoepsilon_2^{2n}\rhoepsilon_3^{2n}\rhoepsilon_4^{2n}$, and  
$C_1$ and $C_2$ are non-vanishing.
Similarly, a direct computation yields
\ba
\tilde T^{({\bf m}),\beta}_\tau\\
&& \hspace{-2cm}=C_1 \det(A) \cdotp
(i\tau)^{n}(1+O(\frac 1\tau))
\vec\rhoepsilon^{\,2\vec n}
(\omega_{43}\omega_{21})^{-1}\rhoepsilon_4^{-4}\rhoepsilon_2^{-4}\rhoepsilon_1^{-4}\rhoepsilon_3^{-4}\cdotp\P\\
&& \hspace{-2cm}=C_2\det(A) \cdotp
(i\tau)^{m}(1+O(\frac 1\tau))
\vec\rhoepsilon^{\,2\vec n}
\rhoepsilon_4^{-4}\rhoepsilon_2^{-4}\rhoepsilon_1^{-4-2}\rhoepsilon_3^{-4-2}\cdotp\P,
\ea
where again, $\P$ is the product of the principal symbols of the waves $u_j$
at zero and $C_1$ and $C_2$ are non-vanishing.

Considering formula (\ref{4th interaction for wave eq B}), we
see that for the wave equation we do not need to consider the
terms  that for the Einstein equations correspond to the
cases when $\vec S^\beta =(Q_0,I)$, $\vec S^\beta =(I,Q_0)$,
or $\vec S^\beta =(I,I)$ as the corresponding terms do not appear
in formula (\ref{4th interaction for wave eq B}).

Let us now consider permutations $\sigma$ of the indexes $(1,2,3,4)$
and compare the terms
\ba
& &L^{({\bf m}),\beta}_{\sigma}=\lim_{\tau\to \infty} \tau^{m}  
T^{({\bf m}),\beta}_{\tau,\sigma},\\
&& \tilde L^{({\bf m}),\beta}_{\sigma}=\lim_{\tau\to \infty}  \tau^{m} \tilde  
T^{({\bf m}),\beta}_{\tau,\sigma}.
\ea
Due to the presence of $\omega_{45}\omega_{12}$ in the above computations,
we observe that  all the terms
  $ \tilde  L^{({\bf m}),\beta}_{\tau,\sigma}/L^{({\bf m}),\beta}_{\tau,id}\to 0$
  as $\vec \rhoepsilon\to 0$, see (\ref{eq: ordering of epsilons wave  
equation}).
Also, if $\sigma\not=(1,2,3,4)$ and  $\sigma\not=\sigma_01=(2,1,3,4)$,
we see that  $ \tilde   
L^{({\bf m}),\beta}_{\tau,\sigma}/L^{({\bf m}),\beta}_{\tau,id}\to 0$
  as $\vec \rhoepsilon\to 0$.
  Also, we observe that $ L^{({\bf m}),\beta}_{\tau,\sigma_1}=
L^{({\bf m}),\beta}_{\tau,id}$.
Thus we see that the equal terms
$ L^{({\bf m}),\beta}_{\tau,\sigma_1}=
  L^{({\bf m}),\beta}_{\tau,id}$
that give the largest contributions as $\vec \rhoepsilon\to 0$
and that when $\P\not =0$ the sum
\ba
S( \vec \rhoepsilon,\P)=\sum_{\sigma\in \Sigma(4)}(
  L^{({\bf m}),\beta}_{\tau,\sigma}+\tilde L^{({\bf m}),\beta}_{\tau,\sigma})
\ea
is non-zero when
$\rhoepsilon_3>0$ is small enough   and
$\rhoepsilon_4,$ $\rhoepsilon_2,$ and $\rhoepsilon_1$ are defined using
$\rhoepsilon_3$ as in (\ref{eq: ordering of epsilons wave equation}).
As the indicator
function is real-analytic, this shows that the indicator function
in non-vanishing in a generic set.   \hfill \Box \medskip

We need also to change the
  the singularity {\it detection condition} ({D}) with light-like  
directions $(\vec x,\vec \xi)$
  as follows:  We define that point $y\in U_{\hat g}$,
   satisfies the singularity {\it detection condition} (${D}^\prime$)  
with light-like directions $(\vec x,\vec \xi)$
   and $t_0,\hat s>0$
  if
  \medskip  
  
($D^\prime$) For any $s,s_0\in (0,\hat s)$  there are  
$(x_j^{\prime},\xi_j^{\prime})\in \W_j(s;x_j,\xi_j)$, $j=1,2,3,4$,
and ${f}_{(j)}\in {\mathcal  
I_{C}}^{n+1}(Y((x_j^{\prime},\xi_j^{\prime});t_0,s_0))$, and such  
that if  $u_{\vec \e}$ of is the solution  of (\ref{PABC eq})
with the source ${f}_{\vec \e}=\sum_{j=1}^4
\e_j{f}_{(j)}$, then
the function
$\p_{\vec \e}^4u_{\vec \e}|_{\vec \e=0}$ is not
$C^\infty$-smooth in any neighborhood of $y$.
  \medskip

When condition (D) is replace by ($D^\prime$), the considerations in  
the Sections
  \ref{sec: normal coordinates} and \ref{subsection combining}  show  
that we can recover
the conformal class of the metric. This proves Theorem \ref{main  
thm3}. \hfill \Box \medskip

}}

\extension{

\section*{Appendix A: Reduced Einstein equation}

In this section we review known results on the Einstein equations and wave maps.

\subsection*{A.1.\ Summary of the used notations}


Let us recall some definitions given in Introduction.
Let $(\hattuM ,\hat g)$ be a $C^\infty$-smooth globally hyperbolic Lorentzian
manifold  and  $\tilde g$ be a $C^\infty$-smooth globally hyperbolic
metric on $M$ such that $\hat g<\tilde g$.
Let us  start by explaining  how one can construct a  $C^\infty$-smooth metric $\tilde g$ 
such that $\hat g<\tilde g$ and  $(M,\tilde g)$ is globally hyperbolic:
When $v(x)$ is an eigenvector corresponding to the negative eigenvalue
of $\hat g(x)$, we can choose a smooth, strictly positive function
$\eta:\hattuM \to \R_+$ such that $$\tilde g^\prime:=\hat g-\eta v\otimes v<\tilde g.$$ Then
$(\hattuM ,\tilde g^\prime)$ is globally hyperbolic, $\tilde g^\prime$ is smooth
and $\hat g<\tilde g^\prime$. Thus we can replace $\tilde g$ by the smooth metric
 $\tilde g^\prime$ having the same properties that are required for $\tilde g$.

Recall that there is an isometry $\Phi:(M,\tilde g)\to (\R\times N,\tilde  h)$,
where $N$
is a 3-dimensional manifold and the metric $\tilde  h$ can be written as
$\tilde  h=-\beta(t,y) dt ^2+\overline h(t,y)$ where $\beta:\R\times  N\to (0,\infty)$ is a smooth function and 
$\overline h(t,\cdotp)$ is a Riemannian metric on $ N$ depending smoothly
on $t\in \R$. 
As in the main text we identify these isometric manifolds and denote $\hattuM =\R\times N$.
%
Also, for $t\in \R$, recall that $\hattuM (t)=(-\infty,t)\times N$. We use parameters $t_1>t_0>0$
and denote $\hattuM _j=\hattuM (t_j)$, $j\in \{0,1\}$.
\modified{We use  the time-like geodesic    $\hat \mu =\mu_{\hat g}$, $\mu_{\hat g}:[-1,1]\to M_0$ on $(M_0,\hat g)$
and  the set $\K_j:=J^+_{\tilde g}(\hat \mu(-1))\cap M_j$ with $\hat \mu(-1)=
\in (-\infty,t_0)\times N$. Then 
$J^+_{\tilde g}(\hat \mu(-1))\cap M_j$ is compact. Also,
there exists $\e_0>0$ such that if 
$g$ is a Lorentz metric in $M_1$ such that 
$\|g-\hat g\|_{C^0_b(\hattuM _1;\hat g^+)}<\e_0$, then
$g|_{\K_1}<\tilde g|_{\K_1}$. In particular, this implies that we have  $J^+_{g}(p)\cap \hattuM _1\subset \K_1$
for all $p\in \K_1$. Later, we use this property to deduce that
when $g$ satisfies the $\hat g$-reduced Einstein equations 
in $M_1$, with a source that is supported in $ \K_1$ and has
small enough norm is a suitable space, then $g$ coincides with $\hat g$ in
$M_1\setminus \K_1$ and satisfies $g<\tilde g$.}

Let us use local coordinates on $\hattuM _1$ and denote by
$\nabla_k=\nabla_{X_k}$  the covariant derivative with respect to the metric $g$
in the direction $X_k=\frac \p{\p x^k}$
 and by
 $\hat \nabla_k=\hat \nabla_{X_k}$ the covariant derivative with respect to the metric $\hat g$
to the direction $X_k$.

\subsection*{A.2.\ Reduced Ricci and Einstein tensors}
Following \cite{FM} we recall that 
\beq\label{q-formula2copy AAA}
\Ric_{\mu\nu}(g)&=& \Ric_{\mu\nu}^{(h)}(g)
+\frac 12 (g_{\mu q}\frac{\p \Gamma^q}{\p x^{\nu}}+g_{\nu q}\frac{\p \Gamma^q}{\p x^{\mu}})
\eeq
where  $\Gamma^q=g^{mn}\Gamma^q_{mn}$,
\beq\label{q-formula2copyBa}
& &\hspace{-1cm}\Ric_{\mu\nu}^{(h)}(g)=
-\frac 12 g^{pq}\frac{\p^2 g_{\mu\nu}}{\p x^p\p x^q}+ P_{\mu\nu},
\\ \nonumber
& &\hspace{-2cm}P_{\mu\nu}=  
g^{ab}g_{ps}\Gamma^p_{\mu b} \Gamma^s_{\nu a}+ 
\frac 12(\frac{\p g_{\mu\nu }}{\p x^a}\Gamma^a  
+ \nonumber
g_{\nu l}  \Gamma^l _{ab}g^{a q}g^{bd}  \frac{\p g_{qd}}{\p x^\mu}+
g_{\mu l} \Gamma^l _{ab}g^{a q}g^{bd}  \frac{\p g_{qd}}{\p x^\nu}).\hspace{-2cm}
\eeq
Note that $P_{\mu\nu}$ is a polynomial of $g_{jk}$ and $g^{jk}$ and first derivatives of $g_{jk}$.
The harmonic  Einstein tensor is  
\beq\label{harmonic Ein}
\Ein^{(h)}_{jk}(g)=
\Ric_{jk}^{(h)}(g)-\frac 12 g^{pq}\Ric_{pq}^{(h)}(g)\, g_{jk}.
\eeq 
The  harmonic  Einstein tensor is extensively used to study
the Einstein equations in local coordinates where one can use the Minkowski
space $\R^4$ as the background space. To do global constructions
with a background space $(M,\hat g)$ one uses  the  reduced Einstein tensor.
The $\hat g$-reduced Einstein tensor 
 $\Ein_{\hat g} (g)$ and the $\hat g$-reduced  Ricci tensor
  $\Ric_{\hat g} (g)$  
are given 
by
\beq\label{Reduced Einstein tensora}
& &(\Ein_{\hat g} (g))_{pq}=(\Ric_{\hat g} (g))_{pq}-\frac 12 (g^{jk}(\Ric_{\hat g} g)_{jk})g_{pq},\\& &
\label{Reduced Ric tensora}
(\Ric_{\hat g} (g))_{pq}=\Ric_{pq} (g)-\frac 12 (g_{pn} \hat \nabla _q\hat F^n+ g_{qn} \hat \nabla _p\hat F^n)
\eeq
where $\hat F^n$ are the harmonicity functions given by
\beq\label{Harmonicity condition CCC}
\hat F^n=\Gamma^n-\hat \Gamma^n,\quad\hbox{where }
\Gamma^n=g^{jk}\Gamma^n_{jk},\quad
\hat \Gamma^n=g^{jk}\hat \Gamma^n_{jk},
\eeq
where $\Gamma^n_{jk}$ and $\hat \Gamma^n_{jk}$ are the Christoffel symbols
for $g$ and $\hat g$, correspondingly.
Note that $\hat \Gamma^n$ depends also on $g^{jk}$.
As $\Gamma^n_{jk}-\hat \Gamma^n_{jk}$ is the difference of two connection
coefficients, it is a tensor. Thus $\hat F^n$ is tensor (actually, a vector field), implying that both
$(\Ric_{\hat g} (g))_{jk}$ and $(\Ein_{\hat g} (g))_{jk}$ are  2-covariant tensors.
A direct calculation shows that  the $\hat g$-reduced Einstein tensor is the sum of the harmonic Einstein tensor 
and a term that is a zeroth order in $g$, 
\beq\label{hat g reduced einstein and reduced einstein}
(\Ein_{\hat g} (g))_{\mu\nu}=\Ein^{(h)}_{\mu\nu} (g)+\frac 12 (g_{\mu q}\frac{\p \hat \Gamma^q}{\p x^{\nu}}+g_{\nu q}\frac{\p \hat \Gamma^q}{\p x^{\mu}}).
\eeq
We also use the wave operator  
\beq\label{eq: wave ope}\\
\nonumber
\square_g\phi =
\sum_{p,q=1}^4(-\det(g(x)))^{-\frac 12}\frac  \p{\p x^p}
\left ((-\det(g(x)))^{\frac12}
g^{pq}(x)\frac \p{\p x^q}\phi(x)\right),
\eeq
can be written  for 
as 
\beq\label{eq: wave ope BB}\square_g \phi =g^{jk}\p_j\p_k \phi-g^{pq} \Gamma^n_{pq} \p_n \phi=
g^{jk}\p_j\p_k \phi- \Gamma^n \p_n \phi.
\eeq

%

%

%

\subsection*{A.3.\ Wave maps and reduced Einstein equations}
Let us consider the manifold $M_1=(-\infty,t_1)\times N$ with a $C^m $-smooth metric $g^{\prime}$,
$m\geq 8$,
 which
is a perturbation of the metric $\hat g$ and satisfies the Einstein equation
\beq\label{Einstein on M_0}
\Ein(g^{\prime})=T^{\prime}\quad \hbox{on }M_1,
\eeq
or equivalently, 
\ba
\Ric(g^{\prime})=\rho^{\prime},\quad \rho^{\prime}_{jk}=T_{jk}^{\prime}-\frac 12 ((g^{\prime})^{  nm}T^{\prime}_{nm})g^{\prime}_{jk}\quad \hbox{on }M_1.
\ea
Assume also that  $g^{\prime}=\hat g$ in the domain $A$,
where $A=M_1\setminus\K_1$ and  $\|g^{\prime}-\hat g\|_{C^2_b(M_1,\hat g^+)}<\e_0$,
so that $(M_1,g^{\prime})$ is globally hyperbolic. Note that then $T^{\prime}=\hat T$ in the set $A$ 
and that the metric $g^{\prime}$ coincides with $\hat g$
in particular in the set $M^-=\R_-\times N$

We recall next the considerations of \cite{ChBook}.
Let us  consider the Cauchy problem for the wave map
$f:(M_1,g^\prime)\to (M,\hat g)$, namely
\beq
\label{C-problem 1 AAA}& &\square_{g^{\prime},\hat g} f=0\quad\hbox{in } M_1,\\
\label{C-problem 2 AAA}& &f=Id,\quad \hbox{in  }\R_-\times N,
\eeq
where $M_1=(-\infty,t_1)\times N\subset M$. In (\ref{C-problem 1 AAA}), 
$\square_{g^{\prime},\hat g} f=g^\prime\,\cdotp\hat \nabla^2 f$ is the wave map operator, where $\hat \nabla$
is the covariant derivative of a map $(M_1,g^\prime)\to (M,\hat g)$, see \cite[Ch. VI, formula (7.32)]{ChBook}.
In  local coordinates
$X:V\to \R^4$ of $V\subset M_1$, denoted by
$X(z)=(x^j(z))_{j=1}^4$ and $Y:W\to \R^4$ of $W\subset \hattuM $, denoted by
$Y(z)=(y^A(z))_{A=1}^4$, 
the wave map $f:M_1\to \hattuM $ has the representation $Y(f(X^{-1}(x)))=(f^A(x))_{A=1}^4$ and
the wave map operator in equation (\ref{C-problem 1 AAA}) is given by
\beq\label{wave maps2}
& &(\square_{g^{\prime},\hat g} f)^A(x)= 
(g^{\prime})^{  jk}(x)\bigg(\frac \p{\p x^j}\frac \p{\p x^k} f^A(x) -\Gamma^{{\prime} n}_{jk}(x)\frac \p{\p x^n}f^A(x)
\\
& &\quad \quad\quad\quad\quad\quad\quad\quad\nonumber
+\hat \Gamma^A_{BC}(f(x))
\,\frac \p{\p x^j}f^B(x)\,\frac \p{\p x^k}f^C(x)\bigg)\eeq
where $\hat \Gamma^A_{BC}$ denotes the Christoffel symbols of metric $\hat g$
and  $\Gamma^{ {\prime} j}_{kl}$ are the Christoffel symbols of metric $g^{\prime}$.
When (\ref{C-problem 1 AAA})
is satisfied, we say that  $f$ is a wave map with respect to the pair $(g^{\prime},\hat g)$.

\extension{It follows from \cite[App.\ III, Thm. 4.2 and sec.\ 4.2.2] {ChBook},
that if $g^\prime \in C^m(M_0)$, $m\geq 5$ 
is sufficiently close to $\hat g$
in $C^m(M_0)$, then (\ref{C-problem 1 pre})-(\ref{C-problem 2 pre}) has a unique 
solution $f\in C^0([0,t_1];H^{m-1}(N))\cap  C^1([0,t_1];H^{m-2}(N))$.
This comes from the fact that the Christoffel symbols of $g^\prime$
are in $ C^{m-1}(M_0)$.
Moreover, when $m$ is even, using 
\cite[Thm.\ 7]{Kato1975} 
for $f$ and $\p_t^p f$,  we see that
the solution $f$ is in $f\in \cap_{p=0}^{m-1}C^p([0,t_1];H^{m-1-p}(N))\subset C^{m-3}(M_0)$ and
$f$ depends in $C^{m-3}(M_0)$ continuously on $g^\prime\in C^m(M_0)$.
We note that these smoothness results for $f$ are not optimal.}  

The wave map operator
$\square_{g^{\prime},\hat g} $ is a coordinate invariant operator.
The important property of the wave maps is that, if 
$f$ is wave map with respect to the pair $(g^{\prime},\hat g)$ and
$g=f_*g^{\prime}$ 
then, as follows from (\ref{wave maps2}), the identity map $Id:x\mapsto x$ is a wave map with respect to the pair $(g,\hat g)$
and, the wave map equation for the identity map is equivalent to (cf.\ \cite[p.\ 162]{ChBook})
\beq\label{wave equation id}
\Gamma^n=\hat\Gamma^n,\quad \hbox{where }\Gamma^n=g^{jk}\Gamma^n_{jk},
\quad \hat \Gamma^n=g^{jk}\hat\Gamma^n_{jk}
\eeq
%
where the Christoffel symbols $\hat\Gamma^n_{jk}$ of the metric $\hat g$ are smooth functions.

Since  $g=g^{\prime}$ outside a compact set $\K_1\subset (0,t_1)\times N$,
we see that   this Cauchy problem is equivalent to the same equation
restricted to the
set $(-\infty,t_1)\times B_0$, where $B_0\subset N$ is an open
relatively compact set such that $\K_1\subset (0,t_1]\times B_0$
with the boundary condition $f=Id$ on  $(0,t_1]\times \p B_0$.
%
 Moreover, by the uniqueness of the wave map, we have
$f|_{M_1\setminus \K_1}=id$ so that
 $f(\K_1)\cap M_0\subset \K_0$.

 As 
 the inverse function of the wave map $f$ depends continuously,
 in
$C^{m-3}_b([0,t_1]\times N,g^+)$, on the metric $g^{\prime} \in C^m(M_0)$
we can also assume that $\e_1$ is  so small that
$\hattuM _0\subset f(M_1)$.

Denote next $g:=f_*g^{\prime}$, $T:=f_*T^{\prime}$, and $\rho:=f_*\rho^{\prime}$ 
and define $\hat \rho=\hat T-\frac 12 (\hbox{Tr}\, \hat T)\hat g.$ 
Then $g$ is $C^{m-6}$-smooth and the equation (\ref{Einstein on M_0}) implies 
\beq\label{Einstein on M1}
\Ein(g)=T\quad \hbox{on }\hattuM _0.
\eeq 
Since $f$ is a wave map and $g=f_*g^\prime$, we have that
the identity map is a $(g,\hat g)$-wave map and thus $g$ satisfies (\ref{wave equation id}) and thus  by the
definition of the reduced Einstein tensor,  (\ref{Reduced Ric tensor})-(\ref{Reduced Einstein tensor}), we have
\ba
\Ein_{pq}(g)=
(\Ein_{\hat g} (g))_{pq}\quad \hbox{on } \hattuM _0.
\ea
This and  (\ref{Einstein on M1}) yield the $\hat g$-reduced Einstein equation
\beq\label{hat g reduced einstein equations}
(\Ein_{\hat g} (g))_{pq}=T_{pq}\quad \hbox{on } \hattuM _0.
\eeq
This equation is useful for our considerations as it is a quasilinear,
hyperbolic equation on $\hattuM _0$. Recall that 
 $g$ coincides with
$\hat g$ in $M_0\setminus \K_0$. The unique solvability of this 
Cauchy problem is studied in e.g.\ \cite[Thm.\ 4.6 and 4.13]{ ChBook}, \cite {HKM} 
and  Appendix B below.


\subsection*{A.4.\ Relation of the reduced Einstein equations and the original Einstein equation}

The metric $g$ which solves the $\hat g$-reduced Einstein
equation $\Ein_{\hat g} (g)=T$ is a solution of the original
Einstein equations  $\Ein (g)=T$ if the harmonicity
functions $\hat F^n$
vanish identically. Next we recall the result that 
the harmonicity functions vanish on $\hattuM _0$
when 
\beq\label{eq: good system}
& &(\Ein_{\hat g}(g))_{jk}=T_{jk},\quad \hbox{on }M_0,\\
\nonumber & &\nabla_pT^{pq}=0,\quad \hbox{on }M_0,\\
\nonumber & &g=\hat g,\quad \hbox{on }M_0\setminus \K_0.
\eeq
To see this, let
us  denote $\Ein_{jk}(g)=S_{jk}$, $S^{jk}=g^{jn}g^{km}S_{nm}$,
and $T^{jk}=g^{jn}g^{km}T_{nm}$. 
Following standard arguments,
see \cite{ ChBook}, we see from (\ref{Reduced Einstein tensor}) that in local coordinates
\ba
S_{jk}-(\Ein_{\hat g}(g))_{jk}=\frac12(g_{jn}\hat \nabla_k \hat F^n+g_{kn}\hat \nabla_j\hat F^n-g_{jk}\hat \nabla_n \hat  F^n).
\ea
Using equations (\ref{eq: good system}), the Bianchi identity  $\nabla_pS^{pq}=0$,
and the basic property of Lorentzian connection,
$\nabla_kg^{nm}=0$,
we obtain 
\ba
0&=&
2\nabla_p(S^{pq}-T^{pq})\\
&=&\nabla_p(g^{qk}\hat  \nabla_k F^p+g^{pm}\hat\nabla_m \hat F^q-  g^{pq}\hat  \nabla_n \hat F^n)
\\
&=&g^{pm}\nabla_p \hat \nabla_m \hat F^q+(g^{qp} \nabla_n \hat \nabla_p \hat F^n-  g^{qp} \nabla_p \hat \nabla_n \hat F^n)\\
&=&g^{pm}\nabla_p \hat \nabla_m \hat F^q+W^q(\hat F)
\ea
where $\hat F=(\hat F^q)_{q=1}^4$ and
the operator 
\ba
W:(\hat F^q)_{q=1}^4\mapsto (g^{qk} (\nabla_p \hat \nabla_k \hat F^p- \nabla_k \hat \nabla_p \hat F^p))_{q=1}^4
\ea
is a linear first order differential operator which coefficients are polynomial functions of
$\hat g_{jk}$, $\hat g^{jk}$, $g_{jk}$, $g^{jk}$ and their first derivatives.

Thus the harmonicity functions $\hat F^q$ satisfy on $\hattuM _0$ the hyperbolic initial 
value problem 
\ba
& &g^{pm}\nabla_p \hat \nabla_m \hat F^q+W^q(\hat F)=0,\quad\hbox{on }M_0,\\
& & \hat F^q=0,\quad \hbox{on }M_0\setminus \K_0,
\ea
and as this initial 
Cauchy problem is uniquely solvable by \cite[Thm.\ 4.6 and 4.13]{ ChBook}
or \cite {HKM}, we see that $\hat F^q=0$ on $\hattuM _0$. Thus
equations (\ref{eq: good system})
yield
that the Einstein equations  $\Ein(g)=T$ hold on $\hattuM _0$.

We note that in the $(g,\hat g)$-wave map coordinates, where $\hat F^q=0$, the
wave operator (\ref{eq: wave ope BB}) has the form
\beq\label{eq: wave operator in wave gauge}
\square_g \phi=g^{jk}\p_j\p_k \phi-g^{pq}\hat \Gamma^n_{pq} \p_n \phi.
\eeq
Thus, the scalar field equation $\square_g \phi-m^2\phi=0$
does not involve derivatives of $g$.

 \subsection*{A.5.\ Linearized Einstein-scalar field equations} 
 
 Next we consider  the linearized equation are obtained as the derivatives
 of the solutions of the non-linear, generalized Einstein-matter field
equations 
(\ref{eq: adaptive model with no source}). 

 
{Observe that if  a family $\mathcal F_\e=(\mathcal F^1_\e,\mathcal F^2_\e)$
of sources and a family  $(g_\e,\phi_\e)$ of functions 
are solutions of  the non-linear reduced Einstein-scalar field  
equations 
(\ref{eq: adaptive model with no source})
that depend smoothly on $\e\in [0,\e_0)$ in $C^{15}(M_0)$
and satisfy $\mathcal F_\e|_{\e=0}=0$, $(g_\e,\phi_\e)|_{\e=0}=(\hat g,\hat \phi)$,
then $\p_\e  \mathcal F_\e|_{\e=0}=f=(f^1,f^2)$, and  $\p_\e  (g_\e,\phi_\e)|_{\e=0}=(\dot g,\dot \phi)$,
satisfy
the linearized version of 
the equation (\ref{eq: adaptive model with no source}) that has the form (in the local coordiantes),
\beq\label{linearized eq: adaptive model with no source BB}
& &\square_{\hat g} \dot g_{jk}+A_{jk}(\dot g,\dot \phi,\p \dot g,\p \dot \phi)=f^1_{jk},\quad
\hbox{in }M_0,\\
\nonumber
& &\square_{\hat g}\dot \phi_\ell +
B_\ell(\dot g,\dot \phi,\p \dot g,\p \dot \phi)
=f^2_\ell,\quad \ell=1,2,3,\dots,L.
\eeq
Here $A_{jk}$ and $B_\ell$ are first order linear differential
operators whose coefficients depend on $\hat g$ and $\hat \phi$.
Let us write these equations in more explicit form. We see that
the linearized reduced Einstein tensor is in local coordinates of the form 
\ba
e_{pq}(\dot g)&:=&\p_\e (\Ein_{\hat g} g_\e)_{pq}|_{\e =0}\\
&=&
-\frac 12 \hat g^{jk}\hat \nabla_j\hat \nabla_k  \dot g_{pq}
+\frac 14(  \hat g^{nm} \hat g^{jk}\hat \nabla_j\hat \nabla_k  \dot g_{nm}) \hat g_{pq}
\\
& &+A_{pq}^{abn}\hat \nabla_n\dot g_{ab}+B_{pq}^{ab}\dot g_{ab},
\ea
where $A_{pq}^{abn}(x)$ and $B_{pq}^{ab}(x)$ depend on  $\hat g_{jk}$
and its derivatives (these terms can be computed explicitly using  (\ref{q-formula2copyBa})).
The linearized scalar field stress-energy tensor is the linear first order differential operator 
\ba
& &t^{(1)}_{pq}(\dot g)+t^{(2)}_{pq}(\dot \phi):=\p_\e ({\bf T}_{jk}(g_\e,\phi_\e)
)
|_{\e =0}\\
&=&
\sum_{\ell=1}^L\bigg(
\p_j\hat \phi_\ell \,\p_k\dot \phi_\ell +
\p_j\dot \phi_\ell \,\p_k\hat  \phi_\ell \\
& &\quad
-\frac 12 \dot g_{jk} \hat  g^{pq}\p_p \hat  \phi_\ell \,\p_q \hat \phi_\ell
-\frac 12 \hat g_{jk} \dot g^{pq}\p_p \hat  \phi_\ell \,\p_q \hat \phi_\ell\\
& &\quad
-\frac 12 \hat g_{jk} \hat g^{pq}\p_p\dot \phi_\ell \,\p_q \hat \phi_\ell
-\frac 12 \hat g_{jk} \hat  g^{pq}\p_p \hat  \phi_\ell \,\p_q\dot  \phi_\ell\\
& &\quad
- m^2\hat  \phi_\ell \dot \phi_\ell\hat g_{jk}
-\frac 12 m^2\hat \phi_\ell^2\dot g_{jk}\bigg).
\ea
Thus when $\mathcal F_\e=(\mathcal F^1_\e,\mathcal F^2_\e)$
is a family 
of sources and   $(g_\e,\phi_\e)$ a family of functions 
that satisfy  the non-linear reduced Einstein-scalar field  
equations 
(\ref{eq: adaptive model with no source}),
the  $\e$-derivatives $\dot u=(\dot g,\dot \phi)$ 
and $\p_\e  \mathcal F_\e|_{\e=0}=f=(f^1,f^2)$
satisfy  in local $(g,\hat g)$-wave coordinates
\beq\label{eq: final linearization}
& &
e_{pq}(\dot g)-t^{(1)}_{pq}(\dot g)-t^{(2)}_{pq}(\dot \phi)=f^1_{pq},\\
\nonumber
& &\square_{\hat g}\dot \phi_\ell
-\hat g^{nm}\hat g^{kj}(\p_n\p_j\hat \phi_\ell+ \hat \Gamma_{nj}^p\p_p
\hat \phi_\ell)\dot g_{mk}-m^2\phi_\ell=f^2_\ell.
\eeq
We call this the linearized Einstein-scalar field equation. 
 }

\subsection*{A.6.\ Linearization of the conservation law.}

Assume that $u=(g,\phi)$ and $\F=(\F^1,\F^2)$  
satisfy equation (\ref{eq: adaptive model with no source}).
Then
  the conservation law (\ref {conservation law0})
 gives for all $j=1,2,3,4$ equations (see \cite[Sect. 6.4.1]{ChBook}) 
 \ba
0&=&\frac 12\nabla_p^g (g^{pk}T_{jk})\\
&=&\frac 12\nabla_p^g (g^{pk}({\bf T}_{jk}(g,\phi)+\F^1_{jk}))\\
&=& \sum_{\ell=1}^L(g^{pk} \nabla^g_p\p_k\phi_\ell )\,\p_j\phi_\ell 
-(m^2_\ell \phi_\ell \p_p \phi_\ell) \delta_j^p +\frac 12 \nabla^g_p (g^{pk}\F^1_{jk})
\\
&=& \sum_{\ell=1}^L(g^{pk} \nabla^g_p\p_k\phi_\ell -m^2_\ell \phi_\ell ) \,\p_j \phi_\ell +\frac 12 \nabla^g_p (g^{pk}\F^1_{jk}).
\ea
This yields by (\ref{eq: adaptive model with no source})
\beq
\label{mod. conservation law aaa}
\frac 12 g^{pk} \nabla^g_p \F^1_{jk}+\sum_{\ell=1}^L \F^2_\ell \,\p_j\phi_\ell=0,
\quad j=1,2,3,4.
\eeq
Next assume that $u_\e=(g_\e,\phi_\e)$ and $\F_\e$
satisfy equation (\ref{eq: adaptive model with no source}) and
$C^1$-smooth functions of $\e\in (-\e_0,\e_0)$ taking values in $H^1(N)$-tensor fields, and $(g_\e,\phi_\e)|_{\e=0}=(\hat g,\hat \phi)$
and  $\F_\e|_{\e=0}=0$. 
 Denote $(f^1,f^2)=\p_\e\F_\e|_{\e=0}$. 
 Then by taking $\e$-derivative of (\ref{mod. conservation law aaa}) at $\e=0$
 we get
\beq
\label{linearized conservation law aaa}
\frac 12 \hat g^{pk} \hat \nabla_p f^1_{jk}+\sum_{\ell=1}^L f^2_\ell \,\p_j\hat \phi_\ell=0,
\quad j=1,2,3,4.
\eeq
We call this the linearized conservation law.

\section*{Appendix B: Stability and existence of the  direct problem}

 Let us next consider existence and stability of the solutions of the
 Einstein-scalar field equations.
Let $t={\bf t}(x)$ be local time so that there is a diffeomorphism
$\Psi:M\to \R\times N$, $\Psi(x)=({\bf t}(x),Y(x))$, and
 $$S(T)=\{x\in \hattuM _0; \ {\bf t}(x)=T\},\quad T\in \R$$ are Cauchy surfaces.
 Let $t_0>0$. Next we identify $M$ and $\R\times N$ via the map $\Psi$ and just
 denote $M=\R\times N$.
Let us denote $$M(t)=(-\infty,t)\times N,\quad M_0=M(t_0).$$
By \cite[Cor.\ A.5.4]{BGP} the set $\K=J^+_{\tilde g}(\hat p^-)\cap \hattuM _0$,
where $\hat p^0\in M_0$,
is compact.
Let $N_1,N_2\subset N$ be such
 open relatively compact sets with smooth boundary so that
 $N_1\subset N_2$ and  $Y(J_{\tilde g}^+(\hat p^0)\cap \hattuM _0)\subset N_1$.  
 
 To simplify citations to the existing literature, 
let us define $\tilde N$ to be a compact manifold without boundary such
that $N_2$ can be considered as a subset of $\tilde N$. Using a construction based on
a suitable partition of unity, the Hopf double of the manifold $N_2$, 
and the Seeley extension of the metric tensor,  we can
endow $\tilde M=(-\infty,t_0)\times \tilde N$ with a smooth Lorentzian
metric $\hat g^e$ (index $e$ is for "extended")
so that $\{t\}\times \tilde N$ are Cauchy surfaces of $\tilde N$
and that $\hat g$ and $\hat g^e$ coincide in  the set $\Psi^{-1}((-\infty,t_0)\times N_1)$ that contains the set $J^+_{\hat g}(\hat p^0)\cap \hattuM _0$. We extend the metric $\tilde g$ to a (possibly non-smooth) globally
hyperbolic metric $\tilde g^e$ on $\tilde M_0=(-\infty,t_0)\times \tilde N)$
such that $\hat g^e<\tilde g^e$.

 To simplify notations below we denote $\hat g^e=\hat g$ and
 $\tilde g^e=\tilde g$ 
 on the whole
 $\tilde M_0$.
 Our aim is to prove the estimate (\ref{eq: Lip estim}).

Let us denote by $t={\bf t}(x)$ the local time. 
Recall that when $(g,\phi)$ is a solution of 
the scalar field-Einstein equation, we denote $u=(g-\hat g,\phi-\hat \phi)$.
We will consider the equation for $u$, and to emphasize that the metric
depends on $u$, we denote $g=g(u)$ and assume below that  both the metric $g$ is  dominated
by $\tilde g$, that is, $ g<\tilde g$.
We use the pairs $${\bf u}(t)=(u(t,\cdotp),\p_t u(t,\cdotp))\in
H^1(\tilde N)\times L^2(\tilde N)$$ the notations
 ${\bf v}(t)=(v(t),\p_t v(t))$ etc. 
We consider sources $\F$ and $F$, cf. Appendix C.
Let us consider a generalization of the system (\ref{eq: notation for hyperbolic system 1}) of the form 
\beq\label{eq: notation for hyperbolic system}
& &\square_{g(u)}u+V(x,D)u+H(u,\p u) =R(x,u,\p u,F,\p F)+\F,\ \ x\in \tilde M_0,\hspace{-1cm}\\
& & \nonumber \supp(u)\subset \K,
\eeq
where $\square_{g(u)}$ is the Lorentzian Laplace operator operating on the sections 
of the bundle $\B^L$ on $M_0$ and 
$\supp(F)\cup \supp(\F)\subset \K$.
Note that above  $u=(g-\hat g,\phi-\hat \phi)$ and $g(u)=g$.
Also, $F\mapsto R(x,u,\p u,F,\p F)$ is a 
 non-linear first order differential operator where $R$ is a smooth function
 that depends on all variables $x$,
 $u(x)$, $\p_j u(x)$, $F(x),$ and $\p_jF(x)$ smoothly and also on the  derivatives 
of $(\hat g,\hat \phi)$ at $x$. Let
\ba 
V(x,D)=V^j(x)\p_j+V(x)
\ea 
be a linear first order differential operator whose coefficients 
at $x$ are depending smoothly on the  derivatives 
of $(\hat g,\hat \phi)$  at $x$, and finally, $H(u,\p u)$ is a polynomial
of $u(x)$ and $\p_j u(x)$ which coefficients 
at $x$ are depending smoothly on the  derivatives 
of $(\hat g,\hat \phi)$   such that
 $$\p^\a_v\p^\beta_w H(v,w)|_{v=0,w=0}=0\hbox { for }
|\a|+|\b|\leq 1.$$ 
By \cite[Lemma 9.7]{Ringstrom}, the equation  (\ref{eq: notation for hyperbolic system})
has at most one solution with given  $C^2$-smooth source functions $F$ and $\F$.
Next we consider the existence of $u$ and its dependency on $F$ and $\F$.


%

Below we use notations, c.f. (\ref{eq: notation for hyperbolic system 1}) and
(\ref{eq: notation for hyperbolic system 1b}) and Appendix C, \ba
{\mathcal R}({\bf u},F)=R(x,u(x),\p u(x),F(x),\p F(x)),\quad
{\mathcal H}({\bf u})=H(u(x),\p u(x)).\ea
Note that $u=0$, i.e., $g=\hat g$ and $\phi=\hat\phi$ satisfies (\ref{eq: notation for hyperbolic system}) with $F=0$
and $\F=0$.
Let us use the same notations as in  \cite {HKM} cf. also \cite[section 16]{Kato1975}, 
 to consider quasilinear wave equation on $[0,t_0]\times \tilde N$.
Let $\H^{(s)}(\tilde N)=H^s(\tilde N)\times H^{s-1}(\tilde N)$ and 
\ba
Z=\H^{(1)}(\tilde N),\quad  Y=\H^{(k+1)}(\tilde N),\quad X=\H^{(k)}(\tilde N).
\ea
The norms on these space are defined invariantly using the smooth Riemannian 
metric $h=\hat g|_{\{0\}\times \tilde N}$ on $\tilde N$.
Note that $\H^{(s)}(\tilde N)$ are in fact the Sobolev spaces of sections on the bundle $\pi:\B_K\to
\tilde N$, where $\B_K$ denotes also the  pull back bundle of $\B_K$ on $\tilde M$ in the map
$id:\{0\}\times \tilde N\to \tilde M_0$,
or on the bundle $\pi:\B_L\to \tilde N$. Below, $\nabla_h$ denotes the  standard connection of the bundle  $\B_K$
or  $\B_L$ associated to the metric $h$.

Let $k\geq 4$ be an even integer.
By definition of $H$ and $R$ we see that  there are  $0<r_0<1$ and $L_1,L_2>0$,
all depending on $\hat g$, $\hat \phi$, $\K$, and $t_0$, such that if 
$0<r\leq r_0$ and
\beq\label{new basic conditions}
& &\|{\bf v} \|_{C([0,t_0];\H^{(k+1)}(\tilde N))}\leq r,\quad \|{\bf v}^\prime\|_{C([0,t_0];\H^{(k+1)}(\tilde N))}\leq r,\\
& &\|F\|_{C([0,t_0];H^{(k+1)}(\tilde N))}\leq r^2,\quad \|\F\|_{C([0,t_0];H^{(k)}(\tilde N))}\leq r^2 \nonumber\\
& &\|F^\prime\|_{C([0,t_0];H^{(k+1)}(\tilde N))}\leq r^2,\quad \|\F^\prime\|_{C([0,t_0];H^{(k)}(\tilde N))}\leq r^2 \nonumber
\eeq
then
\beq\label{new f1f2 conditions}
& &\quad \quad \|g(\cdotp;v)^{-1}\|_{C([0,t_0];H^s(\tilde N))}\leq L_1 ,\\ \nonumber
& &\|{\mathcal H} ({\bf v})\|_{C([0,t_0];H^{s-1}(\tilde N))}\leq L_2r^2,
\quad \|{\mathcal H} ({\bf v}^\prime)\|_{C([0,t_0];H^{s-1}(\tilde N))}\leq L_2r^2,\\ \nonumber
& &\|{\mathcal H} ({\bf v})-{\mathcal H} ({\bf v}^\prime)\|_{C([0,t_0];H^{s-1}(\tilde N))}\leq  L_2r\, \|{\bf v}-{\bf v}^\prime\|_{C([0,t_0];\H^{(s)}(\tilde N))}, 
\\ \nonumber
& &\|\mathcal R ({\bf v}^\prime,F^\prime)\|_{C([0,t_0];H^{s-1}(\tilde N))}\leq L_2r^2,
\quad \|\mathcal R ({\bf v},F)\|_{C([0,t_0];H^{s-1}(\tilde N))}\leq L_2r^2,\hspace{-1cm} \\ \nonumber
& &\|\mathcal R ({\bf v},F)-\mathcal R ({\bf v}^\prime, F^\prime)\|_{C([0,t_0];H^{s-1}(\tilde N))}
\\ \nonumber & &\leq  
L_2r\, \|{\bf v}-{\bf v}^\prime\|_{C([0,t_0];\H^{(s)}(\tilde N))}+L_2\|F-{F}^\prime\|_{C([0,t_0];H^{s+1}(\tilde N))\cap C^1([0,t_0];H^{s}(\tilde N))},\hspace{-2cm}
%
  \eeq
 for all $s\in [1,k+1]$. 
%

%
%

Next we write (\ref{eq: notation for hyperbolic system}) as a first order
system.  To this end,
let $\mathcal A(t,{\bf v}):\H^{(s)}(\tilde N)\to \H^{(s-1)}(\tilde N)$ be the operator 
$$\mathcal A(t,{\bf v})=\mathcal A_0(t,{\bf v})+\mathcal A_1(t,{\bf v})$$
where in local coordinates and in the local trivialization of the bundle $\B^L$
\ba
\mathcal A_0(t,{\bf v})=-\left(\begin{array}{cc} 0 & I\\
 \frac 1{g^{00}(v)}\sum_{j,k=1}^3 g^{jk}(v)\frac \p{\p x^j}\frac\p{ \p x^k}\quad &
 \frac 1{g^{00}(v)}\sum_{m=1}^3 g^{0m}(v)\frac {\p}{\p x^m} \end{array}\right)
\ea 
with $g^{jk}(v)=g^{jk}(t,\cdotp;v)$ is a function on $\tilde N$ and 
\ba
\mathcal A_1(t,{\bf v})= \frac {-1}{g^{00}(v)}\left(\begin{array}{cc} 0 & 0\\
\sum_{j=1}^3 
 B^j(v)
   \frac \p{\p x^j}\quad &
B^{0}(v)\end{array}\right)
\ea 
where $B^j(v)$ depend on $v(t,x)$ and its first derivatives, and 
the connection coefficients (the Christoffel  symbols) corresponding  to $g(v)$. We denote $\S=(F,\F)$ and 
 \ba
& &f_\S(t,{\bf v})=(f_\S^1(t,{\bf v}),f_\S^2(t,{\bf v}))\in \H^{(k)}(\tilde N),
\quad\hbox{where}\\
& & f_\S^1(t,{\bf v})=0,\quad f_\S^2(t,{\bf v})=\mathcal R({\bf v},F)(t,\cdotp)-\mathcal H({\bf v})(t,\cdotp)+\F(t,\cdotp).
 \ea
 
 Note that when (\ref{new basic conditions}) are satisfied with $r<r_0$, inequalities (\ref{new f1f2 conditions}) imply that there exists $C_2>0$ so that
 \beq\label{fF bound r2}
 & &\|f_\S(t,{\bf v})\|_Y+\|f_{\S^\prime}(t,{\bf v}^\prime)\|_Y\leq C_2r^2,\\
 & &\|f_\S(t,{\bf v})-f_\S(t,{\bf v}^\prime)\|_Y\leq C_2
 r\,\|{\bf v}-{\bf v}^\prime\|_{C([0,t_0];Y)}. \nonumber
 \eeq

%
%

Let $U^{\bf v}(t,s)$ be the wave propagator corresponding to metric $g(v)$, that is,
$U^{\bf v}(t,s): {\bf h}\mapsto {\bf w}$, where ${\bf w}(t)=(w(t),\p_t w(t))$  solves
\ba
(\square_{g(v)}+V(x,D))w=0\quad\hbox{for } (t,y)\in [s,t_0]\times \tilde N,\quad\hbox{with }
 {\bf w}(s,y)={\bf h}.
 \ea 
 Let $$S=(\nabla_{h}^*\nabla_{h}+1)^{k/2}:Y\to Z$$ be an isomorphism.
As  $k$ is an even integer, we see using multiplication estimates for Sobolev
spaces, see e.g.\
  \cite [Sec.\ 3.2, point (2)]{HKM}, that there exists $c_1>0$ (depending on
  $r_0,L_1,$ and $L^2$) so that  $$\A(t,{\bf v})S-S\A(t,{\bf v})=C(t,{\bf v}),$$
 where $\|C(t,{\bf v})\|_{Y\to Z}\leq c_1$ for all ${\bf v}$ satisfying
 (\ref{new basic conditions}). This yields that 
 the property (A2) in  \cite {HKM} holds, namely
 that $$S\A(t,{\bf v})S^{-1}=\A(t,{\bf v})+B(t,{\bf v})$$ where
$B(t,{\bf v})$ extends to a bounded operator in $Z$
 for which $$\|B(t,{\bf v})\|_{Z\to Z}\leq c_1$$ for all ${\bf v}$ satisfying
 (\ref{new basic conditions}).
Alternatively, to see the mapping properties of $B(t,{\bf v})$ we could use the fact that 
 $B(t,{\bf v})$ is a zeroth order pseudodifferential operator with
$H^{k}$-symbol.\hiddenfootnote{On mapping properties of such
pseudodifferential operators, see J. Marschall, Pseudodifferential operators with coefficients in Sobolev spaces. Trans. Amer. Math. Soc. 307 (1988), no. 1, 335-361.} 
 
 Thus the proof of  \cite [Lemma 2.6]{HKM} shows\hiddenfootnote{ Alternatively, as $U^{\bf v}(t,s)$ are propagators
  for linear wave equations having finite speed of wave propagation, one can prove 
the  estimates (\ref{U bounds}) by considering the wave equation
in local coordinate neighborhoods $W_j\subset W_j^\prime \subset \tilde N$,
$\Phi_j:  W_j^\prime\to \R^3$ and a partition of unity $\psi_j\in C^\infty_0(W_j)$ and
cut-off functions $\psi^\prime_j\in C^\infty_0(W_j^\prime)$.
Then  \cite [Lemma 2.6]{HKM}, used  in local coordinates, shows
that when $t_k-t_{k-1}$ is small enough then 
\ba 
\| \psi^\prime_j U^{\bf v}(t_k,t_{k-1})(\psi_j {\bf w}) \|_{Y}\leq C_{3,jk}^{\prime}e^{C_4(t-s)}\| \psi_j {\bf w} \|_{Y}.
\ea
Using  sufficiently small time steps $t_k-t_{k-1}$ and combining the estimates
for different $j$:s
together, one obtain estimates (\ref{U bounds}).}
 that 
%
there is a constant  $C_3>0$
   so that 
  \beq
  \label{U bounds}\|U^{\bf v}(t,s)\|_{Z\to Z}\leq C_3\quad\hbox{and}\quad
  \|U^{\bf v}(t,s) \|_{Y\to Y}\leq C_3\hspace{-1.5cm}
  \eeq 
  for $0\leq s<t\leq t_0$.
    By interpolation of estimates (\ref{U bounds}), we see also that 
    \beq
  \label{U bounds2}
  \|U^{\bf v}(t,s)\|_{X\to X}\leq C_3,
  \eeq  for $0\leq s<t\leq t_0$.

Let us next modify the reasoning given in \cite{Kato1975}: let  $r_1\in (0,r_0)$
be a parameter which value will be chosen later,  $C_1>0$
and 
$E$ be the space of functions ${\bf u}\in C([0,t_0];X)$ for which
\beq
\label{1 bounded}& &\|{\bf u}(t)\|_Y\leq r_1\quad\hbox{and}\\
\label{2 Lip}& &\|{\bf u}(t_1)-{\bf u}(t_2)\|_X\leq C_1|t_1-t_2|\eeq
for all $t,t_1,t_2\in [0,t_0]$.
The set $E$ is endowed by the metric of  $C([0,t_0];X)$. We 
note that  by \cite[ Lemma 7.3] 
{Kato1975},
a convex $Y$-bounded, $Y$-closed set is closed also in $X$.
Similarly, functions $G:[0,t_0]\to X$  satisfying (\ref{2 Lip}) 
form a closed subspace of $C([0,t_0];X)$. Thus $E\subset X$ is a closed set implying
that $E$ is a complete metric space.

Let
 \ba 
 & &\hspace{-1cm}W=\{(F,\F)\in C([0,t_0];H^{k+1}(\tilde N))
 \times C([0,t_0];H^{k}(\tilde N));
 \\
 & &\quad\quad
 \ \sup_{t\in [0,t_0]}\|F(t)\|_{H^{k+1}(\tilde N)}+\|\F(t)\|_{H^{k}(\tilde N)}< r_1\}.\hspace{-1cm}
 \ea

Following \cite[p. 44]{Kato1975}, we see that
the solution of
equation (\ref{eq: notation for hyperbolic system}) with the source $\S\in W$  is found as
a fixed point, if it exists, of the map $\Phi_\S:E\to C([0,t_0];Y)$ where 
$\Phi_\S({\bf v})={\bf u}$ is  given by
$$
{\bf u}(t)=\int_0^t U^{\bf v}(t,\tilde t)f_\S(\tilde t,{\bf v})\,d\tilde t,\quad 0\leq t\leq t_0.
$$ 
Below, we denote  ${\bf u}^{\bf v}=\Phi_\S({\bf v})$.

As  $\Phi_{\S_0}(0)=0$ where
$\S_0=(0,0)$, we see using 
the above and the inequality $\|\,\cdotp\|_X\leq \|\,\cdotp\|_Y$ that the
 function ${\bf u}^{\bf v}$ satifies
  \ba
& &  \|{\bf u}^{\bf v}\|_{C([0,t_0];Y)}\leq C_3C_2t_0 r_1^2,\\
& &  \|{\bf u}^{\bf v}(t_2)-{\bf u}^{\bf v}(t_1)\|_{X}\leq C_3C_2r_1^2|t_2-t_1|,\quad t_1,t_2\in [0,t_0].
  \ea
  When $r_1>0$ is so small that $C_3C_2(1+t_0)<r_1^{-1}$ 
  and $C_3C_2r_1^2<C_1$
  we see that 
  \ba
  & & \|\Phi_\S({\bf v})\|_{C([0,t_0];Y)}<r_1,\\
  & & \|\Phi_\S({\bf v})\|_{C^{0,1}([0,t_0];X)}<C_1.
  \ea
  Hence $ \Phi_\S(E)\subset E$
 and we can consider $\Phi_\S$ as a map
 $\Phi_\S:E\to E$.
  
  As  $k>1+\frac 32$, it follows from Sobolev embedding theorem
  that $X=\H^{(k)}(\tilde N)\subset C^1(\tilde N)^2$. This yields that
  by \cite[Thm. 3]{Kato1975}, for the original reference, see Theorems III-IV
  in \cite{Kato1973},
  \beq\label{eq: stab.}
 & & \|(U^{\bf v}(t,s)- U^{\bf v^\prime}(t,s)){\bf h}\|_X \\
 & &\leq  \nonumber
 C_3\bigg(\sup_{t^\prime \in [0,t]}\|\A(t^\prime,{\bf v})-
 \A(t^\prime,{\bf v}^\prime)\|_{Y\to X}  \|U^{\bf v}(t^\prime,0) {\bf h}\|_Y\bigg)\\
 & &\leq  \nonumber
 C_3^2 \|{\bf v}-{\bf v}^\prime\|_{C([0,t_0];X)}  \|{\bf h}\|_Y.
  \eeq
  Thus,
  \ba
&&\|  U^{\bf v}(t,s)f_\S(s,{\bf v})- U^{\bf v^\prime}(t,s)f_{\S}(s,{\bf v}^\prime)\|_X\\
&\leq&\|  (U^{\bf v}(t,s)- U^{\bf v^\prime}(t,s))f_{\S}(s,{\bf v})\|_X+
\|  U^{\bf v^\prime}(t,s)(f_\S(s,{\bf v})- f_{\S}(s,{\bf v}^\prime))\|_X\\
&\leq&
(1+C_3)^2C_2r_1^2 \|{\bf v}-{\bf v}^\prime\|_{C([0,t_0];X)} .
  \ea
This implies that
\ba
& & \| \Phi_\S({\bf v})-\Phi_{\S}({\bf v}^\prime)\|_{C([0,t_0];X)}
\leq t_0(1+C_3)^2C_2r_1^2 \|{\bf v}-{\bf v}^\prime\|_{C([0,t_0];X))}.
\ea

Assume next that  $r_1>0$ is so small that we have also
\ba
t_0(1+C_3)^2C_2r_1^2<\frac 12.
\ea
 cf.\
  Thm.\ I in \cite{HKM} (or
  (9.15)  and (10.3)-(10.5) in \cite{Kato1975}). For $\S\in W$ this 
 implies that $$\Phi_\S:E\to E$$ is a contraction with a contraction constant 
 $C_L\leq \frac 12$,
  and thus 
 $\Phi_\S$ has a unique fixed point ${\bf u}$  in the space $E\subset C^{0,1}([0,t_0];X)$.
 
 Moreover, elementary considerations related to fixed points
 of the map $\Phi_\S$ show that ${\bf u}$  in  $C([0,t_0];X)$ depends in $E\subset C([0,t_0];X)$
 Lipschitz-continuously on
 $\S\in  W\subset C([0,t_0];H^{k+1}(\tilde N)\times H^{k}(\tilde N))$. Indeed, if 
 $$\|\S-\S^\prime\|_{C([0,t_0];H^{k+1}(\tilde N)\times H^{k}(\tilde N))}<\e,$$
 we see that
  \beq\label{fF bound r2 B}
 \|f_\S(t,{\bf v})-f_{\S^\prime}(t,{\bf v})\|_Y\leq C_2\e,\quad t\in [0,t_0],
 \eeq
 and when  (\ref{fF bound r2}) and (\ref{U bounds}) are satisfied 
 with $r=r_1$,
 we have
  \ba
& & \| \Phi_\S({\bf v})-\Phi_{\S^\prime}({\bf v}^\prime)\|_{C([0,t_0];Y)}
\leq C_3C_2t_0 r_1^2.
\ea
 Hence 
 \ba
 \|\Phi_\S({\bf v})-\Phi_{\S^\prime}({\bf v})\|_{C([0,t_0];Y)}\leq  t_0C_3C_2\e.
 \ea
 This and standard estimates for fixed points, yield that when $\e$ is small enough 
 the fixed point ${\bf u}^\prime$ of the map  $\Phi_{\S^\prime}:E\to E$ 
 corresponding to the source ${\S^\prime}$ and
  the fixed point ${\bf u}$ of the map $\Phi_{\S}:E\to E$ 
  corresponding to the source ${\S}$ satisfy
  \beq\label{stability}
  \|{\bf u}-{\bf u}^\prime\|_{C([0,t_0];X)}\leq  \frac 1{1-C_L}t_0C_3C_2\e.
  \eeq
{
Thus the solution ${\bf u}$ 
  depends in $C([0,t_0];X)$ Lipschitz continuously on
 $\S\in  C([0,t_0];H^{k+1}(\tilde N)\times H^{k}(\tilde N))$ (see also \cite[Sect.\ 16]{Kato1975}, 
 and \cite{Ringstrom})}.
In fact,  for analogous systems it is possible to show that $u$ is in $C([0,t_0];Y)$,
but one can not obtain  
 Lipschitz or H\"older  stability for $u$ in the $Y$-norm, see  \cite{Kato1975}, Remark 7.2. 
 
 Finally, we note that the fixed point ${\bf u}$ of $\Phi_\S$ can be found as a limit
 $ {\bf u}=\lim_{n\to \infty}{\bf u}_n$ in $C([0,t_0];X)$, where ${\bf u}_0=0$ and  ${\bf u}_n=\Phi_\S({\bf u}_{n-1})$.
Denote ${\bf u}_n=(g_n-\hat g,\phi_n-\hat \phi)$. We see that if
 $$\supp({\bf u}_{n-1})\subset J_{\hat g}(\supp(\S))$$ 
then also $$\supp(g_{n-1}-\hat g)\subset J_{\hat g}(\supp(\S)).$$
Hence  for all $x\in M_0\setminus J_{\hat g}(\supp(\S))$ we see that $$J_{g_{n-1}}^-(x)
\cap J_{\hat g}(\supp(\S))=\emptyset.$$ Then, using the definition of the map $\Phi_\S$ we see that
 $\supp({\bf u}_n)\subset J_{\hat g}(\supp(\S))$. Using induction we see that this holds
for all $n$ and hence we have that the solution ${\bf u}$ satisfies
 \beq\label{eq: support condition for u}
 \supp({\bf u})\subset J_{\hat g}(\supp(\S)).
 \eeq 
 \hiddenfootnote{REMOVE MATERIAL IN THIS FOOTNOTE OR MOVE IT ELSEWHERE:
  Let us next use the above considerations to study system 
   (\ref{eq: notation for hyperbolic system 1}) on $\R\times N$.
   For this end, we denote below the background metric of  $\tilde M=\R\times \tilde N$
   by $\tilde g$ and the background metric of  $\hattuM _0=\R\times  N$
   by $\hat g$.
 Consider a source 
\ba
\S\in {C^1([0,T];H^{s_0-1}(\tilde N))\cap
C^0([0,t_0];H^{s_0}(\tilde N))}
\ea  that is supported in a compact set $\K$ and
small enough in the norm of this space and let $\tilde u$ 
be the solution of the  
system (\ref{eq: notation for hyperbolic system 1}) on $\R\times \tilde N$.
Assume that  $u$ is some $C^2$-solution of  
(\ref{eq: notation for hyperbolic system 1}) on $\R\times N$
with the same source $\S$. Our aim is to show that it
coincides with $\tilde u$ in its support.
%
For this end, let $\tilde g^\prime$ be a metric tensor which
has the same eigenvectors as $\tilde g$ and the same
eigenvalues except that  the unique negative eigenvalue of $\tilde g$
 is multiplied by $(1-h_1)$, $h_1>0$. 
Assume that $h_1>0$ is so small
that $J^+_{\hattuM _0,\tilde g^\prime}(p^+)\subset  \hattuM _0\subset N_0\times (-\infty,t_0]$
and assume that $\e$ is so small that $g(\tilde u)<\tilde g^\prime$ for all 
sources $\S$ with $\|\S\|_{C^1([0,t_0];H^{s_0-1}(\tilde N))\cap
C^0([0,t_0];H^{s_0}(\tilde N))}\leq \e$. Then
 $J^+_{\tilde M,g(\tilde u)}(p^+)\subset  (-\infty,t_0]\times N_0$.
This  proves the estimate (\ref{eq: Lip estim}).

Next us we consider solution $u$ on $\hattuM _0$ corresponding to source $\S$.
Let  $T(h_1,\S)$
be the supremum of all $T^\prime\leq t_0$ for which
$J^+_{g(u)}(p^+)\cap \hattuM _0(T^\prime)\subset N_1\times (-\infty,T^\prime]$.
Assume that $ T(h_1,\S)<t_0$.
Then for $T^\prime= T(h_1,\S)$ we see that
$u$ can be continued by zero from $ (-\infty,T^\prime]\times N_1$
to a function on $(-\infty,T^\prime]\times \tilde N$ that 
satisfies the system (\ref{eq: notation for hyperbolic system 1}). Thus
$u$ coincides in $ (-\infty,T^\prime]\times N_1$
with $\tilde u$ and vanishes outside this set.
Since  $J^+_{\tilde M,g(\tilde u)}(p^+)\subset  (-\infty,t_0]\times N_0$, this implies that
also $u$ satisfies 
$J^+_{g(u)}(p^+)\cap \hattuM _0(T^\prime)\subset (-\infty,T^\prime]\times N_0$.
As  $u$ is $C^2$-smooth and $N_0\subset \subset N_1$
we see that there is $T^\prime_1> T(h_1,\S)$ so 
that $J^+_{g(u)}(p^+)\cap \hattuM _0(T^\prime_1)\subset  (-\infty,T^\prime_1]
\times N_1$ that is in contradiction with definition of $ T(h_1,\S)$. Thus
we have to have  $T(h_1,\S)=t_0$.


This shows that when the norm of $\S$ is smaller than the above chosen $\e$
then $J^+_{g(u)}(p^+)\cap \hattuM _0\subset N_0\times (-\infty,t_0]$.
Then $u$ vanishes outside $N_1\times (-\infty,t_0]$ by
the support condition in (\ref{eq: notation for hyperbolic system 1}).
Thus $u$ and $\tilde u$ are both supported on  $N_1\times (-\infty,t_0]$ 
and coincide there. Since $\tilde u$ is unique, 
this shows that for a sufficiently small sources $ \S$
 supported on  $\K$   the solution $u$ of
(\ref{eq: notation for hyperbolic system 1})  in $M_0$
exists and is unique.}

\bigskip

\section*{Appendix C: 
An example  satisfying microlocal linearization stability}


\subsection*{C.1.\ Formulation of the direct problem} 
Let us define some physical fields  and introduce a
 model as a system of partial differential
 equations. Later we will motivate this system by discussion
 of the corresponding Lagrangians, but we postpone this
 discussion to  Appendix C.3 as it is not completely rigorous.

We assume that there are $C^\infty$-background fields $\hat g$,  $\hat \phi$,
 on $M$. 
%

We consider 
a Lorentzian metric  $g$ on $\hattuM _0$ and $\phi=(\phi_\ell)_{\ell=1}^L$ where $\phi_\ell$
are scalar fields on $\hattuM _0=(-\infty,t_0)\times N$.

{Let $P=P_{jk}(x)dx^jdx^k$ be a symmetric tensor on $\hattuM _0$, corresponding below
to a direct perturbation to the stress energy tensor, and 
$Q=(Q_{\ell}(x))_{{\ell}=1}^K$ where  $Q_{\ell}(x)$ are real-valued functions on $\hattuM _0$, where $K\geq L+1$.
We denote by ${\mathcal V}( \phi_\ell;S_\ell)$
the potential functions of the fields $\phi_\ell$, 
\beq\label{eq: Slava-aaa}
{\mathcal V}( \phi_\ell;S_\ell)=\frac 12m^2\bigg(\phi_\ell +\frac 1{m^2}S_\ell\bigg)^2.
\eeq 
These potentials depend on the source variables $S_\ell$.
The way how  $S_\ell$,
called below the adaptive source functions, depend on other 
fields is explained later. We assume that there are smooth background
fields $\hat P$ and $\hat Q$. For a while we consider
the case when $\hat P=0$ and $\hat Q=0$, and discuss later
the generalization to non-vanishing background fields.

Using the $\phi$ and $P$ fields, we define the stress-energy tensor
 \beq\label{eq: T tensor1}
\hspace{-5mm} T_{jk}=\sum_{\ell=1}^L(\p_j\phi_\ell \,\p_k\phi_\ell 
-\frac 12 g_{jk}g^{pq}\p_p\phi_\ell \,\p_q\phi_\ell-\mathcal V( \phi_\ell;S_\ell)g_{jk})
+P_{jk}.\hspace{-14mm}
\eeq


We assume that $P$ and $Q$ are supported on
$\K=J^+_{\tilde g}(\hat p^-)\cap \hattuM _0$. 
Let  us represent the stress energy tensor (\ref{eq: T tensor1})
 in the form
  \ba
T_{jk}
&=&P_{jk}+Zg_{jk}+{\bf T}_{jk}(g,\phi),\quad Z=-(\sum_{\ell=1}^L S_\ell\phi_\ell+
\frac 1{2m^2}S_\ell^2),
\\
{\bf T}_{jk}(g,\phi)&=&\sum_{\ell=1}^L(\p_j\phi_\ell \,\p_k\phi_\ell 
-\frac 12 g_{jk}g^{pq}\p_p\phi_\ell \,\p_q\phi_\ell-
\frac 12 m^2\phi_\ell^2g_{jk}),
 \ea
where we call $Z$ the stress energy density caused by the sources $S_\ell$.

\generalizations{We will also add in our considerations linear, first order differential operators
\beq\label{B-interaction}
B_\ell(\phi)=\sum_{\kappa=1}^L a^\kappa_\ell\phi_\kappa(x)
+ \sum_{\kappa,\alpha,\beta =1}^Lb^{\kappa,\a,\beta}_\ell \phi_\kappa (x)
\phi_\alpha (x)\phi_\beta (x)\hspace{-1cm}
\eeq modeling interaction of matter fields,
 where $a^{\kappa}_\ell,b^{\kappa,\a,\beta}_\ell\in \R$ are constants. 
\HOX{ We have  added here some first order terms
to the wave equation that connect $\phi$-variables together. WE NEED TO ADD 
TERMS IN STRESS ENERGY TENSOR AND IN THE APPENDIX B AND THE MAIN THEOREM}}

{Now we are ready to formulate the direct problem for 
the adaptive Einstein-scalar field equations. Let $g$ and $\phi$ satisfy
\beq\label{eq: adaptive model}& &\Ein_{\hat g}(g) =P_{jk}+Zg_{jk}+{\bf T}_{jk}(g,\phi),\quad Z=-(\sum_{\ell=1}^L S_\ell\phi_\ell+
\frac 1{2m^2}S_\ell^2),
\\ \nonumber
& &\square_g\phi_\ell -\mathcal V^{\prime}( \phi_\ell;S_\ell) =0
\generalizations{+B_\ell(\phi)}
\quad
\hbox{in }\hattuM _0,\quad \ell=1,2,3,\dots,L,
\\ \nonumber
& &S_\ell={\mathcal S}_\ell(g,\phi, \nabla
 \phi,Q,\nabla Q,P,\nabla^g P),\quad
\hbox{in }\hattuM _0,
 \\ \nonumber
& & g=\hat g,\quad \phi_\ell=\hat \phi_\ell,\quad
\hbox{in }\hattuM _0\setminus \K.
\eeq
Above, $\mathcal V^{\prime}(\phi ;s)=\p_\phi \mathcal V(\phi;s)$ so that
$
\mathcal V^{\prime}( \phi_\ell;S_\ell)=m^2\phi_\ell+S_\ell.
$
}
We assume that  the background fields $\hat g$,  $\hat \phi$,
 satisfy these equations with  $\hat Q=0$ and  $\hat P=0$.}

We consider  here $P=(P_{jk})_{j,k=1}^4$
and $Q=(Q_\ell)_{\ell=1}^K$   as   fields that we can control and call those
the controlled source fields. 
 Local existence of the solution for small sources $P$ and $Q$ is considered
in Appendix B. 

%
%
To obtain a physically meaningful model,
we need to consider how  the adaptive source functions $\mathcal S_\ell$ should be chosen 
so that the physical conservation law in relativity
\beq\label{conservation law CC}
\nabla_k(g^{kp} T_{pq})=0
\eeq
 is satisfied. Here $\nabla=\nabla^g$ is the connection corresponding to the metric $g$.

  We
note that the conservation law 
is a necessary condition for the equation (\ref{eq: adaptive model})
to have solutions for which
$\Ein_{\hat g}(g)=\Ein(g)$, i.e., that the solutions of
(\ref{eq: adaptive model}) are solutions of  the Einstein field equations.

{The functions ${\mathcal S}_\ell(g,\phi, \nabla
 \phi,Q,\nabla Q,P,\nabla^g P)$ model the devices that we use to perform active
 measurements. Thus, even though the Assumption S below may appear quite technical,
 it can be viewed as the instructions on  how to build a device that can be used 
 to measure the structure of the spacetime far away. Outside the support
 of the measurement device (that contain the union of the supports of $Q$ and $P$) we have just assumed that the standard coupled Einstein-scalar field
 equations hold, c.f. (\ref{S-vanish condition}). 
We can consider them in the form
\ba
& &{\mathcal S}_\ell(g,\phi, \nabla
 \phi,Q,\nabla Q,P,\nabla^g P)=Q_\ell+{\mathcal S}^{2nd}_\ell(g,\phi, \nabla
 \phi,Q,\nabla Q,P,\nabla^g P)
\ea
for $\ell=1,2,\dots,L$
where $Q_\ell$ are the primary sources and ${\mathcal S}^{2nd}_\ell$,
that depend also on $Q_\ell$ with $\ell=L+1,L+2,\dots,K$,
corresponds to the response of the measurement device that forces
the conservation law to be valid.

The solution $(g,\phi)$ of (\ref{eq: adaptive model}) is a solution of the equations
(\ref{eq: adaptive model with no source}) when we denote
 \ba
& &\F^1_{jk}=P_{jk}+Zg_{jk},\\
& &\F^2_\ell=\mathcal V^{\prime}( \phi_\ell;S_\ell)-\mathcal V^{\prime}( \phi_\ell;0)=S_\ell.
\ea
Our next goal is to construct suitable adaptive source
functions ${\mathcal S}_\ell$ and consider what kind of sources $\F^1$ and $\F^2$
of the above form can be obtained by varying $P$ and $Q$.

We will consider  adaptive source
functions ${\mathcal S}_\ell$  satisfying the following conditions:

 {
\medskip 

\noindent{\bf Assumption S}:

\noindent
The adaptive source functions ${\mathcal S}_\ell(g,\phi, \nabla
 \phi,Q,\nabla Q,P,\nabla^g P)$ have the following properties: 
\smallskip

{(i) Denoting  $c=\nabla
 \phi$, $C=\nabla^g P$, and $H=\nabla Q$
 we assume 
  that ${\mathcal S}_\ell(g,\phi, c,Q,H,P,C)$ are smooth {non-linear} functions, 
  of the pointwise values $g_{jk}(x),\phi(x), \nabla
 \phi(x),$ $Q(x),\nabla Q(x),P(x),$ and $\nabla^g P(x)$,
  defined
 near $(g,\phi, c,Q,H,P,C)=(\hat g,\hat \phi,\nabla \hat \phi ,0,0,0,0)$, that
satisfy
\beq\label{S-vanish condition}
 {\mathcal S}_\ell(g,\phi,c,0,0,0,0)=0.
\eeq
 
\smallskip


 \smallskip

(ii)  
We assume that  ${\mathcal S}_\ell$ is independent of $P(x)$ and the dependency of ${\mathcal S}$ on  $\nabla^g P$ 
and $\nabla Q$ is only due to the dependency
in the term $g^{pk}\nabla^g_p( P_{jk}+Zg_{jk})=
g^{pk}\nabla^g_p P_{jk}+\nabla^g_jQ_K$, associated to
the divergence of the perturbation of $T$,
that is,
there exist functions $\tilde {\mathcal S}_\ell$ so that
\ba
{\mathcal S}_\ell(g,\phi, c,Q,H,P,C)=\tilde {\mathcal S}_\ell(g,\phi,c,Q,
R),\quad R=(g^{pk}\nabla^g_p ( P_{jk}+Q_Kg_{jk}))_{j=1}^4.\ea
%
}
 

 \bigskip


Below, denote $\hat R=\hat g^{pk}\hat\nabla_p \hat P_{jk}+\hat\nabla_j\hat Q_K.$
Note that we still are considering the  case when $\hat Q=0$ and $\hat P=0$
so that $\hat R=0$, too. This implies
 that for the background fields that adaptive source functions $\mathcal S_\ell$ vanish.

To simplify notations, we also denote below
$\tilde {\mathcal S}_\ell$ just by ${\mathcal S}_\ell$ and indicate the function
which we use by the used variables in these functions.

Below we will denote $Q=(Q^{\prime},Q_K)$, $Q^{\prime}=(Q_\ell)_{\ell=1}^{K-1}$.
There are examples when the background fields $(\hat g,\hat \phi)$ and 
the adaptive source functions ${\mathcal S}_\ell$ exists and satisfy the Assumption S. 
  This is shown later in the case the following condition is valid for the 
     background fields:
 \medskip

 {\bf Condition A}:
Assume that at any $x\in \overline U_{\hat g}$ there is
a permutation $\sigma:\{1,2,\dots,L\}\to \{1,2,\dots,L\}$, denoted $\sigma_x$, such that the 
$5\times 5$ matrix $[ B_{jk}^\sigma(\hat \phi(x),\nabla \hat \phi(x))]_{j,k\leq 5}$
is invertible, where
  \ba
[ B_{jk}^\sigma(\phi(x),\nabla \phi(x))]_{k,j\leq 5}=\left[\begin{array}{c}
(\,\p_j  \phi_{\sigma(\ell)}(x))_{\ell\leq 5,\ j\leq 4}\\
(\phi_{\sigma(\ell)}(x))_{\ell\leq 5}\end{array}\right].
 \ea
 Below, for a permutation $\sigma:\{1,2,\dots,L\}\to \{1,2,\dots,L\}$ we denote by $U_{\hat g,\sigma}$
 the open set of points $x\in U_{\hat g}$ for which $[ B_{jk}^\sigma(\hat \phi(x),\nabla \hat \phi(x))]_{j,k\leq 5}$
is invertible. So, we assume that the sets $U_{\hat g,\sigma}$,
$\sigma\in \Sigma(L)$ is an open covering of $U_{\hat g}$.
 \medskip

 Our next aim is to prove the following:
  \medskip

 \begin{theorem} \label{thm: Good S functions}
Let $L\geq 5$ and assume that $\hat Q=0$ and $\hat P=0$ so that
 $\hat R=0$.
Moreover, assume that 
 Condition A is valid. Then for all permutations $\sigma:\{1,2,\dots,L\}\to \{1,2,\dots,L\}$ there exists  functions $  {\mathcal S}_{\ell,\sigma}$ satisfying Assumption S
  such that
%
\smallskip

(i) For all
$x\in   U_{\hat g,\sigma}$ 
 the differential of $$  {\mathcal S}_\sigma(\hat g,\hat \phi,  \nabla
 \hat \phi,Q,R)=(  {\mathcal S}_{\ell,\sigma}(\hat g,\hat \phi, \nabla
 \hat \phi,Q,R))_{\ell=1}^L$$ with respect to $Q$ and $R$, that is, the map
  \beq\label{surjective 1}
D_{Q,R}  {\mathcal S}_\sigma(\hat g(x),\hat \phi(x),  \nabla
 \hat \phi(x),Q,R)|_{Q=\hat Q(x), R=\hat R(x)}:\R^{K+4}\to \R^L\hspace{-1cm}
 \eeq
 is surjective.
\smallskip
 
(ii)  
The adaptive source functions ${\mathcal S}_\sigma$ are such that
 for
 $(Q_\ell)_{\ell=1}^K$
and  $(P_{jk})$ that are sufficiently close to $\hat Q=0$ and $\hat P=0$ in the $C^{{3}}_b(M_0)$-topology 
and supported in $  U_{\hat g,\sigma}$
the equations (\ref{eq: adaptive model}) with source functions  ${\mathcal S}_\sigma$ have
a unique solution $(g,\phi)$ and 
the conservation
 law (\ref{conservation law CC}) is valid. 
 \smallskip

 (iii) Under the same assumptions as in (ii), when  $(g,\phi)$ is a solution  of (\ref{eq: adaptive model}) 
 with the controlled source functions $P$ and $Q$, we have
  $Q_K=Z$. This means that  the physical field $Z$  can be directly controlled. 
 
 \end{theorem}

}}

\noindent
{\bf Proof.} As  one can enumerate the $\ell$-indexes 
of the fields $\phi_\ell$ as one wishes,
it is enough to prove the claim with one $\sigma$. We
consider below the case when $\sigma=Id$.

 Consider 
  a symmetric (0,2)-tensor $P$ and a scalar functions $Q_\ell$
  that are $C^3$-smooth and compactly supported in  $U_{\hat g,\sigma}$.
Let $[P_{jk}(x)]_{j,k=1}^4$ be  the coefficients of 
  $P$ in local coordinates and $Q(x)=(Q_\ell(x))_{\ell=1}^L$. 
 

To obtain adaptive required adaptive source functions, let us start implications of the conservation law (\ref {conservation law CC}).
To this end, consider $C^2$-smooth functions 
 $S_\ell(x)$ on $U_{\hat g,\sigma}$.

Note that since $[\nabla_p,\nabla_n]=[\p_p,\p_n]=0$ (see \cite[Sect. III.6.4.1]{ChBook}), 
 \ba
& &\nabla_ p (g^{pj} {\bf T}_{jk}(g,\phi)=
\\ & &\sum_{\ell=1}^L \nabla_ p
g^{pj}  (\p_j\phi_\ell \,\p_k\phi_\ell 
-\frac 12 g_{jk}g^{nm}\p_n\phi_\ell \,\p_m\phi_\ell-
\frac 12 m^2\phi_\ell^2g_{jk})
\\
&=&
\sum_{\ell=1}^L (g^{pj}  \nabla_ p \p_j\phi_\ell) \,\p_k\phi_\ell
+\sum_{\ell=1}^L  (g^{pj}  \p_j\phi_\ell \, \nabla_ p \p_k\phi_\ell)\\
& &-
 \frac 12
 \sum_{\ell=1}^L \delta^p_k 
 (g^{nm}(\nabla_ p \p_n\phi_\ell) \,\p_m\phi_\ell+
 g^{nm}\p_n\phi_\ell \,(\nabla_ p  \p_m\phi_\ell))
  - \sum_{\ell=1}^L m^2  \delta^p_k   \phi_\ell \p_p \phi_\ell
\\
&=&
\sum_{\ell=1}^L (g^{pj}  \nabla_ p \p_j\phi_\ell) \,\p_k\phi_\ell
+\sum_{\ell=1}^L  (g^{pj}  \p_j\phi_\ell \, (\nabla_ p \p_k\phi_\ell))\\
& &-
 \frac 12
 \sum_{\ell=1}^L 
 (g^{nm}\p_m\phi_\ell\,(\nabla_ n \p_k\phi_\ell)+
 g^{nm}\p_n\phi_\ell \,(\nabla_ m  \p_k\phi_\ell))
 - \sum_{\ell=1}^L m^2   \phi_\ell \p_k \phi_\ell
\\
&=&
\sum_{\ell=1}^L (g^{pj}  \nabla_ p \p_j\phi_\ell -m^2   \phi_\ell) \,\p_k\phi_\ell. 
\ea
Thus conservation law (\ref {conservation law CC})
gives for all $j=1,2,3,4$ equations
 \ba
0&=&\nabla_p^g (g^{pk}T_{jk})\\
&=&\nabla_p^g (g^{pk}({\bf T}_{jk}(g,\phi)+P_{jk}+Zg_{jk}))\\
&=& \sum_{\ell=1}^L\bigg((g^{pk} \nabla^g_p\p_k\phi_\ell )\,\p_j\phi_\ell 
-(m^2_\ell \phi_\ell \p_p \phi_\ell) \delta_j^p \\
& &\quad \quad- \nabla^g_p (g^{pk}g_{jk}
(S_\ell\phi_\ell{\color{black}+\frac 1{2m^2}S_\ell^2})+g^{pk}P_{jk})\bigg)
\\
&=& \sum_{\ell=1}^L\bigg(
(g^{pk} \nabla^g_p\p_k\phi_\ell -m^2_\ell \phi_\ell ) \,\p_j \phi_\ell 
\\ & &\quad \quad-\nabla^g_p (g^{pk}g_{jk}(S_\ell\phi_\ell
{\color{black}+\frac 1{2m^2}S_\ell^2})
+g^{pk}P_{jk})\bigg)
\\
&=&\sum_{\ell=1}^L S_\ell \,\p_j\phi_\ell- \nabla^g_p (g^{pk}g_{jk}
(S_\ell \phi_\ell{\color{black}+\frac 1{2m^2}S_\ell^2}))+ g^{pk} \nabla^g_p P_{jk}\\
&=&\left (\sum_{\ell=1}^L S_\ell \,\p_j\phi_\ell\right)-\p_j\left(
\sum_{\ell=1}^L S_\ell \phi_\ell {\color{black}+\frac 1{2m^2}S_\ell^2}\right)+g^{pk}\nabla_p^g P_{jk}.
\ea
Summarizing, the conservation law yields 
\beq
\label{mod. conservation law a}
\left (\sum_{\ell=1}^L S_\ell \,\p_j\phi_\ell\right)-\p_j\left(
\sum_{\ell=1}^L S_\ell \phi_\ell +\frac 1{2m^2}S_\ell^2\right)+g^{pk}\nabla_p^g P_{jk}=0,\hspace{-1cm}
\eeq
for $j=1,2,3,4$.

Recall that the field $Z$ has the definition
\beq\label{meson number}
\sum_{\ell=1}^L S_\ell \phi_\ell+\frac 1{2m^2}S_\ell^2=-Z.
\eeq

Then, the  conservation law  (\ref {conservation law CC}) holds if we have
\beq\label{meson number2}& &
\sum_{\ell=1}^L S_\ell \,\p_j\phi_\ell
=- g^{pk}\nabla^g_p  {\color{black}V}_{jk},\quad
 {\color{black}V}_{jk}=(P_{jk}+g_{jk}Z),
\eeq
for $j=1,2,3,4.$

Equations (\ref{meson number}) and (\ref{meson number2}) give together  five point-wise equations
 for the  functions $S_1,\dots,S_L$.

%
%
%
Recall that we consider here the case when $\sigma=Id$.
By Condition A, at any $x\in U_{\hat g,\sigma}$ that the 
$5\times 5$ matrix $( B_{jk}^\sigma(\hat \phi(x),\nabla \hat \phi(x)))_{j,k\leq 5}$
is invertible, where
  \ba
( B_{jk}^\sigma(\phi(x),\nabla \phi(x)))_{j,k\leq 5}=\left(\begin{array}{c}
( \,\p_j  \phi_{\cell}(x))_{j\leq 4,\ \ell\leq 5}\\
(   \phi_{\cell}(x)) _{\ell\leq 5}\end{array}\right).
 \ea

 We  consider
a $\R^K$ valued function  $Q(x)=(Q^\prime(x),Q_K(x))$,
where  $$Q^\prime=(Q_\ell)_{\ell=1}^{K-1}.$$

Also, below $ {\color{black}R}_j=g^{pk}\nabla^g_p  {\color{black}V} _{jk},$ $ {\color{black}V}_{jk}=P_{jk}+g_{jk}Z$ and  we require that identity
\beq
& &Q_K=Z
\label{S R Z equations 3bbb}
\eeq
holds.

Motivated by equations (\ref{meson number}), (\ref{meson number2}),
and (\ref{S R Z equations 3bbb}),
our next aim is to consider a point $x\in U_{\hat g,\sigma}$,
and construct  functions
${\mathcal S_{\sigma,\ell}}(\phi,\nabla \phi, Q^\prime,Q_K,R,g)$, $\ell=1,2,\dots,L$ that satisfy
\beq\label{S R Z equations 1bbb}
& &\sum_{\ell=1}^5 {\mathcal S_{\sigma,\ell}}(\phi,\nabla \phi, Q^\prime,Q_K,R,g)
 \,\p_j\phi_{\cell}=- R_{j} -
\sum_{\ell=6}^{L} Q_{\sigma,\ell} \,\p_j\phi_{\cell},\hspace{-1cm} \\
& &\sum_{\ell=1}^5 {\mathcal S_{\sigma,\ell}}(\phi,\nabla \phi, Q^\prime,Q_K,R,g)\, \phi_{\cell}=-\bigg(Q_K+\sum_{\ell=6}^L Q_{\sigma,\ell} \phi_{\cell}
+\label{S R Z equations 2bbb}\\ \nonumber\\ & &\quad
+\sum_{\ell=1}^L\frac 1{2m^2} {\mathcal S_{\sigma,\ell}}(\phi,\nabla \phi, Q^\prime,Q_K,R,g)^2\bigg)
.\hspace{-1cm}
\eeq
%
%
Let
\ba
& &(Y_\sigma(\phi,\nabla \phi))(x)=\psi(x)(B^\sigma(\phi,\nabla \phi))^{-1},\quad \hbox{for
$x\in U_{\hat g,\sigma}$,}
\\
& &(Y_\sigma(\phi,\nabla \phi))(x)=0,\quad \hbox{for
$x\not \in U_{\hat g,\sigma}$,}
\ea 
where $\psi\in C^\infty_0(U_{\hat g,\sigma})$ has value 1 in $\supp(Q)\cup\supp(P)$.

Then we define 
${\mathcal S_{\sigma,\ell}}={\mathcal S_{\sigma,\ell}}(g,\phi,\nabla \phi, Q^\prime,Q_K,R)$, $\ell=1,2,\dots,L,$ to be the solution of the system 
\beq\label{S sigma formulas}
& &(S_{\sigma,\ell})_{\ell\leq 5}= Y_\sigma(\phi,\nabla \phi) \left(\begin{array}{c}
(-  {\color{black}R}_j -
\sum_{\ell=6}^L Q_{\sigma,\ell} \,\p_j\phi_{\cell})_{j\leq 4}
\\
-Q_K-\sum_{\ell=6}^L Q_{\sigma,\ell}\phi_{\cell}-\sum_{\ell=1}^L\frac 1{2m^2}S_{\sigma,\ell}^2 \end{array}\right)
\hspace{-1cm}\\ \nonumber
%
\hspace{-1cm}\\
& & \nonumber
(S_{\sigma,\ell})_{\ell\geq 6}= (Q_{\cell})_{\ell\geq 6}.
\eeq  
{When $Q$ and $R$ are sufficiently small,
 this equation can be solved point-wisely, 
 at each point $x\in U_{\hat g,\sigma}$, using iteration by the Banach fixed point theorem.}


Let
 \ba
( K_{jk}^\sigma(\phi(x),\nabla \phi(x)))_{j,k\leq 5}=\left(\begin{array}{c}
( \,\p_j  \phi_{\cell}(x))_{j\leq 4,\ 6\leq \ell\leq L}\\
(   \phi_{\cell}(x)) _{6\leq \ell\leq L}\end{array}\right).
 \ea
Then
we see that
 the differential of 
$ {\mathcal S}_\sigma=( {\mathcal S}_{\sigma,\ell})_{\ell=1}^L$
with respect to $(Q^{\prime},Q_{K}, {R})$ at $(Q,R)=(0,0)$, that is,
 \beq\label{eq: deriv1}\\
 \nonumber
& &D_{Q^{\prime}, Q_{K}, {R}} {\mathcal S}_\sigma(\hat g,\hat \phi, \nabla
 \hat \phi,Q^{\prime},Q_{K}, {R})|_{Q=0,R=0}:\R^{K+4}\to \R^L,\\
 & &\nonumber (Q^{\prime},Q_{K}, R)
\mapsto 
 -\left(\begin{array}{cc}
Y_\sigma(\hat \phi,\nabla\hat\phi)   &Y_\sigma(\hat\phi,\nabla\hat\phi)K(\hat\phi,\nabla\hat\phi) 
\\
0 &I_\sigma  \end{array}\right)
 \left(\begin{array}{c}
\left(\begin{array}{c}
  R\\
 Q_K \end{array}\right) \\
Q^\prime\\
\end{array}\right),
 \eeq
 is surjective, where $I_\sigma=[\delta_{k,j+5}]_{k\leq K-1,\ j\leq L-5,}\in \R^{(K-1)\times (L-5)}$. Hence (i) is valid.




%

%
%
%

By their construction, the functions ${\mathcal S}_\sigma=( {\mathcal S}_{\sigma,\ell})_{\ell=1}^L$ satisfy the equations
(\ref{meson number}) and (\ref{meson number2})
for all $x\in U_{\hat g,\sigma}$ and also equation  (\ref{S R Z equations 3bbb}) holds.
\hiddenfootnote{ 
We make the following note related 
the case considered in the main part of the paper when $Q,P,R\in \I^m(Y)$, where $Y$ is 2-dimensional sufrace,
and we  need to use the principal symbol ${\bf r}$ of $R$
as an independent variable:
Let us consider also a point $x_0\in \hattuM _0$ and $\eta$
be a light-like covector choose
coordinates so that $g=\diag(-1,1,1,1)$ and $\eta=(1,1,0,0).$
When $c=(c^k)_{k=1}^4\in \R^4$ and $P^{jk}=C^{jk}(x\,\cdotp \eta)^a_+$, where $C^{jk}$ is such a symmetric
matrix that
$C^{11}=c_1$, $
C^{12}=C^{21}=\frac  12 c_2$,
 $C^{13}=C^{31}=c_3$, and  $C^{14}=C^{41}=c_4$ and other $C^{jk}$ are zeros. 
Then we have
\ba
\nabla^g_j (C^{jk}(x\,\cdotp \eta)_+^a)&=&\eta_jP^{jk}=(C^{1k}+ C^{2k})(x\,\cdotp \eta)_+^{a-1}=c_k(x\,\cdotp \eta)_+^{a-1}.
\ea
As we can always choose coordinates
that $\hat g$ and $\eta$ have at a given point the above forms,
we can obtain arbitrary vector ${\bf r}$, as
principal symbol of $R$, by considering
$Q^{\prime},Q_{L+1},P\in \I^m(Y)$, where $Y$ is 2-submanifold, with principal symbols of
$\tilde {\bf p}$ and $\tilde {\bf z}$ satisfing equations corresponding to
equations
$\hat g^{jk}\hat \nabla_j ({\bf p}_{kn}+ {\bf z}\hat g_{kn})\in \I^m(Y)$, 
and the sub-principal symbols of ${\bf p}_{kn}$ and ${\bf z}$
vary arbitrarily.} 
Hence (iii) is valid.

 Above, the equation (\ref{meson number2}) is valid by construction of the
 functions $( {\mathcal S}_\ell)_{\ell=1}^L$. Thus the conservation law is valid.
 This proves (ii). \hfill \Box \medskip



%

%

Note that as  the adaptive source functions
 $\mathcal S_\ell$  were constructed in Theorem \ref{thm: Good S functions}
 using the inverse function theorem, the results of Theorem \ref{thm: Good S functions} are valid also if
 $\hat Q$  and $\hat P$ are sufficiently small non-vanishing fields
 and $\hat g$ and $\hat \phi$ satisfy the Einstein scalar field equations  (\ref{eq: adaptive model})
 with these background fields. Next we return to the case when $\hat P=0$  and $\hat Q=0$.

 \subsection*{C.2.\  Microlocal linearization stabililty}
 Below we consider the case when $\hat P=0$  and $\hat Q=0$ and use the adaptive source functions
 $\mathcal S_\ell$  constructed in Theorem \ref{thm: Good S functions} and
 its proof.
 
 Assume that $Y\subset \hattuM _0$
is a 2-dimensional space-like submanifold and consider local coordinates defined 
in  $V\subset \hattuM _0$. Moreover, assume that 
in these local coordinates $Y\cap V\subset \{x\in \R^4;\ x^jb_j=0,\  x^jb^\prime_j=0\}$,
where $b^\prime_j\in \R$  and let ${\bf p}\in \I^{n}(Y)$, $n\leq n_0=-17$,  be defined by 
\beq\label{eq: b b-prime BBB}
{\bf p}_{jk}(x^1,x^2,x^3,x^4)=
\re \int_{\R^2}e^{i(\theta_1b_m+\theta_2b^\prime_m)x^m}\symbolP_{jk}(x,\theta_1,\theta_2)\,
d\theta_1d\theta_2.\hspace{-2cm}
\eeq
Here, we assume that 
 $\symbolP_{jk}(x,\theta)$, $\theta=(\theta_1,\theta_2)$ 
 are  classical symbols and we denote their principal symbols by
 $\sigma_p({\bf p}_{jk})(x,\theta)$. When 
$x\in Y$ and $\xi=(\theta_1b_m+\theta_2b^\prime_m)dx^m$ so that $(x,\xi)\in N^*Y$,
we denote the value of the principal symbol $\sigma_p({\bf p})$ at $(x,\theta_1,\theta_2)$ by 
 $\tilde \symbolP ^{(a)}_{jk}(x,\xi)$, that is,  $ \tilde \symbolP ^{(a)}_{jk} (x,\xi)=\sigma_p({\bf p}_{jk})(x,\theta_1,\theta_2)$,
 and say that it is the principal symbol of ${\bf p}_{jk}$ at $(x,\xi)$, associated to the phase function $\phi(x,\theta_1,\theta_2)=(\theta_1b_m+\theta_2b^\prime_m)x^m$.
The above defined principal symbols can be defined invariantly, see \cite{GuU1}.

We assume that also ${\bf q}^\prime,{\bf z}\in \I^{n}(Y)$ have 
representations (\ref{eq: b b-prime BBB}) with classical symbols. Below
we consider symbols in local coordinates.
{Let
 us denote the principal symbols of ${\bf p},{\bf q^{\prime}},{\bf z}\in \I^{n}(Y)$ 
 by  $\tilde \symbolP ^{(a)}(x,\xi),$ $\tilde \symbolQ_{1}^{(a)}(x,\xi)$,
 $\tilde \symbolQ_2^{(a)}(x,\xi)$, respectively
  and let $\tilde \symbolP ^{(b)}(x,\xi)$ and $\tilde \symbolQ_2^{(b)}(x,\xi)$ denote the sub-principal
symbols of ${\bf p}$ and ${\bf z}$, correspondingly, at $(x,\xi)\in N^*Y$.


We will below consider what happens when  $({\bf p}_{jk}+{\bf z}\hat g_{jk})\in \I^{n}(Y)$
 satisfies 
\beq\label{eq: leading order lineariz. conse. law}
\hat g^{lk} \nabla_l^{\hat g} ({\bf p}_{jk}+{\bf z}\hat g_{jk})\in \I^{n}(Y),\quad j=1,2,3,4.
\eeq
Note that a priori this function is only in $\I^{n+1}(Y)$, so the
assumption (\ref{eq: leading order lineariz. conse. law}) means that
$\hat g^{lk} \nabla_l^{\hat g} ({\bf p}_{jk}+{\bf z}\hat g_{jk})$
is one degree smoother than it should be a priori.

When (\ref{eq: leading order lineariz. conse. law}) is valid, we say
that  {\it  the leading order of singularity of the wave satisfies
the linearized conservation
law}. This  corresponds to the assumption that the principal
symbol of the sum of divergence of the first two terms appearing in the stress energy tensor
on the right hand side of (\ref{eq: adaptive model}) vanishes.
}


{By \cite{GuU1},  
the identity  (\ref{eq: leading order lineariz. conse. law}) is equivalent to the vanishing
of the principal symbol on $N^*Y$,} that is,
\beq\label{eq: lineariz. conse. law symbols mLS}
& &\hat g^{lk}\xi_l(\tilde \symbolP ^{(a)}_{kj} (x,\xi)
+\hat g_{kj}(x)\tilde \symbolQ_2^{(a)}(x,\xi))=0,\ \hbox{for 
$j\leq 4$ and $\xi\in N^*_x Y$}.
\eeq
We say that this is the {\it linearized conservation law for principal symbol of $R$}.

  Let us consider source fields
that have the form $Q^{\prime}_\e=((Q_\e)_\ell)_{\ell=1}^{K-1}=\e{\bf q^{\prime}}$, $(Q_\e)_K=\e{\bf z}$ and 
 $P_\e=\e {\bf p}$.
We denote 
${\bf q}=({\bf q^{\prime}},{\bf z})$. We assume that 
${\bf q^{\prime}}$, ${\bf z}$, and ${\bf p}$
%
%
are supported in $\hat V\subset \subset U_{\hat g}$. 

Let $u_\e=(g_\e,\phi_\e)$ be the  solution of
 (\ref{eq: adaptive model}) with source  $P_\e$  and $Q_\e$.
 Then $u_\e$ depends $C^4$-smoothly on $\e$ and
$(g_\e,\phi_\e)|_{\e=0}=(\hat g,\hat \phi)$.
Denote $\p_\e  (g_\e,\phi_\e)|_{\e=0}=(\dot g,\dot \phi)$.
When $\e_0$ is small enough,
$P_\e$ and $Q_\e$ are supported in $U_{g_\e}$ for all $\e\in (0,\e_0)$.

Let $$R_\e=
g_\e ^{pk}\nabla_p^{g_\e} ((P_\e)_{jk}+ g^\e_{jk}(Q_\e)_K)$$ and 
$$
({S_\e})_\ell=
 {\mathcal S}_\ell(g_\e,\phi_\e, \nabla
 \phi_\e ,Q_\e ^{\prime},(Q_\e)_{K},R_\e),
$$
where
 $\mathcal S_\ell$  are   the adaptive source functions constructed in Theorem \ref{thm: Good S functions} and
 its proof.

Then $S_\e|_{\e=0}=0$ and $\p_\e  S_\e|_{\e=0}={\dot S}$
satisfy
\beq\label{eq: formula for sources 1}\\
\nonumber
{\dot S}_\ell=
D_{Q',Q_K, R} {\mathcal S}_\ell(\hat g,\hat \phi, \nabla
 \hat \phi,Q^{\prime},Q_{K},R)\bigg|_{Q'=0,Q_K=0,R=0}
 \left(\begin{array}{c}
{\bf q'}
\\
{\bf z}
\\
{\bf r}
 \end{array}\right),
 \eeq
where 
${\bf r}=
\hat g^{pk}\hat \nabla_p ({\bf p}_{jk}+\hat g_{jk}{\bf z})$.

Functions $\dot u=(\dot g,\dot \phi)$ satisfy
the  linearized Einstein-scalar field equation (\ref{eq: final linearization}).
The  linearized Einstein-scalar field equation (\ref{eq: final linearization})
is
\ba
& &e_{pq}(\dot g)-t^{(1)}_{pq}(\dot g)-t^{(2)}_{pq}(\dot \phi)={\bf f}^1_{pq}
\\
& &\square_{\hat g}\dot \phi^\ell
-\hat g^{nm}\hat g^{kj}(\p_n\p_j\hat \phi_\ell+ \hat \Gamma_{nj}^p\p_p
\hat \phi_\ell)\dot g_{mk}
-m^2\dot \phi^\ell= {\bf f}^2_\ell,
\ea
where 
\beq\label{eq: formula for sources 2}
& &{\bf f}^1_{pq}=
{\bf p}_{pq}-(\sum_{\ell=1}^L \dot S_\ell \hat \phi_\ell)\hat g_{pq},\quad
\hbox{and }
-(\sum_{\ell=1}^L \dot S_\ell \hat \phi_\ell)={\bf z},\\
\nonumber 
& & {\bf f}^2_\ell=\dot S_\ell. \nonumber 
\eeq
By Theorem \ref{thm: Good S functions} (ii),
  $u_\e=(g_\e,\phi_\e)$ satisfy  the conservation
 law (\ref{conservation law CC}).
This yields that  $\dot u=(\dot g,\dot \phi)$ satisfies the
  linearized Einstein-scalar field equation (\ref{eq: final linearization})
 and  the linearized conservation law (\ref{eq: lineariz. conse. law PRE}) is valid, too.
The linearized  conservation law (\ref{eq: lineariz. conse. law PRE}) gives,
by the considerations before (\ref {eq: formula for sources 2}),
\beq\label{lin cons 2B}
& &\hspace{-1.5cm}\left (\sum_{\ell=1}^L {\bf f}^2_\ell \, \p_j\hat \phi_\ell\right)
+ \hat g^{pk}\hat \nabla_p {\bf f}^1_{kj}=0.\hspace{-1cm}
\eeq

Below,
we use the adaptive source functions
 $\mathcal S_\ell$  constructed in Theorem \ref{thm: Good S functions} and
 its proof.
 
%

{We see that 
\beq\label{eq: source f}
{\bf f}=F(x;{\bf p},{\bf q})=(F^{(1)}(x;{\bf p},{\bf q}),F^{(2)}(x;{\bf p},{\bf q}))
\eeq
has by formulas (\ref{eq: formula for sources 1}) and (\ref{eq: formula for sources 2})
and Assumption S 
the form
\beq
\label{lin wave eq source A1a} F^{(1)}_{jk}(x;{\bf p},{\bf q})={\bf p}_{jk}+{\bf z}(x) \hat g_{jk}(x)\eeq
  and
\beq
\label{lin wave eq source B}
F^{(2)}(x;{\bf p},{\bf q})
&=& M_{(2)} {\bf q^{\prime}}+ L_{(2)}{\bf z}+N^j_{(2)} \hat g^{lk}\hat \nabla_l ( {\bf p}_{jk}+{\bf z}\hat g_{jk}),
  \eeq
  where 
  \ba
  & &M_{(2)}=M_{(2)}(\hat \phi(x),\hat \nabla  \hat \phi(x),\hat g(x)),\\
  & & L_{(2)}=L_{(2)}(\hat \phi(x),\hat \nabla  \hat \phi(x),\hat g(x)),\\
  & &N^j_{(2)}=N^j_{(2)}(\hat \phi(x),\hat \nabla  \hat \phi(x),\hat g(x))
  \ea
  are, in local coordinates, matrices whose elements  are smooth functions 
of  $\hat \phi(x),\hat \nabla  \hat \phi(x),$ and $\hat g(x)$.
By Thm.\ \ref{thm: Good S functions} (i), the union of the image spaces of  the matrices 
$M_{(2)}(x)$ and
$L_{(2)}(x) $,  
and $N^j_{(2)}(x) $,  $j=1,2,3,4$,
span the space $\R^L$  for all $x\in U_{\hat g,\sigma}$.
}

Consider $n\in \Z$, $t_0,s_0>0$, $Y=Y(x_0,\zeta_0;t_0,s_0)$, $K=K(x_0,\zeta_0;t_0,s_0)$
and $ (x,\xi)\in N^*Y$ (to recall the definitions of these notations, see formula (\ref{associated submanifold})
and definitions below it).
Let  $\mathcal Z=\mathcal Z (x,\xi)$ be the set of the values of 
the principal symbol $\tilde f(x,\xi)=(\tilde f_1(x,\xi),\tilde f_2(x,\xi))$,
at $(x,\xi)$,  
of the source ${\bf f}=(f_1(x),f_2(x))\in \I^{n}(Y)$ that
 satisfy  the linearized conservation law for principal symbols (\ref{mu linearized conservation law for symbols}). 
 
 We use the following auxiliary result:


\begin{lemma}\label{lem: Lagrangian} 
Assume the the Condition A is satisfied and $\hat Q=0$  and $\hat P=0$.
Let $k_0\geq 8$,  $s_1\geq k_0+5$, and 
$Y\subset W_{\hat g}$
be a 2-dimensional space-like submanifold and $y\in Y$, $\xi\in N^*_yY$,
and let $\W$ be a conic neighborhood of $(y,\xi)$ in $T^*M$. Also, let $y\in U_{\hat g,\sigma}$
with some permutation $\sigma$.
Let us consider an open, relatively compact local coordinate neighborhood $V\subset  U_{\hat g,\sigma}$ of $y$ such that
in the coordinates $X:V\to \R^4$, $X^j(x)=x^j$,
we have $X(Y\cap V)\subset \{x\in \R^4;\ x^jb^1_j=0,\  x^jb^2_j=0\}$.
Let  $n_1\in \Z_+$  be sufficiently large and  $n\leq -n_1$.
 Let us consider
 ${\bf p},{\bf q^{\prime}},{\bf z}\in \I^{n}(Y)$,  supported in $V$, 
 that  have classical symbols 
 with principal symbols $\tilde \symbolP ^{(a)}(x,\xi),$ $\tilde \symbolQ_{1}^{(a)}(x,\xi)$,
 $\tilde \symbolQ_2^{(a)}(x,\xi)$, correspondingly, at $(x,\xi)\in N^*Y$.
 Moreover, assume that the principal symbols of ${\bf p}$ and ${\bf z}$ 
 satisfy the linearized conservation law 
 for the principal symbols, that is, (\ref{eq: lineariz. conse. law symbols mLS}),  at 
 all $N^*Y\cap N^*K$ and assume that they vanish outside the conic neighborhood $\W$ of
 $(y,\xi)$ in $T^*M$. Let ${\bf f}=(f^1,f^2) \in\I^{n}(Y)$ be given by
(\ref{eq: formula for sources 1}) and (\ref{eq: formula for sources 2}).

Then  the principal symbol $\tilde f(y,\xi)=(\tilde f_1(y,\xi),\tilde f_2(y,\xi))$ 
of the source ${\bf f}$ at $(y,\xi)$ is the set $\mathcal Z=\mathcal Z(y,\xi)$.
 Moreover,  by varying ${\bf p},{\bf q^{\prime}},{\bf z}$  so that the linearized conservation law  (\ref{eq: lineariz. conse. law symbols mLS})  for principal symbols 
is satisfied, {the principal symbol  $\tilde f(y,\xi)$} at $(y,\xi)$
 achieves all values in the $(L+6)$ dimensional space $\mathcal Z$.

\end{lemma}

\noindent 
{\bf Proof.} 
{
Let us use local coordinates $X:V\to \R^4$ where $V\subset M_0$ is a neighborhood of $x$.
In these coordinates, let $\tilde \symbolP ^{(b)}(x,\xi)$ and $\tilde \symbolQ_{2}^{(b)}(x,\xi)$ denote the sub-principal
symbols of ${\bf p}$ and ${\bf z}$, respectively, at $(x,\xi)\in N^*Y$. 
 Moreover, let $\tilde \symbolP ^{(c)}_j(x,\xi)=\frac \p {\p x^j}\tilde \symbolP ^{(a)}(x,\xi)$ and $\tilde d^{(c)}_j(x,\xi)=\frac \p {\p x^j}\tilde d^{(a)}_2(x,\xi)$, $j=1,2,3,4$ be the $x$-derivatives of the principal symbols
and let us denote \ba
\tilde \symbolP ^{(c)}(x,\xi)=(\tilde \symbolP ^{(c)}_j(x,\xi))_{j=1}^4,\quad 
\tilde d^{(c)}(x,\xi)=(\tilde d^{(c)}_{j}(x,\xi))_{j=1}^4.
\ea


 Let ${\bf f}=({\bf f}_1,{\bf f}_2)=F(x;{\bf p},{\bf q})$ 
 be defined by (\ref{lin wave eq source A1a}) and
 (\ref{lin wave eq source B}).
{When the principal symbols of  ${\bf p},{\bf q^{\prime}},{\bf z}\in \I^{n}(Y)$ 
 are such that  the linearized conservation law  (\ref{eq: lineariz. conse. law symbols mLS})  for principal symbols 
is satisfied, we see that  ${\bf f}\in \I^{n}(Y)$ has the principal symbol
  $\tilde f(x,\xi)=(\tilde f_{1}(x,\xi),\tilde f_{2}(x,\xi))$ at $(x,\xi)$, given by  
  \ba
  & &\tilde  f_{1}(x,\xi) = s_1(x,\xi),\\
\nonumber 
& & \tilde  f_{2}(x,\xi)= s_2(x,\xi),
\ea
where   
     \beq\label{lin wave eq source symbols}
& &s_{1}(x,\xi) = (\tilde \symbolP ^{(a)}+\hat g\tilde \symbolQ_2^{(a)})(x,\xi),\\
\nonumber 
& & s_{2}(x,\xi)= \bigg(
 M_{(2)}(x)\tilde \symbolQ_{1}^{(a)} + J_{(2)}(\tilde \symbolP ^{(c)}+\hat g\tilde d^{(c)}) +\\
 \nonumber
 & &\quad +
L_{(2)}\tilde \symbolQ_2^{(a)} +N^j_{(2)} \, \hat g^{lk}\xi_l (\tilde \symbolP_{1}^{(b)}+\hat g\tilde \symbolQ_2^{(b)})_{jk}\bigg)(x,\xi).
  \eeq
Here, roughly speaking, the $J_{(2)}$ term appears when the $\nabla$-derivatives in $R$ hit to the symbols of the conormal distributions having the form
(\ref{eq: b b-prime BBB}).  
We emphasize that here the symbols $s_1(x,\xi)$  and $s_2(x,\xi)$ are 
well defined objects (in fixed local coordinates) also when   the linearized conservation law  (\ref{eq: lineariz. conse. law symbols mLS})  for principal symbols is not valid.  When (\ref{eq: lineariz. conse. law symbols mLS})
 is valid, ${\bf f}\in \I^{n}(Y)$ and $s_1(x,\xi)$  and $s_1(x,\xi)$
coincide with the principal symbols of ${\bf f}_1$ and ${\bf f}_2$.}

   Observe that  the map $(c^{(b)}_{jk})\mapsto (\hat g^{lk}\xi_lc^{(b)}_{jk})_{j=1}^4$, 
 defined as
   $\hbox{Symm}(\R^{4\times 4})\to \R^4$, is surjective.
  Denote 
  \ba
& &  \tilde m^{(a)}=(\tilde \symbolP ^{(a)}+\hat g\tilde \symbolQ_2^{(a)})(x,\xi),
 \\
 & &\tilde m^{(b)}=(\tilde \symbolP ^{(b)}+\hat g\tilde \symbolQ_2^{(b)})(x,\xi),
\\
 & &\tilde m^{(c)}=(\tilde \symbolP ^{(c)}+\hat g\tilde d^{(c)})(x,\xi).
 \ea
 As noted above,
by {(\ref{eq: deriv1})},   the union of the image spaces of  the matrices 
$M_{(2)}(x)$ and
$L_{(2)}(x) $,  
and $N^j_{(2)}(x) $,  $j=1,2,3,4$,
span the space $\R^L$  for all $x\in U_{\hat g}$.  
Hence the map 
  \ba
  {\bf A}:(\tilde m^{(a)},\tilde m^{(b)},\tilde m^{(c)},
  \tilde \symbolQ_{1}, \tilde \symbolQ_2^{(a)})|_{(x,\xi)}\mapsto
{(s_1(x,\xi),s_2(x,\xi)),}
  \ea 
  given by (\ref{lin wave eq source symbols}), considered as a 
  map ${\bf A}:\mathcal Y=(\hbox{Symm}(\R^{4\times 4}))^{1+1+4}\times \R^{K}\times \R\to
  \hbox{Symm}(\R^{4\times 4})\times \R^L$,
  is surjective. 
  Let $\mathcal X$ be  the set of
  elements $(\tilde m^{(a)}(x,\xi),\tilde m^{(b)}(x,\xi),\tilde m^{(c)}(x,\xi),
  \tilde \symbolQ_{1}^{(a)}(x,\xi),$ $\tilde \symbolQ_2^{(a)}(x,\xi))\in \mathcal Y$
  where 
   $\tilde m^{(a)}(x,\xi)=(\tilde \symbolP ^{(a)}+\hat g\tilde \symbolQ_2^{(a)})(x,\xi)$ is such that the  pair  $(\tilde \symbolP ^{(a)}(x,\xi),\tilde \symbolQ_2^{(a)}(x,\xi))$  satisfies
 the linearized conservation law for principal symbols, see (\ref{eq: lineariz. conse. law symbols mLS}). 
 Then    $\mathcal X$ has codimension 4 in $Y$ and we see
 that the image  ${\bf A}(\mathcal X)$ has  in $ \hbox{Symm}(\R^{4\times 4})\times \R^L$ co-dimension less or equal
 to 4.
 
By (\ref{lin cons 2B})}) and the considerations above it, we have that $\bf f$ satisfies the linearized
 conservation law 
(\ref{eq: lineariz. conse. law PRE}). This implies that  its
principal symbol ${\bf A}((\tilde m^{(a)}(x,\xi),\tilde m^{(b)}(x,\xi),$ $ \tilde m^{(c)}(x,\xi),
  \tilde \symbolQ_{1}^{(a)}(x,\xi),$ $\tilde \symbolQ_2^{(a)}(x,\xi)))$ 
has to satisfy the linearized conservation law for principal symbols 
(\ref{mu linearized conservation law for symbols})
and  hence ${\bf A}(\mathcal X)\subset \mathcal Z$. 
As  $\mathcal Z$ has codimension 4, this and the above
prove that    ${\bf A}(\mathcal X)=\mathcal Z$.
\hfill \Box \medskip

Now we are ready to prove the microlocal stability result for 
the Einstein-scalar field equation (\ref{eq: adaptive model with no source}).
Note that the claim of the following theorem does not 
involve the  adaptive source functions constructed in Theorem \ref{thm: Good S functions}
as these functions are needed only as an auxiliary tool in the proof.

 \begin{theorem}\label{microloc linearization stability thm Einstein} ($\mu$-LS i.e., Microlocal linearization stability) 
Let $k_0\geq 8$,  $s_1\geq k_0+5$, and 
$Y\subset W_{\hat g}$
be a 2-dimensional space-like submanifold, and $y\in Y$
and let $\W$ be a conic neighborhood of $(y,\eta)$ in $T^*M$. 
Also, 
consider in an open local coordinate neighborhood $V\subset W_{\hat g}$ of $y\in Y$ such that
in the coordinates $X:V\to \R^4$, $X^j(x)=x^j$,
we have $X(Y\cap V)\subset \{x\in \R^4;\ x^jb^1_j=0,\  x^jb^2_j=0\}$.
Then there there is  $n_1\in \Z_+$ such that if $(y,\eta)\in N^*Y$ is light-like, $n\leq -n_1$ and  $\tilde c_{jk}(x,\theta)$ and $\tilde d_{\ell}(x,\theta)$
 are positively $n$-homogeneous in the $\theta$-variable in the domain $ \{(x,\theta)\in V\times \R^2;\ |\theta|>1\}$  
and satisfy 
\beq\label{mu linearized conservation law for symbols BB}
& &\hat g^{lk}(\theta_1  b^1_l
+\theta_2  b^2_l)
\tilde c_{kj} (x,\theta_1,\theta_2)
=0,\quad |\theta|>1,\quad j=1,2,3,4,
\eeq
then there are $f^1\in \I^{n}(Y)$
and $f^2\in \I^{n}(Y)$
 supported in $V$ such that the principal symbols of these distributions vanish outside $\W$ and are equal to 
$\tilde c_{jk}(y,\eta)$  and $\tilde d_{\ell}(y,\eta)$  at $(y,\eta)$, respectively.
Moreover, 
 $f=(f^1,f^2)$  satisfies the linearized conservation law (\ref{eq: lineariz. conse. law PRE})
and  
 $f$  has the LS-property (\ref{eq: LSp}) 
 in $C^{s_1}(M_0)$
  with a family  $\F_\e$, $\e\in [0,\e_0)$
  such that all functions   $\F_\e$,$\e\in [0,\e_0)$ are supported in $V$.
 \end{theorem}

\noindent {\bf Proof.}
Let $\sigma\in \Sigma(K)$ be such that $y\in U_{\hat g,\sigma}$.
Let
 ${\bf p}$ and ${\bf q}$  be the functions constructed in
 Lemma \ref{lem: Lagrangian} such that
they are supported in $V$ and their principal symbols vanish outside $\W$.  We can assume 
 that these functions are supported in $W_0= V\cap U_{\hat g,\sigma}$. 
Let  $P_\e=\e {\bf p}$ and $Q_\e=\e {\bf q}$   be sources depending on $\e\in \R$
and $u_\e=(g_\e,\phi_\e)$ be the  solution of
 (\ref{eq: adaptive model}) with the sources  $P_\e$  and $Q_\e$.
 Also, let
\ba
& &\F_\e^1=P_\e+Z_\e g_\e ,\quad Z_\e=-(\sum_{\ell=1}^L S^\e_\ell\phi^\e_\ell+
\frac 1{2m^2}(S^\e_\ell)^2),\\
& &(\F_\e^2)_\ell=S^\e_\ell,
\ea
where 
\ba
& &S^\e_\ell= \mathcal S_\ell(g_\e,\phi_\e,\nabla \phi_\e,Q_\e,\nabla Q_\e,P_\e,
\nabla^{g_\e} P_\e),
\ea
where 
 $\mathcal S_\ell$ are the adaptive source functions constructed in Theorem \ref{thm: Good S functions} and
 its proof.

By (\ref{S-vanish condition}) also $S^\e_\ell$ and 
the family $\F_\e$, $\e\in [0,\e_0]$ of non-linear sources are supported in  $V_0$
and we have shown that  $u_\e=(g_\e,\phi_\e)$ and $\F_\e$
satisfy
the reduced Einstein-scalar field equation (\ref{eq: adaptive model with no source})
and  the conservation
 law (\ref{conservation law CC}).
This proves Theorem \ref{microloc linearization stability thm Einstein}.
\hfill \Box \medskip

  {\extension 
  
  \subsection*{C.3. A special case when the whole metric is determined}
  Finally, we consider the case when  $\hat Q$ and $\hat P$ are non-zero, and we define
\ba
{\cal D}^{mod}(\hat g,\hat \phi,\e)=\{[(U_g,g|_{U_g},\phi|_{U_g},F|_{U_g})]&;&(g,\phi,F)\hbox{ are smooth 
solutions,}\\
& & 
\hspace{-6cm}\hbox{ of  
(\ref{eq: adaptive model}) with $F=(P-\hat P,Q-\hat Q)$, }F\in C^{19}_0(W_{g};\B^K),
\\
& &\hspace{-6cm}\hbox{$J_g^+(\supp(F))\cap J_g^-(\supp(F))\subset 
W_g, $ $\mathcal N^{({ {17}})}(F)<\e$, $\mathcal N^{({{17}})}_{\hat g}(g)<\e$}\}.
\ea
Next we consider the case when $\hat Q^{(j)}$ and  $\hat P^{(j)}$ are not assumed to be zero,
see Fig.\ A11. 

\begin{corollary}\label{coro of main thm original Einstein}
Assume that $(\hattuM^{(1)} ,\hat g^{(1)})$ and $(\hattuM^{(2)} ,\hat g^{(2)})$ are  globally hyperbolic
manifolds and $\hat \phi^{(j)}$, $\hat Q^{(j)}$,  and $\hat P^{(j)}$ are background 
fields satisfying (\ref{eq: adaptive model}).
Also, assume that there are neighborhoods $U_{\hat g^{(j)}}$, $j=1,2$ of time-like geodesics $\mu_j\subset \hattuM^{(j)}$ where the  Condition $\mu$-LS is valid,
and points $p_j^-=\mu_j(s_-)$ and $p_j^+=\mu_j(s_+)$.
Moreover,
 assume  
also that for $j=1,2$  there are sets $W_j\subset M_j$ such that 
$\hat \phi^{(j)}$, $\hat Q^{(j)}$ and $\hat P^{(j)}$ are zero (and thus
the metric tensors $\hat g_j$ have vanishing Ricci curvature)
in $W_j$ and that  $I_{\hat g_j}(\hat p^-_j,\hat p^+_j)\subset W_j \cup U_{\hat g}$. 
If there is $\e>0$  such that 
\ba
{\cal D}^{mod}(\hat g^{(1)},\hat \phi^{(1)},\e)={\cal D}^{mod}(\hat g^{(2)},\hat \phi^{(2)},\e)
\ea
then the metric $\Psi^*\hat g_2$ is isometric to $\hat g_1$ in $I_{\hat g_1}(\hat p^-_1,\hat p^+_1)$.
\end{corollary}
\medskip

\begin{center}\label{Fig-4}

\psfrag{1}{}
\psfrag{2}{}
\psfrag{3}{}
\psfrag{4}{$W_{\hat g}$}
\psfrag{5}{$\mu_{z,\eta}$}
\includegraphics[width=4.5cm]{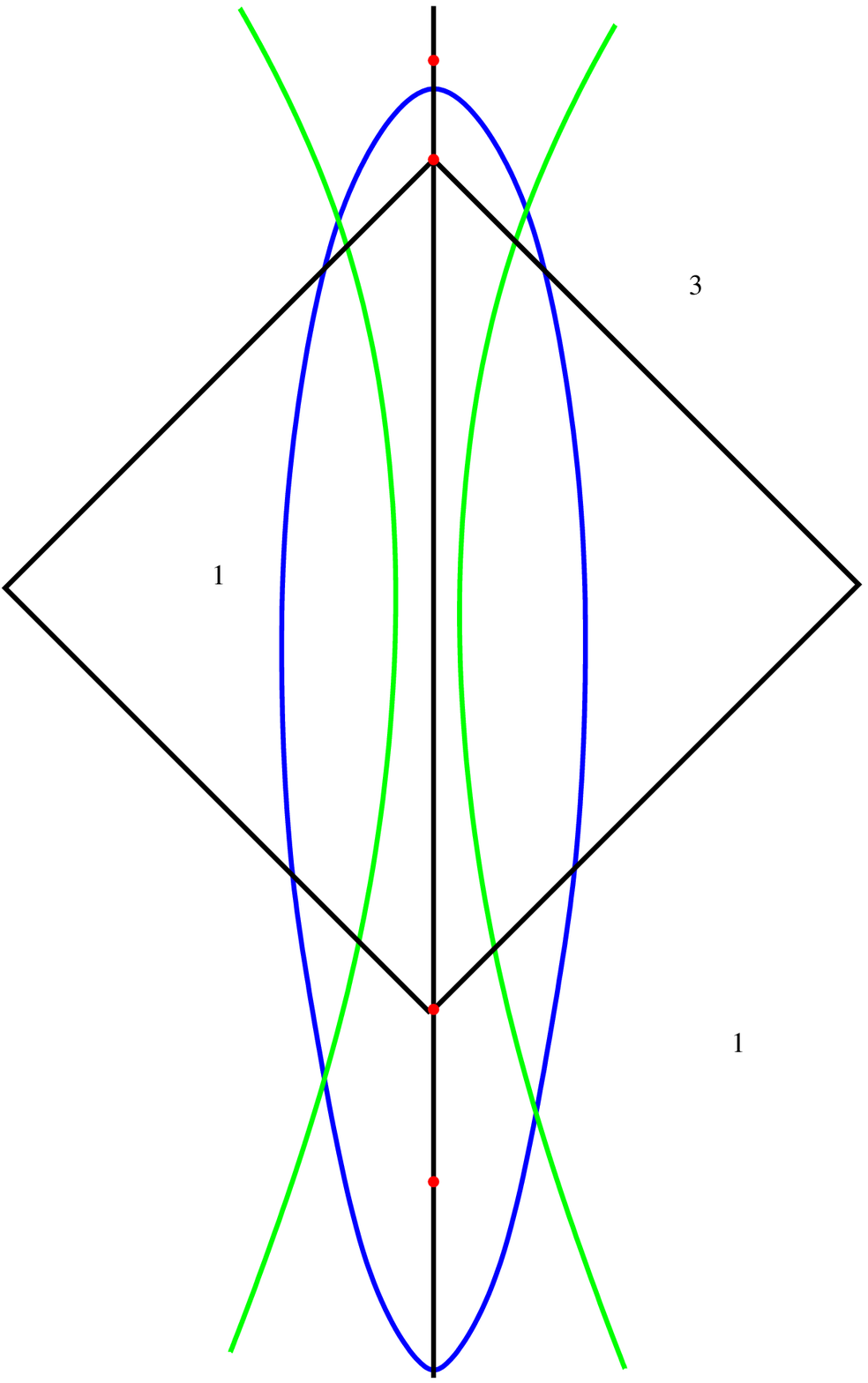}
\end{center}
{\it FIGURE A11: A schematic figure where the space-time is represented as the  2-dimensional set $\R^{1+1}$ on the setting in Corollary \ref{coro of main thm original Einstein}.
The set $J_{\hat g}(\hat p^-,\hat p^+)$, i.e., the diamond with the black boundary, is contained in the union of the blue set $U_{\hat g}$
and the set $W$. The set $W$ is in the figure the area outside of the green curves. 
The sources are controlled in the set $U_{\hat g}\setminus W$ and  the set $W$ consists of vacuum.
}

\medskip

\noindent
{\bf Proof.}\  (of Corollary \ref{coro of main thm original Einstein})  
{In the above proof of  Theorem \ref{alternative main thm Einstein}, we used the  
 assumption that $\hat Q=0$ and $\hat P=0$ to obtain  
 equations (\ref{w1-3 solutions}). In the   
setting of Cor.\ \ref{coro of main thm original Einstein}  
 where the background source fields  
 $\hat Q$ and $\hat P$  
  are not zero,  
we need to assume   
in the computations related to sources (\ref{eq: f vec e sources})  
that there are  
neighborhoods $V_j$ of the geodesics $\gamma_j$ that satisfy  
$\supp({\bf f}_j)\subset V_j$,  
the linearized waves $u_j={\bf Q}_{\hat g}{\bf f}_j$ satisfy  singsupp$(u_j)\subset V_j$,  
and     
$V_i\cap V_j\cap (\supp(\hat Q)\cup \supp(\hat P))=\emptyset$ for $i\not =j$.}  
To this end, we have first consider  measurements for the linearized waves   
and check for given $(\vec x,\vec\xi)$ that no two geodesic $\gamma_{x_j,\xi_j}(\R_+)$  
do intersect at $U_{\hat g}$ and restrict all considerations for such $(\vec x,\vec\xi)$.  
Notice that such $(\vec x,\vec\xi)$ form an open and dense set in   
$(TU_{\hat g})^4$. If then the sources are supported in balls so small
that each ball is in some set $U_{\hat g,\sigma}$ and the width $\hat s$ of the used  
distorted plane  waves is chosen to be small enough, we see that condition   
$V_i\cap V_j\cap (\supp(\hat Q)\cup \supp(\hat P))=\emptyset$ is satisfied.  
  
The above restriction causes only minor modifications in the above
proof and thus, mutatis mutandis, we see that we can determine  
the conformal  
type of the  metric in relatively compact sets 
$I_{\hat g}(\hat\mu(s^\prime),\hat\mu(s^{\prime\prime}))\setminus (\supp(\hat Q)\cup \supp(\hat P))$ for all $s_-<s^\prime<s^{\prime\prime}<s_+$.  
By glueing these  manifolds and $U_{\hat g}$ together, we find the conformal  
type of the  metric in $I_{\hat g}(\hat p^-,\hat p^+)$.  
After this the claim follows similarly to the  
 proof of  Theorem \ref{alternative main thm Einstein} using 
 Corollary 1.3 of \cite{Paper-part-2}.  
\hfill \Box \medskip

In the setting of Corollary \ref{coro of main thm original Einstein} the set
 $W$ is such $I(\hat p^-,\hat p^+)\cap (M\setminus W)\subset U$. This means that if we restrict to the domain $I(\hat p^-,\hat p^+)$
then we have the Vacuum Einstein equations in the unknown domain $I(\hat p^-,\hat p^+)\setminus U$
and have matter only in the domain $U$ where we implement our measurement (c.f.\ a space
ship going around in a system of black holes, see \cite{ChM}). This could be considered as an ``Inverse problem
for the vacuum Einstein equations''. 

  }

\subsection*{C.4. Motivation using Lagrangian formulation}

To motivate the system (\ref{eq: adaptive model}) of partial differential equations, we 
give in this subsection a non-rigorous discussion.

Following  \cite[Ch. III, Sect. 6.4, 7.1, 7.2, 7.3]{ChBook}
and  \cite[p. 36]{AnderssonComer}
we start by considering  the 
Lagrangians, associated to gravity, scalar fields $\phi=(\phi_\ell)_{\ell=1}^L$ and non-interacting fluid fields, that is, the number density four-currents ${\bf n}=({\bf n}_\kappa(x))_{\kappa=1}^J$
(where each ${\bf n}_\kappa$ is a vector field,  see \cite[p. 33]{AnderssonComer}).
We consider also products of  vector fields ${\bf n}_\kappa(x)$ 
and $\frac 12$-density $|\det(g)|^{1/2}$ that denote ${\bf p}_\kappa$,
\beq\label{half density}
{\bf p}_\kappa^j(x)\frac\p{\p x^j}={\bf n}^j_\kappa(x)|\det(g)|^{1/2}\frac\p{\p x^j}.
\eeq
see  \cite[p. 53]{Dirac}. Also, $\rho=(-g_{jk}{\bf n}^j_\kappa{\bf n}^k_\kappa)^{1/2}$ 
corresponds to the energy density of the fluid.
Below, we use the variation of density with respect to the metric,
\beq\label{eee}
\frac {\delta}
{\delta g_{jk}}( 
\sum_{\kappa=1}^J  (-g_{nm}{\bf p}^n_\kappa{\bf p}^m_\kappa)^{1/2})
&=&- \sum_{\kappa=1}^J   \frac 12 (-g_{nm}{\bf p}^n_\kappa{\bf p}^m_\kappa)^{-1/2}
{\bf p}_\kappa^j
{\bf p}_\kappa^k
\\ \nonumber
&=& - \sum_{\kappa=1}^J   \frac 12 \rho\,{\bf n}_\kappa^j
{\bf n}_\kappa^k\, |\det(g)|^{1/2}. 
\eeq
Due to this, we denote  
\beq\label{P formula}\\
\nonumber
P&=& \sum_{\kappa=1}^J   \frac 12 \rho\,{\bf n}^\kappa_j
{\bf n}^\kappa_k dx^j\otimes dx^k,\quad\hbox{where }
{\bf n}^\kappa_k =g_{ki}  {\bf n}_\kappa^i=g_{ki}  {\bf p}_\kappa^i|\det(g)|^{-1/2}.
\eeq

Below, we consider a model for $g$, $\phi$, and ${\bf p}$.
We also add in to the model a Lagrangian  associated with
some scalar valued source fields $S=(S_\ell)_{\ell=1}^L$  and
$Q=(Q_k)_{k=1}^K$. We consider the action corresponding to
the coupled Lagrangian\hiddenfootnote{
An alternative Lagrangian would be the one where the term
$P_{pq}g^{pq}$ is replaced by $\frac 12 P_{pq}P_{jk}g^{pj}g^{qk}$.
Then in the stress energy tensor $P_{jk}$ is replaced by
${\mathcal T}^{(1)}_{jk}= F_{jp}F_{kq}g^{qp}$, i.e. ${\mathcal T}^{(1)}=\hat F\hat G\hat F.$
By [Duistermaat, Prop. 1.3.2], 
its derivative with respect to $P_{jk}$, that is,
$((D_F{\mathcal T}^{(1)}_{jk}|_{\hat F})^{ab})$ (COMPUTE THIS)
is a surjective pointwise map is surjective at points $x\in M$ if
$\hat F_{qp}$ is positive definite symmetric (and thus invertible) matrix (CHECK THIS!!).
}
%
\ba
& &\mathcal A=\int_{M}\bigg(L_{grav}(x)+L_{fields}(x)+L_{source}(x)\bigg)\,dV_g(x),\\
& &L_{grav}=\frac 12 R(g),\\
& &L_{fields}=\sum_{\ell=1}^L\bigg(-\frac 12 g^{jk}\p_j\phi_\ell \,\p_k\phi_\ell 
-{\mathcal V}( \phi_\ell;S_\ell)\bigg)
+  \\
& &\quad\quad\quad\quad\quad+\sum_{\kappa=1}^J\bigg(- \frac 12 (-g_{jk}{\bf p}^j_\kappa {\bf p}^k_\kappa)^{\frac 12}
\bigg)|\det (g)|^{-\frac 12},
\\
& &L_{source}=\e \mathcal H_\e(g,S,Q, {\bf p},\phi),
\ea
where $R(g)$ is the scalar curvature, $dV_g=(-\det (g))^{1/2}dx$ is the volume form on $(M,g)$, 
\beq\label{eq: Slava-a}
{\mathcal V}( \phi_\ell;S_\ell)=\frac 12m^2\bigg(\phi_\ell +\frac 1{m^2}S_\ell\bigg)^2
\eeq are energy potentials
of the scalar fields $\phi_\ell$ that depend on $S_\ell$,
 and  $\mathcal H_\e(g,S,Q, {\bf p},\phi)$ is a function modeling the measurement
 device we use.
 We assume that $\mathcal H_\e$ is bounded and   its derivatives with respect
 to $S,Q,{\bf p}$  are very large (like of order $O((\e)^{-2})$) and its derivatives
with respect of $g$ and $\phi$ are bounded when $\e>0$ is small.
We note that the  above Lagrangian for the fluid fields is the sum of the single fluid Lagrangians.
where for all fluids the master function  $\Lambda(s)=s^{1/2}$, that is, the energy density
of each fluid is given by
$\rho=\Lambda(-g_{jk}{\bf n}^j{\bf n}^k)$. On fluid Lagrangians, see the discussions in
 \cite[p. 33-37]{AnderssonComer},  \cite[Ch. III, Sect. 8]{ChBook}, \cite[p. 53]{Dirac}, and
 \cite{Taub} and \cite[p.\ 196]{Felice}.

%
When  we compute the critical points of the Lagrangian $L$ and 
neglect the  $O(\e)$-terms, the equation  $\frac {\delta  \mathcal A}{\delta g}=0$,
together with formulas (\ref{eee}) and  (\ref{P formula}), 
give the Einstein equations  with a stress-energy tensor $T_{jk}$
defined in (\ref{eq: T tensor1}).
The equation  $\frac {\delta  \mathcal A}{\delta \phi}=0$
gives the wave equations with sources $S_\ell$. We assume that $O(\e^{-1})$ order equations
obtained from 
the equation $(\frac {\delta \mathcal A}{\delta S},\frac {\delta  \mathcal A}{\delta Q},
\frac {\delta  \mathcal A}{\delta {\bf p}})=0$ 
fix the values of the scalar functions $Q$ and 
the fields ${\bf p}^\kappa$, $\kappa=1,2,\dots,J$, and moreover, 
  yield for the sources $S=(S_\ell)_{\ell=1}^L$ equations of the form
\beq\label{S! formula}
S_\ell=\mathcal S_\ell(g,\phi, \nabla
 \phi,Q,\nabla Q,P,\nabla^g P)
 \eeq
where $P$ is given by (\ref{P formula}).
Let us aslo write (\ref{S! formula}) using different notations, as
\ba
S_\ell=Q_\ell+\mathcal S_\ell^{2nd}(g,\phi, \nabla
 \phi,Q,\nabla Q,P,\nabla^g P).
 \ea 
 Summarizing, we have obtained, up to the above used approximations,
the model  (\ref{eq: adaptive model}). However, note that above the field $P$  is not directly controlled 
but instead, we control ${\bf p}$ and the value of the field $P$ is determined by
the solution $n$ and  formula  (\ref{P formula}). In this sense $P$  is not controlled,
but an observed field. 

 Above, the function
 $\mathcal H_\e$ models the way the measurement device works. 
Due to this we will assume that $\mathcal H_\e$ and thus functions ${\mathcal S}_\ell$ may
be quite complicated.   The interpretation of the above is that in  each measurement event we use a device that
fixes the values of the scalar functions $Q$, ${\bf p}$, and gives the equations
for $S^{2nd}$ that 
 tell how the sources of the $\phi$-fields adapt to these changes
 so that 
the physical conservation laws are satisfied.

\section*{Appendix D: An inverse problem for a non-linear wave equation}

In this appendix  we explain how a problem for a scalar wave equation can
be solved with the same techniques that we used for the Einstein equations.

Let $(M_j,g_j)$, $j=1,2$ be two globally hyperbolic  $(1+3)$  
dimensional Lorentzian
manifolds represented using
  global smooth time functions as $M_j=\R\times N_j$, $\mu_j=\mu_j([-1,1])
\subset M_j$ be
a time-like geodesic and $U_j\subset M_j$ be open, relatively compact  
neighborhood
of $\mu_j([s_-,s_+])$, $-1<s_-<s_+<1$. Let $M_j^0=(-\infty,T_0)\times  
N_j$ where $T_0>0$ is such
that $U_j\subset  M_j^0$.
{\mltext Consider the
non-linear wave equation
\beq\nonumber
& &\hspace{-.5cm}  \square_{g_j}u(x)+a_j(x)\,u(x)^2=f(x)
\quad\hbox{on }M_j^0,\hspace{-.5cm}\\
\label{eq: wave-eq general} & &\quad \supp(u)\subset J^+_{g_j}(\supp(f)),
\eeq
where $\supp(f)\subset U_j,$
\ba
\square_gu=
\sum_{p,q=1}^4(-\det (g))^{-1/2}\frac  \p{\p x^p}
\left ((-\det(g))^{1/2}
g^{pq}\frac \p{\p x^q}u(x)\right),
\ea
$\det(g)=\det((g_{pq}(x))_{p,q=1}^4)$,
  $f\in C^6_0(U_j)$ is a controllable source, and $a_j$ is
a non-vanishing $C^\infty$-smooth  function.}
Our goal is to prove the following result:

\begin{theorem}\label{main thm3}
Let $(M_j,g_j)$, $j=1,2$ be two
open,  smooth, globally hyperbolic    Lorentzian manifolds of  
dimension $(1+3)$.
Let
$p^+_j=\mu_j(s_+), p^-_j=\mu_j(s_-)\in M_j$ the points of a  time-like  
geodesic  $\mu_j=
\mu_j([-1,1])\subset M_j$, $-1<s_-<s_+<1$,
and let  $U_j\subset M_j$ be an open  relatively compact neighborhood
of $\mu_j([s_-,s_+])$  given in (\ref{eq: Def Wg with hat}).  Let  
$a_j:M_j\to \R$, $j=1,2$ be $C^\infty$-smooth functions that are
non-zero on $M_j$.

Let   $L_{U_j}$, $j=1,2$ be measurement operators defined
in an open set $\mathcal W_j\subset C^6_0(U_j)$ containing the zero  
function by setting
\beq\label{measurement operator}
L_{U_j}: f\mapsto u|_{U_j},\quad f\in C^6_0(U_j),
\eeq
where $u$ satisfies the wave equation (\ref{eq: wave-eq general}) on  
$(M_j^0,g_j)$.

Assume that there is a diffeomorphic isometry
$\Phi:U_1\to U_2$ so that $\Phi(p^-_1)=p^-_2$ and
$\Phi(p^+_1)=p^+_2$ and the measurement maps
satisfy
\ba
((\Phi^{-1})^*\circ L_{U_1}\circ \Phi^*) f =L_{U_2}f
\ea
for all $f\in \W$ where $\W$ is  some  neighborhood of the zero function in $C^6_0(U_2)$.

Then there is a diffeomorphism $\Psi:I(p^-_1,p^+_1)\to I(p^-_2,p^+_2)$,
and the metric $\Psi^*g_2$ is conformal to $g_1$ in  
$I(p^-_1,p^+_1)\subset M_1$,
that is, there is $\beta(x)$ such that $g_1(x)=\beta(x)(\Psi^*g_2)(x)$ in  
$I(p^-_1,p^+_1)$.
\end{theorem}

We note that the  smoothness assumptions assumed above on the functions
$a$ and the source $f$ are not optimal.
The proof, presented below, is based on using the interaction
of singular waves.  The techniques used can be modified to study
different non-linearities, such as the equations
$\square_{\hat g} u+a(x)u^3=f$, $\square_{\hat g} u+a(x)u_t^2=f$,
or $\square_{g(x,u(x))} u=f$, but these considerations are outside the scope
of this paper.

Theorem \ref{main thm3} can be applicable
  for example in the mathematical analysis of  non-destructive
testing or imaging in non-linear medium e.g,
in imaging  the non-linearity  of the acoustic material parameter  inside a
given body  when it
is under large, time-varying, possibly periodic, changes of the external
pressure and at the same time the body is probed with small-amplitude fields.
  Such acoustic measurements are analogous to
  the recently developed Ultrasound Elastography imaging technique where the interaction
  of the elastic shear and pressure
waves is used for medical imaging, see  
e.g.\ \cite{Hoskins,McLaughlin1,McLaughlin2,Ophir}. There, the slowly  
progressing
shear wave is imaged using a pressure wave and the
image of the shear wave inside the body is used to determine
approximately the material parameters. In other words,
the changes which the elastic wave causes in the medium
are imaged using the interaction of the s-wave and p-wave
components of the  elastic wave.

Let us also consider some implications of theorem \ref{main thm3}  for inverse
problems for a non-linear equation involving a time-independent metric
\beq\label{product metric} & &g(t,y)=-dt^2+
\sum_{\alpha,\beta=1}^3 h_{\alpha\beta}(y)dy^\alpha dy^\beta,\quad
(t,y)\in \R\times N.
\eeq
The metric (\ref{product metric}) corresponds to
the hyperbolic operator $\p_t^2- \Delta_{h}$,
with a time-independent Riemannian metric  
$h=(h_{\alpha\beta}(y))_{\alpha,\beta=1}^3,$
$y\in N$, where $N$ is a 3-dimensional manifold and
\ba
  \Delta_{h}u(y,t)=\sum_{\alpha,\beta=1}^3\frac \p{\p y^\alpha}
\left( h^{\alpha\beta}(y) \frac \p{\p y^\beta}u(t,y)\right).
\ea

\begin{corollary}\label{coro thm3}
Let $(M_j,g_j)$, $M_j=\R\times N_j$, $j=1,2$ be two
open,  smooth, globally hyperbolic    Lorentzian manifolds of  
dimension $(1+3)$.
Assume that $g_j$ is the product metric of the type (\ref{product metric}),
$g_j=-dt^2+h_j(y)$, $j=1,2$.
Assume that $\mu_j$ is a time-like geodesic  $\mu_j([-1,1])\subset  \R\times \{p_j\}$, where
$p_j\in M_j$.

Let
$p^+_j=\mu_j(s_+),p^-_j=\mu_j(s_-)\in (0,T_0)\times N_j$, $-1<s_-<s_+<1$ and 
 assume that $\mu_j(s_-)=(1,p_j)$.
Moreover, 
Let  $U_j\subset (0,T_0)\times N_j$ be an open  relatively compact neighborhood
of $\mu_j([s_-,s_+])$ given in (\ref{eq: Def Wg with hat}). 
  Let $a_j:M_j\to \R$, $j=1,2$ be  
$C^\infty$-smooth functions that are
non-zero on $M_j$ and $x=(t,y)\in \R\times N$.

For $j=1,2$, consider the
non-linear wave equations
\beq\nonumber
& &\hspace{-.5cm}  (\frac {\p^2}{\p  
t^2}-\Delta_{{h_j}})u(t,y)+a_j(y,t)(u(t,y))^2=f(t,y)
\quad\hbox{on }(0,T_0)\times N_j,\hspace{-.5cm}\\
\label{eq: wave-eq} & &\quad \supp(u)\subset J^+_{g_j}(\supp(f)),
\eeq
where $f\in C^6_0(U_j)$, $j=1,2$.
Let $L_{U_j}:f\mapsto u|_{U_j}$ be the measurement operator (\ref{measurement operator})
for the wave equation (\ref{eq: wave-eq}) with the  Riemannian metric $h_j(x)$ and the coefficient
$a_j(x,t)$ for $j=1,2$,  defined
in some $C^6_0(U_j)$ neighborhood of the zero function.

Assume that there is a diffeomorphism
$\Phi:U_1\to U_2$
of the form $\Phi(t,y)=(t,\phi(y))$  so that
\ba
((\Phi^{-1})^*\circ L_{U_1}\circ \Phi^*)f=L_{U_2} f
\ea
for all $f\in \W$ where $\W$ is  some  neighborhood of the zero function in $C^6_0(U_2)$.

Then there is a diffeomorphism $\Psi:I^+(p^-_1)\cap I^-(p^+_1)\to  
I^+(p^-_2)\cap I^-(p^+_2)$ of the form
$\Psi(t,y)=(\psi(y),t)$,
the metric $\Psi^*g_2$ is isometric to $g_1$ in $I^+(p^-_1)\cap I^-(p^+_1)$,
and $a_1(t,y)=a_2(t,\psi(y))$ in $I^+(p^-_1)\cap I^-(p^+_1)$.
\end{corollary}

\medskip

Next we consider the proofs.

\medskip

 {\bf Proof.} (of Theorem \ref{main thm3}).
We will explain how the proof of
Theorem \ref{alternative main thm Einstein} for the Einstein equations  needs to
be modified to obtain the similar result for the non-linear wave equation.

Let  $(M,\hat g)$ be a  smooth
globally hyperbolic Lorentzian manifold that we represent
using a global smooth time function as $M=(-\infty,\infty)\times N$,
and consider   $M^0=(-\infty,T)\times N\subset M$.
Assume that
the set $U$, where the sources are supported and where
we observe the waves, satisfies
$U\subset [0,T]\times N$.

The results of  section \ref{subsec: Direct problem} concerning the  
direct problem
for Einstein equations can be modified for the wave equation
\beq\label{PABC eq}
& &\square_{\hat g}u+au^2=f,\quad\hbox{in }M^0=(-\infty,T)\times N,\\
& &u|_{(-\infty,0)\times N}=0,\nonumber
\eeq
where
$a=a(x)$ is a smooth, non-vanishing function. Here we denote
the metric by $\hat g$ to emphasize the fact that it is independent
on the solution $u$.  Below,
let  $Q$ be the causal inverse  operator of $\square_{\hat g}$.

When $f$ in $C_0([0,t_0];H^6_0(B))\cap
C_0^1([0,t_0];H^5_0(B))$
is small enough,
we see by using  \cite[Prop.\ 9.17]{Ringstrom} and   \cite [Thm.\ III]{HKM},
see also (\ref{stability}) in Appendix B, that  the
equation (\ref{PABC eq}) has a unique solution
$u\in  C([0,t_0];H^{5}(N))\cap C^1([0,t_0];H^{4}(N))$.
Moreover, 
we  can consider the case when $f=\e f_0$ where $\e>0$
is small.
Then, we can write
\ba
u=\e w_1+\e^2 w_2+\e^3 w_3+\e^4 w_4+E_\e
\ea
where $w_j$ and the reminder term $E_\e$ satisfy
\ba
w_1&=&Qf,\\
w_2&=&-Q(a\,w_1\,w_1),\\
w_3&=&-2Q(a\, w_1\,w_2)\\
&=&2Q(a\, w_1\,Q(a\,w_1\,w_1))
,\\
w_4
&=&-Q(a\, w_2\,w_2)-2Q(a\, w_1\,w_3)\\
&=&-Q(a\, Q(a\,w_1\,w_1)\,Q(a\,w_1\,w_1))\\
& &+4Q(a\, w_1\,Q(a\, w_1\,w_2))
\\
&=&-Q(a\, Q(a\,w_1\,w_1)\,Q(a\,w_1\,w_1))\\
& &-4Q(a\, w_1\,Q(a\, w_1\,Q(a\,w_1\,w_1))),\\
& &\hspace{-1.5cm}\|E_\e\|_{C([0,t_0];H^{4}_0(N))\cap  
C^1([0,t_0];H^{3}_0(N))}\leq C\e^5.
\ea

If we consider sources  $f_{\vec\e}(x)=\sum_{j=1}^4\e_j f_{(j)}(x)$,
$\vec\e=(\e_1,\e_2,\e_3,\e_4),$
and the corresponding solution $u_{\vec \e}$ of (\ref{PABC eq}), we see that
\beq \nonumber
\M^{(4)}&=&\p_{\vec \e}^4u_{\vec \e}|_{\vec\e=0}\\
&=& \nonumber
\p_{\e_1}\p_{\e_2}\p_{\e_3}\p_{\e_4}u_{\vec \e}|_{\vec\e=0}\\
\label{4th interaction for wave eq B}
&=&-\sum_{\sigma\in \Sigma(4)}\bigg(Q(a\,  
Q(a\,u_{(\sigma(1))}\,u_{(\sigma(2))})\,Q(a\,u_{(\sigma(3))}\,u_{(\sigma(4))}))\\  
\nonumber
& &\quad+ 4Q(a\, u_{(\sigma(1))}\,Q(a\,  
u_{(\sigma(2))}\,Q(a\,u_{(\sigma(3))}\,u_{(\sigma(4))})))\bigg),
\eeq
where $u_{(j)}=Qf_{(j)}$ and $\cell$ is the set
of permutations of the set $\{1,2,3,\dots,\ell\}$.

The results of Lemma \ref{lem: Lagrangian 1} can
be replaced by the results of \cite[Prop.\ 2.1]{GU1} as follows.
Using the same notations as in  Lemma \ref{lem: Lagrangian 1}, let
$Y=Y(x_0,\zeta_0;t_0,s_0)$, $K=K(x_0,\zeta_0;t_0,s_0)$, and  
$\Lambda_1=\Lambda(x_0,\zeta_0;t_0,s_0)$, and consider a source $f\in  
\I^{n+1}(Y)$.  Then $u=Qf$ satisfies
$u|_{M_0\setminus Y}\in
  \I^{n-1/2} ( M_0\setminus Y;\Lambda_1)$. Assume that  
$(x,\xi),(y,\eta)\in L^+M$  are on the same bicharacteristics of  
$\square_{\hat g}$,
  and $x<y$, that is, $((x,\xi),(y,\eta))\in \Lambda_{\hat g}^\prime$.  
Moreover, assume
  that $(x,\xi)\in N^*Y$.
Let   $\tilde b(x,\xi)$ be the principal
   symbol of $f$ at $(x,\xi)$ and
    $\tilde a(y,\eta)$ be the principal
   symbol of $u$ at $(y,\eta)$. Then $\tilde a(y,\eta)$
   depends linearly on $ \tilde f(x,\xi)$ and
    $\tilde a(y,\eta)$ vanishes if and only if
   $ \tilde f(x,\xi)$ vanishes.

Analogously to the Einstein equations,
we consider the indicator function
\beq\label{test sing 2}
\Theta_\tau^{(4)}=\bra F_{\tau},\M^{(4)}\cet_{L^2(U)},
\eeq
where
$\M^{(4)}$ is given by (\ref{4th interaction for wave eq B})
with  $u_{(j)}=Qf_{(j)}$, $j=1,2,3,4$, where  $f_{(j)}\in  
\I^{n+1}(Y(x_j,\xi_j;t_0,s_0))$, $n\leq -n_1$,
  and $F_\tau$ is the source producing a gaussian beam $Q^*F_\tau$
that propagates to the past along the geodesic $\gamma_{x_5,\xi_5}(\R_-)$,
see (\ref{Ftau source}).

Similar results to the ones given in Proposition \ref{lem:analytic limits A}
are
valid. Let us consider next the case when $(x_5,\xi_5)$ comes from the  
4-intersection
of  rays corresponding to $(\vec x,\vec \xi)=((x_j,\xi_j))_{j=1}^4$ and
$q$ is the corresponding intersection point, that is, $q=\gamma_{x_j,\xi_j}(t_j)$ for
all $j=1,2,3,4,5$.
Then 
\beq\label{indicator 2}
\Theta^{(4)}_\tau\sim 
\sum_{k=m}^\infty s_{k}\tau^{-k}
  \eeq
as $\tau\to \infty$ where   $m=-4n+4$.
Moreover,
let $b_j=(\dot\gamma_{x_j,\xi_j}(t_j))^\flat$ and
  $\bsequence=(b_{j})_{j=1}^5\in (T^*_q\hattuM _0)^5$,
  $w_j$ be the principal symbols of the waves $u_{(j)}$
  at $(q,b_j)$, and ${\bf w}=(w_j)_{j=1}^5$.
Then we see as in Proposition \ref{lem:analytic limits A}  that there is
  a real-analytic function $\mathcal G(\bsequence,{\bf w})$ such that
  the leading order term in (\ref{indicator})  satisfies
  \beq\label{definition of G 2}
s_{m}=
\mathcal  G(\bsequence,{\bf w}).
\eeq

The proof of Prop.\ \ref{singularities in Minkowski space} dealing
with the Einstein equations needs significant changes and we need to prove  
the following:

\begin{proposition}\label{singularities in Minkowski space for wave equation}
The  function $\ \mathcal  G(\bsequence,{\bf w})
$ given in (\ref{definition of G 2}) for the non-linear wave equation
is a non-identically vanishing real-analytic function.
\end{proposition}

\noindent{\bf Proof.}
Let us use the notations introduced in Prop.\ \ref{singularities in  
Minkowski space}.

As for the Einstein equations, we consider light-like vectors
\ba
b_5=(1,1,0,0),\quad b_j=(1,1-\frac  
12\rhoepsilon_j^2,\rhoepsilon_j+O(\rhoepsilon_j^3),\rhoepsilon_j^3),\quad  
j=1,2,3,4,
\ea
  in the Minkowski space $\R^{1+3}$, endowed with the standard metric $g=\diag(-1,1,1,1)$, where the terms
  $O(\rhoepsilon_k^3)$ are such that the vectors $b_j$, $j\leq 5$,
are  light-like. Then
\ba
g(b_5,b_j)= -\frac 12 \rhoepsilon_j^2,\quad
g(b_k,b_j)=-\frac 12 \rhoepsilon_k^2-\frac 12  
\rhoepsilon_j^2+O(\rhoepsilon_k\rhoepsilon_j).
\ea
Below, we denote $\omega_{kj}=g(b_k,b_j)$.
Note that if $\rhoepsilon_j<\rhoepsilon_k^4$, we have
$\omega_{kj}=-\frac 12 \rhoepsilon_k^2+O(\rhoepsilon_k^3).$

{For the wave equation,
we use different parameters $\rhoepsilon_j$ than for the Einstein  
equations, and
define (so, we use here the "unordered" numbering 4-2-1-3)
\beq\label{eq: ordering of epsilons wave equation}
\rhoepsilon_4=\rhoepsilon_2^{100},\  
\rhoepsilon_2=\rhoepsilon_1^{100},\hbox{ and  
}\rhoepsilon_1=\rhoepsilon_3^{100}.
\eeq
Below in this proof, we denote $\vec\rhoepsilon\to 0$ when
$\rhoepsilon_3\to 0$ and
$\rhoepsilon_4,$ $\rhoepsilon_2,$ and $\rhoepsilon_1$ are defined using
$\rhoepsilon_3$ as in (\ref{eq: ordering of epsilons wave equation}).

Let us next consider in Minkowski space
the coordinates $(x^j)_{j=1}^4$ such that
$K_j=\{x^j=0\}$ are light-like hyperplanes and  the waves $u_j=u_{(j)}$ that
satisfy in the Minkowski space $\square u_j=0$ and can be written
as
\ba
u_j(x)=\int_{\R}e^{i x^j \theta }a_j(x,\theta)\,d\theta,
\quad
a_j(x,\theta^{\prime})\in S^{n}(\R^4;\R\setminus 0),\quad j\leq 4,
\ea
and
\ba
u^\tau(x) =\chi(x^0)w_{(5)} \exp(i\tau b^{(5)}\,\cdotp x).
\ea
Note that the singular supports of the waves $u_j$, $j=1,2,3,4,$  
intersect then at the point
$\cap_{j=1}^4 K_j=\{0\}$.
Analogously to the definition  (\ref{definition of G}) we considered
  for  the Einstein equations, we
   define  the (Minkowski) indicator function
\ba
\mathcal G^{({\bf m})}(v,{\bf b})=
\lim_{\tau\to\infty} \tau^{m}(\sum_{\b\leq n_1}
\sum_{\sigma\in \Sigma(4)}
T^{({\bf m}),\b}_{\tau,\sigma}+\tilde T^{({\bf m}),\b}_{\tau,\sigma}),
\ea
where
\ba
T^{({\bf m}),\beta}_{\tau,\sigma}
&=&\bra Q_0(u^\tau\, \cdotp \a u_{\sigma(4)}), h\,\cdotp  \a  
u_{\sigma(3)}\,\cdotp Q_0(\a u_{\sigma(2)}\,\cdotp u_{\sigma(1)})\cet,\\
\tilde T^{({\bf m}),\beta}_{\tau,\sigma}
&=&\bra u^\tau,h\a\,Q_0(\a u_{\sigma(4)}\,\cdotp  \,u_{\sigma(3)})\,\cdotp
  Q_0(\a u_{\sigma(2)}\,\cdotp \, u_{\sigma(1)})\cet.
\ea

As for the Einstein equations, we
see that when $\alpha$ is equal to the value of the function $a(t,y)$
at the intersection point $q=0$ of the waves,
we have
  $\mathcal G^{({\bf m})}(v,{\bf b})=\mathcal G(v,{\bf b})$.

Similarly to the Lemma \ref{lem:analytic limits A}  we analyze next  
the functions
\ba
\Theta_\tau^{({\bf m})}=\sum_{\beta\in J_\ell}\sum_{\sigma\in  
\Sigma(4)}(T_{\tau,\sigma}^{({\bf m}),,\beta}+\tilde  
T_{\tau,\sigma}^{({\bf m}),\beta}).
\ea
Here $({\bf m})$ refers to ``Minkowski''.
We denote $T_{\tau}^{({\bf m}),\beta}=T_{\tau,id}^{({\bf m}),\beta}$
and  $\tilde T_{\tau}^{({\bf m}),\beta}=\tilde T_{\tau,id}^{({\bf m}),\beta}$.

Let us first consider the case when the permutation
$\sigma=id$. Then, as in the proof of Prop.\ \ref{singularities in  
Minkowski space},
in the case  when $\vec S^\beta =(Q_0,Q_0)$, we have
\ba
T^{({\bf m}),\beta}_\tau\\
&& \hspace{-2cm}=C_1 \det(A) \cdotp
(i\tau)^{m}(1+O(\frac 1\tau))
\vec\rhoepsilon^{\,2\vec n}
(\omega_{45}\omega_{12})^{-1}\rhoepsilon_4^{-4}\rhoepsilon_2^{-4}\rhoepsilon_1^{-4}\rhoepsilon_3^{2}\cdotp\P\\
&& \hspace{-2cm}=C_2\det(A) \cdotp
(i\tau)^{m}(1+O(\frac 1\tau))
\vec\rhoepsilon^{\,2\vec n}
\rhoepsilon_4^{-4-2}\rhoepsilon_2^{-4}\rhoepsilon_1^{-4-2}\rhoepsilon_3^{-2}\cdotp\P
\ea
where $\P$ is the product of the principal symbols of the waves $u_j$
at zero, $\vec\rhoepsilon^{\,2\vec n}=
\rhoepsilon_1^{2n}\rhoepsilon_2^{2n}\rhoepsilon_3^{2n}\rhoepsilon_4^{2n}$, and  
$C_1$ and $C_2$ are non-vanishing.
Similarly, a direct computation yields
\ba
\tilde T^{({\bf m}),\beta}_\tau\\
&& \hspace{-2cm}=C_1 \det(A) \cdotp
(i\tau)^{n}(1+O(\frac 1\tau))
\vec\rhoepsilon^{\,2\vec n}
(\omega_{43}\omega_{21})^{-1}\rhoepsilon_4^{-4}\rhoepsilon_2^{-4}\rhoepsilon_1^{-4}\rhoepsilon_3^{-4}\cdotp\P\\
&& \hspace{-2cm}=C_2\det(A) \cdotp
(i\tau)^{m}(1+O(\frac 1\tau))
\vec\rhoepsilon^{\,2\vec n}
\rhoepsilon_4^{-4}\rhoepsilon_2^{-4}\rhoepsilon_1^{-4-2}\rhoepsilon_3^{-4-2}\cdotp\P,
\ea
where again, $\P$ is the product of the principal symbols of the waves $u_j$
at zero and $C_1$ and $C_2$ are non-vanishing.

Considering formula (\ref{4th interaction for wave eq B}), we
see that for the wave equation we do not need to consider the
terms  that for the Einstein equations correspond to the
cases when $\vec S^\beta =(Q_0,I)$, $\vec S^\beta =(I,Q_0)$,
or $\vec S^\beta =(I,I)$ as the corresponding terms do not appear
in formula (\ref{4th interaction for wave eq B}).

Let us now consider permutations $\sigma$ of the indexes $(1,2,3,4)$
and compare the terms
\ba
& &L^{({\bf m}),\beta}_{\sigma}=\lim_{\tau\to \infty} \tau^{m}  
T^{({\bf m}),\beta}_{\tau,\sigma},\\
&& \tilde L^{({\bf m}),\beta}_{\sigma}=\lim_{\tau\to \infty}  \tau^{m} \tilde  
T^{({\bf m}),\beta}_{\tau,\sigma}.
\ea
Due to the presence of $\omega_{45}\omega_{12}$ in the above computations,
we observe that  all the terms
  $ \tilde  L^{({\bf m}),\beta}_{\tau,\sigma}/L^{({\bf m}),\beta}_{\tau,id}\to 0$
  as $\vec \rhoepsilon\to 0$, see (\ref{eq: ordering of epsilons wave  
equation}).
Also, if $\sigma\not=(1,2,3,4)$ and  $\sigma\not=\sigma_01=(2,1,3,4)$,
we see that  $ \tilde   
L^{({\bf m}),\beta}_{\tau,\sigma}/L^{({\bf m}),\beta}_{\tau,id}\to 0$
  as $\vec \rhoepsilon\to 0$.
  Also, we observe that $ L^{({\bf m}),\beta}_{\tau,\sigma_1}=
L^{({\bf m}),\beta}_{\tau,id}$.
Thus we see that the equal terms
$ L^{({\bf m}),\beta}_{\tau,\sigma_1}=
  L^{({\bf m}),\beta}_{\tau,id}$
that give the largest contributions as $\vec \rhoepsilon\to 0$
and that when $\P\not =0$ the sum
\ba
S( \vec \rhoepsilon,\P)=\sum_{\sigma\in \Sigma(4)}(
  L^{({\bf m}),\beta}_{\tau,\sigma}+\tilde L^{({\bf m}),\beta}_{\tau,\sigma})
\ea
is non-zero when
$\rhoepsilon_3>0$ is small enough   and
$\rhoepsilon_4,$ $\rhoepsilon_2,$ and $\rhoepsilon_1$ are defined using
$\rhoepsilon_3$ as in (\ref{eq: ordering of epsilons wave equation}).
Since the indicator
function is real-analytic, this shows that the indicator function
in non-vanishing in a generic set.   \hfill \Box \medskip

We need also to change the
  the singularity {\it detection condition} ({D}) with light-like  
directions $(\vec x,\vec \xi)$
  as follows:  We define that point $y\in U_{\hat g}$,
   satisfies the singularity {\it detection condition} (${D}^\prime$)  
with light-like directions $(\vec x,\vec \xi)$
   and $t_0,\hat s>0$
  if
  \medskip  
  
($D^\prime$) For any $s,s_0\in (0,\hat s)$  there are  
$(x_j^{\prime},\xi_j^{\prime})\in \W_j(s;x_j,\xi_j)$, $j=1,2,3,4$,
and ${f}_{(j)}\in {\mathcal  
I_{S}}^{n+1}(Y((x_j^{\prime},\xi_j^{\prime});t_0,s_0))$, and such  
that if  $u_{\vec \e}$ of is the solution  of (\ref{PABC eq})
with the source ${f}_{\vec \e}=\sum_{j=1}^4
\e_j{f}_{(j)}$, then
the function
$\p_{\vec \e}^4u_{\vec \e}|_{\vec \e=0}$ is not
$C^\infty$-smooth in any neighborhood of $y$.
  \medskip

When condition (D) is replace by ($D^\prime$), the considerations in  
the Sections
  \ref{sec: normal coordinates} and \ref{subsection combining}  show  
that we can recover
the conformal class of the metric. This proves Theorem \ref{main  
thm3}. \hfill \Box \medskip

\noindent
{\bf Proof.} (of Corollary \ref{coro thm3}). 
Let us denote $W_j=I^+(p^-_j)\cap I^-(p^+_j)\subset M_j$.
By Theorem \ref{main thm3},  there is a map $\Psi:W_1\to W_2$
such that  the product type metrics $g_1=-dt^2+h_1(y)$
and $g_2=-dt^2+h_2(y)$ are conformal. Note that the above methods
to determine the conformal class of the metric but not the metric itself. 

To construct the metric, let us consider the linearized waves. 
Let $B_j(y,r)$ denote the Riemannian ball of $(N_j,h_j)$ with center $y$ and radius $r$.
Consider a set 
$K_j=(s_j^\prime,s_j^{\prime\prime})\times B_j(p_j,r)\subset W_j$.
On $(s_j^\prime,s_j^{\prime\prime})\times  (N_j\setminus B_j(p_j,r)) \subset M_j$ we can consider
the linear wave equation
\beq\label{lin wave again 1}
& &\square_{g_j}w^{(j)}=0,\quad \hbox{on }(s_j^\prime,s_j^{\prime\prime})\times  (N_j\setminus B_j(p_j,r)) ,
\\ & &w^{(j)}|_{(s_j^\prime,s_j^{\prime\prime})\times  \p B_j(p_j,r))}=h 
\nonumber \\
\nonumber & & w^{(j)}|_{t= s_j^\prime}=0,\quad  \p_t w^{(j)}|_{t= s_j^\prime}=0.
\eeq
For this wave equation we define the Dirichlet-to-Neumann operator
$\Lambda^{(j)}:h\mapsto \p_\nu w^{(j)}|_{(s_j^\prime,s_j^{\prime\prime})\times  \p B_j(p_j,r))}$.
We observe that any solution $w^{(j)}$ of (\ref{lin wave again 1}) can
be continued to a solution of 
\beq\label{lin wave again 2}
& &\square_{g_j}u^{(j)}=F_j,\quad \hbox{on }(s_j^\prime,s_j^{\prime\prime})\times  N_j ,\\
\nonumber & &u^{(j)}|_{t= s_j^\prime}=0,\quad  \p_t u^{(j)}|_{t= s_j^\prime}=0.
\eeq
where $F_j$ is some source supported in $K_j$, that is, for all $ w^{(j)}$ and $h$ there
exists $ u^{(j)}$ and $F_j$ such that $ u^{(j)}= w^{(j)}$ on $(s_j^\prime,s_j^{\prime\prime})\times  (N_j\setminus B_j(p_j,r)).$ 

Consider next a given $h$.
We see that if $F_j$ is such that it satisfies
\beq\label{cond NL}
h=\p_\e\bigg (L_{U_j}(\e F_j)|_{(s_j^\prime,s_j^{\prime\prime})\times  \p B_j(p_j,r))}
\bigg)\bigg|_{\e=0},
\eeq
then
  \ba
\Lambda_jh&=& \p_\nu w^{(j)}|_{(s_j^\prime,s_j^{\prime\prime})\times  \p B_j(p_j,r))}\\
&=& \p_\e\bigg (\p_\nu(L_{U_j} (\e F_j))|_{(s_j^\prime,s_j^{\prime\prime})\times  \p B_j(p_j,r))}
\bigg)\bigg|_{\e=0},
\ea
and by the above considerations, for any $h$ there exists some $F_j$ for which
(\ref{cond NL}) is valid.

This means that by linearizing the map  $L_{U_j}$ 
we can determine the map $\Lambda_j$. Now, using
 \cite{KaKu2,KKL}, see also \cite{BelKur,KrKL}, we see that if  $L_{U_1}= L_{U_2}$ then $(W_1,g_1)$ and $(W_2,g_2)$ are isometric.


 As 
$g_1$ and $g_2$ are independent
of $t$, we see that there is a diffeomorphism $\Psi:W_1\to W_2$
of the form $\Psi(t,y)=(t,\psi(y))$ such that $g_1=\Psi^*g_2$.
Note also that if $\pi_2:(t,y)\mapsto y$,
then $h_1=\psi^*h_2$ on $\pi_2(W_1)$.
Thus, the metric tensors $h_1$ and $h_2$
are isometric. 

As  the linearized waves $u_{(j)}=Qf_{(j)}$ depend only on $W_j$ and the
metric $g_j$, using the proof of Theorem \ref{main thm3}
we see that the indicator functions
$\mathcal  G(\bsequence,{\bf w})$ for $(U_1,g_1,a_1)$
and $(U_1,g_1,a_1)$ coincide for all $\bsequence$ and ${\bf w}$.
This implies that $a_1(t,y)^3=a_2(t,y)^3$ for all $(t,y)\in W$.
Hence $a_1(t,y)=a_2(t,y)$ for all $(t,y)\in W$.
  \hfill \Box \medskip
}
}

\medskip

\noindent{\bf Acknowledgements.} The authors express their gratitude to
 MSRI, the Newton  
Institute, the Fields Institute and  
the Mittag-Leffler Institute, where parts of this work have been done.

YK was partly supported by EPSRC and the AXA professorship.
ML was partly supported by the Finnish Centre of
Excellence in Inverse Problems Research 2012-2017.
GU was partly supported  by NSF, a Clay Senior Award at MSRI,  a Chancellor Professorship at UC Berkeley,  
a Rothschild Distinguished Visiting Fellowship at the Newton Institute, the Fondation de Sciences Math\'ematiques de Paris,
and a Simons Fellowship.


 \end{document}